\documentclass[10pt]{amsart}
\usepackage[english,russian]{babel}
\usepackage[cp1251]{inputenc}
\usepackage{amsmath}
\usepackage{amssymb}
\usepackage{amsfonts}

\usepackage[linktocpage=true, colorlinks=true, linkcolor=blue, citecolor=blue, urlcolor=blue]{hyperref}

\begin{document}
\renewcommand{\refname}{References}
\renewcommand\contentsname{Contents}

\thispagestyle{empty}

\title[Four New Forms of the Taylor--Ito and Taylor--Stratonovich Expansions]
{Four New Forms of the Taylor--Ito and Taylor--Stratonovich Expansions and its
Application to the High-Order Strong Numerical Methods
for Ito Stochastic Differential Equations}
\author[D.F. Kuznetsov]{Dmitriy F. Kuznetsov}
\address{Dmitriy Feliksovich Kuznetsov
\newline\hphantom{iii} Peter the Great Saint-Petersburg Polytechnic University,
\newline\hphantom{iii} Polytechnicheskaya ul., 29,
\newline\hphantom{iii} 195251, Saint-Petersburg, Russia}%
\email{sde\_kuznetsov@inbox.ru}

\maketitle {\small
\begin{quote}
\vspace{5mm}
\noindent{\sc Abstract.} 
The problem of the Taylor--Ito and Taylor--Stratonovich 
expansions of the Ito stochastic 
processes in a neighborhood of a fixed moment of time is considered.
The classical forms of the Taylor--Ito and Taylor--Stratonovich 
expansions are transformed to 
the four new representations, which includes
the minimal sets of different types of iterated Ito and Stratonovich
stochastic integrals. 
Therefore, these representations
(the so-called unified Taylor--Ito and Taylor--Stratonovich expansions)
are more convenient for constructing 
of high-order strong numerical 
methods for Ito stochastic differential equations.
Explicit one-step strong numerical schemes with the orders
of convergence 1.0, 1.5, 2.0, 2.5, and 3.0
based on the unified Taylor--Ito and Taylor--Stratonovich expansions 
are derived. Effective mean-square approximations 
of iterated Ito and Stratonovich 
stochastic integrals from these numerical schemes
are constructed on the base of the multiple Fourier--Legendre
series with multiplicities 1 to 6.

\medskip
\end{quote}
}

\vspace{15mm}

\linespread{1.75}

\tableofcontents

\linespread{1.0}

\vspace{15mm}

\section{Introduction}

\vspace{5mm}

Let $(\Omega,$ ${\rm F},$ ${\sf P})$ be a complete probability space, let 
$\{{\rm F}_t, t\in[0, T]\}$ be a nondecreasing right-con\-ti\-nu\-o\-us family 
of $\sigma$-subfields of ${\rm F},$
and let ${\bf f}_t$ be a standard $m$-dimensional Wiener stochastic 
process, which is
${\rm F}_t$-measurable for any $t\in[0, T].$ We assume that the components
${\bf f}_{t}^{(i)}$ $(i=1,\ldots,m)$ of this process are independent. Consider
an Ito stochastic differential equation (SDE) in the integral form

\begin{equation}
\label{1.5.2}
{\bf x}_t={\bf x}_0+\int\limits_0^t {\bf a}({\bf x}_{\tau},\tau)d\tau+
\int\limits_0^t B({\bf x}_{\tau},\tau)d{\bf f}_{\tau},\ \ \
{\bf x}_0={\bf x}(0,\omega),\ \ \ \omega\in\Omega.
\end{equation}

\vspace{3mm}
\noindent
Here ${\bf x}_t$ is some $n$-dimensional stochastic process 
satisfying to Ito SDE (\ref{1.5.2}). 
The nonrandom functions ${\bf a}: \mathbb{R}^n\times[0, T]\to\mathbb{R}^n$,
$B: \mathbb{R}^n\times[0, T]\to\mathbb{R}^{n\times m}$
guarantee the existence and uniqueness up to stochastic 
equivalence of a solution
of (\ref{1.5.2}) \cite{1}. The second integral on 
the right-hand side of (\ref{1.5.2}) is 
interpreted as an Ito stochastic integral.
Let ${\bf x}_0$ be an $n$-dimensional random variable, which is 
${\rm F}_0$-measurable and 
${\sf M}\{\left|{\bf x}_0\right|^2\}<\infty$ 
(${\sf M}$ denotes a mathematical expectation).
Also we assume that
${\bf x}_0$ and ${\bf f}_t-{\bf f}_0$ are independent when $t>0.$

It is well known \cite{KlPl2}-\cite{Mi3}
that Ito SDEs are 
adequate mathematical models of dynamic systems of 
different physical origin that are affected by random
perturbations. For example, Ito SDEs are used as 
mathematical models in stochastic mathematical finance,
hydrology, seismology, geophysics, chemical kinetics, 
population dynamics, electrodynamics, medicine and other fields 
\cite{KlPl2}-\cite{Mi3}.
Numerical integration of Ito SDEs based on the 
strong convergence criterion of approximations \cite{KlPl2} is widely 
used for the numerical simulation of sample trajectories of 
solutions to Ito SDEs (which is required for constructing new 
mathematical models on the basis of such equations and for the 
numerical solution of different mathematical problems
connected with Ito SDEs).
Among these problems, we note the following:
filtering of signals under influence of random noises 
in various statements (linear Kalman--Bucy filtering, 
nonlinear optimal filtering, filtering of continuous time Markov 
chains with a finite space of states, etc.), optimal stochastic 
control (including incomplete data control), 
testing estimation procedures of parameters of stochastic systems, 
stochastic stability and bifurcations analysis
\cite{KlPl2}, \cite{KPS}.

Exact solutions of Ito SDEs are known in rather rare cases. 
It is for this reason that it becomes necessary 
to construct numerical procedures for solving these equations.

In this paper, a promising approach \cite{KlPl2}-\cite{Mi3} to the numerical 
integration of Ito SDEs based on the stochastic analogues of the 
Taylor formula (Taylor--Ito and Taylor--Stratonovich 
expansions) \cite{PW1}-\cite{k1} is used. 
This approach uses a finite discretization of the time variable 
and the numerical simulation of the solution to the Ito SDE at 
discrete instants of time using the stochastic analogues of the Taylor 
formula mentioned above. 
A number of works (e.g., \cite{KlPl2}-\cite{Mi3}) describe numerical schemes 
with the strong orders of convergence $1.5, 2.0$, $2.5$, and $3.0$ for 
the Ito SDEs; however, they do not contain efficient procedures 
of the mean-square approximation of iterated stochastic 
integrals involved in these schemes
for the case of non-commutative noise.

In this paper we consider the unified Taylor--Ito 
and Taylor--Stratonovich expansions \cite{kk1}, \cite{k1} 
which makes 
it possible (in contrast with its classical analogues \cite{PW1}, \cite{KlPl1}) 
to use the 
minimal sets of iterated Ito and Stratonovich stochastic integrals; this is a 
simplifying factor 
for the numerical methods implementation. 
We prove 
the unified Taylor--Ito expansion \cite{kk1}
with using the slightly different approach (which is taken from \cite{k1})
in comparison with the approach from \cite{kk1}.
Moreover, we obtain another (second) version of 
the unified Taylor--Ito expansion \cite{20a}, \cite{20b}.
In addition, we construct two new forms of the Taylor--Stratonovich expansion
(the so-called unified Taylor--Stratonovich expansions \cite{k1}).
Futhermore, in this paper
we study methods \cite{arxiv-1}-\cite{new-art-1xxy}
of numerical simulation for iterated Ito and Stratonovich
stochastic integrals of multiplicities $1, 2, 3, 4, 5, 6, \ldots $
used in the 
strong numerical methods for Ito SDEs 
\cite{KlPl2}-\cite{Mi3}, \cite{7}-\cite{11}, \cite{19}-\cite{12aa-afterxxx}.
To approximate the iterated Ito and Stratonovich stochastic 
integarls appearing in the numerical schemes with the strong 
orders of convergence $1.0, 1.5,$ $2.0, 2.5, 3.0$ etc., 
the method of generalized multiple Fourier series
and especially method of multiple Fourier--Legendre 
series are studied in \cite{arxiv-1}-\cite{new-art-1xxy}. 
It is important to note that the method of 
generalized multiple Fourier series \cite{arxiv-1}-\cite{new-art-1xxy} 
does not lead to the partitioning of the 
integration interval of the iterated Ito and Stratonovich 
stochastic integrals under consideration; 
this interval is the integration step of the numerical methods used 
to solve Ito SDEs; therefore, it is already fairly small and does not 
need to be partitioned. 
Computational experiments \cite{7} show that the application of numerical 
simulation for iterated stochastic integrals
(in which the interval of integration is partitioned) leads to 
unacceptably high computational cost and accumulation of 
computation errors. 
Also note that the Legendre polynomials have
essential advantage over the trigonomentric functions
(see \cite{arxiv-12}, \cite{5-004}) in the framework of the 
method of generalized multiple Fourier series \cite{arxiv-1}-\cite{new-art-1xxy}
for the mean-square approximation of iterated Ito and
Stratonovich stochastic integrals.

The rest of the article is organized as follows.
In the introduction (below) we 
consider a brief review of the literature on the 
problem of construction of the Taylor--Ito and Taylor--Stratonovich
expansions for the solutions of Ito SDEs.
Sect. 2 is devoted 
to the integration order replacement technique for 
iterated Ito stochastic integrals.
In Sect. 3, we consider the classical Taylor--Ito expansion
while Sect. 4 and Sect. 5 are devoted to the 
first and second forms of the so-called unified
Taylor--Ito expansion correspondingly.
The classical Taylor--Stratonovich expansion is considered in Sect. 6.
The first and second forms of the
unified Taylor--Stratonovich expansion 
are derived in Sect. 7--9.
In Sect. 10, we give a comparative analysis 
of the unified Taylor--Ito and Taylor--Stratonovich expansions
with the classical Taylor--Ito and Taylor--Stratonovich expansions.
Application of the first form of the unified 
Taylor--Ito expansion to the high-order
strong numerical methods for Ito SDEs
is considered in Sect. 11.
In Sect. 12, we construct 
the high-order
strong numerical methods for Ito SDEs
on the base of the first form of the 
unified Taylor--Stratonovich expansion.
Sect. 13 is devoted to the 
effective method 
of the mean-square approximation
of iterated Ito and Stratonovich 
stochastic integrals based on generalized multiple 
Fourier series. In Sect.~14, we discuss
the connection between Theorems 7--16 (see Sect.~13) 
and Wong--Zakai approximation.

Let us consider the following 
iterated Ito and Stratonovich stochastic integrals

\vspace{-1mm}
\begin{equation}
\label{ito}
J[\psi^{(k)}]_{s,t}=\int\limits_t^s\psi_k(t_k) \ldots \int\limits_t^{t_{2}}
\psi_1(t_1) d{\bf w}_{t_1}^{(i_1)}\ldots
d{\bf w}_{t_k}^{(i_k)},
\end{equation}
\begin{equation}
\label{str}
J^{*}[\psi^{(k)}]_{s,t}=
{\int\limits_t^{*}}^s \psi_k(t_k) \ldots {\int\limits_t^{*}}^{t_2}
\psi_1(t_1) d{\bf w}_{t_1}^{(i_1)}\ldots
d{\bf w}_{t_k}^{(i_k)},
\end{equation}

\vspace{2mm}
\noindent
where $0\le t<s\le T,$ every $\psi_l(\tau)$ $(l=1,\ldots,k)$ is a nonrandom
function 
at the interval $[t,T],$ ${\bf w}_{\tau}^{(i)}={\bf f}_{\tau}^{(i)}$
for $i=1,\ldots,m$ and
${\bf w}_{\tau}^{(0)}=\tau,$\ \
$i_1,\ldots,i_k = 0, 1,\ldots,m,$

\vspace{-1mm}
$$
\int\limits\ \ \hbox{and}\ \  \int\limits^{*}
$$ 

\vspace{2mm}
\noindent
denote Ito and 
Stratonovich stochastic integrals,
respectively.
In this paper we use the definition of the Stratonovich 
stochastic integral from \cite{KlPl2} (also see \cite{20aa}, Chapter~2).

It sould be noted that one of the main problems 
when constructing the high-order strong numerical methods
for Ito SDEs on the base of the Taylor--Ito and
Taylor--Stratonovich expansions is the mean-square approximation of
iterated Ito and Stratonovich stochastic integrals.
Obviously, in the absence of procedures for the numerical simulation 
of stochastic integrals, the mentioned numerical methods 
are unrealizable in practice.
For this reason, in Sect. 13 
we give a brief overview to the effective method
of the mean-square approximation of
iterated Ito and Stratonovich stochastic integrals
(\ref{ito}) and (\ref{str}) of arbitrary multiplicity $k$ ($k\in\mathbb{N}$),
which is proposed and developed by the author of this article in
a number of publications \cite{arxiv-1}-\cite{new-art-1xxy}. This method is based
on the generalized multiple Fourier series converging 
in the mean-square sense.
The extensive practical material on expansions
and mean-square approximations 
of iterated Ito and Stratonovich stochastic integrals
of multiplicities 1 to 6 from the Taylor--Ito and 
Taylor--Stratonovich expansions is given in Sect. 13.
In the mentioned section, the main focus is on approximations based on 
multiple Fourier--Legendre series. Such approximations 
is more effective in comparison with the trigonometric 
approximations \cite{arxiv-12}, \cite{5-004} at least 
for the numerical methods with the
strong order 1.5 of convergence and higher \cite{arxiv-12}, \cite{5-004}.

Let us give a brief review of the literature on the 
problem of construction of the Taylor--Ito and Taylor--Stratonovich
expansions for the solutions of Ito SDEs.
A few variants of a stochastic analog of the Taylor formula 
have been obtained in \cite{KlPl2}-\cite{KlPl1}
for the stochastic processes
in the form $R({\bf x}_s, s)$, where ${\bf x}_s$ is a solution 
of the Ito SDE (\ref{1.5.2}) and $R: \mathbb{R}^n\times [0, T]
\to \mathbb{R}^1$ is a nonrandom
sufficiently smooth function.

The first result in this direction called the Ito--Taylor expansion 
has been obtained in \cite{PW1}, \cite{KlPl1}. This result gives
an expansion of the process $R({\bf x}_s, s)$ into a series 
such that every term (if $k>0$) contains an iterated Ito stochastic
integral 

\begin{equation}
\label{re11}
\int\limits_t^s\ldots \int\limits_t^{t_{2}}
d{\bf w}_{t_1}^{(i_1)}\ldots
d{\bf w}_{t_k}^{(i_k)}
\end{equation}

\vspace{3mm}
\noindent
as a factor, where $0\le t<s\le T,$  $i_1,\ldots,i_k=0,1,\ldots,m.$
Obviously that the iterated Ito stochastic integral (\ref{re11})
is a particular case of (\ref{ito}) for
$\psi_1(\tau),\dots,\psi_k(\tau)\equiv 1.$

In \cite{KlPl1}, another expansion of 
the stochastic process $R({\bf x}_s, s)$ in a series has been
derived. Instead of the Ito integrals, the iterated Stratonovich 
stochastic integrals

\begin{equation}
\label{str11}
{\int\limits_t^{*}}^s \ldots {\int\limits_t^{*}}^{t_2}
d{\bf w}_{t_1}^{(i_1)}\ldots
d{\bf w}_{t_k}^{(i_k)}
\end{equation}

\vspace{3mm}
\noindent
were used; the corresponding expansion was called the 
Stratonovich--Taylor expansion. In the formula (\ref{str11}), the indices
$i_1,\ldots,i_k$ take values $0, 1,\ldots,m$.

In \cite{kk1}
the Ito--Taylor expansion from \cite{PW1}, \cite{KlPl1} 
is reduced to the
interesting and unexpected
form
(called the unified Taylor--Ito expansion)
with the help of special transformations \cite{3000}
(also see \cite{arxiv-22}, \cite{444}, \cite{3001}). 
Every term of this
expansion (if $k>0$) contains an iterated Ito stochastic
integral of the form
\begin{equation}
\label{ll1}
\int\limits^ {s} _ {t} (s-t _
{k}) ^ {l_ {k}} 
\ldots \int\limits^ {t _ {2}} _ {t} (s-t _ {1}) ^ {l_ {1}} d
{\bf f} ^ {(i_ {1})} _ {t_ {1}} \ldots 
d {\bf f} _ {t_ {k}} ^ {(i_ {k})},
\end{equation}

\vspace{2mm}
\noindent
where
$0\le t<s\le T,$ $l_1,\ldots,l_k=0, 1, 2,\ldots $ and $i_1,\ldots,i_k = 1,\ldots,m$.

It is worth to mention another form of the 
unified Taylor--Ito expansion \cite{20a}, \cite{20b}, \cite{444} (also see
\cite{7}-\cite{11}, \cite{15}-\cite{12aa-afterxxx}). 
Terms of the latter expansion contain iterated 
Ito stochastic integrals of the form
\begin{equation}
\label{ll11}
\int\limits^ {s} _ {t} (t-t _{k}) ^ {l_ {k}} 
\ldots \int\limits^ {t _ {2}} _ {t} (t-t _ {1}) ^ {l_ {1}} 
d{\bf f} ^ {(i_ {1})} _ {t_ {1}} \ldots 
d {\bf f} _ {t_ {k}} ^ {(i_ {k})},
\end{equation}

\vspace{2mm}
\noindent
where
$l_1,\ldots,l_k=0, 1, 2,\ldots $ and $i_1,\ldots,i_k = 1,\ldots,m$.

In this paper, we derive 
two new forms of the Taylor--Ito expansions (the so-called
unified Taylor--Ito expansions 
\cite{20a}, \cite{20b}, \cite{444} (also see
\cite{7}-\cite{11}, \cite{15}-\cite{12aa-afterxxx}))
using an approach which is taken
from \cite{k1}.
Obviously that some of iterated Ito stochastic integrals of the form 
(\ref{re11}) or (\ref{str11}) are connected by linear relations, 
while
this is not the case for integrals of the form 
(\ref{ll1}), (\ref{ll11}). In this sense, the total quantity
of stochastic integrals of the
form (\ref{ll1}) or (\ref{ll11}) is minimal. 
Futhermore, in this article we construct two new forms of the 
Taylor--Stratonovich expansion (the so-called unified
Taylor--Stratonovich expansions \cite{k1}) such that every term
(if $k>0$)
contains as a multiplier 
an iterated Stratonovich stochastic integral
of one of two types
\begin{equation}
\label{a111}
{\int\limits_t^{*}}^{s}(t-t_k)^{l_{k}}\ldots
{\int\limits_t^{*}}^{t _ {2}}(t-t _ {1}) ^ {l_ {1}} 
d{\bf f} ^ {(i_ {1})} _ {t_ {1}} \ldots
d{\bf f}_{t_ {k}}^{(i_ {k})},
\end{equation}
\begin{equation}
\label{a112}
{\int\limits_t^{*}}^{s}(s-t_k)^{l_{k}}\ldots
{\int\limits_t^{*}}^{t _ {2}}(s-t _ {1}) ^ {l_ {1}} 
d{\bf f} ^ {(i_ {1})} _ {t_ {1}} \ldots
d{\bf f}_{t_ {k}}^{(i_ {k})},
\end{equation}

\vspace{2mm}
\noindent
where 
$l_1,\ldots,l_k=0, 1, 2,\ldots $, $i_1,\ldots,i_k = 1,\ldots,m$,
and $k=1, 2,\ldots $

It is not difficult to see that for the sets of iterated Stratonovich
stochastic integrals
(\ref{a111}) and (\ref{a112}) 
the property of minimality (see above) 
also holds as for 
the sets of iterated Ito 
stochastic integrals (\ref{ll1}), (\ref{ll11}).

As we noted above, the main problem in implementation
of high-order strong numerical methods 
for Ito SDEs is the mean-square approximation of iterated stochastic integrals
(\ref{re11})--(\ref{a112}). 
Obviously, these stochastic integrals are particular cases
of the stochastic integrals 
(\ref{ito}), (\ref{str}).

Taking into account the results of \cite{arxiv-1}-\cite{322} 
and the minimality
of the sets of stochastic integrals of the forms (\ref{ll1})--(\ref{a112}), 
we conclude that the unified Taylor--Ito and Taylor--Stratonovich
expansions based on the iterated stochastic
integrals (\ref{ll1})--(\ref{a112})  may be useful for
constructing of high-order strong numerical
methods with the orders of convergence 1.5, 
2.0, 2.5, 3.0, 3.5, 4.0, 4.5, $\ldots$ 
for Ito SDEs.

\vspace{5mm}

\section{Integration Order Replacement Technique for 
Iterated Ito Stochastic Integrals}

\vspace{5mm}

Let $f_{\tau}, \tau\in [0, T]$ be a scalar standard Wiener process 
that is ${\rm F}_{\tau}$-measurable for every $\tau\in [0, T]$. We
introduce a class $\mathfrak{M}_2([t, T])$ $(t\ge 0)$ of random functions 
$\xi(\tau,\omega)\stackrel{\sf def}{=}
\xi_{\tau}:$ $[t, T]\times\Omega \to \mathbb{R}^1$ 
having the following properties:
these functions are measurable with respect to the pair $(\tau,\omega)$ 
of variables, ${\rm F}_{\tau}$-measurable for every $\tau\in [t, T]$,
and satisfy the conditions

\vspace{-1mm}
$$
\int\limits_t^T
{\sf M}\left\{\xi_{\tau}^2\right\}d\tau<\infty
$$

\vspace{2mm}
\noindent
and ${\sf M}\left\{\xi_{\tau}^2\right\}<\infty$
for any $\tau\in [t, T]$.

On the class $\mathfrak{M}_2([t, T])$ $(t\ge 0)$, we introduce the Hilbert norm

\vspace{-1mm}
$$
\Vert \xi \Vert_{2,T,t}=\left(
\int\limits_t^T{\sf M}\left\{\xi_{\tau}^2\right\}d\tau
\right)^{1/2}.
$$

\vspace{2mm}

Let $\bigl\{\tau_j^{(N)}\bigr\}_{j=1}^N$ 
be a partition of the interval $[t, T]$ such that

\begin{equation}
\label{usl}
t=\tau_0^{(N)}<\tau_1^{(N)}<\ldots <\tau_N^{(N)}=T,\ \ \
\Delta_N=
\hbox{\vtop{\offinterlineskip\halign{
\hfil#\hfil\cr
{\rm max}\cr
$\stackrel{}{{}_{0\le j\le N-1}}$\cr
}} }\left|\tau_{j+1}^{(N)}-\tau_j^{(N)}\right|\to 0\ \ \ \hbox{as}\ \ \ N\to \infty.
\end{equation}

\vspace{3mm}

Let $\xi^{(N)}(\tau,\omega)$
be a sequence of step functions from the space 
$\mathfrak{M}_2([t, T])$ defined as follows

$$
\xi^{(N)}(\tau,\omega)=\xi_{j}(\omega)\ \ \
{\rm w.\ p.\ 1\ \ for}\ \ \tau\in\left[\tau_j^{(N)},\tau_{j+1}^{(N)}\right),\ \ \
j=0, 1,\ldots,N-1,
$$

\vspace{3mm}
\noindent
where here and further w.~p.~1 means with probability 1.

It is known \cite{1} that
for any function $\xi_{\tau}\in\mathfrak{M}_2([t,T])$
there exists a sequence 
$\xi^{(N)}(\tau,\omega)\in\mathfrak{M}_2([t,T]),$  which
converges to the function $\xi_{\tau}$ in the sence of norm 
$\Vert \cdot \Vert_{2,T,t}.$

The mean-square limit

\vspace{-1mm}
\begin{equation}
\label{10.1}
\hbox{\vtop{\offinterlineskip\halign{
\hfil#\hfil\cr
{\rm l.i.m.}\cr
$\stackrel{}{{}_{N\to \infty}}$\cr
}} }\sum_{j=0}^{N-1}\xi^{(N)}(\tau_j^{(N)},\omega)
\left(f(\tau_{j+1}^{(N)},\omega)-
f(\tau_j^{(N)},\omega)\right)
\stackrel{\sf def}{=}\int\limits_t^T\xi_\tau df_\tau
\end{equation}

\vspace{3mm}
\noindent
is called \cite{1} the Ito stochastic integral of a function 
$\xi_{\tau}\in \mathfrak{M}_2([t, T])$. Here
$\xi^{(N)}(\tau,\omega)$ is an arbitrary sequence of step functions
from the class
$\mathfrak{M}_2([t,T])$ converging to the function 
$\xi(\tau,\omega)$ in the sense
of norm
$\Vert\cdot\Vert_{2,T,t}$, i.e.

\vspace{-1mm}
\begin{equation}
\label{jjj}
\hbox{\vtop{\offinterlineskip\halign{
\hfil#\hfil\cr
{\rm lim}\cr
$\stackrel{}{{}_{N\to \infty}}$\cr
}} }\int\limits_t^T{\sf M}\left\{\left|
\xi^{(N)}(\tau,\omega)-\xi(\tau,\omega)\right|^2\right\}d\tau=0.
\end{equation}

\vspace{3mm}

We introduce the class $\mathfrak{Q}_m([t, T])$ $(t\ge 0)$ of Ito processes 
$\eta_{\tau},$ $\tau \in [t, T]$ of the form

\begin{equation}
\label{z900}
\eta_{\tau} = \eta_t + \int\limits_t^{\tau} a_sds +
\int\limits_t^{\tau} b_sdf_s, 
\end{equation}

\vspace{2mm}
\noindent
where $(a_{\tau})^m, (b_{\tau})^m \in  \mathfrak{M}_2([t, T])$ and

$$
\lim\limits_{s\to\tau} {\sf M}\left\{\left\vert b_s - b_{\tau}\right\vert ^4\right\}=0
$$

\vspace{3mm}
\noindent
for all $\tau \in [t, T]$.

Let $C^{2,1}(\mathbb{R}^1\times [t, T])$ $(t\ge 0)$ be the space of functions
$F(x,\tau): \mathbb{R}^1 \times [t, T] \to \mathbb{R}^1$ such that

\vspace{1mm}
$$
\left|\frac{\partial F}{\partial x}(x,\tau)\right|\le K,\ \ 
\left|\frac{\partial^2 F}{\partial x^2}(x,\tau) \right|\le K,\ \ 
\left|\frac{\partial F}{\partial \tau}(x,\tau)\right|\le K,\ \ 
\left|\frac{\partial^2 F}{\partial \tau \partial x}(x,\tau) \right|\le K
$$ 

\vspace{5mm}
\noindent
for all $x\in \mathbb{R}^1$ and $\tau\in [t, T],$ where constant $K$ does not depend on $x,\tau.$

The mean-square limit

\vspace{-1mm}
\begin{equation}
\label{123321.2}
\hbox{\vtop{\offinterlineskip\halign{
\hfil#\hfil\cr
{\rm l.i.m}\cr
$\stackrel{}{{}_{N\to \infty}}$\cr
}} }\sum_{j=0}^{N-1}F\left(\frac{1}{2}\left(
\eta_{\tau_j^{(N)}}+\eta_{\tau_{j+1}^{(N)}}\right),\tau_j^{(N)}\right)
\left(f_{\tau_{j+1}^{(N)}}-
f_{\tau_j^{(N)}}\right)
\stackrel{\sf def}{=}{\int\limits_t^{*}}^T F(\eta_{\tau},\tau)df_\tau
\end{equation}

\vspace{3mm}
\noindent
is called \cite{str} the Stratonovich stochastic integral 
of the process $F(\eta_{\tau}, \tau)$, $\tau\in [t, T]$ $(t\ge 0)$,
where $F(x,\tau)\in C^{2,1}(\mathbb{R}^1\times [t, T]).$
We apply in the formula
(\ref{123321.2}) the same notations as in the formula (\ref{10.1}).

It is known \cite{str} (also see \cite{KlPl2}) 
that under proper conditions, the following 
relation holds

\begin{equation}
\label{d11}
{\int\limits_{t}^{*}}^T F(\eta_{\tau},\tau)df_{\tau}=
\int\limits_t^T F(\eta_{\tau},\tau)df_{\tau}+
\frac{1}{2}\int\limits_t^T \frac{\partial F}{\partial x}(\eta_{\tau},\tau)
b_{\tau}d\tau\ \ \ \hbox{w.\ p.\ 1}.
\end{equation}

\vspace{3mm}

If the Wiener processes in the formulas (\ref{z900}) and 
(\ref{123321.2}) are independent, then

\begin{equation}
\label{d11a}
{\int\limits_{t}^{*}}^T F(\eta_{\tau},\tau)df_{\tau}=
\int\limits_t^T F(\eta_{\tau},\tau)df_{\tau}\ \ \ \hbox{w.\ p.\ 1}.
\end{equation}

\vspace{3mm}

Note that a possible variant of conditions providing
the correctness of the formulas 
(\ref{d11}) and (\ref{d11a}) consists of the following 
conditions

\vspace{1mm}
$$
\eta_{\tau}\in \mathfrak{Q}_4([t,T]),\ \ \
F(\eta_{\tau},\tau)\in\mathfrak{M}_2([t,T]),\ \ \ \hbox{and}\ \ \ 
F(x,\tau)\in C^{2,1}(\mathbb{R}^1\times [t, T]).
$$

\vspace{4mm}

Note that if $F(x,\tau)=F_1(x)F_2(\tau),$ then it suffices to require
that $F(x,\tau)$ be twice differentiable with respect to $x$ 
$($with bounded derivatives$)$ and continuous with respect to $\tau$
$($instead of the condition $F(x,\tau)\in C^{2,1}(\mathbb{R}^1\times [t, T])).$

A theorem allowing the change of order of integration in 
the iterated Ito stochastic integrals has been
proved in \cite{3000}-\cite{3001} (also see \cite{arxiv-22}). 
In what follows, we apply this theorem; let us 
cite its exact formulation and the notations.

It is well known that
the Ito stochastic integral exists
in the mean-square sense (see (\ref{10.1})), if the stochastic process
$\xi(\tau,\omega)\in \mathfrak{M}_2([t,T]),$ 
that is, perhaps this process does not satisfy 
the property of the mean-square continuity on the interval 
$[t,T].$ Let us formulate 
the theorem on integration order replacement for the special 
class of iterated Ito stochastic integrals.
At the same time, the condition of the mean-square continuity 
of integrand in the innermost
stochastic integral will be significant.

Let $\mathfrak{S}_2([t,T])$ $(t\ge 0)$ be the class of functions 
$\xi:$ $[t,T]\times\Omega\rightarrow
\mathbb{R}^1,$ which satisfy the conditions:

\vspace{2mm}

1. $\xi_{\tau} \in \mathfrak{M}_2([t,T])$.

\vspace{2mm}

2. $\xi_{\tau}$
is the mean-square continuous stochastic process at the interval
$[t,T].$

\vspace{2mm}

Let us introduce the following class
of iterated Ito stochastic integrals

\vspace{-1mm}
$$
J[\phi,\psi^{(k)}]_{T,t}=\int\limits_{t}^{T}\psi_1(t_1)\ldots
\int\limits_t^{t_{k-1}}\psi_k(t_k)\int\limits_t^{t_k}
\phi_{\tau}dw^{(k+1)}_{\tau}dw^{(k)}
_{t_k}
\ldots dw^{(1)}_{t_1},
$$

\vspace{2mm}
\noindent
where
$\phi_\tau\in\mathfrak{S}_2([t,T]),$ every
$\psi_l(\tau)$ $(l=1,\ldots,k)$ is a continuous nonrandom function
at the interval $[t, T]$,
here and further 
$w_\tau^{(l)}=f_\tau$ or $w_\tau^{(l)}=\tau$
for $\tau\in[t,T]$
$(l=1,\ldots,k+1),$
$(\psi_1,\ldots,\psi_k)\stackrel{\rm def}{=}\psi^{(k)},$
$\psi^{(1)}\stackrel{\rm def}{=}\psi_1.$

In \cite{str} Stratonovich introduced the definition 
of the so-called 
combined stochastic integral for the specific class of 
integrated processes.
Taking this definition as a foundation, let us consider the 
following construction of stochastic integral 

\vspace{-2mm}
\begin{equation}
\label{1.5000000}
\hbox{\vtop{\offinterlineskip\halign{
\hfil#\hfil\cr
{\rm l.i.m.}\cr
$\stackrel{}{{}_{N\to \infty}}$\cr
}} }
\sum^{N-1}_{j=0} \phi_{\tau_{j}}\left(
f_{\tau_{j+1}} - f_{\tau_{j}}
\right)\theta_{\tau_{j+1}} \stackrel {{\rm def}}{=}
\int\limits_{t}^{T}\phi_{\tau}df_{\tau}\theta_{\tau},
\end{equation}

\vspace{3mm}
\noindent
where $\phi_{\tau},$ $\theta_\tau\in\mathfrak{S}_2([t,T]),$
$\{\tau_j\}_{j=0}^{N}$ is the partition 
of the interval $[t, T],$ which satisfies the condition (\ref{usl})
(for simplicity we write here and sometimes further 
$\tau_j$ instead of $\tau_j^{(N)}$).

Further, we will use integrals of the type 
(\ref{1.5000000}) ($\phi_{\tau}\in\mathfrak{S}_2([t,T])$
and $\theta_\tau$ from a little bit narrower class of stochastic processes
than $\mathfrak{S}_2([t,T])$) 
for formulation 
the theorem on integration order replacement for iterated Ito
stochastic integrals $J[\phi,\psi^{(k)}]_{T,t},$ $k \ge 1.$

Note that under the appropriate conditions the following 
properties of stochastic integrals
defined by the formula (\ref{1.5000000}) can be proved

$$
\int\limits_{t}^{T} \phi_\tau df_\tau g(\tau)=
\int\limits_{t}^{T} \phi_\tau g(\tau) df_\tau\ \ \ \hbox{w.\ p.\ 1}, 
$$

\vspace{3mm}
\noindent
where $g(\tau)$ is a continuous nonrandom function at the
interval $[t,T]$,

\vspace{1mm}
$$
\int\limits_{t}^{T} 
\left(\alpha\phi_{\tau}+\beta\psi_{\tau}\right)df_\tau\theta_\tau=
\alpha\int\limits_{t}^{T} \phi_\tau df_\tau\theta_\tau+
\beta\int\limits_{t}^{T} \psi_{\tau} df_\tau\theta_\tau\ \ \ \hbox{w.\ p.\ 1},
$$

\vspace{2mm}
$$
\int\limits_{t}^{T} \phi_\tau df_\tau
\left(\alpha\theta_\tau+\beta\psi_\tau\right)=
\alpha\int\limits_{t}^{T} \phi_\tau df_\tau\theta_\tau+
\beta\int\limits_{t}^{T} \phi_\tau df_\tau\psi_\tau\ \ \ \hbox{w.\ p.\ 1},
$$

\vspace{4mm}
\noindent
where $\alpha,$ $\beta\in\mathbb{R}^1.$
At that, we suppose that the stochastic processes 
$\phi_\tau,$ $\theta_\tau$, and $\psi_\tau$ are such that
all integrals (included in the mentioned 
properties) exist.

Let us define the stochastic integrals 
$\hat I[\psi^{(k)}]_{T,t},$ $k\ge 1$ of the form

\vspace{1mm}
$$
\hat I[\psi^{(k)}]_{T,t}=\int\limits_t^T\psi_k(t_k)dw_{t_k}^{(k)}
\int\limits_{t_k}^T\psi_{k-1}(t_{k-1})dw_{t_{k-1}}^{(k-1)}
\ldots \int\limits_{t_{2}}^T \psi_1(t_1)dw_{t_1}^{(1)}
$$

\vspace{4mm}
\noindent
in accordance with the definition (\ref{1.5000000}) by the following 
recurrence relation

\vspace{1mm}
\begin{equation}
\label{2.1000000}
\hat I[\psi^{(k)}]_{T,t}
\stackrel{\rm def}{=}
\hbox{\vtop{\offinterlineskip\halign{
\hfil#\hfil\cr
{\rm l.i.m.}\cr
$\stackrel{}{{}_{N\to \infty}}$\cr
}} }
\sum^{N-1}_{l=0} \psi_k(\tau_{l})
\left(w_{\tau_{l+1}}^{(k)}-w_{\tau_l}^{(k)}\right)
\hat I[\psi^{(k-1)}]_{T,\tau_{l+1}},
\end{equation}

\vspace{4mm}
\noindent
where $k\ge 1,$\ \ 
$\hat I[\psi^{(0)}]_{T,s}\stackrel{\rm def}{=}1,$ and
$[s,T] \subseteq [t,T].$

Then, we will define the iterated stochastic 
integral $\hat J[\phi,\psi^{(k)}]_{T,t},$\ $k\ge 1$

\vspace{1mm}
$$
\hat J[\phi,\psi^{(k)}]_{T,t}=\int\limits_{t}^T \phi_s dw_s^{(k+1)}
\hat I[\psi^{(k)}]_{T,s}
$$

\vspace{4mm}
\noindent
similarly in accordance with the definition (\ref{1.5000000})

\vspace{1mm}
$$
\hat J[\phi,\psi^{(k)}]_{T,t}
\stackrel{\rm def}{=}
\hbox{\vtop{\offinterlineskip\halign{
\hfil#\hfil\cr
{\rm l.i.m.}\cr
$\stackrel{}{{}_{N\to \infty}}$\cr
}} }
\sum^{N-1}_{l=0} \phi_{\tau_{l}}
\left(w_{\tau_{l+1}}^{(k+1)}-w_{\tau_l}^{(k+1)}\right)
\hat I[\psi^{(k)}]_{T,\tau_{l+1}}.
$$

\vspace{4mm}

Let us formulate the theorem on
integration order replacement for
iterated Ito stochastic integrals. 

\vspace{2mm}

{\bf Theorem 1}\ \cite{3000}-\cite{3001} (also see \cite{arxiv-22}). 
{\it Suppose that 
$\phi_\tau\in\mathfrak{S}_2([t,T])$ and every $\psi_l(\tau)$
$(l=1,\ldots,k)$ is a continuous nonrandom function
at the interval
$[t,T]$. Then, the stochastic integral
$\hat J[\phi,\psi^{(k)}]_{T,t}$ $(k\ge 1)$ exists and 

\vspace{-2mm}
$$
J[\phi,\psi^{(k)}]_{T,t}=\hat J[\phi,\psi^{(k)}]_{T,t}\ \ \ 
\hbox{w.\ p.\ {\rm 1.}}
$$
}

\vspace{2mm}

Let us consider some propositions related to 
Theorem 1.

\vspace{2mm}

{\bf Proposition 1}\ \cite{3000}-\cite{3001} (also see \cite{arxiv-22}). 
{\it Let the 
conditions of Theorem {\rm 1} are fulfilled and $h(\tau)$  
is a continuous nonrandom function at the interval $[t,T]$.
Then

\begin{equation}
\label{2.19000000}
\int\limits_{t}^T \phi_\tau dw_\tau^{(k+1)}h(\tau) 
\hat I[\psi^{(k)}]_{T,\tau}=
\int\limits_{t}^T 
\phi_\tau h(\tau) dw_\tau^{(k+1)}\hat I[\psi^{(k)}]_{T,\tau}\ \ \
\hbox{{\rm w.\ p.\ 1}},
\end{equation}

\vspace{3mm}
\noindent
and the integrals on the left-hand side of {\rm (\ref{2.19000000})} 
as well as on the right-hand side of {\rm (\ref{2.19000000})} 
exist.} 

\vspace{2mm}

{\bf Proposition 2}\ \cite{3000}-\cite{3001} (also see \cite{arxiv-22}). 
{\it Under the conditions of Theorem {\rm 1}
the following equality is satisfied

\vspace{-2mm}
$$
\int\limits_{t}^T h(t_1)\int\limits_{t}^{t_1}\phi_\tau dw_\tau^{(k+2)}
dw_{t_1}^{(k+1)} \hat I[\psi^{(k)}]_{T,t_1}=
$$

\begin{equation}
\label{2.21000000}
=\int\limits_{t}^T \phi_\tau dw_\tau^{(k+2)}\int\limits_{\tau}^T
h(t_1)dw_{t_1}^{(k+1)}\hat I[\psi^{(k)}]_{T,t_1}\ \ \ 
\hbox{w.\ p.\ {\rm 1}}.
\end{equation}

\vspace{3mm}
\noindent
Moreover, the stochastic integrals in {\rm (\ref{2.21000000})} exist.}

\vspace{2mm}

Using the integration order replacement technique
for iterated Ito stochastic integrals (Theorem 1), we can obtain 
different equalities for iterated Ito stochastic integrals.
At that, the mentioned technique
is essentially simpler (for the specific class
of Ito processes which are the iterated Ito stochastic integrals)
in application than the Ito formula.
Let us consider two examples on application of 
the integration order replacement technique
for iterated Ito stochastic integrals.

\vspace{2mm}

{\bf Example 1.} {\it Using Theorem {\rm 1} and Proposition {\rm 1,} 
we obtain}

\vspace{1mm}
$$
\int\limits_t^T\int\limits_t^{t_3}\int\limits_t^{t_2} df_{t_1}df_{t_2}dt_3=
\int\limits_t^T df_{t_1}\int\limits_{t_1}^T df_{t_2} \int\limits_{t_2}^T dt_3=
$$

\vspace{1mm}
$$
=
\int\limits_t^T df_{t_1}\int\limits_{t_1}^T df_{t_2}(T-t_2)=
\int\limits_t^T df_{t_1}\int\limits_{t_1}^T (T-t_2) df_{t_2}=
$$

\vspace{1mm}
$$
=
\int\limits_t^T(T-t_2)\int\limits_t^{t_2}df_{t_1}df_{t_2}\ \ \ 
\hbox{w.\ p.\ 1.}
$$
 
\vspace{5mm}

{\bf Example 2.} {\it Using Theorem {\rm 1} and Proposition {\rm 1,} 
we obtain}

\vspace{1mm}
$$
\int\limits_t^T\int\limits_t^{t_4}
\int\limits_t^{t_3}\int\limits_t^{t_2}df_{t_1}dt_2df_{t_3}dt_4=
\int\limits_t^T df_{t_1}\int\limits_{t_1}^T dt_2 \int\limits_{t_2}^T df_{t_3}
\int\limits_{t_3}^T dt_4
=
$$

\vspace{1mm}
$$
=\int\limits_t^T df_{t_1}\int\limits_{t_1}^T 
dt_2 \int\limits_{t_2}^T df_{t_3}(T-t_3)=
\int\limits_t^T df_{t_1}
\int\limits_{t_1}^T dt_2 \int\limits_{t_2}^T(T-t_3)df_{t_3}=
$$

\vspace{1mm}
$$
=
\int\limits_t^T(T-t_3)
\int\limits_t^{t_3}\int\limits_t^{t_2} df_{t_1}dt_2df_{t_3}=
\int\limits_t^T(T-t_3)\left(\int\limits_t^{t_3}
\int\limits_t^{t_2} df_{t_1}dt_2\right)df_{t_3}=
$$

\vspace{1mm}
$$
=
\int\limits_t^T(T-t_3)\left(
\int\limits_t^{t_3}df_{t_1}\int\limits_{t_1}^{t_3}dt_2\right)df_{t_3}=
\int\limits_t^T(T-t_3)\left(
\int\limits_t^{t_3}df_{t_1}(t_3-t_1)\right)df_{t_3}=
$$

\vspace{1mm}
$$
=
\int\limits_t^T(T-t_3)\left(
\int\limits_t^{t_3}(t_3-t_1)df_{t_1}\right)df_{t_3}
=\int\limits_t^T(T-t_2)\int\limits_t^{t_2}(t_2-t_1)
df_{t_1}df_{t_2}\ \ \ \hbox{w.\ p.\ 1.}
$$

\vspace{5mm}

Let us apply Theorem 1 to deriving of one propetry for Ito
stochastic integrals.

\vspace{2mm}

{\bf Lemma 1.}\ {\it Let $h(\tau), 
g(\tau), G(\tau): [t,s]\to\mathbb{R}^1$ 
be continuous nonrandom functions at the interval $[t, s]$ and let $G(\tau)$         
be a antiderivative of the function
$g(\tau).$ Furthermore, let $\xi_{\tau}\in \mathfrak{S}_2([t,s]).$ 
Then

\vspace{1mm}
\begin{equation}
\label{a1}
\int\limits_t^s g(\tau)\int\limits_t^{\tau}h(\theta)
\int\limits_t^{\theta}\xi_u d{\bf f}_u^{(i)}
d{\bf f}_{\theta}^{(j)}
d\tau=\int\limits_t^s(G(s)-G(\theta))h(\theta)
\int\limits_t^{\theta}\xi_u d{\bf f}_u^{(i)}d{\bf f}_{\theta}^{(j)}
\end{equation}

\vspace{4mm}
\noindent
w.\ p.\ {\rm 1}, 
where $i, j =1, 2$ and ${\bf f}_{\tau}^{(1)}$, ${\bf f}_{\tau}^{(2)}$
are independent standard Wiener processes that are
${\rm F}_{\tau}$--measurable for all $\tau\in[t,s]$.}

\vspace{2mm}

{\bf Proof.} Applying Theorem 1 twice and Proposition 1, 
we get the following relations

$$
\int\limits_t^s g(\tau)\int\limits_t^{\tau}h(\theta)
\int\limits_t^{\theta}\xi_u d{\bf f}_u^{(i)}
d{\bf f}_{\theta}^{(j)}
d\tau=\int\limits_t^s \xi_u d{\bf f}_u^{(i)}
\int\limits_u^s h(\theta)d{\bf f}_{\theta}^{(j)}\int\limits_{\theta}^s
g(\tau)d\tau=
$$

\vspace{1mm}

$$
=G(s)\int\limits_t^s \xi_u d{\bf f}_u^{(i)}
\int\limits_u^s h(\theta)d{\bf f}_{\theta}^{(j)}-
\int\limits_t^s \xi_u d{\bf f}_u^{(i)}
\int\limits_u^s G(\theta)h(\theta)d{\bf f}_{\theta}^{(j)}=
$$

\vspace{1mm}
$$
=G(s)\int\limits_t^{s}h(\theta)
\int\limits_t^{\theta}\xi_u d{\bf f}_u^{(i)}d{\bf f}_{\theta}^{(j)}
-
\int\limits_t^{s}G(\theta)h(\theta)
\int\limits_t^{\theta}\xi_u d{\bf f}_u^{(i)}d{\bf f}_{\theta}^{(j)}=
$$

\begin{equation}
\label{a3}
=\int\limits_t^s(G(s)-G(\theta))h(\theta)
\int\limits_t^{\theta}\xi_u d{\bf f}_u^{(i)}d{\bf f}_{\theta}^{(j)}\ \ \ 
\hbox{w.\ p.\ 1.}
\end{equation}

\vspace{5mm}

The proof of Lemma 1 is completed.
Let us consider an analogue of Lemma 1 for
Stratonovich stochastic integrals.

\vspace{2mm}

{\bf Lemma 2}\ \cite{k1}. {\it Let $h(\tau), 
g(\tau), G(\tau): [t,s]\to\mathbb{R}^1$ 
be continuous nonrandom functions at the interval $[t, s]$ and let $G(\tau)$         
be a antiderivative of the function
$g(\tau).$ Let $\xi_{\tau}^{(l)}\in \mathfrak{Q}_4([t,s])$ 
and

$$
\xi_{\tau}^{(l)}=\int\limits_t^{\tau}a_u du+
\int\limits_t^{\tau}b_u d{\bf f}_u^{(l)},\ \ \ l=1, 2.
$$

\vspace{2mm}

Then
\begin{equation}
\label{a1xxx}
\int\limits_t^s g(\tau){\int\limits_t^{*}}^{\tau}h(\theta)
{\int\limits_t^{*}}^{\theta}\xi^{(l)}_u d{\bf f}_u^{(i)}
d{\bf f}_{\theta}^{(j)}
d\tau={\int\limits_t^{*}}^s(G(s)-G(\theta))h(\theta)
{\int\limits_t^{*}}^{\theta}\xi_u^{(l)} d{\bf f}_u^{(i)}d{\bf f}_{\theta}^{(j)}
\end{equation}

\vspace{2mm}
\noindent
w.\ p.\ {\rm 1}, 
where $i, j, l=1, 2$ and ${\bf f}_{\tau}^{(1)}$, ${\bf f}_{\tau}^{(2)}$
are independent standard Wiener processes that are
${\rm F}_{\tau}$--measurable for all $\tau\in[t,s]$.}

\vspace{2mm}

{\bf Proof.} Under the conditions of Lemma 2, we can 
apply equalities (\ref{d11}) and (\ref{d11a}) with
$F(x,\theta)\equiv x h(\theta),$
$$
\eta_{\theta}={\int\limits_t^{*}}^{\theta}\xi_u^{(l)} d{\bf f}_u^{(i)},
$$

\vspace{2mm}
\noindent
since the function $x h(\theta)$ is sufficiently smooth and the following obvious inclusions hold:
$\eta_{\theta}\in \mathfrak{Q}_4([t,s])$
and
$\eta_{\theta}h(\theta)\in 
\mathfrak{M}_2([t,s]).$ Thus, we have the equalities

\vspace{-1mm}
\begin{equation}
\label{a2}
{\int\limits_t^{*}}^{\tau}h(\theta)
{\int\limits_t^{*}}^{\theta}\xi^{(l)}_u d{\bf f}_u^{(i)}d{\bf f}_{\theta}^{(j)}
=\int\limits_t^{\tau}h(\theta)
{\int\limits_t^{*}}^{\theta}\xi_u^{(l)} 
d{\bf f}_u^{(i)}d{\bf f}_{\theta}^{(j)}+
\frac{1}{2}{\bf 1}_{\{i=j\}}\int\limits_t^{\tau}h(\theta)\xi_{\theta}^{(l)}
d\theta,
\end{equation}

\begin{equation}
\label{a2.hoho}
{\int\limits_t^{*}}^{\theta}\xi^{(l)}_u d{\bf f}_u^{(i)}
=\int\limits_t^{\theta}\xi^{(l)}_u d{\bf f}_u^{(i)}+
\frac{1}{2}{\bf 1}_{\{l=i\}}\int\limits_t^{\theta}b_u du
\end{equation}

\vspace{3mm}
\noindent
w.\ p.\ {\rm 1}, where  ${\bf 1}_A$ 
is the indicator of a set $A$. 
Substituting formulas (\ref{a2}) and (\ref{a2.hoho}) into the left-hand
side of equality (\ref{a1xxx}) and applying Theorem 1 twice
and Proposition 1, 
we get the following relations

$$
\int\limits_t^s g(\tau){\int\limits_t^{*}}^{\tau}h(\theta)
{\int\limits_t^{*}}^{\theta}\xi_u^{(l)} d{\bf f}_u^{(i)}d{\bf f}_{\theta}^{(j)}
d\tau=
$$

\vspace{1mm}
$$
=\int\limits_t^s \xi_u^{(l)} d{\bf f}_u^{(i)}
\int\limits_u^s h(\theta)d{\bf f}_{\theta}^{(j)}\int\limits_{\theta}^s
g(\tau)d\tau+
$$

\vspace{1mm}
$$
+\frac{1}{2}{\bf 1}_{\{l=i\}}\int\limits_t^s b_u du
\int\limits_u^s h(\theta)d{\bf f}_{\theta}^{(j)}
\int\limits_{\theta}^s g(\tau)d\tau
+\frac{1}{2}{\bf 1}_{\{i=j\}}\int\limits_t^s h(\theta)\xi_{\theta}^{(l)}
d\theta \int\limits_{\theta}^s g(\tau)d\tau=
$$

\vspace{1mm}
$$
=G(s)\left(\int\limits_t^s \xi_u^{(l)} d{\bf f}_u^{(i)}
\int\limits_u^s h(\theta)d{\bf f}_{\theta}^{(j)}+
\frac{1}{2}{\bf 1}_{\{i=j\}}\int\limits_t^{s}h(\theta)\xi_{\theta}^{(l)}
d\theta
+\right.
$$

\vspace{1mm}
$$
\left.+\frac{1}{2}{\bf 1}_{\{l=i\}}\int\limits_t^s b_u du
\int\limits_u^s h(\theta)d{\bf f}_{\theta}^{(j)}
\right)-
$$

\vspace{1mm}
$$
-\left(\int\limits_t^s \xi_u^{(l)} d{\bf f}_u^{(i)}
\int\limits_u^s G(\theta)h(\theta)d{\bf f}_{\theta}^{(j)}+
\frac{1}{2}{\bf 1}_{\{i=j\}}\int\limits_t^{s}
G(\theta)h(\theta)\xi_{\theta}^{(l)}d\theta
+\right.
$$

\vspace{1mm}
$$
\left.+\frac{1}{2}{\bf 1}_{\{l=i\}}\int\limits_t^s b_u du
\int\limits_u^s h(\theta)G(\theta)d{\bf f}_{\theta}^{(j)}
\right)=
$$

\vspace{1mm}
$$
=G(s)\left(\int\limits_t^{s}h(\theta)
\int\limits_t^{\theta}\xi_u^{(l)} d{\bf f}_u^{(i)}d{\bf f}_{\theta}^{(j)}+
\frac{1}{2}{\bf 1}_{\{i=j\}}\int\limits_t^{s}h(\theta)\xi_{\theta}^{(l)}
d\theta+\right.
$$

\vspace{1mm}
$$
\left.+\frac{1}{2}{\bf 1}_{\{l=i\}}\int\limits_t^s
h(\theta)\int\limits_t^{\theta}b_u du d{\bf f}_{\theta}^{(j)}
\right)-
$$

\vspace{1mm}
$$
-\left(\int\limits_t^{s}G(\theta)h(\theta)
\int\limits_t^{\theta}\xi_u^{(l)} d{\bf f}_u^{(i)}d{\bf f}_{\theta}^{(j)}+
\frac{1}{2}{\bf 1}_{\{i=j\}}\int\limits_t^{s}
G(\theta)h(\theta)\xi_{\theta}^{(l)}d\theta
+\right.
$$

\vspace{1mm}
\begin{equation}
\label{a3xxx}
\left.+\frac{1}{2}{\bf 1}_{\{l=i\}}\int\limits_t^s
h(\theta)G(\theta)\int\limits_t^{\theta}b_u du d{\bf f}_{\theta}^{(j)}
\right)
\end{equation}

\vspace{5mm}
\noindent
w.\ p.\ 1.
Applying successively the formulas (\ref{a2}), (\ref{a2.hoho}) 
together with the 
formula (\ref{a2}) in which $h(\theta)$ replaced 
by $G(\theta)h(\theta)$ as well as
the relation (\ref{a3xxx}), we obtain the equality (\ref{a1xxx}).

\vspace{5mm}

\section{The Taylor--Ito Expansion}

\vspace{5mm}

In this section, we cite the Taylor-Ito expansion \cite{KlPl1} 
and introduce some necessary notations.
At that, we will use the original notation introduced by the author
of this paper.

Let $\mathfrak{L}$ be the class of functions $R({\bf x}, t): 
\mathbb{R}^n\times [0, T] \to \mathbb{R}^1$ with the 
following property: these functions are twice
continuously differentiable in ${\bf x}$ and have one continuous 
derivative in $t$. We consider the following operators on
the space $\mathfrak{L}$

\begin{equation}
\label{2.3}
L= \frac{\partial}{\partial t}
+ \sum^ {n} _ {i=1} a^{(i)} ({\bf x},  t) 
\frac{\partial}{\partial  {\bf  x}^{(i)}} +
\frac{1}{2}\sum^ {m} _ {j=1} \sum^ {n} _ {l,i=1}
B^ {(lj)} ({\bf x}, t) B^ {(ij)} ({\bf x}, t) 
\frac{\partial^{2}}{\partial{\bf x}^{(l)}\partial{\bf x}^{(i)}},
\end{equation}

\vspace{1mm}

\begin {equation}
\label{2.4}
G^ {(i)} _ {0} = \sum^ {n} _ {j=1} B^ {(ji)} ({\bf x}, t)
\frac{\partial}{\partial {\bf x} ^ {(j)}},\ \ \
i=1,\ldots,m.
\end {equation}

\vspace{4mm}

By the Ito formula, we have the equality

\begin{equation}
\label{sa}
R({\bf x}_s,s)= R({\bf x}_t,t) + 
\int\limits_t^s LR({\bf x}_\tau,\tau)d\tau
+\sum^m_{i=1}\int\limits_t^s G_0^{(i)}R({\bf x}_\tau,\tau)d{\bf f}_\tau^{(i)}
\end{equation}

\vspace{3mm}
\noindent
w.~p.~1, where  $0 \le t<s \le T.$ In the formula (\ref{sa})
it is assumed that the functions ${\bf a}({\bf x}, t)$, 
$B({\bf x}, t)$, and
$R({\bf x}, t)$ satisfy the following condition: 
$LR({\bf x}_{\tau},\tau)$, $G_0^{(i)}R({\bf x}_{\tau},\tau)
\in \mathfrak{M}_2([0, T])$
for $i = 1,\ldots,m$.

\vspace{2mm}

Introduce the following notation

\begin{equation}
\label{999.003}
{}^{(k)}A =\Biggl\| A^{(i_{1}\ldots  i_{k})}
\Biggr\|^{m_{1}~\ldots~~ m_{k}}_{i_{1}=1,\ldots ,i_{k}=1},\ \ \ 
m_1,\ldots,m_k\ge 1,
\end{equation}

\vspace{5mm}

$$
{ }^{(k+l)}A \stackrel{l}{\cdot}{}^{(l)}B^{(k)}
=\begin{cases} 
\Biggl\|\sum\limits^{m_{1}}_{i_1=1}\ldots\sum\limits_{i_l=1}^{m_l}
A^{(i_{1}\ldots i_{k+l})}B^{(i_{1}
\ldots  i_{l})}\Biggr\|^{m_{l+1}~\ldots ~~m_{l+k}}_{i_{l+1}=1,
\ldots ,i_{l+k}=1}\ \ \ &\hbox{for}\ \ \ k\ge 1\cr\cr
\sum\limits^{m_{1}}_{i_1=1}\ldots\sum\limits_{i_l=1}^{m_l}
A^{(i_{1}\ldots i_{l})}B^{(i_{1}
\ldots  i_{l})}\ \ \ &\hbox{for}\ \ \ k= 0
\end{cases},
$$

\vspace{5mm}

\begin{equation}
\label{sa901}
\Biggl\| A_{k+1}D_{k}^{(i_k)}A_k\ldots A_{2}D_{1}^{(i_1)}A_1
R({\bf x},t)
\Biggr\|^{m_{1}~\ldots~~ m_{k}}_{i_{1}=1,\ldots ,i_{k}=1}=
{}^{(k)}A_{k+1}D_{k}A_k\ldots A_{2}D_{1}A_1
R({\bf x},t),
\end{equation}

\vspace{6mm}
\noindent
where $A_p$ and $D_q^{(i_q)}$ are operators defined on the 
space $\mathfrak{L}$ for $p = 1,\ldots,k+1$, 
$q = 1,\ldots,k$, and $i_q = 1,\ldots,m_q$. It
is assumed that the left-hand side of (\ref{sa901}) exists. 
The symbol $\stackrel{0}{\cdot}$
is treated as the usual multiplication. If $m_l = 0$
in (\ref{999.003}) for some $l\in\{1,\ldots,k\}$, 
then the right-hand side of (\ref{999.003}) 
is treated as

$$
\Biggl\| A^{(i_{1}\ldots i_{l-1}i_{l+1}\ldots i_{k})}
\Biggr\|^{m_{1}~\ldots~~ m_{l-1}~~m_{l+1}~\ldots~~ m_{k}}
_{i_{1}=1,\ldots, i_{l-1}=1, i_{l+1}=1,\ldots,i_{k}=1},
$$

\vspace{3mm}
\noindent
(shortly, ${}^{(k-1)}A$).

We also introduce the following notation

$$
\Biggl\Vert Q_{\lambda_l}^{(i_l)}\ldots Q_{\lambda_1}^{(i_1)}
R({\bf x},t)\Biggr\Vert_{i_1=\lambda_1,\ldots,i_l=\lambda_l}^{m\lambda_1
~\ldots~m\lambda_l}\stackrel{\sf def}{=}{}^{(p_l)}Q_{\lambda_l}
\ldots Q_{\lambda_1}R({\bf x},t),
$$

\vspace{5mm}

$$
{}^{(p_{k})}J_{(\lambda_{k}\ldots \lambda_1)s,t}
=\Biggl\Vert J_{(\lambda_{k}\ldots \lambda_1)s,t}^{(i_k\ldots
i_1)}\Biggr\Vert_
{i_1=\lambda_1,\ldots,i_k=\lambda_k}^{m\lambda_1~\ldots~m\lambda_k},
$$

\vspace{4mm}

$$
{M}_k=\biggl\{(\lambda_k,\ldots,\lambda_1):
\lambda_l=1\ \hbox{or}\ \lambda_l=0;\ l=1,\ldots,k\biggr\},\ \ \ k\ge 1,
$$

\vspace{4mm}

$$
J_{(\lambda_{k}\ldots \lambda_1)s,t}^{(i_k\ldots
i_1)}=
\int\limits_t^s\ldots
\int\limits_t^{t_{2}}
d{\bf w}_{t_{1}}^{(i_k)}\ldots
d{\bf w}_{t_k}^{(i_1)},\ \ \ k\ge 1,
$$

\vspace{5mm}
\noindent
where $\lambda_l=1$ or $\lambda_l=0$,
$Q_{\lambda_l}^{(i_l)}={L}$
and $i_l=0$ for $\lambda_l=0,$  $Q_{\lambda_l}^{(i_l)}=G_0^{(i_l)}$
and $i_l=1,\ldots,m$ for $\lambda_l=1,$ 
$$
p_l=\sum\limits_{j=1}^l \lambda_j\ \ \ \hbox{for}\ \ \ l=1,\ldots, r+1,\ \ \
r\in\mathbb{N},
$$
${\bf w}_{\tau}^{(i)}$ $(i=1,\ldots,m)$ are ${\rm F}_{\tau}$-measurable 
for all $\tau\in [0, T]$
independent standard Wiener processes and
${\bf w}_{\tau}^{(0)}=\tau.$

Applying the formula (\ref{sa}) to the process 
$R({\bf x}_s, s)$ repeatedly, we obtain the following 
Taylor--Ito expansion
\cite{KlPl1}

\vspace{-1mm}
\begin{equation}
\label{5.7.11}
R({\bf x}_s, s)=R({\bf x}_t, t)+\sum_{k=1}^r \sum_{(\lambda_{k},\ldots,
\lambda_1)\in {M}_k}\
{}^{(p_{k})}Q_{\lambda_{k}}\ldots Q_{\lambda_1}
R({\bf x}_t, t)\stackrel{p_k}{\cdot}
{}^{(p_{k})}J_{(\lambda_{k}\ldots \lambda_1)s,t}
+\left(D_{r+1}\right)_{s,t}
\end{equation}

\vspace{4mm}
\noindent
w.\ p.\ 1, where

\vspace{-1mm}
\begin{equation}
\label{5.7.12}
\left(D_{r+1}\right)_{s,t}
=\sum_{(\lambda_{r+1},\ldots,\lambda_1)\in{M}_{r+1}}\
\int\limits_t^s\ldots
\left(
\int\limits_t^{t_2}
{}^{(p_{r+1})}Q_{\lambda_{r+1}}\ldots
Q_{\lambda_1}
R({\bf x}_{t_1}, t_1)\stackrel{\lambda_{r+1}}{\cdot}
d{\bf w}_{t_{1}}\right)
\ldots
\stackrel{\lambda_1}{\cdot}d{\bf w}_{t_{r+1}}.
\end{equation}

\vspace{6mm}
\noindent
It is assumed that the right-hand sides of (\ref{5.7.11}), 
(\ref{5.7.12}) exist.

A possible variant of the conditions under which the right-hand 
sides of (\ref{5.7.11}), (\ref{5.7.12}) exist is as follows

\vspace{2mm}

(i)\ $Q_{\lambda_l}^{(i_l)}\ldots Q_{\lambda_1}^{(i_1)}R({\bf x},t)
\in \mathfrak{L}$ for all 
$(\lambda_l,\ldots,\lambda_1)\in\bigcup\limits_{g=1}^{r} M_g$;

\vspace{1mm}

(ii)\ $Q_{\lambda_l}^{(i_l)}\ldots Q_{\lambda_1}^{(i_1)}R({\bf x}_{\tau},\tau)
\in \mathfrak{M}_2([0,T])$
for all 
$(\lambda_l,\ldots,\lambda_1)\in\bigcup\limits_{g=1}^{r+1}
{M}_g.$ 

\vspace{3mm}

Let us write 
the expansion (\ref{5.7.11}) in the another form

$$
R({\bf x}_s,s)=R({\bf x}_t,t)
+\sum_{k=1}^r \sum_{(\lambda_{k},\ldots,\lambda_1)
\in{M}_k}\
\sum_{i_1=\lambda_1}^{m\lambda_1}
\ldots 
\sum_{i_k=\lambda_k}^{m\lambda_k}
Q_{\lambda_k}^{(i_k)}\ldots Q_{\lambda_1}^{(i_1)}
R({\bf x}_t,t)\
{J}_{(\lambda_{k}\ldots \lambda_1)s,t}^{(i_k\ldots
i_1)}
+
$$

\vspace{1mm}
$$
+\left(D_{r+1}\right)_{s,t}\ \ \ \hbox{w.\ p.\ 1}.
$$

\vspace{5mm}

Denote
$$
{G}_{rk}=\biggl\{(\lambda_k,\ldots,\lambda_1):\ r+1\le 
2k-\lambda_1-\ldots-\lambda_k\le 2r\biggr\},
$$

\vspace{1mm}

$$
{E}_{qk}=\biggl\{(\lambda_k,\ldots,\lambda_1):\ 2k-\lambda_1-\ldots-
\lambda_k=q\biggr\},
$$

\vspace{5mm}
\noindent
where $\lambda_l=1$ or $\lambda_l=0$ $(l=1,\ldots,k).$

The Taylor--Ito expansion ordered according to the 
order of smallness (in the mean-square sense when $s\downarrow  t$) 
of its terms has the form

$$
R({\bf x}_s,s)=R({\bf x}_t,t)+
\sum_{q,k=1}^r \sum_{(\lambda_{k},\ldots,\lambda_1)
\in{\rm E}_{qk}}\
\sum_{i_1=\lambda_1}^{m\lambda_1}
\ldots 
\sum_{i_k=\lambda_k}^{m\lambda_k}
Q_{\lambda_k}^{(i_k)}\ldots Q_{\lambda_1}^{(i_1)}
R({\bf x}_t,t)\
{J}_{(\lambda_{k}\ldots \lambda_1)s,t}^{(i_k\ldots
i_1)}+
$$

\begin{equation}
\label{5.6.1rrr}
+\left(H_{r+1}\right)_{s,t}\ \ \ \hbox{w.\ p.\ 1},
\end{equation}

\vspace{4mm}
\noindent
where

\vspace{-2mm}
$$
\left(H_{r+1}\right)_{s,t}=
\sum_{k=1}^r \sum_{(\lambda_{k},\ldots,\lambda_1)
\in{\rm G}_{rk}}\
\sum_{i_1=\lambda_1}^{m\lambda_1}
\ldots 
\sum_{i_k=\lambda_k}^{m\lambda_k}
Q_{\lambda_k}^{(i_k)}\ldots Q_{\lambda_1}^{(i_1)}
R({\bf x}_t,t)\
{J}_{(\lambda_{k}\ldots \lambda_1)s,t}^{(i_k\ldots
i_1)}+
$$

\vspace{1mm}
\begin{equation}
\label{5.6.1rrrh}
+\left(D_{r+1}\right)_{s,t}.
\end{equation}

\vspace{5mm}

\section{The First Form of the Unified Taylor--Ito Expansion}

\vspace{5mm}

In this section, we transform the right-hand side of (\ref{5.7.11}) with 
the help of Theorem 1 and Lemma 1 to a
representation including iterated 
Ito stochastic integrals of the form (\ref{ll11}).

Denote

\vspace{-2mm}
\begin{equation}
\label{opr1}
I^{(i_1\ldots i_k)} _ {{l_1 \ldots l_k}_{s, t}} 
=\int\limits_t^{s}(t-t_k)^{l_{k}}\ldots\
\int\limits_t^{t _ {2}}(t-t _ {1}) ^ {l_ {1}} 
d{\bf f} ^ {(i_ {1})} _ {t_ {1}} \ldots
d{\bf f}_{t_ {k}}^{(i_ {k})}\ \ \ \hbox{for}\ \ \ k\ge 1
\end{equation}

\vspace{2mm}
\noindent
and
$$
I^{(i_1\ldots i_k)} _ {{l_1 \ldots l_k}_{s, t}}=1\ \ \ 
\hbox{for}\ \ \ k=0,
$$

\vspace{5mm}
\noindent
where $i_1,\ldots,i_k=1,\ldots,m.$ Moreover, let

\vspace{1mm}
$$
{}^{(k)}I_{{l_1\ldots l_k}_{s,t}}=\Biggl\Vert 
I^{(i_1\ldots i_k)} _ {{l_1 \ldots l_k}_{s, t}}\Biggr\Vert
_{i_1,\ldots,i_k=1}^{m},
$$

\vspace{4mm}

\begin{equation}
\label{a9}
G_p^{(i)}\stackrel{\sf def}{=}\frac{1}{p}\left(
G_{p-1}^{(i)}L-LG_{p-1}^{(i)}\right),\ \ \
p=1, 2,\ldots,\ \ \ i=1,\ldots,m,
\end{equation}

\vspace{4mm}
\noindent
where $L$ and $G_0^{(i)},$ $i=1,\ldots,m,$
are determined by the equalities
(\ref{2.3}), (\ref{2.4}). Denote

\vspace{3mm}

$$
{A}_q\stackrel{\sf def}{=}
\Biggl\{
(k,j,l_1,\ldots,l_k):\ k+j+\sum_{p=1}^k l_p=q;\ k,j,l_1,\ldots,l_k=0, 1,\ldots
\Biggr\},
$$

\vspace{4mm}
$$
\Biggl\Vert 
G_{l_1}^{(i_1)}\ldots  G_{l_k}^{(i_k)} L^j
R({\bf x},t) \Biggr\Vert
_{i_1,\ldots,i_k=1}^
{m}\stackrel{\sf def}{=}{}^{(k)}
G_{l_1}\ldots G_{l_k} L^j
R({\bf x},t),
$$

\vspace{4mm}

$$
L^j R({\bf x},t)\stackrel{\sf def}{=}
\begin{cases}\underbrace{L\ldots L}_j
R({\bf x},t)\ &\hbox{for}\ j\ge 1\cr\cr
R({\bf x},t)\ &\hbox{for}\ j=0
\end{cases}.
$$

\vspace{6mm}

{\bf Theorem 2.}\ {\it Let conditions {\rm (i), (ii)}
be satisfied. Then for any $s, t \in [0, T]$ such that $s>t$ 
and for any positive
integer $r$, the following expansion takes place w.\ p.\ {\rm 1}

$$
R({\bf x}_s,s)=
R({\bf x}_t,t)+
\sum_{q=1}^r
\sum_{(k,j,l_1,\ldots,l_k) \in {\rm A}_q}\ 
\frac{(s-t)^j}{j!}
\sum_{i_1,\ldots,i_k=1}^m G_{l_1}^{(i_1)}\ldots
G_{l_k}^{(i_k)}L^j R({\bf x}_t,t)\ 
I_{{l_1\ldots l_k}_{s,t}}^{(i_1\ldots i_k)}+
$$

\vspace{2mm}
\begin{equation}
\label{razl4}
+\left(D_{r+1}\right)_{s,t},
\end{equation}

\vspace{5mm}
\noindent
where $\left(D_{r+1}\right)_{s,t}$ has the form {\rm (\ref{5.7.12})}.
}

\vspace{2mm}

{\bf Proof.} We claim that

\vspace{2mm}
$$
\sum_{(\lambda_{q},\ldots,
\lambda_1)\in {M}_q}\
{}^{(p_{q})}Q_{\lambda_{q}}\ldots Q_{\lambda_1}
R({\bf x}_t,t)\stackrel{p_q}{\cdot}
{}^{(p_{q})}J_{(\lambda_{q}\ldots \lambda_1)s,t}=
$$

\vspace{2mm}
\begin{equation}
\label{a22}
=
\sum_{(k,j,l_1,\ldots,l_k) \in {\rm A}_q}\ 
\frac{(s-t)^j}{j!}
\sum_{i_1,\ldots,i_k=1}^m G_{l_1}^{(i_1)}\ldots
G_{l_k}^{(i_k)}L^j R({\bf x}_t,t)\ 
I_{{l_1\ldots l_k}_{s,t}}^{(i_1\ldots i_k)}
\end{equation}

\vspace{7mm}
\noindent
w.\ p.\ 1. The equality (\ref{a22}) is valid for $q = 1$. Assume 
that (\ref{a22}) is valid for some $q > 1$. 
In this case, using the induction hypothesis, we obtain

\vspace{2mm}
$$
\sum_{(\lambda_{q+1},\ldots,
\lambda_1)\in {M}_{q+1}}\
{}^{(p_{q+1})}Q_{\lambda_{1}}\ldots Q_{\lambda_{q+1}}
R({\bf x}_t,t)\stackrel{p_{q+1}}{\cdot}
{}^{(p_{q+1})}J_{(\lambda_{1}\ldots \lambda_{q+1})s,t}=
$$

\vspace{3mm}
$$
=\sum_{\lambda_{q+1}\in\{1,\ 0\}}
\int\limits_t^s
\sum_{(\lambda_{q},\ldots,
\lambda_1)\in {M}_{q}}
\Biggl({}^{(p_{q+1})}Q_{\lambda_{1}}
\ldots
Q_{\lambda_{q+1}}
R({\bf x}_t,t)\stackrel{p_{q}}{\cdot}
{}^{(p_{q})}J_{(\lambda_{1}\ldots \lambda_{q})\theta,t}\Biggr)
\stackrel{\lambda_{q+1}}{\cdot}d{\bf w}_{\theta}=
$$

\vspace{3mm}
$$
=\sum_{\lambda_{q+1}\in\{1,\ 0\}}
\int\limits_t^s
\sum_{(k,j,l_1,\ldots,l_k)\in{A}_q}
\frac{(\theta-t)^j}{j!}\times
$$

\vspace{3mm}
$$
\times
\Biggl({}^{(k+\lambda_{q+1})}G_{l_1}\ldots G_{l_k} L^j
Q_{\lambda_{q+1}}R({\bf x}_t,t)
\stackrel{k}{\cdot}
{}^{(k)}I_{{l_1\ldots l_k}_{s,t}}\Biggr)
\stackrel{\lambda_{q+1}}{\cdot}d{\bf w}_{\theta}=
$$

\vspace{3mm}
$$
=\sum_{(k,j,l_1,\ldots,l_k)\in{A}_q}\left(
{}^{(k)}G_{l_1}\ldots G_{l_k} L^{j+1}
R({\bf x}_t,t)
\stackrel{k}{\cdot}
\int\limits_t^s\frac{(\theta-t)^j}{j!}
{}^{(k)}I_{{l_1\ldots l_k}_{\theta,t}}d\theta+\right.
$$

\vspace{3mm}
\begin{equation}
\label{a30}
\left.+\left(
{}^{(k+1)}G_{l_1}\ldots G_{l_k} L^{j} G_0
R({\bf x}_t,t)
\stackrel{k}{\cdot}
\int\limits_t^s\frac{(\theta-t)^j}{j!}
{}^{(k)}I_{{l_1\ldots l_k}_{\theta,t}}\right)
\stackrel{1}{\cdot}d{\bf f}_{\theta}\right)
\end{equation}

\vspace{6mm}
\noindent
w.\ p.\ 1.

Using Lemma 1, we obtain

\vspace{-2mm}
$$
\int\limits_t^s\frac{(\theta-t)^j}{j!}
{}^{(k)}I_{{l_1\ldots l_k}_{\theta,t}}d\theta=
$$
\begin{equation}
\label{a31}
=\frac{1}{(j+1)!}
\begin{cases}(s-t)^{j+1}\ &\hbox{for}\ k=0 \cr \cr 
(s-t)^{j+1} \cdot {}^{(k)}I_{{l_1\ldots l_k}_{s,t}}-
(-1)^{j+1} \cdot
{}^{(k)}I_{{l_1\ldots l_{k-1}\ l_k+j+1}_{s,t}}\ &\hbox{for}\ k>0
\end{cases}
\end{equation}

\vspace{6mm}
\noindent
w.\ p.\ 1. In addition (see (\ref{opr1})), we get

\vspace{1mm}
\begin{equation}
\label{a32}
\int\limits_t^s\frac{(\theta-t)^j}{j!}
I^{(i_1\ldots i_k)}_{{l_1\ldots l_k}_{\theta,t}}
d{\bf f}_{\theta}^{(i_{k+1})}=
\frac{(-1)^j}{j!} I_{{l_1\ldots l_k j}_{s,t}}^{(i_1\ldots i_k i_{k+1})}
\end{equation}

\vspace{4mm}
\noindent
in the notations just introduced.
Substitute the relations (\ref{a31}) and (\ref{a32}) 
into the formula (\ref{a30}). Grouping summands of the 
obtained expression with
equal lower indices at iterated Ito stochastic
integrals and using (\ref{a9}) and the equality

\vspace{1mm}
\begin{equation}
\label{a33}
G_p^{(i)}R({\bf x},t)=\frac{1}{p!}
\sum_{q=0}^p(-1)^q C_p^q L^q G_0^{(i)} L^{p-q}
R({\bf x},t),\ \ \ \hbox{where}\ \ \ C_p^q=\frac{p!}{q!(p-q)!},
\end{equation}

\vspace{5mm}
\noindent
(this equality follows from (\ref{a9})), we note that the obtained 
expression is equal to

\vspace{1mm}
$$
\sum_{(k,j,l_1,\ldots,l_k)\in{A}_{q+1}}
\frac{(s-t)^j}{j!}
{}^{(k)} G_{l_1}\ldots G_{l_k} L^j\{\eta_t\} 
\stackrel{k}{\cdot}
{}^{(k)}I_{{l_1\ldots l_k}_{s,t}}
$$

\vspace{5mm}
\noindent
w.\ p.\ 1. Summing the equalities (\ref{a22}) for $q = 1, 2,\ldots,r$ 
and applying the formula (\ref{5.7.11}), we obtain the expression
(\ref{razl4}). The proof is completed.

Let us order terms of the expansion (\ref{razl4}) according to 
their smallness orders as $s \downarrow t$ in the mean-square sense

\vspace{1mm}
$$
R({\bf x}_s,s)=
R({\bf x}_t,t)+
\sum_{q=1}^r
\sum_{(k,j,l_1,\ldots,l_k) \in {\rm D}_q}\ 
\frac{(s-t)^j}{j!}
\sum_{i_1,\ldots,i_k=1}^m  G_{l_1}^{(i_1)}\ldots
G_{l_k}^{(i_k)} L^j R({\bf x}_t,t)\
I_{{l_1\ldots l_k}_{s,t}}^{(i_1\ldots i_k)}+
$$

\vspace{2mm}
\begin{equation}
\label{t100}
+\left(H_{r+1}\right)_{s,t}\ \ \ \hbox{w.\ p.\ 1},
\end{equation}

\vspace{6mm}
\noindent
where

$$
\left(H_{r+1}\right)_{s,t}=
\sum_{(k,j,l_1,\ldots,l_k) \in {\rm U}_r}\ 
\frac{(s-t)^j}{j!}
\sum_{i_1,\ldots,i_k=1}^m  G_{l_1}^{(i_1)}\ldots
G_{l_k}^{(i_k)} L^j R({\bf x}_t,t)\ 
I_{{l_1\ldots l_k}_{s,t}}^{(i_1\ldots i_k)}+
\left(D_{r+1}\right)_{s,t},
$$

\vspace{3mm}
\begin{equation}
\label{w1}
{\rm D}_{q}=\left\{
(k, j, l_ {1},\ldots, l_ {k}): k + 2\left(j +
\sum\limits_{p=1}^k l_p\right)= q;\ k, j, l_{1},\ldots, 
l_{k} =0,1,\ldots\right\},
\end{equation}

\vspace{3mm}

\begin{equation}
\label{w2}
{\rm U}_{r}=\left\{
(k, j, l_ {1},\ldots, l_ {k}):
k + j +
\sum\limits_{p=1}^k l_p\le r,\  k + 2\left(j +
\sum\limits_{p=1}^k l_p\right)\ge r+1;\ k, j, l_
{1},\ldots, l_ {k} =0,1,\ldots\right\},
\end{equation}

\vspace{5mm}
\noindent
and $\left(D_{r+1}\right)_{s,t}$ has the form (\ref{5.7.12}).
Note that the remainder term $\left(H_{r+1}\right)_{s,t}$ in 
(\ref{t100}) has a higher order of smallness in the mean-square sense as
$s \downarrow t$ than the terms of the main part of expansion (\ref{t100}).

\vspace{5mm}

\section{The Second Form of the Unified Taylor--Ito Expansion}

\vspace{5mm}

Consider iterated Ito stochastic integrals of the form

$$
J^{(i_1\ldots i_k)} _ {{l_1 \ldots l_k}_{s, t}} 
=\int\limits_t^{s}(s-t_k)^{l_{k}}\ldots\
\int\limits_t^{t_ {2}}(s-t _ {1}) ^ {l_ {1}} 
d{\bf f} ^ {(i_ {1})} _ {t_ {1}} \ldots
d{\bf f}_{t_ {k}}^{(i_ {k})}\ \ \ \hbox{for}\ \ \ k\ge 1
$$

\vspace{2mm}
\noindent
and
$$
J^{(i_1\ldots i_k)} _ {{l_1 \ldots l_k}_{s, t}}=1\ \ \ \hbox{for}\ \ \
k=0,
$$

\vspace{3mm}
\noindent
where $i_1,\ldots,i_k=1,\ldots,m.$ 

The additive property of stochastic integrals and the 
Newton binomial formula imply the following
equality

\vspace{-1mm}
\begin{equation}
\label{70}
I^{(i_1\ldots i_k)} _ {{l_1 \ldots l_k}_{s, t}}=
\sum_{j_1=0}^{l_1}\ldots \sum_{j_k=0}^{l_k}
\prod_{g=1}^k C_{l_g}^{j_g}(t-s)^{l_1+\ldots+l_k-j_1-\ldots-j_k}\
J^{(i_1\ldots i_k)} _ {{j_1 \ldots j_k}_{s, t}}\ \ \ \hbox{w.\ p.\ 1},
\end{equation}

\vspace{3mm}
\noindent
where 
$$
C_l^k=\frac{l!}{k!(l-k)!}
$$

\vspace{3mm}
\noindent
is the binomial coefficient. Thus, the Taylor--Ito expansion 
of the process $\eta_s = R({\bf x}_s, s)$ 
can be constructed either using the iterated 
stochastic integrals 
$I^{(i_1\ldots i_k)} _ {{l_1 \ldots l_k}_{s, t}}$
similarly to the
previous section or using the iterated stochastic integrals 
$J^{(i_1\ldots i_k)} _ {{l_1 \ldots l_k}_{s, t}}.$
This is the main subject of this section.

Denote

\vspace{-1mm}
$$
\Biggl\Vert 
J^{(i_1\ldots i_k)} _ {{l_1 \ldots l_k}_{s, t}}\Biggr\Vert
_{i_1,\ldots,i_k=1}^{m}\stackrel{\sf def}
{=}{}^{(k)}J_{{l_1\ldots l_k}_{s,t}},
$$

\vspace{4mm}

$$
\Biggl\Vert 
L^j  G_{l_1}^{(i_1)}\ldots  G_{l_k}^{(i_k)}
R({\bf x},t) \Biggr\Vert
_{i_1,\ldots,i_k=1}^
{m}\stackrel{\sf def}{=}{}^{(k)}
L^j  G_{l_1}\ldots  G_{l_k}
R({\bf x},t).
$$

\vspace{5mm}

{\bf Theorem 3.}\ {\it Let conditions {\rm (i), (ii)} be satisfied. 
Then for any $s, t\in [0, T]$ such that $s>t$ and for any positive
integer $r$, the following expansion is valid w.\ p.\ {\rm 1}

\vspace{1mm}
$$
R({\bf x}_s,s)=
R({\bf x}_t,t)+
\sum_{q=1}^r
\sum_{(k,j,l_1,\ldots,l_k) \in {\rm A}_q}\ 
\frac{(s-t)^j}{j!}
\sum_{i_1,\ldots,i_k=1}^m  L^j G_{l_1}^{(i_1)}\ldots
G_{l_k}^{(i_k)}R({\bf x}_t,t)\ 
J_{{l_1\ldots l_k}_{s,t}}^{(i_1\ldots i_k)}+
$$

\vspace{1mm}
\begin{equation}
\label{razl44}
+\left(D_{r+1}\right)_{s,t},
\end{equation}

\vspace{5mm}
\noindent
where $\left(D_{r+1}\right)_{s,t}$ has the form {\rm (\ref{5.7.12})}.
}

\vspace{2mm}

{\bf Proof.} To prove the theorem, we check the equalities

\vspace{2mm}
$$
\sum_{(k,j,l_1,\ldots,l_k) \in {\rm A}_q}\ 
\frac{(s-t)^j}{j!}
\sum_{i_1,\ldots,i_k=1}^m  L^j G_{l_1}^{(i_1)}\ldots
G_{l_k}^{(i_k)}R({\bf x}_t,t)\ 
J_{{l_1\ldots l_k}_{s,t}}^{(i_1\ldots i_k)}=
$$

\vspace{1mm}
\begin{equation}
\label{f11}
\sum_{(k,j,l_1,\ldots,l_k) \in {\rm A}_q}\ 
\frac{(s-t)^j}{j!}
\sum_{i_1,\ldots,i_k=1}^m  G_{l_1}^{(i_1)}\ldots
G_{l_k}^{(i_k)} L^j R({\bf x}_t,t)\ 
I_{{l_1\ldots l_k}_{s,t}}^{(i_1\ldots i_k)}\ \ \ \hbox{w.\ p.\ 1}
\end{equation}

\vspace{6mm}
\noindent
for $q=1,2,\ldots,r.$ To check (\ref{f11}), substitute the expression 
(\ref{70}) 
into the right-hand side of (\ref{f11}) and then use the formulas
(\ref{a9}), (\ref{a33}).

Let us rank terms of the expansion (\ref{razl44}) according to their orders of 
smallness in the mean-square sense as $s \downarrow t$

\vspace{-2mm}
$$
R({\bf x}_s,s)=
R({\bf x}_t,t)+
\sum_{q=1}^r
\sum_{(k,j,l_1,\ldots,l_k) \in {\rm D}_q}\ 
\frac{(s-t)^j}{j!}
\sum_{i_1,\ldots,i_k=1}^m  L^j  G_{l_1}^{(i_1)}\ldots
G_{l_k}^{(i_k)}R({\bf x}_t,t)\ 
J_{{l_1\ldots l_k}_{s,t}}^{(i_1\ldots i_k)}+
$$

\vspace{2mm}
$$
+\left(H_{r+1}\right)_{s,t}\ \ \ \hbox{w.\ p.\ 1},
$$

\vspace{5mm}
\noindent
where

\vspace{-1mm}
$$
\left(H_{r+1}\right)_{s,t}=
\sum_{(k,j,l_1,\ldots,l_k) \in {\rm U}_r}\ 
\frac{(s-t)^j}{j!}
\sum_{i_1,\ldots,i_k=1}^m  L^j G_{l_1}^{(i_1)}\ldots
G_{l_k}^{(i_k)}R({\bf x}_t,t)\ 
J_{{l_1\ldots l_k}_{s,t}}^{(i_1\ldots i_k)}+
\left(D_{r+1}\right)_{s,t}.
$$

\vspace{4mm}

The term $\left(D_{r+1}\right)_{s,t}$ 
has the form (\ref{5.7.12}); the terms ${\rm D}_q$ and ${\rm U}_r$ have the 
forms (\ref{w1}) and (\ref{w2}), respectively.
Finally, we note that the convergence w.\ p.\ 1 of the 
truncated Taylor--Ito expansion (\ref{5.7.11}) (without
the remainder term $\left(D_{r+1}\right)_{s,t}$) to the process 
$R({\bf x}_s, s)$ as $r\to\infty$  
for all $s, t \in [0, T]$ such that $s>t$ and $T < \infty$
has been proved in \cite{KlPl2} (Proposition 5.9.2).
Since expansions (\ref{razl4}) and 
(\ref{razl44}) are obtained from the Taylor--Ito 
expansion (\ref{5.7.11}) without any additional
conditions, the truncated expansions (\ref{razl4}) and 
(\ref{razl44}) (without the 
reminder term $\left(D_{r+1}\right)_{s,t}$) under the
conditions of Proposition 5.9.2 \cite{KlPl2} converge to the 
process $R({\bf x}_s, s)$ w.\ p.\ 1 as $r\to\infty$ for all
$s, t \in [0, T]$ such that $s>t$ and $T < \infty.$

\vspace{5mm}

\section{The Taylor--Stratonovich Expansion}

\vspace{5mm}

In this section, we cite the Taylor--Stratonovich expansion \cite{KlPl1} 
and introduce some necessary no\-ta\-ti\-ons.
At that, we will use the original no\-ta\-ti\-ons introduced by the author
of this paper.

Let us consider two classic results.

\vspace{2mm}

{\bf Proposition~3}\ \cite{1}.\ {\it Suppose that the following conditions
are satisfied.

\vspace{1mm}

{\rm AI.}\ The functions ${\bf a}({\bf x},t),\
B_j({\bf x},t):\ \mathbb{R}^n\times[0, T]\to \mathbb{R}^{n}$ $(j=1,\ldots,m)$ 
are measurable for all
$({\bf x},t)\in \mathbb{R}^n \times [0,T],$
where
$B_j({\bf x},t)$ is the $j$th column
of the matrix $B({\bf x},t)$ {\rm(}see {\rm (\ref{1.5.2}))}.

\vspace{1mm}

{\rm AII.}\ For all ${\bf x},\ {\bf y} \in \mathbb{R}^n$
there exists a constant $K<\infty$ such that

$$
\left|{\bf a}({\bf x},t)-{\bf a}({\bf y},t)\right| + 
\sum_{j=1}^m \left|B_{j}({\bf x},t)-B_{j}({\bf y},t)\right|
\le K \left|{\bf x}-{\bf y}\right|,
$$

\vspace{1mm}
$$
\left|{\bf a}({\bf x},t)\right|^2 + 
\sum_{j=1}^m\left|B_j({\bf x},t)\right|^2 \le 
K^2\bigl(1+\left|{\bf x}\right|^2 \bigr),
$$

\vspace{3mm}
\noindent
where $\left|\cdot\right|$ is the Euclidean norm of the vector.

\vspace{1mm}

{\rm AIII.}\ A random variable ${\bf x}_0$ is ${\rm F}_0$-measurable and
${\sf M}\bigl\{\left|{\bf x}_0\right|^2\bigr\}<\infty.$

\vspace{1mm}

Then there exists a unique {\rm(}up to stochastic equivalence{\rm)}
and continuous w.~p.~{\rm 1} strong solution of the Ito
SDE {\rm (\ref{1.5.2})}.}

\vspace{2mm}

{\bf Proposition~4}\ \cite{KlPl2}.\ {\it Suppose that the conditions
{\rm AI--AIII} {\rm(}see Proposition~{\rm 3)}
are satisfied and
${\sf M}\bigl\{\left|{\bf x}_{t_0}\right|^{2n}\bigr\}<\infty$ $(n\ge 1).$
Then

$$
{\sf M}\bigl\{\left|{\bf x}_t\right|^{2n}\bigr\}\le \bigl(1+{\sf M}\bigl\{\left|{\bf x}_{t_0}\right|
^{2n}\bigr\}\bigr)e^{C(t-t_0)},
$$

\vspace{1mm}
$$
{\sf M}\bigl\{\left|{\bf x}_t - {\bf x}_{t_0}\right|^{2n}\bigr\}\le
C_1\bigl(1+{\sf M}\bigl\{\left|{\bf x}_{t_0}\right|^{2n}\bigr\}\bigr)(t-t_0)^n e^{C(t-t_0)},
$$

\vspace{5mm}
\noindent
where ${\bf x}_t$ is the solution of the Ito
SDE {\rm (\ref{1.5.2}),}
$t\in[t_0,T],$ $T<\infty,$
constant $C_1$ $(C_1\in (0, \infty))$ depends only on 
$n, K, T-t_0,$ $C=2n(2n+1)K^2,$ $K<\infty$
is a constant.}

\vspace{2mm}

Assume that
$R({\bf x},t)\in \mathfrak{L},$
$LR({\bf x}_{\tau},\tau)$, $G_0^{(i)}R({\bf x}_{\tau},\tau)
\in \mathfrak{M}_2([0, T])$
for $i = 1,\ldots,$ $m$ and
consider the Ito formula (\ref{sa}).

In addition, suppose that the function $G_0^{(i)}R({\bf x},t)$
$(i = 1,\ldots,m)$ is such that the formulas
(\ref{d11}) and (\ref{d11a}) can be applied. 
For example, assume that

\vspace{2mm}

{\rm 1}.\ $G_0^{(i)}R({\bf x},t)\in \mathfrak{L},$ $i=1,\ldots,m.$   

\vspace{1mm}

{\rm 2}.\ For all ${\bf x}, {\bf y}\in\mathbb{R}^n,$ $t, s\in[0,T],$
$i_1, i_2=1,\ldots,m$ and for some $\nu>0$

$$
\left|G_0^{(i_2)}G_0^{(i_1)}R({\bf x},t)-G_0^{(i_2)}G_0^{(i_1)}R({\bf y},t)\right|
\le K_1\left|{\bf x}-{\bf y}\right|,
$$

\vspace{1mm}

$$
\left|G_0^{(i_2)}G_0^{(i_1)}R({\bf x},t)\right|+
\left|LG_0^{(i_1)}R({\bf x},t)\right|\le K_1\left(1+\left|{\bf x}\right|\right),
$$

\vspace{1mm}

$$
\left|G_0^{(i_2)}G_0^{(i_1)}R({\bf x},s)-
G_0^{(i_2)}G_0^{(i_1)}R({\bf x},t)\right|
\le K_1\left|s-t\right|^{\nu}
\left(1+\left|{\bf x}\right|\right),
$$

\vspace{3mm}
\noindent
where $K_1<\infty$ is a constant.

\vspace{1mm}

{\rm 3}.\ Conditions {\rm AI, AII} are fulfilled {\rm (see Proposition~3)}.

\vspace{1mm}

{\rm 4}.\ ${\sf M}\{\left|{\bf x}_0\right|^8\}<\infty.$

\vspace{3mm}

Indeed, using the above conditions, Proposition~4 and the 
elementary inequality
$(a+b)^2\le 2a^2+2b^2,$ we obtain 

$$
{\sf M}\left\{\left|
G_0^{(i_2)}G_0^{(i_1)}R({\bf x}_{s},s)-
G_0^{(i_2)}G_0^{(i_1)}R({\bf x}_{t},t)
\right|^4\right\}\le
$$

\vspace{1mm}

$$
\le 8{\sf M}\left\{\left|
G_0^{(i_2)}G_0^{(i_1)}R({\bf x}_{s},s)-
G_0^{(i_2)}G_0^{(i_1)}R({\bf x}_{t},s)
\right|^4\right\}+
$$

\vspace{1mm}

$$
+ 8{\sf M}\left\{\left|
G_0^{(i_2)}G_0^{(i_1)}R({\bf x}_{t},s)
-G_0^{(i_2)}G_0^{(i_1)}R({\bf x}_{t},t)
\right|^4\right\}\le
$$

\vspace{1mm}
$$
\le 8K_1^4{\sf M}\bigl\{\left|{\bf x}_{s}-{\bf x}_{t}\right|^4 \bigr\}+
8K_1^4 \left|s-t\right|^{4\nu}{\sf M}\bigl\{\bigl(1+\left|{\bf x}_{t}\right|\bigr)^4\bigr\}\le
$$

\begin{equation}
\label{2024str1}
\le C_2\left|s-t\right|^2+C_3\left|s-t\right|^{4\nu}\to 0\ \ \ \hbox{if}\ \ \ 
s-t\to 0,
\end{equation}

\vspace{3mm}

$$
{\sf M}\left\{\left|
G_0^{(i_2)}G_0^{(i_1)}R({\bf x}_{s},s)\right|^8\right\}
\le K_1^8 {\sf M}\bigl\{\bigl(1+\left|{\bf x}_{s}\right|\bigr)^8\bigr\}
\le
$$

\vspace{1mm}

\begin{equation}
\label{2024str2}
\le C_4\bigl(1+{\sf M}\bigl\{\left|{\bf x}_{s}\right|^8\bigr\}\bigr)\le
C_5\bigl(1+\bigl(1+{\sf M}\bigl\{\left|{\bf x}_{0}\right|^{8}\bigr\}\bigr)e^{Cs}\bigr)<\infty,
\end{equation}

\vspace{5mm}
\noindent
where $C_2,\ldots,C_5<\infty$ are constants, $t, s\in[0,T].$

Analogously, we get

\vspace{-1mm}
\begin{equation}
\label{2024str3}
{\sf M}\left\{\left|
LG_0^{(i_1)}R({\bf x}_s,s)\right|^8\right\}<\infty,\ \ \ s\in[0,T].
\end{equation}

\vspace{3.5mm}

Applying the Ito formula, we obtain w.~p.~1

\vspace{-1mm}
\begin{equation}
\label{8888.01}
R({\bf x}_s,s)=R({\bf x}_t,t)+\int\limits_t^s
LR({\bf x}_{\tau},\tau)d\tau+\sum_{i_1=1}^m
\int\limits_t^s
G_0^{(i_1)}R({\bf x}_{\tau},\tau)d{\bf f}_{\tau}^{(i_1)},
\end{equation}

\vspace{2mm}

$$
G_0^{(i_1)}R({\bf x}_s,s)=G_0^{(i_1)}R({\bf x}_t,t)+\int\limits_t^s
LG_0^{(i_1)}R({\bf x}_{\tau},\tau)d\tau+
$$

\begin{equation}
\label{8888.02}
+\sum_{i_2=1}^m
\int\limits_t^s
G_0^{(i_2)}G_0^{(i_1)}R({\bf x}_{\tau},\tau)d{\bf f}_{\tau}^{(i_2)},
\end{equation}

\vspace{3mm}
\noindent
where $i_1, i_2=1,\ldots,m.$

Thus, using (\ref{2024str1})--(\ref{2024str3}),
(\ref{d11}) and (\ref{d11a}),
we have 

\begin{equation}
\label{8888.03}
\int\limits_t^s G_0^{(i)}R({\bf x}_{\tau},\tau)d{\bf f}_{\tau}^{(i)}=
{\int\limits_t^{*}}^s 
G_0^{(i)}R({\bf x}_{\tau},\tau)d{\bf f}_{\tau}^{(i)}-
\frac{1}{2}\int\limits_t^s G_0^{(i)}G_0^{(i)}R({\bf x}_{\tau},\tau)d\tau
\end{equation}

\vspace{3mm}
\noindent
w.~p.~1, where $s, t\in [0, T],$ $s>t,$ $i = 1,\ldots,m.$

Using the relation (\ref{8888.03}), let us write 
(\ref{8888.01}) in the following form

\begin{equation}
\label{sa1}
R({\bf x}_s,s)= R({\bf x}_t,t) + 
\int\limits_t^s {\bar L}R({\bf x}_\tau,\tau)d\tau
+\sum^m_{i=1}{\int\limits_t^{*}}^s 
G_0^{(i)}R({\bf x}_\tau,\tau)d{\bf f}_\tau^{(i)}\ \ \ \hbox{w.\ p.\ 1,}
\end{equation}

\noindent
where
\begin{equation}
\label{2.4a}
{\bar L}R({\bf x},t)=LR({\bf x},t)-
\frac{1}{2}\sum^m_{i=1}G_0^{(i)}G_0^{(i)}R({\bf x},t).
\end{equation}

\vspace{3mm}

Introduce the following notations

$$
\Biggl\Vert D_{\lambda_l}^{(i_l)}\ldots D_{\lambda_1}^{(i_1)}
R({\bf x},t)\Biggr\Vert_{i_1=\lambda_1,\ldots,i_l=\lambda_l}^{m\lambda_1
~\ldots~m\lambda_l}\stackrel{\sf def}{=}{}^{(p_l)}D_{\lambda_l}
\ldots D_{\lambda_1}R({\bf x},t),
$$

\vspace{4mm}

$$
{}^{(p_{k})}J^{*}_{(\lambda_{k}\ldots \lambda_1)s,t}
=\Biggl\Vert J_{(\lambda_{k}\ldots \lambda_1)s,t}^{*(i_k\ldots
i_1)}\Biggr\Vert_
{i_1=\lambda_1,\ldots,i_k=\lambda_k}^{m\lambda_1~\ldots~m\lambda_k},
$$

\vspace{3mm}

$$
{M}_k=\biggl\{(\lambda_k,\ldots,\lambda_1):
\lambda_l=1\ \hbox{or}\ \lambda_l=0;\ l=1,\ldots,k\biggr\},\ \ \ k\ge 1,
$$

\vspace{3mm}

$$
J_{(\lambda_{k}\ldots \lambda_1)s,t}^{*(i_k\ldots
i_1)}=
{\int\limits_t^{*}}^s\ldots
{\int\limits_t^{*}}^{t_{2}}
d{\bf w}_{t_{1}}^{(i_k)}\ldots
d{\bf w}_{t_k}^{(i_1)},\ \ \ k\ge 1,
$$

\vspace{3mm}
\noindent
where $\lambda_l=1$ or $\lambda_l=0$,
$D_{\lambda_l}^{(i_l)}={\bar L}$
and $i_l=0$ for $\lambda_l=0,$  $D_{\lambda_l}^{(i_l)}=G_0^{(i_l)}$
and $i_l=1,\ldots,m$ for $\lambda_l=1,$ 
$$
p_l=\sum\limits_{j=1}^l \lambda_j\ \ \ \hbox{for}\ \ \ l=1,\ldots, r+1,\ \ \
r\in\mathbb{N},
$$
${\bf w}_{\tau}^{(i)}$ $(i=1,\ldots,m)$ are ${\rm F}_{\tau}$-measurable 
for all $\tau\in [0, T]$
independent standard Wiener processes and
${\bf w}_{\tau}^{(0)}=\tau.$

Applying the formula (\ref{sa1}) to the process 
$R({\bf x}_s, s)$ repeatedly, we obtain the following 
Taylor--Stra\-to\-no\-vich expansion
\cite{KlPl1}

\begin{equation}
\label{5.7.11xxx}
R({\bf x}_s, s)=R({\bf x}_t, t)+\sum_{k=1}^r \sum_{(\lambda_{k},\ldots,
\lambda_1)\in {M}_k}\
{}^{(p_{k})}D_{\lambda_{k}}\ldots D_{\lambda_1}
R({\bf x}_t, t)\stackrel{p_k}{\cdot}
{}^{(p_{k})}J^{*}_{(\lambda_{k}\ldots \lambda_1)s,t}
+\left(D_{r+1}\right)_{s,t}
\end{equation}

\vspace{6mm}
\noindent
w.\ p.\ 1, where

\begin{equation}
\label{5.7.12xxx}
\left(D_{r+1}\right)_{s,t}
=\sum_{(\lambda_{r+1},\ldots,\lambda_1)\in{M}_{r+1}}\
{\int\limits_t^{*}}^s\ldots
\Biggl(
{\int\limits_t^{*}}^{t_2}
{}^{(p_{r+1})}D_{\lambda_{r+1}}\ldots
D_{\lambda_1}
R({\bf x}_{t_1}, t_1)\stackrel{\lambda_{r+1}}{\cdot}
d{\bf w}_{t_{1}}\Biggr)
\ldots
\stackrel{\lambda_1}{\cdot}d{\bf w}_{t_{r+1}}.
\end{equation}

\vspace{5mm}
\noindent
It is assumed that the right-hand sides of (\ref{5.7.11xxx}),
(\ref{5.7.12xxx}) exist.

Let us write 
the expansion (\ref{5.7.11xxx}) in the another form 

\vspace{-1mm}
$$
R({\bf x}_s,s)=R({\bf x}_t,t)
+\sum_{k=1}^r \sum_{(\lambda_{k},\ldots,\lambda_1)
\in{M}_k}\
\sum_{i_1=\lambda_1}^{m\lambda_1}
\ldots 
\sum_{i_k=\lambda_k}^{m\lambda_k}
D_{\lambda_k}^{(i_k)}\ldots D_{\lambda_1}^{(i_1)}
R({\bf x}_t,t)\
{J}_{(\lambda_{k}\ldots \lambda_1)s,t}^{*(i_k\ldots
i_1)}
+
$$

\vspace{1mm}
$$
+\left(D_{r+1}\right)_{s,t}\ \ \ \hbox{w.\ p.\ 1}.
$$

\vspace{4mm}

Denote
$$
{G}_{rk}=\biggl\{(\lambda_k,\ldots,\lambda_1):\ r+1\le 
2k-\lambda_1-\ldots-\lambda_k\le 2r\biggr\},
$$

\vspace{1mm}

$$
{E}_{qk}=\biggl\{(\lambda_k,\ldots,\lambda_1):\ 2k-\lambda_1-\ldots-
\lambda_k=q\biggr\},
$$

\vspace{5mm}
\noindent
where $\lambda_l=1$ or $\lambda_l=0$ $(l=1,\ldots,k).$

The Taylor--Stratonovich expansion ordered according to the 
order of smallness (in the mean-square sense when $s\downarrow  t$) 
of its terms has the form

$$
R({\bf x}_s,s)=R({\bf x}_t,t)+
\sum_{q,k=1}^r \sum_{(\lambda_{k},\ldots,\lambda_1)
\in{\rm E}_{qk}}\
\sum_{i_1=\lambda_1}^{m\lambda_1}
\ldots 
\sum_{i_k=\lambda_k}^{m\lambda_k}
D_{\lambda_k}^{(i_k)}\ldots D_{\lambda_1}^{(i_1)}
R({\bf x}_t,t)\
{J}_{(\lambda_{k}\ldots \lambda_1)s,t}^{*(i_k\ldots
i_1)}+
$$

\begin{equation}
\label{5.6.1rrrx}
+\left(H_{r+1}\right)_{s,t}\ \ \ \hbox{w.\ p.\ 1},
\end{equation}

\vspace{2mm}
\noindent
where
$$
\left(H_{r+1}\right)_{s,t}=
\sum_{k=1}^r \sum_{(\lambda_{k},\ldots,\lambda_1)
\in{\rm G}_{rk}}\
\sum_{i_1=\lambda_1}^{m\lambda_1}
\ldots 
\sum_{i_k=\lambda_k}^{m\lambda_k}
D_{\lambda_k}^{(i_k)}\ldots D_{\lambda_1}^{(i_1)}
R({\bf x}_t,t)\
{J}_{(\lambda_{k}\ldots \lambda_1)s,t}^{*(i_k\ldots
i_1)}+
$$

\begin{equation}
\label{5.6.1rrrxh}
+\left(D_{r+1}\right)_{s,t}.
\end{equation}

\vspace{3mm}

The following two questions seem interesting.

\vspace{2mm}

{\bf 1.}\ Under what conditions do the right-hand sides
of the formulas (\ref{5.6.1rrrx})
and (\ref{5.6.1rrrxh})
exist for $r\ge 2$?

\vspace{2mm}

{\bf 2.}\ Is it possible to obtain another
representation of the remainder term (\ref{5.6.1rrrxh})
for $r\ge 2$?

\vspace{2mm}

Below we will provide compelling 
arquments in favor of the following two facts.

\vspace{2mm}

{\bf (A)}\ First, one can construct the Taylor--Stratonovich expansion
(\ref{5.6.1rrrx}) $(r\ge 2)$ in such a way that its remainder term 
will coincide w.~p.~1 with the remainder term (\ref{5.6.1rrrh}) $(r\ge 2)$
of the Taylor--Ito expansion (\ref{5.6.1rrr}) $(r\ge 2)$.

\vspace{2mm}

{\bf (B)}\ Second, the truncated Taylor--Stratonovich expansion (\ref{5.6.1rrrx}) $(r\ge 2)$
(without the remainder term (\ref{5.6.1rrrxh}) $(r\ge 2)$)
will coincide w.~p.~1 with the truncated Taylor--Ito
expansion (\ref{5.6.1rrr}) $(r\ge 2)$
(without the remainder term (\ref{5.6.1rrrh}) $(r\ge 2)$).

\vspace{2mm}

This means that the right-hand side of (\ref{5.6.1rrrx}) $(r\ge 2)$
(in which the remainder term will have the form (\ref{5.6.1rrrh}) $(r\ge 2)$)
will exist under the conditions (i), (ii) (see Sect.~3).

Let us begin our reasoning with Proposition~5 (see below).
This proposition allows us to represent 
the iterated Stratonovich stochastic integral
of multiplicity $k$ $(k\in\mathbb{N})$
as a sum of iterated Ito stochastic integrals
and its mathematical expectation.

Let us introduce the following notations

$$
J[\psi^{(k)}]_{T,t}^{s_l,\ldots,s_1}\ \stackrel{\rm def}{=}\ 
\prod_{p=1}^l {\bf 1}_{\{i_{s_p}=
i_{s_{p}+1}\ne 0\}}\ \times
$$

\vspace{1mm}
$$
\times
\int\limits_t^T\psi_k(t_k)\ldots \int\limits_t^{t_{s_l+3}}
\psi_{s_l+2}(t_{s_l+2})
\int\limits_t^{t_{s_l+2}}\psi_{s_l}(t_{s_l+1})
\psi_{s_l+1}(t_{s_l+1}) \times
$$

\vspace{1mm}
$$
\times
\int\limits_t^{t_{s_l+1}}\psi_{s_l-1}(t_{s_l-1})
\ldots
\int\limits_t^{t_{s_1+3}}\psi_{s_1+2}(t_{s_1+2})
\int\limits_t^{t_{s_1+2}}\psi_{s_1}(t_{s_1+1})
\psi_{s_1+1}(t_{s_1+1}) \times
$$

\vspace{1mm}
$$
\times
\int\limits_t^{t_{s_1+1}}\psi_{s_1-1}(t_{s_1-1})
\ldots \int\limits_t^{t_2}\psi_1(t_1)
d{\bf w}_{t_1}^{(i_1)}\ldots d{\bf w}_{t_{s_1-1}}^{(i_{s_1-1})}
dt_{s_1+1}d{\bf w}_{t_{s_1+2}}^{(i_{s_1+2})}\ldots
$$

\vspace{1mm}
\begin{equation}
\label{30.1}
\ldots\
d{\bf w}_{t_{s_l-1}}^{(i_{s_l-1})}
dt_{s_l+1}d{\bf w}_{t_{s_l+2}}^{(i_{s_l+2})}\ldots d{\bf w}_{t_k}^{(i_k)},
\end{equation}

\vspace{2mm}
\noindent
where 
\begin{equation}
\label{30.5550001}
{\rm A}_{k,l}
=\bigl\{(s_l,\ldots,s_1):\
s_l>s_{l-1}+1,\ldots,s_2>s_1+1,\ s_l,\ldots,s_1=1,\ldots,k-1\bigr\},
\end{equation}

\vspace{1mm}
$$
(s_l,\ldots,s_1)\in{\rm A}_{k,l},\ \ \ 
l=1,\ldots,\left[k/2\right],\ \ \
i_s=0, 1,\ldots,m,\ \ \
s=1,\ldots,k,
$$

\vspace{3mm}
\noindent
$[x]$ is an
integer
part of a real number $x,$
and ${\bf 1}_A$ is the indicator of the set $A$.

Let us formulate the statement on connection 
between
iterated 
Stra\-to\-no\-vich and Ito stochastic integrals 
$J^{*}[\psi^{(k)}]_{T,t},$ $J[\psi^{(k)}]_{T,t}$ 
of fixed multiplicity $k,$ $k\in\mathbb{N}$ (see (\ref{str}), (\ref{ito})).

\vspace{2mm}

{\bf Proposition~5} \cite{20aa} (Theorem~2.12), \cite{322} (1997). {\it Suppose that
every $\psi_l(\tau)$ $(l=1,\ldots,k)$ is a continuous
function at the interval $[t, T]$.
Then$,$ the following relation between iterated
Stra\-to\-no\-vich and Ito stochastic integrals 

\begin{equation}
\label{30.4}
J^{*}[\psi^{(k)}]_{T,t}=J[\psi^{(k)}]_{T,t}+
\sum_{r=1}^{\left[k/2\right]}\frac{1}{2^r}
\sum_{(s_r,\ldots,s_1)\in {\rm A}_{k,r}}
J[\psi^{(k)}]_{T,t}^{s_r,\ldots,s_1}\ \ \ \hbox{{\rm w.\ p.\ 1}}
\end{equation}

\vspace{4mm}
\noindent
is correct, 
where $\sum\limits_{\emptyset}$ is supposed to be equal to zero{\rm .}
}

\vspace{2mm}

For example, from Proposition~5 for $k=1, 2, 3, 4$ we obtain
the following well known equalities, which are fulfilled
w.~p.~1

$$
{\int\limits_t^{*}}^T\psi_1(t_1)d{\bf w}_{t_1}^{(i_1)}=
\int\limits_t^T\psi_1(t_1)d{\bf w}_{t_1}^{(i_1)},
$$

\vspace{2mm}

$$
{\int\limits_t^{*}}^T\psi_2(t_2)
{\int\limits_t^{*}}^{t_2}\psi_1(t_1)d{\bf w}_{t_1}^{(i_1)}
d{\bf w}_{t_2}^{(i_2)}=
\int\limits_t^T\psi_2(t_2)
\int\limits_t^{t_2}\psi_1(t_1)d{\bf w}_{t_1}^{(i_1)}
d{\bf w}_{t_2}^{(i_2)}+
$$

\vspace{1mm}
$$
+ \frac{1}{2}{\bf 1}_{\{i_1=i_2\ne 0\}}
\int\limits_t^T\psi_2(t_2)\psi_1(t_2)dt_2,
$$

\vspace{3mm}
$$
{\int\limits_t^{*}}^T\psi_3(t_3)\ldots
{\int\limits_t^{*}}^{t_2}\psi_1(t_1)d{\bf w}_{t_1}^{(i_1)}\ldots
d{\bf w}_{t_3}^{(i_3)}=\hspace{-1mm}
\int\limits_t^T\psi_3(t_3)\ldots\hspace{-0.5mm}
\int\limits_t^{t_2}\psi_1(t_1)d{\bf w}_{t_1}^{(i_1)}\ldots
d{\bf w}_{t_3}^{(i_3)}+
$$

\vspace{1mm}
$$
+ \frac{1}{2}{\bf 1}_{\{i_1=i_2\ne 0\}}
\int\limits_t^T\psi_3(t_3)
\int\limits_t^{t_3}\psi_2(t_2)\psi_1(t_2)dt_2
d{\bf w}_{t_3}^{(i_3)}+
$$

\vspace{1mm}
$$
+ \frac{1}{2}{\bf 1}_{\{i_2=i_3\ne 0\}}
\int\limits_t^T\psi_3(t_3)\psi_2(t_3)
\int\limits_t^{t_3}\psi_1(t_1)d{\bf w}_{t_1}^{(i_1)}
dt_3,
$$

\vspace{3mm}
$$
{\int\limits_t^{*}}^T\psi_4(t_4)\ldots
{\int\limits_t^{*}}^{t_2}\psi_1(t_1)d{\bf w}_{t_1}^{(i_1)}\ldots
d{\bf w}_{t_4}^{(i_4)}=
\int\limits_t^T\psi_4(t_4)\ldots
\int\limits_t^{t_2}\psi_1(t_1)d{\bf w}_{t_1}^{(i_1)}\ldots
d{\bf w}_{t_4}^{(i_4)}+
$$

\vspace{1mm}
$$
+ \frac{1}{2}{\bf 1}_{\{i_1=i_2\ne 0\}}
\int\limits_t^T\psi_4(t_4)\int\limits_t^{t_4}
\psi_3(t_3)\int\limits_t^{t_3}\psi_1(t_2)\psi_2(t_2)dt_2
d{\bf w}_{t_3}^{(i_3)}d{\bf w}_{t_4}^{(i_4)} +
$$

\vspace{1mm}
$$
+ \frac{1}{2}{\bf 1}_{\{i_2=i_3\ne 0\}}
\int\limits_t^T\psi_4(t_4)\int\limits_t^{t_4}
\psi_3(t_3)\psi_2(t_3)\int\limits_t^{t_3}\psi_1(t_1)
d{\bf w}_{t_1}^{(i_1)}dt_3d{\bf w}_{t_4}^{(i_4)} +
$$

\vspace{1mm}
$$
+ \frac{1}{2}{\bf 1}_{\{i_3=i_4\ne 0\}}
\int\limits_t^T\psi_4(t_4)\psi_3(t_4)\int\limits_t^{t_4}
\psi_2(t_2)\int\limits_t^{t_2}\psi_1(t_1)
d{\bf w}_{t_1}^{(i_1)}d{\bf w}_{t_2}^{(i_2)}dt_4 +
$$

\vspace{1mm}
$$
+ \frac{1}{4}{\bf 1}_{\{i_1=i_2\ne 0\}}{\bf 1}_{\{i_3=i_4\ne 0\}}
\int\limits_t^T\psi_4(t_4)\psi_3(t_4)\int\limits_t^{t_4}
\psi_2(t_2)\psi_1(t_2)dt_2 dt_4.
$$

\vspace{4mm}

It is obvious that it is possible to obtain an inverse
formula that will express the iterated Ito stochastic integral (\ref{ito})
as a sum of iterated Stratonovich stochastic integrals (\ref{str}).
Below we present the corresponding proposition.

\vspace{2mm}

{\bf Proposition~6}\ \cite{20aa} (Proposition~4.3).\ {\it Suppose that
every $\psi_l(\tau)$ $(l=1,\ldots,k)$ is a continuous
function at the interval $[t, T]$.
Then$,$ the following relation between iterated
Ito and Stra\-to\-no\-vich stochastic integrals 

\begin{equation}
\label{2024str11}
J[\psi^{(k)}]_{T,t}=J^{*}[\psi^{(k)}]_{T,t}+
\sum_{r=1}^{\left[k/2\right]}\frac{(-1)^r}{2^r}
\sum_{(s_r,\ldots,s_1)\in {\rm A}_{k,r}}
J^{*}[\psi^{(k)}]_{T,t}^{s_r,\ldots,s_1}\ \ \ \hbox{{\rm w.\ p.\ 1}}
\end{equation}

\vspace{4mm}
\noindent
is correct, 
where $\sum\limits_{\emptyset}$ is supposed to be equal to zero$,$
$J[\psi^{(k)}]_{T,t}$ and $J^{*}[\psi^{(k)}]_{T,t}$
are defined by {\rm (\ref{ito})} and {\rm (\ref{str}),}
respectively$,$

\vspace{-2mm}
$$
J^{*}[\psi^{(k)}]_{T,t}^{s_l,\ldots,s_1}\ \stackrel{\rm def}{=}\ 
\prod_{p=1}^l {\bf 1}_{\{i_{s_p}=
i_{s_{p}+1}\ne 0\}}\ \times
$$

\vspace{1mm}
$$
\times
{\int\limits_t^{*}}^T
\psi_k(t_k)\ldots 
{\int\limits_t^{*}}^{t_{s_l+3}}
\psi_{s_l+2}(t_{s_l+2})
\int\limits_t^{t_{s_l+2}}\psi_{s_l}(t_{s_l+1})
\psi_{s_l+1}(t_{s_l+1}) \times
$$

\vspace{1mm}
$$
\times
{\int\limits_t^{*}}^{t_{s_l+1}}
\psi_{s_l-1}(t_{s_l-1})
\ldots
{\int\limits_t^{*}}^{t_{s_1+3}}
\psi_{s_1+2}(t_{s_1+2})
\int\limits_t^{t_{s_1+2}}\psi_{s_1}(t_{s_1+1})
\psi_{s_1+1}(t_{s_1+1}) \times
$$

\vspace{1mm}
$$
\times
{\int\limits_t^{*}}^{t_{s_1+1}}
\psi_{s_1-1}(t_{s_1-1})
\ldots 
{\int\limits_t^{*}}^{t_2}
\psi_1(t_1)
d{\bf w}_{t_1}^{(i_1)}\ldots d{\bf w}_{t_{s_1-1}}^{(i_{s_1-1})}
dt_{s_1+1}d{\bf w}_{t_{s_1+2}}^{(i_{s_1+2})}\ldots
$$

\vspace{2mm}
$$
\ldots\
d{\bf w}_{t_{s_l-1}}^{(i_{s_l-1})}
dt_{s_l+1}d{\bf w}_{t_{s_l+2}}^{(i_{s_l+2})}\ldots d{\bf w}_{t_k}^{(i_k)},
$$

\vspace{2mm}
\noindent
where 
$$
{\rm A}_{k,l}
=\bigl\{(s_l,\ldots,s_1):\
s_l>s_{l-1}+1,\ldots,s_2>s_1+1,\ s_l,\ldots,s_1=1,\ldots,k-1\bigr\},
$$

$$
(s_l,\ldots,s_1)\in{\rm A}_{k,l},\ \ \ 
l=1,\ldots,\left[k/2\right],\ \ \
i_s=0, 1,\ldots,m,\ \ \
s=1,\ldots,k,
$$

\vspace{3mm}
\noindent
$[x]$ is an integer part of a real number $x,$
${\bf 1}_A$ is the indicator of the set $A.$}

\vspace{3mm}

For example, from Proposition~6 for $k=1, 2, 3, 4$ we obtain
the following equalities w.~p.~1

$$
\int\limits_t^T\psi_1(t_1)d{\bf w}_{t_1}^{(i_1)}=
{\int\limits_t^{*}}^T\psi_1(t_1)d{\bf w}_{t_1}^{(i_1)},
$$

\vspace{2mm}

$$
\int\limits_t^T\psi_2(t_2)
\int\limits_t^{t_2}\psi_1(t_1)d{\bf w}_{t_1}^{(i_1)}
d{\bf w}_{t_2}^{(i_2)}=
{\int\limits_t^{*}}^T\psi_2(t_2)
{\int\limits_t^{*}}^{t_2}\psi_1(t_1)d{\bf w}_{t_1}^{(i_1)}
d{\bf w}_{t_2}^{(i_2)}-
$$

\vspace{1mm}
$$
-\frac{1}{2}{\bf 1}_{\{i_1=i_2\ne 0\}}
\int\limits_t^T\psi_2(t_2)\psi_1(t_2)dt_2,
$$

\vspace{3mm}
$$
\int\limits_t^T\psi_3(t_3)\ldots\hspace{-0.5mm}
\int\limits_t^{t_2}\psi_1(t_1)d{\bf w}_{t_1}^{(i_1)}\ldots
d{\bf w}_{t_3}^{(i_3)}=
\hspace{-1mm}
{\int\limits_t^{*}}^T\psi_3(t_3)\ldots
{\int\limits_t^{*}}^{t_2}\psi_1(t_1)d{\bf w}_{t_1}^{(i_1)}\ldots
d{\bf w}_{t_3}^{(i_3)}-
$$

\vspace{1mm}
$$
-\frac{1}{2}{\bf 1}_{\{i_1=i_2\ne 0\}}
{\int\limits_t^{*}}^T
\psi_3(t_3)
\int\limits_t^{t_3}\psi_2(t_2)\psi_1(t_2)dt_2
d{\bf w}_{t_3}^{(i_3)}-
$$

\vspace{1mm}
$$
-\frac{1}{2}{\bf 1}_{\{i_2=i_3\ne 0\}}
\int\limits_t^T\psi_3(t_3)\psi_2(t_3)
{\int\limits_t^{*}}^{t_3}
\psi_1(t_1)d{\bf w}_{t_1}^{(i_1)}
dt_3,
$$

\vspace{3mm}
$$
\int\limits_t^T\psi_4(t_4)\ldots\hspace{-0.5mm}
\int\limits_t^{t_2}\psi_1(t_1)d{\bf w}_{t_1}^{(i_1)}\ldots
d{\bf w}_{t_4}^{(i_4)}=
\hspace{-1mm}
{\int\limits_t^{*}}^T
\psi_4(t_4)\ldots
{\int\limits_t^{*}}^{t_2}\psi_1(t_1)d{\bf w}_{t_1}^{(i_1)}\ldots
d{\bf w}_{t_4}^{(i_4)}
-
$$

\vspace{1mm}
$$
-\frac{1}{2}{\bf 1}_{\{i_1=i_2\ne 0\}}
{\int\limits_t^{*}}^T
\psi_4(t_4)
{\int\limits_t^{*}}^{t_4}
\psi_3(t_3)\int\limits_t^{t_3}\psi_1(t_2)\psi_2(t_2)dt_2
d{\bf w}_{t_3}^{(i_3)}d{\bf w}_{t_4}^{(i_4)} -
$$

\vspace{1mm}
$$
-\frac{1}{2}{\bf 1}_{\{i_2=i_3\ne 0\}}
{\int\limits_t^{*}}^T\psi_4(t_4)\int\limits_t^{t_4}
\psi_3(t_3)\psi_2(t_3)
{\int\limits_t^{*}}^{t_3}
\psi_1(t_1)
d{\bf w}_{t_1}^{(i_1)}dt_3d{\bf w}_{t_4}^{(i_4)} -
$$

\vspace{1mm}
$$
- \frac{1}{2}{\bf 1}_{\{i_3=i_4\ne 0\}}
\int\limits_t^T\psi_4(t_4)\psi_3(t_4)
{\int\limits_t^{*}}^{t_4}
\psi_2(t_2)
{\int\limits_t^{*}}^{t_2}
\psi_1(t_1)
d{\bf w}_{t_1}^{(i_1)}d{\bf w}_{t_2}^{(i_2)}dt_4 +
$$

\vspace{1mm}
$$
+ \frac{1}{4}{\bf 1}_{\{i_1=i_2\ne 0\}}{\bf 1}_{\{i_3=i_4\ne 0\}}
\int\limits_t^T\psi_4(t_4)\psi_3(t_4)\int\limits_t^{t_4}
\psi_2(t_2)\psi_1(t_2)dt_2 dt_4.
$$

\vspace{4mm}

Further, using Proposition~6, we obtain for $r\ge 2$

$$
R({\bf x}_t,t)+
\sum_{q,k=1}^r \sum_{(\lambda_{k},\ldots,\lambda_1)
\in{\rm E}_{qk}}
\sum_{i_1=\lambda_1}^{m\lambda_1}
\ldots 
\sum_{i_k=\lambda_k}^{m\lambda_k}
Q_{\lambda_k}^{(i_k)}\ldots Q_{\lambda_1}^{(i_1)}
R({\bf x}_t,t){J}_{(\lambda_{k}\ldots \lambda_1)s,t}^{(i_k\ldots
i_1)}=
$$

\vspace{1mm}
\begin{equation}
\label{2024str201}
=R({\bf x}_t,t)+
\sum_{q,k=1}^r \sum_{(\lambda_{k},\ldots,\lambda_1)
\in{\rm E}_{qk}}
\sum_{i_1=\lambda_1}^{m\lambda_1}
\ldots 
\sum_{i_k=\lambda_k}^{m\lambda_k}
D_{\lambda_k}^{(i_k)}\ldots D_{\lambda_1}^{(i_1)}
R({\bf x}_t,t)
{J}_{(\lambda_{k}\ldots \lambda_1)s,t}^{*(i_k\ldots
i_1)}
\end{equation}
        
\vspace{3mm}
\noindent
w.~p.~1, where notations are the same as in 
(\ref{5.6.1rrr}), (\ref{5.6.1rrrx}).
Thus, (A) and (B) take place.

\vspace{5mm}

\section{The First Form of the Unified Taylor--Stratonovich Expansion}

\vspace{5mm}

In this section, we transform the right-hand side of (\ref{5.7.11xxx}) with 
the help of Theorem 1 and Lemma 2 to a
representation including iterated 
Stratonovich stochastic integrals of the form (\ref{a111}).

Denote
\begin{equation}
\label{opr1xxx}
I^{*(i_1\ldots i_k)} _ {{l_1 \ldots l_k}_{s, t}} 
={\int\limits_t^{*}}^{s}(t-t_k)^{l_{k}}\ldots\
{\int\limits_t^{*}}^{t _ {2}}(t-t _ {1}) ^ {l_ {1}} 
d{\bf f} ^ {(i_ {1})} _ {t_ {1}} \ldots
d{\bf f}_{t_ {k}}^{(i_ {k})}\ \ \ \hbox{for}\ \ \ k\ge 1\ \ \ 
\left(I^{*(i_1\ldots i_k)} _ {{l_1 \ldots l_k}_{s, t}}=1\ \ \ 
\hbox{for}\ \ \ k=0\right),
\end{equation}

\vspace{3mm}
\noindent
where $i_1,\ldots,i_k=1,\ldots,m.$ Moreover, let

$$
{}^{(k)}I^{*}_{{l_1\ldots l_k}_{s,t}}=\Biggl\Vert 
I^{*(i_1\ldots i_k)} _ {{l_1 \ldots l_k}_{s, t}}\Biggr\Vert
_{i_1,\ldots,i_k=1}^{m},
$$

\vspace{2mm}

\begin{equation}
\label{a9x}
\bar G_p^{(i)}\stackrel{\sf def}{=}\frac{1}{p}\left(
\bar G_{p-1}^{(i)}\bar L-\bar L\bar G_{p-1}^{(i)}\right),\ \ \
p=1, 2,\ldots,\ \ \ i=1,\ldots,m,
\end{equation}

\vspace{4mm}
\noindent
where $\bar G_0^{(i)}\stackrel{\sf def}{=}G_0^{(i)},$
$i=1,\ldots,m.$
The operators $\bar L$ and $G_0^{(i)},$ $i=1,\ldots,m,$
are determined by the equalities
(\ref{2.3}), (\ref{2.4}), and (\ref{2.4a}). Denote

\vspace{1mm}

$$
{A}_q\stackrel{\sf def}{=}
\Biggl\{
(k,j,l_1,\ldots,l_k):\ k+j+\sum_{p=1}^k l_p=q;\ k,j,l_1,\ldots,l_k=0, 1,\ldots
\Biggr\},
$$

\vspace{3mm}
$$
\Biggl\Vert 
\bar G_{l_1}^{(i_1)}\ldots \bar G_{l_k}^{(i_k)}\bar L^j
R({\bf x},t) \Biggr\Vert
_{i_1,\ldots,i_k=1}^
{m}\stackrel{\sf def}{=}{}^{(k)}
\bar G_{l_1}\ldots \bar G_{l_k}\bar L^j
R({\bf x},t),
$$

\vspace{3mm}

$$
\bar L^j R({\bf x},t)\stackrel{\sf def}{=}
\begin{cases}\underbrace{\bar L\ldots\bar L}_j
R({\bf x},t)\ &\hbox{for}\ j\ge 1\cr\cr
R({\bf x},t)\ &\hbox{for}\ j=0
\end{cases}.
$$

\vspace{5mm}

{\bf Theorem 4.}\ {\it Suppose that sufficient conditions are 
satisfied under which the right-hand 
sides of {\rm(\ref{5.7.11xxx}), (\ref{5.7.12xxx})} exist. 
Then for any $s, t \in [0, T]$ such that $s>t$ 
and for any positive
integer $r$, the following expansion takes place w.\ p.\ {\rm 1}
\begin{equation}
\label{razl4x}
R({\bf x}_s,s)=
R({\bf x}_t,t)+
\sum_{q=1}^r
\sum_{(k,j,l_1,\ldots,l_k) \in {\rm A}_q}\ 
\frac{(s-t)^j}{j!}
\sum_{i_1,\ldots,i_k=1}^m \bar G_{l_1}^{(i_1)}\ldots
\bar G_{l_k}^{(i_k)}\bar L^j R({\bf x}_t,t)\
I_{{l_1\ldots l_k}_{s,t}}^{*(i_1\ldots i_k)}+\left(D_{r+1}\right)_{s,t},
\end{equation}

\vspace{4mm}
\noindent
where $\left(D_{r+1}\right)_{s,t}$ has the form {\rm (\ref{5.7.12xxx})}.
}

\vspace{3mm}

{\bf Proof.} We claim that

\vspace{1mm}
$$
\sum_{(\lambda_{q},\ldots,
\lambda_1)\in {M}_q}\
{}^{(p_{q})}D_{\lambda_{q}}\ldots D_{\lambda_1}
R({\bf x}_t,t)\stackrel{p_q}{\cdot}
{}^{(p_{q})}J^{*}_{(\lambda_{q}\ldots \lambda_1)s,t}=
$$

\vspace{2mm}
\begin{equation}
\label{a22x}
=
\sum_{(k,j,l_1,\ldots,l_k) \in {\rm A}_q}\ 
\frac{(s-t)^j}{j!}
\sum_{i_1,\ldots,i_k=1}^m \bar G_{l_1}^{(i_1)}\ldots
\bar G_{l_k}^{(i_k)}\bar L^j R({\bf x}_t,t)\ 
I_{{l_1\ldots l_k}_{s,t}}^{*(i_1\ldots i_k)}
\end{equation}

\vspace{6mm}
\noindent
w.\ p.\ 1. The equality (\ref{a22x}) is valid for $q = 1$. Assume 
that (\ref{a22x}) is valid for some $q > 1$. In 
this case, using the induction hypothesis, we obtain

\vspace{1mm}
$$
\sum_{(\lambda_{q+1},\ldots,
\lambda_1)\in {M}_{q+1}}\
{}^{(p_{q+1})}D_{\lambda_{1}}\ldots D_{\lambda_{q+1}}
R({\bf x}_t,t)\stackrel{p_{q+1}}{\cdot}
{}^{(p_{q+1})}J^{*}_{(\lambda_{1}\ldots \lambda_{q+1})s,t}=
$$

\vspace{3mm}
$$
=\sum_{\lambda_{q+1}\in\{1,\ 0\}}
{\int\limits_t^{*}}^s
\sum_{(\lambda_{q},\ldots,
\lambda_1)\in {M}_{q}}
\Biggl({}^{(p_{q+1})}D_{\lambda_{1}}
\ldots
D_{\lambda_{q+1}}
R({\bf x}_t,t)\stackrel{p_{q}}{\cdot}
{}^{(p_{q})}J^{*}_{(\lambda_{1}\ldots \lambda_{q})\theta,t}\Biggr)
\stackrel{\lambda_{q+1}}{\cdot}d{\bf w}_{\theta}=
$$

\vspace{3mm}
$$
=\sum_{\lambda_{q+1}\in\{1,\ 0\}}
{\int\limits_t^{*}}^s
\sum_{(k,j,l_1,\ldots,l_k)\in{A}_q}
\frac{(\theta-t)^j}{j!}\times
$$

\vspace{3mm}
$$
\times
\Biggl({}^{(k+\lambda_{q+1})}\bar G_{l_1}\ldots \bar G_{l_k}\bar L^j
D_{\lambda_{q+1}}R({\bf x}_t,t)
\stackrel{k}{\cdot}
{}^{(k)}I^{*}_{{l_1\ldots l_k}_{s,t}}\Biggr)
\stackrel{\lambda_{q+1}}{\cdot}d{\bf w}_{\theta}=
$$

\vspace{3mm}
$$
=\sum_{(k,j,l_1,\ldots,l_k)\in{A}_q}\left(
{}^{(k)}\bar G_{l_1}\ldots \bar G_{l_k}\bar L^{j+1}
R({\bf x}_t,t)
\stackrel{k}{\cdot}
\int\limits_t^s\frac{(\theta-t)^j}{j!}
{}^{(k)}I^{*}_{{l_1\ldots l_k}_{\theta,t}}d\theta+\right.
$$

\vspace{3mm}
\begin{equation}
\label{a30x}
\left.+\left(
{}^{(k+1)}\bar G_{l_1}\ldots \bar G_{l_k}\bar L^{j}\bar G_0
R({\bf x}_t,t)
\stackrel{k}{\cdot}
{\int\limits_t^{*}}^s
\frac{(\theta-t)^j}{j!}
{}^{(k)}I^{*}_{{l_1\ldots l_k}_{\theta,t}}\right)
\stackrel{1}{\cdot}d{\bf f}_{\theta}\right)
\end{equation}

\vspace{6mm}
\noindent
w.\ p.\ 1.

Using Lemma 2, we obtain

\vspace{-2mm}
$$
\int\limits_t^s\frac{(\theta-t)^j}{j!}
{}^{(k)}I^{*}_{{l_1\ldots l_k}_{\theta,t}}d\theta=
$$
\begin{equation}
\label{a31x}
=\frac{1}{(j+1)!}
\begin{cases}(s-t)^{j+1}\ &\hbox{for}\ k=0 \cr \cr 
(s-t)^{j+1} \cdot {}^{(k)}I^{*}_{{l_1\ldots l_k}_{s,t}}-
(-1)^{j+1} \cdot
{}^{(k)}I^{*}_{{l_1\ldots l_{k-1}\ l_k+j+1}_{s,t}}\ &\hbox{for}\ k>0
\end{cases}
\end{equation}

\vspace{4mm}
\noindent
w.\ p.\ 1. In addition (see (\ref{opr1xxx})), we get

\begin{equation}
\label{a32x}
{\int\limits_t^{*}}^s\frac{(\theta-t)^j}{j!}
I^{*(i_1\ldots i_k)}_{{l_1\ldots l_k}_{\theta,t}}
d{\bf f}_{\theta}^{(i_{k+1})}=
\frac{(-1)^j}{j!}I_{{l_1\ldots l_k j}_{s,t}}^{*(i_1\ldots i_k i_{k+1})}
\end{equation}

\vspace{2mm}
\noindent
in the notations just introduced.
Substitute the relations (\ref{a31x}) and (\ref{a32x}) 
into the formula (\ref{a30x}). Grouping summands of the 
obtained expression with
equal lower indices at iterated Stratonovich stochastic
integrals and using (\ref{a9x}) and the equality

\begin{equation}
\label{a33x}
\bar G_p^{(i)}R({\bf x},t)=\frac{1}{p!}
\sum_{q=0}^p(-1)^q C_p^q\bar L^q\bar G_0^{(i)}\bar L^{p-q}
R({\bf x},t),\ \ \ \hbox{where}\ \ \ C_p^q=\frac{p!}{q!(p-q)!},
\end{equation}

\vspace{3mm}
\noindent
(this equality follows from (\ref{a9x})), we note that the obtained 
expression is equal to

$$
\sum_{(k,j,l_1,\ldots,l_k)\in{A}_{q+1}}
\frac{(s-t)^j}{j!}
{}^{(k)}\bar G_{l_1}\ldots \bar G_{l_k}\bar L^j\{\eta_t\} 
\stackrel{k}{\cdot}
{}^{(k)}I^{*}_{{l_1\ldots l_k}_{s,t}}
$$

\vspace{4mm}
\noindent
w.\ p.\ 1. Summing the equalities (\ref{a22x}) for $q = 1, 2,\ldots,r$ 
and applying the formula (\ref{5.7.11xxx}), we obtain the expression
(\ref{razl4x}). The proof is completed.

Let us order terms of the expansion (\ref{razl4x}) according to 
their smallness orders as $s \downarrow t$ in the mean-square sense

\vspace{-1mm}
$$
R({\bf x}_s,s)=
R({\bf x}_t,t)+
\sum_{q=1}^r
\sum_{(k,j,l_1,\ldots,l_k) \in {\rm D}_q}\ 
\frac{(s-t)^j}{j!}
\sum_{i_1,\ldots,i_k=1}^m \bar G_{l_1}^{(i_1)}\ldots
\bar G_{l_k}^{(i_k)}\bar L^j R({\bf x}_t,t)\
I_{{l_1\ldots l_k}_{s,t}}^{*(i_1\ldots i_k)}+
$$

\vspace{1mm}
\begin{equation}
\label{t100x}
+\left(H_{r+1}\right)_{s,t}\ \ \ \hbox{w.\ p.\ 1},
\end{equation}

\vspace{5mm}
\noindent
where 

\vspace{-1mm}
$$
\left(H_{r+1}\right)_{s,t}=
\sum_{(k,j,l_1,\ldots,l_k) \in {\rm U}_r}\ 
\frac{(s-t)^j}{j!}
\sum_{i_1,\ldots,i_k=1}^m \bar G_{l_1}^{(i_1)}\ldots
\bar G_{l_k}^{(i_k)}\bar L^j R({\bf x}_t,t)\ 
I_{{l_1\ldots l_k}_{s,t}}^{*(i_1\ldots i_k)}+
\left(D_{r+1}\right)_{s,t},
$$

\vspace{1mm}
\begin{equation}
\label{w1x}
{\rm D}_{q}=\left\{
(k, j, l_ {1},\ldots, l_ {k}): k + 2\left(j +
\sum\limits_{p=1}^k l_p\right)= q;\ k, j, l_{1},\ldots, 
l_{k} =0,1,\ldots\right\},
\end{equation}

\vspace{2mm}

\begin{equation}
\label{w2x}
{\rm U}_{r}=\left\{
(k, j, l_ {1},\ldots, l_ {k}):
k + j +
\sum\limits_{p=1}^k l_p\le r,\  k + 2\left(j +
\sum\limits_{p=1}^k l_p\right)\ge r+1;\ k, j, l_
{1},\ldots, l_ {k} =0,1,\ldots\right\},
\end{equation}

\vspace{5mm}
\noindent
and $\left(D_{r+1}\right)_{s,t}$ has the form (\ref{5.7.12xxx}).
Note that the remainder term $\left(H_{r+1}\right)_{s,t}$ in 
(\ref{t100x}) has a higher order of smallness in the mean-square sense as
$s \downarrow t$ than the terms of the main part of expansion (\ref{t100x}).

\vspace{5mm}

\section{The Second Form of the Unified Taylor--Stratonovich Expansion}

\vspace{5mm}

Consider iterated Stratonovich stochastic integrals of the form

$$
J^{*(i_1\ldots i_k)} _ {{l_1 \ldots l_k}_{s, t}} 
={\int\limits_t^{*}}^{s}(s-t_k)^{l_{k}}\ldots\
{\int\limits_t^{*}}^{t _ {2}}(s-t_ {1}) ^ {l_ {1}} 
d{\bf f} ^ {(i_ {1})} _ {t_ {1}} \ldots
d{\bf f}_{t_ {k}}^{(i_ {k})}\ \ \ \hbox{for}\ \ \ k\ge 1
$$

\vspace{2mm}
\noindent
and
$$
J^{*(i_1\ldots i_k)} _ {{l_1 \ldots l_k}_{s, t}}=1\ \ \ \hbox{for}\ \ \
k=0,
$$

\vspace{4mm}
\noindent
where $i_1,\ldots,i_k=1,\ldots,m.$ 

The additive property of stochastic integrals and the 
Newton binomial formula imply the following
equality

\vspace{-1mm}
\begin{equation}
\label{70x}
I^{*(i_1\ldots i_k)} _ {{l_1 \ldots l_k}_{s, t}}=
\sum_{j_1=0}^{l_1}\ldots \sum_{j_k=0}^{l_k}
\prod_{g=1}^k C_{l_g}^{j_g}(t-s)^{l_1+\ldots+l_k-j_1-\ldots-j_k}\
J^{*(i_1\ldots i_k)} _ {{j_1 \ldots j_k}_{s, t}}\ \ \ \hbox{w.\ p.\ 1},
\end{equation}

\vspace{2mm}
\noindent
where 
$$
C_l^k=\frac{l!}{k!(l-k)!}
$$

\vspace{2mm}
\noindent
is the binomial coefficient. Thus, the Taylor--Stratonovich expansion 
of the process $\eta_s = R({\bf x}_s, s)$, $s\in [0, T]$
can be constructed either using the iterated 
stochastic integrals 
$I^{*(i_1\ldots i_k)} _ {{l_1 \ldots l_k}_{s, t}}$
similarly to the
previous section or using the iterated stochastic integrals 
$J^{*(i_1\ldots i_k)} _ {{l_1 \ldots l_k}_{s, t}}.$
This is the main subject of this section.

Denote

\vspace{-2mm}
$$
\Biggl\Vert 
J^{*(i_1\ldots i_k)} _ {{l_1 \ldots l_k}_{s, t}}\Biggr\Vert
_{i_1,\ldots,i_k=1}^{m}\stackrel{\sf def}
{=}{}^{(k)}J^{*}_{{l_1\ldots l_k}_{s,t}},
$$

\vspace{3mm}

$$
\Biggl\Vert 
\bar L^j \bar G_{l_1}^{(i_1)}\ldots \bar G_{l_k}^{(i_k)}
R({\bf x},t) \Biggr\Vert
_{i_1,\ldots,i_k=1}^
{m}\stackrel{\sf def}{=}{}^{(k)}
\bar L^j \bar G_{l_1}\ldots \bar G_{l_k}
R({\bf x},t).
$$

\vspace{5mm}

{\bf Theorem 5.}\ {\it Suppose that sufficient conditions are 
satisfied under which the right-hand 
sides of {\rm(\ref{5.7.11xxx}), (\ref{5.7.12xxx})} exist. 
Then for any $s, t \in [0, T]$ such that $s>t$ 
and for any positive
integer $r$, the following expansion takes place w.\ p.\ {\rm 1}
$$
R({\bf x}_s,s)=
R({\bf x}_t,t)+
\sum_{q=1}^r
\sum_{(k,j,l_1,\ldots,l_k) \in {\rm A}_q}\ 
\frac{(s-t)^j}{j!}
\sum_{i_1,\ldots,i_k=1}^m \bar L^j\bar G_{l_1}^{(i_1)}\ldots
\bar G_{l_k}^{(i_k)}R({\bf x}_t,t)\ 
J_{{l_1\ldots l_k}_{s,t}}^{*(i_1\ldots i_k)}+
$$

\vspace{1mm}
\begin{equation}
\label{razl44x}
+\left(D_{r+1}\right)_{s,t},
\end{equation}

\vspace{5mm}
\noindent
where $\left(D_{r+1}\right)_{s,t}$ has the form {\rm (\ref{5.7.12xxx})}.
}

\vspace{2mm}

{\bf Proof.} To prove the theorem, we check the equalities

\vspace{2mm}
$$
\sum_{(k,j,l_1,\ldots,l_k) \in {\rm A}_q}\ 
\frac{(s-t)^j}{j!}
\sum_{i_1,\ldots,i_k=1}^m \bar L^j\bar G_{l_1}^{(i_1)}\ldots
\bar G_{l_k}^{(i_k)}R({\bf x}_t,t)\ 
J_{{l_1\ldots l_k}_{s,t}}^{*(i_1\ldots i_k)}=
$$

\vspace{1mm}
\begin{equation}
\label{f11x}
\sum_{(k,j,l_1,\ldots,l_k) \in {\rm A}_q}\ 
\frac{(s-t)^j}{j!}
\sum_{i_1,\ldots,i_k=1}^m \bar G_{l_1}^{(i_1)}\ldots
\bar G_{l_k}^{(i_k)}\bar L^j R({\bf x}_t,t)\ 
I_{{l_1\ldots l_k}_{s,t}}^{*(i_1\ldots i_k)}\ \ \ \hbox{w.\ p.\ 1}
\end{equation}

\vspace{6mm}
\noindent
for $q=1,2,\ldots,r.$ 
To check (\ref{f11x}), substitute the expression 
(\ref{70x}) 
into the right-hand side of (\ref{f11x}) and then use the formulas
(\ref{a9x}), (\ref{a33x}).

Let us rank terms of the expansion (\ref{razl44x}) according to their orders of 
smallness in the mean-square sense as $s \downarrow t$

\vspace{-2mm}
$$
R({\bf x}_s,s)=
R({\bf x}_t,t)+
\sum_{q=1}^r
\sum_{(k,j,l_1,\ldots,l_k) \in {\rm D}_q}\ 
\frac{(s-t)^j}{j!}
\sum_{i_1,\ldots,i_k=1}^m \bar L^j \bar G_{l_1}^{(i_1)}\ldots
\bar G_{l_k}^{(i_k)}R({\bf x}_t,t)\ 
J_{{l_1\ldots l_k}_{s,t}}^{*(i_1\ldots i_k)}+
$$

\vspace{2mm}
$$
+\left(H_{r+1}\right)_{s,t}\ \ \ \hbox{w.\ p.\ 1},
$$

\vspace{5mm}
\noindent
where

\vspace{-1mm}
$$
\left(H_{r+1}\right)_{s,t}=
\sum_{(k,j,l_1,\ldots,l_k) \in {\rm U}_r}\ 
\frac{(s-t)^j}{j!}
\sum_{i_1,\ldots,i_k=1}^m \bar L^j\bar G_{l_1}^{(i_1)}\ldots
\bar G_{l_k}^{(i_k)}R({\bf x}_t,t)\ 
J_{{l_1\ldots l_k}_{s,t}}^{*(i_1\ldots i_k)}+
\left(D_{r+1}\right)_{s,t}.
$$

\vspace{5mm}

The term $\left(D_{r+1}\right)_{s,t}$ 
has the form (\ref{5.7.12xxx}); the terms ${\rm D}_q$ and ${\rm U}_r$ have the 
forms (\ref{w1x}) and (\ref{w2x}), respectively.
Finally, we note that the convergence w.\ p.\ 1 of the 
truncated Taylor--Stratonovich expansion (\ref{5.7.11xxx}) (without
the remainder term $\left(D_{r+1}\right)_{s,t}$) to the process 
$R({\bf x}_s, s)$ as $r\to\infty$  
for all $s, t \in [0, T]$ such that $s>t$ and $T < \infty$
has been proved in \cite{KlPl2} (Proposition 5.10.2).
Since the expansions (\ref{razl4x}) and 
(\ref{razl44x}) are obtained from the Taylor--Stratonovich 
expansion (\ref{5.7.11xxx}) without any additional
conditions, the truncated expansions (\ref{razl4x}) and 
(\ref{razl44x}) (without the 
reminder term $\left(D_{r+1}\right)_{s,t}$) under the
conditions of Proposition 5.10.2 \cite{KlPl2} converge to the 
process $R({\bf x}_s, s)$ w.\ p.\ 1 as $r\to\infty$ for all
$s, t \in [0, T]$ such that $s>t$ and $T < \infty.$

\vspace{5mm}

\section{A Remark on Theorems~4 and 5}

\vspace{5mm}

Note that when proving Theorems~4 and 5
we established the following equalities w.~p.~1

\vspace{1mm}
$$
R({\bf x}_t,t)+
\sum_{k=1}^r \sum_{(\lambda_{k},\ldots,\lambda_1)
\in{M}_k}
\sum_{i_1=\lambda_1}^{m\lambda_1}
\ldots 
\sum_{i_k=\lambda_k}^{m\lambda_k}
D_{\lambda_k}^{(i_k)}\ldots D_{\lambda_1}^{(i_1)}
R({\bf x}_t,t){J}_{(\lambda_{k}\ldots \lambda_1)s,t}^{*(i_k\ldots i_1)}
=
$$

\vspace{2mm}
\begin{equation}
\label{2024str202}
=R({\bf x}_t,t)+
\sum_{q=1}^r
\sum_{(k,j,l_1,\ldots,l_k) \in {\rm A}_q}
\frac{(s-t)^j}{j!}
\sum_{i_1,\ldots,i_k=1}^m \bar G_{l_1}^{(i_1)}\ldots
\bar G_{l_k}^{(i_k)}\bar L^j R({\bf x}_t,t)
I_{{l_1\ldots l_k}_{s,t}}^{*(i_1\ldots i_k)},
\end{equation}

\vspace{4mm}
$$
R({\bf x}_t,t)+
\sum_{k=1}^r \sum_{(\lambda_{k},\ldots,\lambda_1)
\in{M}_k}
\sum_{i_1=\lambda_1}^{m\lambda_1}
\ldots 
\sum_{i_k=\lambda_k}^{m\lambda_k}
D_{\lambda_k}^{(i_k)}\ldots D_{\lambda_1}^{(i_1)}
R({\bf x}_t,t){J}_{(\lambda_{k}\ldots \lambda_1)s,t}^{*(i_k\ldots i_1)}
=
$$

\vspace{3mm}
\begin{equation}
\label{2024str203}
=R({\bf x}_t,t)+
\sum_{q=1}^r
\sum_{(k,j,l_1,\ldots,l_k) \in {\rm A}_q}
\frac{(s-t)^j}{j!}
\sum_{i_1,\ldots,i_k=1}^m \bar L^j\bar G_{l_1}^{(i_1)}\ldots
\bar G_{l_k}^{(i_k)}R({\bf x}_t,t)
J_{{l_1\ldots l_k}_{s,t}}^{*(i_1\ldots i_k)}.
\end{equation}

\vspace{5mm}

It is easy to see that by analogy with 
(\ref{2024str202}) and (\ref{2024str203})
the following equalities can be obtained w.~p.~1

\vspace{1mm}
$$
R({\bf x}_t,t)+
\sum_{q,k=1}^r \sum_{(\lambda_{k},\ldots,\lambda_1)
\in{\rm E}_{qk}}
\sum_{i_1=\lambda_1}^{m\lambda_1}
\ldots 
\sum_{i_k=\lambda_k}^{m\lambda_k}
D_{\lambda_k}^{(i_k)}\ldots D_{\lambda_1}^{(i_1)}
R({\bf x}_t,t)
{J}_{(\lambda_{k}\ldots \lambda_1)s,t}^{*(i_k\ldots
i_1)}=
$$

\vspace{3mm}
\begin{equation}
\label{2024str204}
=R({\bf x}_t,t)+
\sum_{q=1}^r
\sum_{(k,j,l_1,\ldots,l_k) \in {\rm D}_q}
\frac{(s-t)^j}{j!}
\sum_{i_1,\ldots,i_k=1}^m \bar G_{l_1}^{(i_1)}\ldots
\bar G_{l_k}^{(i_k)}\bar L^j R({\bf x}_t,t)
I_{{l_1\ldots l_k}_{s,t}}^{*(i_1\ldots i_k)},
\end{equation}

\vspace{5mm}
$$
R({\bf x}_t,t)+
\sum_{q,k=1}^r \sum_{(\lambda_{k},\ldots,\lambda_1)
\in{\rm E}_{qk}}
\sum_{i_1=\lambda_1}^{m\lambda_1}
\ldots 
\sum_{i_k=\lambda_k}^{m\lambda_k}
D_{\lambda_k}^{(i_k)}\ldots D_{\lambda_1}^{(i_1)}
R({\bf x}_t,t)
{J}_{(\lambda_{k}\ldots \lambda_1)s,t}^{*(i_k\ldots
i_1)}=
$$

\vspace{3mm}
\begin{equation}
\label{2024str205}
=R({\bf x}_t,t)+
\sum_{q=1}^r
\sum_{(k,j,l_1,\ldots,l_k) \in {\rm D}_q}\ 
\frac{(s-t)^j}{j!}
\sum_{i_1,\ldots,i_k=1}^m \bar L^j \bar G_{l_1}^{(i_1)}\ldots
\bar G_{l_k}^{(i_k)}R({\bf x}_t,t)
J_{{l_1\ldots l_k}_{s,t}}^{*(i_1\ldots i_k)}.
\end{equation}

\vspace{5mm}

Recall the equality (\ref{2024str201})

\vspace{1mm}
$$
R({\bf x}_t,t)+
\sum_{q,k=1}^r \sum_{(\lambda_{k},\ldots,\lambda_1)
\in{\rm E}_{qk}}
\sum_{i_1=\lambda_1}^{m\lambda_1}
\ldots 
\sum_{i_k=\lambda_k}^{m\lambda_k}
Q_{\lambda_k}^{(i_k)}\ldots Q_{\lambda_1}^{(i_1)}
R({\bf x}_t,t){J}_{(\lambda_{k}\ldots \lambda_1)s,t}^{(i_k\ldots
i_1)}=
$$

\vspace{3mm}
\begin{equation}
\label{2024str206}
=R({\bf x}_t,t)+
\sum_{q,k=1}^r \sum_{(\lambda_{k},\ldots,\lambda_1)
\in{\rm E}_{qk}}
\sum_{i_1=\lambda_1}^{m\lambda_1}
\ldots 
\sum_{i_k=\lambda_k}^{m\lambda_k}
D_{\lambda_k}^{(i_k)}\ldots D_{\lambda_1}^{(i_1)}
R({\bf x}_t,t)
{J}_{(\lambda_{k}\ldots \lambda_1)s,t}^{*(i_k\ldots
i_1)}
\end{equation}

\vspace{5mm}
\noindent
w.~p.~1, where $r\ge 2$.

Combining (\ref{2024str204})--(\ref{2024str206}), we obtain

\vspace{1mm}
$$
R({\bf x}_t,t)+
\sum_{q,k=1}^r \sum_{(\lambda_{k},\ldots,\lambda_1)
\in{\rm E}_{qk}}
\sum_{i_1=\lambda_1}^{m\lambda_1}
\ldots 
\sum_{i_k=\lambda_k}^{m\lambda_k}
Q_{\lambda_k}^{(i_k)}\ldots Q_{\lambda_1}^{(i_1)}
R({\bf x}_t,t){J}_{(\lambda_{k}\ldots \lambda_1)s,t}^{(i_k\ldots
i_1)}=
$$

\vspace{3mm}
$$
=R({\bf x}_t,t)+
\sum_{q=1}^r
\sum_{(k,j,l_1,\ldots,l_k) \in {\rm D}_q}
\frac{(s-t)^j}{j!}
\sum_{i_1,\ldots,i_k=1}^m \bar G_{l_1}^{(i_1)}\ldots
\bar G_{l_k}^{(i_k)}\bar L^j R({\bf x}_t,t)
I_{{l_1\ldots l_k}_{s,t}}^{*(i_1\ldots i_k)}=
$$

\vspace{3mm}
\begin{equation}
\label{2024str207}
=R({\bf x}_t,t)+
\sum_{q=1}^r
\sum_{(k,j,l_1,\ldots,l_k) \in {\rm D}_q}\ 
\frac{(s-t)^j}{j!}
\sum_{i_1,\ldots,i_k=1}^m \bar L^j \bar G_{l_1}^{(i_1)}\ldots
\bar G_{l_k}^{(i_k)}R({\bf x}_t,t)
J_{{l_1\ldots l_k}_{s,t}}^{*(i_1\ldots i_k)}
\end{equation}

\vspace{5mm}
\noindent
w.~p.~1, where $r\ge 2$.

The equality (\ref{2024str207}) means that we have the following theorem.

\vspace{2mm}

{\bf Theorem 6}\ \cite{20aa}.\ {\it Let conditions {\rm (i), (ii) (}see
Sect.~{\rm 3)}
be satisfied. Then for any $s, t \in [0, T]$ such that $s>t$ 
the following unified 
Taylor--Stratonovich expansions take place w.~p.~{\rm 1}

\vspace{1mm}
$$
R({\bf x}_s,s)=R({\bf x}_t,t)+
$$

\vspace{2mm}
$$
+
\sum_{q=1}^r
\sum_{(k,j,l_1,\ldots,l_k) \in {\rm D}_q}
\frac{(s-t)^j}{j!}
\sum_{i_1,\ldots,i_k=1}^m \bar G_{l_1}^{(i_1)}\ldots
\bar G_{l_k}^{(i_k)}\bar L^j R({\bf x}_t,t)
I_{{l_1\ldots l_k}_{s,t}}^{*(i_1\ldots i_k)}+\left(H_{r+1}\right)_{s,t},
$$

\vspace{5mm}
$$
R({\bf x}_s,s)=R({\bf x}_t,t)+
$$

\vspace{2mm}
$$
+
\sum_{q=1}^r
\sum_{(k,j,l_1,\ldots,l_k) \in {\rm D}_q}\ 
\frac{(s-t)^j}{j!}
\sum_{i_1,\ldots,i_k=1}^m \bar L^j \bar G_{l_1}^{(i_1)}\ldots
\bar G_{l_k}^{(i_k)}R({\bf x}_t,t)
J_{{l_1\ldots l_k}_{s,t}}^{*(i_1\ldots i_k)}
+\left(H_{r+1}\right)_{s,t},
$$

\vspace{5mm}
\noindent
where $r\ge 2,$ the reainder term $\left(H_{r+1}\right)_{s,t}$
is definded by the relations {\rm(\ref{5.6.1rrrh})} and {\rm (\ref{5.7.12});}
another notations are the same as in the previous sections.}

\vspace{5mm}

\section{Comparison of the Unified Taylor--Ito and
Taylor--Stra\-to\-no\-vich Expansions
With the Classical Taylor--Ito and Taylor--Stratonovich Expansions}

\vspace{5mm}

Note that the 
truncated 
unified 
Taylor--Ito and Taylor--Stratonovich expansions contain 
the less number of various 
iterated Ito and Stratonovich stochastic integrals (moreover, their major part 
will have 
less multiplicity) in comparison with 
the classical Taylor--Ito and
Taylor--Stratonovich expansions \cite{KlPl1}.

It is easy to notice that the stochastic integrals from the families 
(\ref{re11}), (\ref{str11}) are connected by linear relations. 
However, the stochastic integrals from the families 
(\ref{ll1}), (\ref{ll11}) cannot be connected by linear relations. 
This holds for the 
stochastic integrals from the families
(\ref{a111}), (\ref{a112}).
Therefore, we will call the families 
(\ref{ll1})--(\ref{a112}) as the {\it stochastic 
bases.}

Let us 
call
the numbers ${\rm rank}_{\rm A}(r)$ and
${\rm rank}_{\rm D}(r)$ of various iterated 
Ito and Stratonovich stochastic integrals which are included in the 
families (\ref{ll1})--(\ref{a112}) as the 
{\it ranks of stochastic bases} when 
summation in the stochastic expansions is performed using the 
sets 
${\rm A} _ {q}$ ($q=1,\ldots,r)$
and ${\rm D} _ {q}$ ($q=1,\ldots,r$) correspondently.
Here $r$ is a fixed natural number.

At the beginning, let us analyze several examples
related to the Taylor--Ito expansions (obviously,
the same conclusions will hold for the Taylor--Stratonovich expansions).

Assume that summation in the unified  
Taylor--Ito expansions is performed using 
the sets ${\rm D} _ {q}$ ($q=1,\ldots,r$).
It is easy to see that the 
truncated 
unified Taylor--Ito expansion (\ref{t100}), 
where summation is performed with respect to the sets ${\rm D}_q$ when $r=3$  
includes 4 (${\rm rank}_{\rm D}(3)=4$) various iterated Ito
stochastic 
integrals

\vspace{-2mm}
$$
I_{0_{s,t}}^{(i_1)},\ \ \ I_{{00}_{s,t}}^{(i_1 i_2)},\ \ \
I_{1_{s,t}}^{(i_1)},\ \ \
I_{{000}_{s,t}}^{(i_1 i_2 i_3)}.
$$

\vspace{3mm}

The same truncated 
classical
Taylor--Ito expansion (\ref{5.6.1rrr}) \cite{KlPl2} contains 
5 various iterated Ito stochastic integrals

\vspace{-2mm}
$$
{J}_{(1)s,t}^{(i_1)},\ \ \ {J}_{(11)s,t}^{(i_1 i_2)},\ \ \ 
{J}_{(10)s,t}^{(i_1 0)},\ \ \ {J}_{(01)s,t}^{(0 i_1)},\ \ \
{J}_{(111)s,t}^{(i_1 i_2 i_3)}.
$$

\vspace{4mm}

For $r=4$ we have
7 (${\rm rank}_{\rm D}(4)=7$) integrals

$$
I_{0_{s,t}}^{(i_1)},\ \ \ I_{{00}_{s,t}}^{(i_1 i_2)},\ \ \
I_{1_{s,t}}^{(i_1)},\ \ \
I_{{000}_{s,t}}^{(i_1 i_2 i_3)},\ \ \
I_{{01}_{s,t}}^{(i_1 i_2)},\ \ \ I_{{10}_{s,t}}^{(i_1 i_2)},\ \ \
I_{{0000}_{s,t}}^{(i_1 i_2 i_3 i_4)}
$$

\vspace{4mm}
\noindent
against 9 stochastic integrals

$$
{J}_{(1)s,t}^{(i_1)},\ \ \ {J}_{(11)s,t}^{(i_1 i_2)},\ \ \
{J}_{(10)s,t}^{(i_1 0)},\ \ \ {J}_{(01)s,t}^{(0 i_1)},\ \ \
{J}_{(111)s,t}^{(i_1 i_2 i_3)},
$$

\vspace{2mm}
$$
{J}_{(101)s,t}^{(i_1 0 i_3)},\ \ \
{J}_{(110)s,t}^{(i_1 i_2 0)},\ \ \ {J}_{(011)s,t}^{(0 i_1 i_2)},\ \ \
{J}_{(1111)s,t}^{(i_1 i_2 i_3 i_4)}.
$$

\vspace{4mm}

For $r=5$ (${\rm rank}_{\rm D}(5)=12$) we get 
12 integrals against 17 integrals
and for 
$r=6$ and $r=7$ we have
20 against 29 and 33 against 50 
correspondently.

We will obtain the same results when compare the unified Taylor--Stratonovich 
expansions \cite{k1}, \cite{7}-\cite{11}, \cite{15}-\cite{12aa-afterxxx}
with their classical analogues \cite{KlPl2}, \cite{KlPl1}
(see previous sections).

Note that summation according to the sets ${\rm D}_q$ is usually used 
while constructing strong numerical methods (built according to 
the mean-square criterion of convergence) for Ito SDEs
\cite{KlPl2}, \cite{Mi2}, \cite{7}-\cite{11}, \cite{19}-\cite{12aa-afterxxx}.
Summation according to the sets ${\rm A}_q$ is usually used when building 
weak numerical methods (built in accordance with the weak 
criterion of convergence) for Ito SDEs
\cite{KlPl2}, \cite{Mi2}.
For example, ${\rm rank}_{\rm A}(4)=15$, while the total number of various 
iterated Ito stochastic integrals (included in the 
classical Taylor--Ito 
expansion \cite{KlPl2} when $r=4$) equals to 26.

Let us show that \cite{9}-\cite{11}, \cite{19}-\cite{12aa-afterxxx}

$$
{\rm rank}_{\rm A}(r)=2^r-1.
$$ 

\vspace{3mm}

Let $(l_1,\ldots,l_k)$ be an ordered set such that 
$l_1,\ldots,$ $l_k$ $=$ $0, 1,\ldots$ and
$k=1, 2,\ldots $ Consider $S(k)\stackrel{\sf def}{=}l_1+\ldots+l_k=p$
($p$ is a fixed natural number or zero).
Let $N(k,p)$ be a number of all ordered
combinations
$(l_1,\ldots,l_k)$ such that $l_1,\ldots,l_k=0, 1,\ldots,$
$k=1, 2,\ldots,$ and 
$S(k)=p$. First let us show that

\vspace{-2mm}
$$
N(k,p)=C_{p+k-1}^{k-1},
$$

\vspace{4mm}
\noindent
where 
$$
C_n^m=\frac{n!}{m!(n-m)!}
$$ 

\vspace{3mm}
\noindent
is a binomial coefficient.

It is not difficult to see that

$$
N(1,p)=1=C_{p+1-1}^{1-1},
$$

\vspace{2mm}
$$
N(2,p)=p+1=C_{p+2-1}^{2-1},
$$

\vspace{2mm}
$$
N(3,p)=\frac{(p+1)(p+2)}{2}=C_{p+3-1}^{3-1}.
$$ 

\vspace{3mm}

Moreover,

\vspace{-2mm}
$$
N(k+1,p)=\sum\limits_{l=0}^p N(k,l)=\sum_{l=0}^p C_{l+k-1}^{k-1}=C_{p+k}^{k},
$$

\vspace{3mm}
\noindent
where we used the induction assumption and the well known 
property of binomial coefficients.

Then

\vspace{-2mm}
$$
{\rm rank}_{A}(r)=
$$

$$
=N(1,0)+(N(1,1)+N(2,0))+
(N(1,2)+N(2,1)+N(3,0))+\ldots
$$

$$
\ldots
+(N(1,r-1)+N(2,r-2)+\ldots+N(r,0))=
$$

$$
=C_0^0+(C_1^0+C_1^1)+(C_2^0+C_2^1+C_2^2)+ \ldots
$$

$$
\ldots +(C_{r-1}^0+C_{r-1}^1+C_{r-1}^2+\ldots+C_{r-1}^{r-1})=
$$

$$
=2^0+2^1+2^2+\ldots+2^{r-1}=2^r-1.
$$ 

\vspace{4mm}

Let $n_{{\rm M}}(r)$ be the total number of various iterated stochastic 
integrals included in the 
classical
Taylor--Ito expansion (\ref{5.7.11}) \cite{KlPl2},
where summation is performed with respect to the set  
$$
\bigcup\limits_{k=1}^r{\rm M}_k.
$$

If we exclude from the consideration 
the integrals which are equal to
$(s-t)^j/j!,$ then 

\vspace{-1mm}
$$
n_{{M}}(r)=
$$

$$
=(2^1-1)+(2^2-1)+(2^3-1)+\ldots+(2^r-1)=
$$

$$
=2(1+2+2^2+\ldots+2^{r-1})-r=2(2^r-1)-r.
$$

\vspace{3mm}

It means that
$$
\lim\limits_{r\to\infty}
\frac{n_{{M}}(r)}{{\rm rank}_{A}(r)}=2.
$$

\vspace{4mm}

In Table 1 we can see the numbers

\vspace{-2mm}
$$
{\rm rank}_{\rm A}(r),\ \ \ n_{{\rm M}}(r),\ \ \ 
f(r)=n_{{\rm M}}(r)/{\rm rank}_{\rm A}(r)
$$

\vspace{2mm}
\noindent
for various values $r$.

Let us show that \cite{9}-\cite{11}, \cite{19}-\cite{12aa-afterxxx}

\vspace{1mm}
\begin{equation}
\label{gg1aa}
{\rm rank}_{\rm D}(r)=
\begin{cases}
\sum\limits_{s=0}^{r-1}~~
\sum\limits_{l=s}^{(r-1)/2+[s/2]}~ C_l^s\ \ \  &\hbox{for}\ \ \ 
r=1,\ 3,\ 5,\ldots\cr\cr\cr
\sum\limits_{s=0}^{r-1}~~ \sum\limits_{l=s}^{
r/2-1+[(s+1)/2]}~ C_l^s\ \ \  &\hbox{for}\ \ \ 
r=2,\ 4,\ 6,\ldots
\end{cases},
\end{equation}

\vspace{6mm}
\noindent
where $[x]$ is an integer part of a number $x,$
and $C_n^m$ is a binomial coefficient.

For proving (\ref{gg1aa}) we write
the condition

\vspace{-2mm}
$$
k + 2(j + S(k))\le r,
$$ 

\vspace{2mm}
\noindent
where 
$S(k)\stackrel{\sf def}{=}l_1+\ldots+l_k$
($k, j, l_1,\ldots, l_k =0,1,\ldots$)
in the form
$j+S(k)\le (r-k)/2,$ 
and perform the consideration of all possible combinations
with respect to
$k=1,\ldots,r$. Moreover, we take into account the above 
reasoning.

Let us calculate the number
$n_{{E}}(r)$ of all different iterated Ito
stochastic integrals from the classical Taylor--Ito
expansion (\ref{5.6.1rrr}) \cite{KlPl2} 
if the summation in this expansion is performed
with respect to the set
$$
\bigcup\limits_{q,k=1}^r{E}_{qk}.
$$

\vspace{1mm}

The summation condition can be written in this case
in the form: $0\le p+2q\le r,$ where $q$ 
is a total number of integrations
with respect to time while
$p$ is a total number of integrations
with respect to the Wiener processes in the selected iterated
stochastic integral from the Taylor--Ito expansion 
(\ref{5.6.1rrr}) \cite{KlPl2}.
At that, the multiplicity of the mentioned
stochastic integral equals to $p+q$ and it is not more than
$r.$ Let us write the above condition ($0\le p+2q\le r$)
in the form:
$0\le q\le (r-p)/2$ $\Leftrightarrow$ $0\le q\le [(r-p)/2]$,
where $[x]$ means an integer part of a real number $x$.
Then, performing the consideration of all possible combinations
with respect to $p=1,\ldots,r$ and using the
combinatorial reasoning, we obtain
the formula

\vspace{-2mm}
\begin{equation}
\label{123u}
n_{{E}}(r)=\sum\limits_{s=1}^r~~ \sum\limits_{l=0}^{[(r-s)/2]}~
C_{[(r-s)/2]+s-l}^s,
\end{equation}

\vspace{3mm}
\noindent
where $[x]$ means an integer part of a real number $x$.

\begin{table}
\centering
\caption{Numbers ${\rm rank}_{\rm A}(r),$ $n_{{\rm M}}(r),$ 
$f(r)=n_{{\rm M}}(r)/{\rm rank}_{\rm A}(r)$}
\label{tab:100}      
\begin{tabular}{p{1.3cm}p{0.2cm}p{0.9cm}p{0.9cm}p{0.9cm}p{0.9cm}p{1cm}p{1cm}p{1cm}p{1cm}p{1cm}}
\hline\noalign{\smallskip}
$ r$ & 1 & 2& 3& 4 & 5 & 6 & 7 & 8 & 9 & 10\\
\noalign{\smallskip}\hline\noalign{\smallskip}
${\rm rank}_{\rm A}(r)$ &  1 & 3 & 7& 15 & 31 & 63& 127 & 255& 511& 1023\\
${n}_{\rm M}(r)$ &  1 & 4 & 11& 26 & 57 & 120& 247 & 502 & 1013& 2036\\
$f(r)$ & 1& 1.3333& 1.5714& 1.7333 & 1.8387 & 1.9048 & 1.9449 & 1.9686 & 
1.9824 & 1.9902\\
\noalign{\smallskip}\hline\noalign{\smallskip}
\end{tabular}
\vspace{7mm}
\end{table}

\begin{table}
\centering
\caption{Numbers ${\rm rank}_{\rm D}(r),$ $n_{{\rm E}}(r),$ 
$g(r)=n_{{\rm E}}(r)/{\rm rank}_{\rm D}(r)$}
\label{tab:101}      
\begin{tabular}{p{1.3cm}p{0.2cm}p{0.3cm}p{0.9cm}p{1cm}p{1cm}p{1cm}p{1cm}p{1cm}p{1cm}p{1cm}}
\hline\noalign{\smallskip}
$ r$ & 1 & 2& 3& 4 & 5 & 6 & 7 & 8 & 9 & 10\\
\noalign{\smallskip}\hline\noalign{\smallskip}
${\rm rank}_{\rm D}(r)$ &  1 & 2 & 4& 7 & 12 & 20& 33 & 54& 88 & 143\\
${n}_{\rm E}(r)$ &  1 & 2 & 5& 9 & 17 & 29& 50 & 83 & 138 & 261\\
$g(r)$ & 1& 1& 1.2500& 1.2857 & 1.4167 & 1.4500 & 1.5152 & 1.5370 & 
1.5682 & 1.8252\\
\noalign{\smallskip}\hline\noalign{\smallskip}
\end{tabular}
\vspace{15mm}
\end{table}

In Table 2 we can see the numbers

\vspace{-2mm}
$$
{\rm rank}_{\rm D}(r),\ \ \ n_{{\rm E}}(r),\ \ \
g(r)=n_{{\rm E}}(r)/{\rm rank}_{\rm D}(r)
$$

\vspace{3mm}
\noindent
for various values $r$.

\vspace{5mm}

\section{Application of First Form of the Unified 
Taylor--Ito Expansion to the High-Order
Strong Numerical Methods for Ito SDEs}

\vspace{5mm}

Let us write (\ref{t100}) for all $s, t\in[0,T]$ such that $s>t$
in the following form

\vspace{1mm}
$$
R({\bf x}_s,s)=R({\bf x}_t,t)+\sum_{q=1}^r\sum_{(k,j,l_1,\ldots,l_k)\in
{\rm D}_q}
\frac{(s-t)^j}{j!} \sum\limits_{i_1,\ldots,i_k=1}^m
G_{l_1}^{(i_1)}\ldots G_{l_k}^{(i_k)}
L^j R({\bf x}_t,t)\
I^{(i_1\ldots i_k)}_{{l_1\ldots l_k}_{s,t}}+
$$

\vspace{2mm}
\begin{equation}
\label{15.001}
+{\bf 1}_{\{r=2d-1,
d\in N\}}\frac{(s-t)^{(r+1)/2}}{\left(
(r+1)/2\right)!}L^{(r+1)/2}R({\bf x}_t,t)
+\left(\bar H_{r+1}\right)_{s,t}\ \ \ \hbox{w.\ p.\ 1},
\end{equation}

\vspace{5mm}
\noindent
where

\vspace{-3mm}
$$
\left(\bar H_{r+1}\right)_{s,t}= \left(H_{r+1}\right)_{s,t}
-{\bf 1}_{\{r=2d-1,
d\in N\}}\frac{(s-t)^{(r+1)/2}}{\left(
(r+1)/2\right)!}L^{(r+1)/2}R({\bf x}_t,t).
$$

\vspace{6mm}

Consider the partition 
$\{\tau_p\}_{p=0}^N$ of the interval
$[0,T]$ such that

\vspace{-1mm}
$$
0=\tau_0<\tau_1<\ldots<\tau_N=T,\ \ \
\Delta_N=
\max\limits_{0\le j\le N-1}\left|\tau_{j+1}-\tau_j\right|.
$$

\vspace{3mm}

From (\ref{15.001}) for $s=\tau_{p+1},$
$t=\tau_p$ we obtain the following representation of
explicit one-step strong numerical scheme for Ito SDE (\ref{1.5.2}),
which is based on the first form of the unified Taylor--Ito
expansion

\vspace{1mm}

$$
{\bf y}_{p+1}={\bf y}_{p}+
\sum_{q=1}^r\sum_{(k,j,l_1,\ldots,l_k)\in{\rm D}_q}
\frac{(\tau_{p+1}-\tau_p)^j}{j!} \sum\limits_{i_1,\ldots,i_k=1}^m
G_{l_1}^{(i_1)}\ldots G_{l_k}^{(i_k)}
L^j\hspace{0.4mm}{\bf y}_{p}\
{\hat I}^{(i_1\ldots i_k)}_{{l_1\ldots l_k}_{\tau_{p+1},\tau_p}}+
$$

\vspace{2mm}
\begin{equation}
\label{15.002}
+{\bf 1}_{\{r=2d-1,
d\in N\}}\frac{(\tau_{p+1}-\tau_p)^{(r+1)/2}}{\left(
(r+1)/2\right)!}L^{(r+1)/2}\hspace{0.4mm}{\bf y}_{p},
\end{equation}

\vspace{5mm}
\noindent
where
$\hat I^{(i_1\ldots i_k)}_{{l_1\ldots l_k}_{\tau_{p+1},\tau_p}}$ 
is an approximation of the iterated Ito
stochastic integral 
$I^{(i_1\ldots i_k)}_{{l_1\ldots l_k}_{\tau_{p+1},\tau_p}}$
defined as
$$
I_{{l_1\ldots l_k}_{s,t}}^{(i_1\ldots i_k)}=
\int\limits_t^s
(t-t_{k})^{l_{k}}\ldots 
\int\limits_t^{t_2}
(t-t _ {1}) ^ {l_ {1}} d
{\bf f} ^ {(i_ {1})} _ {t_ {1}} \ldots 
d {\bf f} _ {t_ {k}} ^ {(i_ {k})}.
$$

\vspace{2mm}

Note that we understand the equality (\ref{15.002}) componentwise
with respect to the components ${\bf y}_p^{(i)}$ of the column
${\bf y}_p.$
Also for simplicity we put 
$\tau_p=p\Delta$,
$\Delta=T/N,$ $T=\tau_N,$ $p=0,1,\ldots,N.$

It is known \cite{KlPl2} that under the appropriate conditions
the numerical scheme (\ref{15.002}) has strong order of convergence $r/2$
($r\in\mathbb{N}$).

Let $B_j({\bf x},t)$ is the $j$-th column of the matrix
function $B({\bf x},t)$ (see (\ref{1.5.2})).

Below we consider particular cases of the numerical scheme
(\ref{15.002}) for $r=2,3,4,5,$ and $6,$ i.e. 
explicit one-step strong numerical schemes for Ito SDE (\ref{1.5.2})
with orders $1.0, 1.5, 2.0, 2.5,$ and $3.0$ of convergence.
At that, for simplicity 
we will write ${\bf a},$ $L {\bf a},$ 
$B_i,$ $G_0^{(i)}B_{j}$ etc.
instead of ${\bf a}({\bf y}_p,\tau_p),$ 
$L {\bf a}({\bf y}_p,\tau_p),$ 
$B_i({\bf y}_p,\tau_p),$ 
$G_0^{(i)}B_{j}({\bf y}_p,\tau_p)$ etc. correspondingly.
Moreover, the operators $L$ and $G_0^{(i)},$ $i=1,\ldots,m,$
are determined by the equalities
(\ref{2.3}), (\ref{2.4})  as before.

\vspace{8mm}

\centerline{\bf Scheme with strong order 1.0}

\vspace{2mm}
\begin{equation}
\label{al1}
{\bf y}_{p+1}={\bf y}_{p}+\sum_{i_{1}=1}^{m}B_{i_{1}}
\hat I_{0_{\tau_{p+1},\tau_p}}^{(i_{1})}+\Delta{\bf a}
+\sum_{i_{1},i_{2}=1}^{m}G_0^{(i_{2})}
B_{i_{1}}\hat I_{00_{\tau_{p+1},\tau_p}}^{(i_{2}i_{1})}.
\end{equation}

\vspace{10mm}

\centerline{\bf Scheme with strong order 1.5}

\vspace{5mm}
$$
{\bf y}_{p+1}={\bf y}_{p}+\sum_{i_{1}=1}^{m}B_{i_{1}}
\hat I_{0_{\tau_{p+1},\tau_p}}^{(i_{1})}+\Delta{\bf a}
+\sum_{i_{1},i_{2}=1}^{m}G_0^{(i_{2})}
B_{i_{1}}\hat I_{00_{\tau_{p+1},\tau_p}}^{(i_{2}i_{1})}+
$$

\vspace{2mm}
$$
+
\sum_{i_{1}=1}^{m}\left[G_0^{(i_{1})}{\bf a}\left(
\Delta \hat I_{0_{\tau_{p+1},\tau_p}}^{(i_{1})}+
\hat I_{1_{\tau_{p+1},\tau_p}}^{(i_{1})}\right)
- LB_{i_{1}}\hat I_{1_{\tau_{p+1},\tau_p}}^{(i_{1})}\right]+
$$

\vspace{2mm}
$$
+\sum_{i_{1},i_{2},i_{3}=1}^{m} G_0^{(i_{3})}G_0^{(i_{2})}
B_{i_{1}}\hat I_{000_{\tau_{p+1},\tau_p}}^{(i_{3}i_{2}i_{1})}+
$$

\vspace{2mm}
\begin{equation}
\label{al2}
+
\frac{\Delta^2}{2}L{\bf a}.
\end{equation}

\vspace{10mm}

\centerline{\bf Scheme with strong order 2.0}

\vspace{4mm}
$$
{\bf y}_{p+1}={\bf y}_{p}+\sum_{i_{1}=1}^{m}B_{i_{1}}
\hat I_{0_{\tau_{p+1},\tau_p}}^{(i_{1})}+\Delta{\bf a}
+\sum_{i_{1},i_{2}=1}^{m}G_0^{(i_{2})}
B_{i_{1}}\hat I_{00_{\tau_{p+1},\tau_p}}^{(i_{2}i_{1})}+
$$

\vspace{2mm}
$$
+
\sum_{i_{1}=1}^{m}\left[G_0^{(i_{1})}{\bf a}\left(
\Delta \hat I_{0_{\tau_{p+1},\tau_p}}^{(i_{1})}+
\hat I_{1_{\tau_{p+1},\tau_p}}^{(i_{1})}\right)
- LB_{i_{1}}\hat I_{1_{\tau_{p+1},\tau_p}}^{(i_{1})}\right]+
$$

\vspace{2mm}
$$
+\sum_{i_{1},i_{2},i_{3}=1}^{m} G_0^{(i_{3})}G_0^{(i_{2})}
B_{i_{1}}\hat I_{000_{\tau_{p+1},\tau_p}}^{(i_{3}i_{2}i_{1})}+
\frac{\Delta^2}{2} L{\bf a}+
$$

\vspace{2mm}
$$
+\sum_{i_{1},i_{2}=1}^{m}
\left[G_0^{(i_{2})} LB_{i_{1}}\left(
\hat I_{10_{\tau_{p+1},\tau_p}}^{(i_{2}i_{1})}-
\hat I_{01_{\tau_{p+1},\tau_p}}^{(i_{2}i_{1})}
\right)
-LG_0^{(i_{2})}
B_{i_{1}}\hat I_{10_{\tau_{p+1},\tau_p}}^{(i_{2}i_{1})}
+\right.
$$

\vspace{2mm}
$$
\left.+G_0^{(i_{2})}G_0^{(i_{1})}{\bf a}\left(
\hat I_{01_{\tau_{p+1},\tau_p}}
^{(i_{2}i_{1})}+\Delta \hat I_{00_{\tau_{p+1},\tau_p}}^{(i_{2}i_{1})}
\right)\right]+
$$

\vspace{2mm}
\begin{equation}
\label{al3}
+
\sum_{i_{1},i_{2},i_{3},i_{4}=1}^{m}G_0^{(i_{4})}G_0^{(i_{3})}G_0^{(i_{2})}
B_{i_{1}}\hat I_{0000_{\tau_{p+1},\tau_p}}^{(i_{4}i_{3}i_{2}i_{1})}.
\end{equation}

\vspace{10mm}

\centerline{\bf Scheme with strong order 2.5}

\vspace{4mm}
$$
{\bf y}_{p+1}={\bf y}_{p}+\sum_{i_{1}=1}^{m}B_{i_{1}}
\hat I_{0_{\tau_{p+1},\tau_p}}^{(i_{1})}+\Delta{\bf a}
+\sum_{i_{1},i_{2}=1}^{m}G_0^{(i_{2})}
B_{i_{1}}\hat I_{00_{\tau_{p+1},\tau_p}}^{(i_{2}i_{1})}+
$$

\vspace{2mm}
$$
+
\sum_{i_{1}=1}^{m}\left[G_0^{(i_{1})}{\bf a}\left(
\Delta \hat I_{0_{\tau_{p+1},\tau_p}}^{(i_{1})}+
\hat I_{1_{\tau_{p+1},\tau_p}}^{(i_{1})}\right)
- LB_{i_{1}}\hat I_{1_{\tau_{p+1},\tau_p}}^{(i_{1})}\right]+
$$

\vspace{2mm}
$$
+\sum_{i_{1},i_{2},i_{3}=1}^{m} G_0^{(i_{3})}G_0^{(i_{2})}
B_{i_{1}}\hat I_{000_{\tau_{p+1},\tau_p}}^{(i_{3}i_{2}i_{1})}+
\frac{\Delta^2}{2} L{\bf a}+
$$

\vspace{2mm}
$$
+\sum_{i_{1},i_{2}=1}^{m}
\left[G_0^{(i_{2})} LB_{i_{1}}\left(
\hat I_{10_{\tau_{p+1},\tau_p}}^{(i_{2}i_{1})}-
\hat I_{01_{\tau_{p+1},\tau_p}}^{(i_{2}i_{1})}
\right)
- LG_0^{(i_{2})}
B_{i_{1}}\hat I_{10_{\tau_{p+1},\tau_p}}^{(i_{2}i_{1})}
+\right.
$$

\vspace{2mm}
$$
\left.+G_0^{(i_{2})}G_0^{(i_{1})}{\bf a}\left(
\hat I_{01_{\tau_{p+1},\tau_p}}
^{(i_{2}i_{1})}+\Delta \hat I_{00_{\tau_{p+1},\tau_p}}^{(i_{2}i_{1})}
\right)\right]+
$$

$$
+
\sum_{i_{1},i_{2},i_{3},i_{4}=1}^{m}G_0^{(i_{4})}G_0^{(i_{3})}G_0^{(i_{2})}
B_{i_{1}}\hat I_{0000_{\tau_{p+1},\tau_p}}^{(i_{4}i_{3}i_{2}i_{1})}+
$$

\vspace{2mm}
$$
+\sum_{i_{1}=1}^{m}\Biggl[G_0^{(i_{1})} L{\bf a}\left(\frac{1}{2}
\hat I_{2_{\tau_{p+1},\tau_p}}
^{(i_{1})}+\Delta \hat I_{1_{\tau_{p+1},\tau_p}}^{(i_{1})}+
\frac{\Delta^2}{2}\hat I_{0_{\tau_{p+1},\tau_p}}^{(i_{1})}\right)\Biggr.+
$$

\vspace{2mm}
$$
+\frac{1}{2} L LB_{i_{1}}\hat I_{2_{\tau_{p+1},\tau_p}}^{(i_{1})}-
 LG_0^{(i_{1})}{\bf a}\Biggl.
\left(\hat I_{2_{\tau_{p+1},\tau_p}}^{(i_{1})}+
\Delta \hat I_{1{\tau_{p+1},\tau_p}}^{(i_{1})}\right)\Biggr]+
$$

\vspace{2mm}
$$
+
\sum_{i_{1},i_{2},i_{3}=1}^m\left[
G_0^{(i_{3})} LG_0^{(i_{2})}B_{i_{1}}
\left(\hat I_{100_{\tau_{p+1},\tau_p}}
^{(i_{3}i_{2}i_{1})}-\hat I_{010_{\tau_{p+1},\tau_p}}
^{(i_{3}i_{2}i_{1})}\right)
\right.+
$$

\vspace{2mm}
$$
+G_0^{(i_{3})}G_0^{(i_{2})} LB_{i_{1}}\left(
\hat I_{010_{\tau_{p+1},\tau_p}}^{(i_{3}i_{2}i_{1})}-
\hat I_{001_{\tau_{p+1},\tau_p}}^{(i_{3}i_{2}i_{1})}\right)+
$$

\vspace{2mm}
$$
+G_0^{(i_{3})}G_0^{(i_{2})}G_0^{(i_{1})} {\bf a}
\left(\Delta \hat I_{000_{\tau_{p+1},\tau_p}}^{(i_{3}i_{2}i_{1})}+
\hat I_{001_{\tau_{p+1},\tau_p}}^{(i_{3}i_{2}i_{1})}\right)-
$$

\vspace{2mm}

$$
\left.- LG_0^{(i_{3})}G_0^{(i_{2})}B_{i_{1}}
\hat I_{100_{\tau_{p+1},\tau_p}}^{(i_{3}i_{2}i_{1})}\right]+
$$

\vspace{2mm}
$$
+\sum_{i_{1},i_{2},i_{3},i_{4},i_{5}=1}^m
G_0^{(i_{5})}G_0^{(i_{4})}G_0^{(i_{3})}G_0^{(i_{2})}B_{i_{1}}
\hat I_{00000_{\tau_{p+1},\tau_p}}^{(i_{5}i_{4}i_{3}i_{2}i_{1})}+
$$

\vspace{1mm}
\begin{equation}
\label{al4}
+
\frac{\Delta^3}{6}LL{\bf a}.
\end{equation}

\vspace{10mm}

\centerline{\bf Scheme with strong order 3.0}

\vspace{4mm}

$$
{\bf y}_{p+1}={\bf y}_{p}+\sum_{i_{1}=1}^{m}B_{i_{1}}
\hat I_{0_{\tau_{p+1},\tau_p}}^{(i_{1})}+\Delta{\bf a}
+\sum_{i_{1},i_{2}=1}^{m}G_0^{(i_{2})}
B_{i_{1}}\hat I_{00_{\tau_{p+1},\tau_p}}^{(i_{2}i_{1})}+
$$

\vspace{2mm}
$$
+
\sum_{i_{1}=1}^{m}\left[G_0^{(i_{1})}{\bf a}\left(
\Delta \hat I_{0_{\tau_{p+1},\tau_p}}^{(i_{1})}+
\hat I_{1_{\tau_{p+1},\tau_p}}^{(i_{1})}\right)
- LB_{i_{1}}\hat I_{1_{\tau_{p+1},\tau_p}}^{(i_{1})}\right]+
$$

\vspace{2mm}
$$
+\sum_{i_{1},i_{2},i_{3}=1}^{m} G_0^{(i_{3})}G_0^{(i_{2})}
B_{i_{1}}\hat I_{000_{\tau_{p+1},\tau_p}}^{(i_{3}i_{2}i_{1})}+
\frac{\Delta^2}{2} L{\bf a}+
$$

\vspace{2mm}
$$
+\sum_{i_{1},i_{2}=1}^{m}
\left[G_0^{(i_{2})} LB_{i_{1}}\left(
\hat I_{10_{\tau_{p+1},\tau_p}}^{(i_{2}i_{1})}-
\hat I_{01_{\tau_{p+1},\tau_p}}^{(i_{2}i_{1})}
\right)
- LG_0^{(i_{2})}
B_{i_{1}}\hat I_{10_{\tau_{p+1},\tau_p}}^{(i_{2}i_{1})}
+\right.
$$

\vspace{2mm}
$$
\left.+G_0^{(i_{2})}G_0^{(i_{1})}{\bf a}\left(
\hat I_{01_{\tau_{p+1},\tau_p}}
^{(i_{2}i_{1})}+\Delta \hat I_{00_{\tau_{p+1},\tau_p}}^{(i_{2}i_{1})}
\right)\right]+
$$

\vspace{2mm}
\begin{equation}
\label{al5}
+
\sum_{i_{1},i_{2},i_{3},i_{4}=1}^{m}G_0^{(i_{4})}G_0^{(i_{3})}G_0^{(i_{2})}
B_{i_{1}}\hat I_{0000_{\tau_{p+1},\tau_p}}^{(i_{4}i_{3}i_{2}i_{1})}+
{\bf q}_{p+1,p}+{\bf r}_{p+1,p},
\end{equation}

\vspace{6mm}
\noindent
where

\vspace{1mm}

$$
{\bf q}_{p+1,p}=
\sum_{i_{1}=1}^{m}\Biggl[G_0^{(i_{1})} L{\bf a}\left(\frac{1}{2}
\hat I_{2_{\tau_{p+1},\tau_p}}
^{(i_{1})}+\Delta \hat I_{1_{\tau_{p+1},\tau_p}}^{(i_{1})}+
\frac{\Delta^2}{2}\hat I_{0_{\tau_{p+1},\tau_p}}^{(i_{1})}\right)\Biggr.+
$$

\vspace{2mm}
$$
+\frac{1}{2} L LB_{i_{1}}\hat I_{2_{\tau_{p+1},\tau_p}}^{(i_{1})}-
LG_0^{(i_{1})}{\bf a}\Biggl.
\left(\hat I_{2_{\tau_{p+1},\tau_p}}^{(i_{1})}+
\Delta \hat I_{1{\tau_{p+1},\tau_p}}^{(i_{1})}\right)\Biggr]+
$$

\vspace{2mm}
$$
+
\sum_{i_{1},i_{2},i_{3}=1}^m\left[
G_0^{(i_{3})} LG_0^{(i_{2})}B_{i_{1}}
\left(\hat I_{100_{\tau_{p+1},\tau_p}}
^{(i_{3}i_{2}i_{1})}-\hat I_{010_{\tau_{p+1},\tau_p}}
^{(i_{3}i_{2}i_{1})}\right)
\right.+
$$

\vspace{2mm}
$$
+G_0^{(i_{3})}G_0^{(i_{2})} LB_{i_{1}}\left(
\hat I_{010_{\tau_{p+1},\tau_p}}^{(i_{3}i_{2}i_{1})}-
\hat I_{001_{\tau_{p+1},\tau_p}}^{(i_{3}i_{2}i_{1})}\right)+
$$

\vspace{2mm}
$$
+G_0^{(i_{3})}G_0^{(i_{2})}G_0^{(i_{1})} {\bf a}
\left(\Delta \hat I_{000_{\tau_{p+1},\tau_p}}^{(i_{3}i_{2}i_{1})}+
\hat I_{001_{\tau_{p+1},\tau_p}}^{(i_{3}i_{2}i_{1})}\right)-
$$

\vspace{2mm}

$$
\left.- LG_0^{(i_{3})}G_0^{(i_{2})}B_{i_{1}}
\hat I_{100_{\tau_{p+1},\tau_p}}^{(i_{3}i_{2}i_{1})}\right]+
$$

\vspace{2mm}
$$
+\sum_{i_{1},i_{2},i_{3},i_{4},i_{5}=1}^m
G_0^{(i_{5})}G_0^{(i_{4})}G_0^{(i_{3})}G_0^{(i_{2})}B_{i_{1}}
\hat I_{00000_{\tau_{p+1},\tau_p}}^{(i_{5}i_{4}i_{3}i_{2}i_{1})}+
$$

\vspace{2mm}
$$
+
\frac{\Delta^3}{6}LL {\bf a},
$$

\vspace{6mm}
\noindent
and

\vspace{2mm}

$$
{\bf r}_{p+1,p}=\sum_{i_{1},i_{2}=1}^{m}
\Biggl[G_0^{(i_{2})}G_0^{(i_{1})} L {\bf a}\Biggl(
\frac{1}{2}\hat I_{02_{\tau_{p+1},\tau_p}}^{(i_{2}i_{1})}
+
\Delta \hat I_{01_{\tau_{p+1},\tau_p}}^{(i_{2}i_{1})}
+
\frac{\Delta^2}{2}
\hat I_{00_{\tau_{p+1},\tau_p}}^{(i_{2}i_{1})}\Biggr)+\Biggr.
$$

\vspace{2mm}
$$
+
\frac{1}{2} L LG_0^{(i_{2})}B_{i_{1}}
\hat I_{20_{\tau_{p+1},\tau_p}}^{(i_{2}i_{1})}
$$

\vspace{2mm}
$$
+G_0^{(i_{2})} LG_0^{(i_{1})} {\bf a}\left(
\hat I_{11_{\tau_{p+1},\tau_p}}
^{(i_{2}i_{1})}-\hat I_{02_{\tau_{p+1},\tau_p}}^{(i_{2}i_{1})}+
\Delta\left(\hat I_{10_{\tau_{p+1},\tau_p}}
^{(i_{2}i_{1})}-\hat I_{01_{\tau_{p+1},\tau_p}}^{(i_{2}i_{1})}
\right)\right)+
$$

\vspace{2mm}
$$
+ LG_0^{(i_{2})} LB_{i_1}\left(
\hat I_{11_{\tau_{p+1},\tau_p}}
^{(i_{2}i_{1})}-\hat I_{20_{\tau_{p+1},\tau_p}}^{(i_{2}i_{1})}\right)+
$$

\vspace{2mm}
$$
+G_0^{(i_{2})} L LB_{i_1}\Biggl(
\frac{1}{2}\hat I_{02_{\tau_{p+1},\tau_p}}^{(i_{2}i_{1})}+
\frac{1}{2}\hat I_{20_{\tau_{p+1},\tau_p}}^{(i_{2}i_{1})}-
\hat I_{11_{\tau_{p+1},\tau_p}}^{(i_{2}i_{1})}\Biggr)-
$$

\vspace{2mm}
$$
\Biggl.- LG_0^{(i_{2})}G_0^{(i_{1})}{\bf a}\left(
\Delta \hat I_{10_{\tau_{p+1},\tau_p}}
^{(i_{2}i_{1})}+\hat I_{11_{\tau_{p+1},\tau_p}}^{(i_{2}i_{1})}\right)
\Biggr]+
$$

\vspace{2mm}
$$
+
\sum_{i_{1},i_2,i_3,i_{4}=1}^m\Biggl[
G_0^{(i_{4})}G_0^{(i_{3})}G_0^{(i_{2})}G_0^{(i_{1})}{\bf a}
\left(\Delta \hat I_{0000_{\tau_{p+1},\tau_p}}
^{(i_4i_{3}i_{2}i_{1})}+\hat I_{0001_{\tau_{p+1},\tau_p}}
^{(i_4i_{3}i_{2}i_{1})}\right)
+\Biggr.
$$

\vspace{2mm}
$$
+G_0^{(i_{4})}G_0^{(i_{3})} LG_0^{(i_{2})}B_{i_1}
\left(\hat I_{0100_{\tau_{p+1},\tau_p}}
^{(i_4i_{3}i_{2}i_{1})}-\hat I_{0010_{\tau_{p+1},\tau_p}}
^{(i_4i_{3}i_{2}i_{1})}\right)-
$$

\vspace{2mm}
$$
- LG_0^{(i_{4})}G_0^{(i_{3})}G_0^{(i_{2})}B_{i_1}
\hat I_{1000_{\tau_{p+1},\tau_p}}
^{(i_4i_{3}i_{2}i_{1})}+
$$

\vspace{2mm}
$$
+G_0^{(i_{4})} LG_0^{(i_{3})}G_0^{(i_{2})}B_{i_1}
\left(\hat I_{1000_{\tau_{p+1},\tau_p}}
^{(i_4i_{3}i_{2}i_{1})}-\hat I_{0100_{\tau_{p+1},\tau_p}}
^{(i_4i_{3}i_{2}i_{1})}\right)+
$$

\vspace{2mm}
$$
\Biggl.+G_0^{(i_{4})}G_0^{(i_{3})}G_0^{(i_{2})}LB_{i_1}
\left(\hat I_{0010_{\tau_{p+1},\tau_p}}
^{(i_4i_{3}i_{2}i_{1})}-\hat I_{0001_{\tau_{p+1},\tau_p}}
^{(i_4i_{3}i_{2}i_{1})}\right)\Biggr]+
$$

\vspace{2mm}
$$
+\sum_{i_{1},i_2,i_3,i_4,i_5,i_{6}=1}^m
G_0^{(i_{6})}G_0^{(i_{5})}
G_0^{(i_{4})}G_0^{(i_{3})}G_0^{(i_{2})}B_{i_{1}}
\hat I_{000000_{\tau_{p+1},\tau_p}}^{(i_6i_{5}i_{4}i_{3}i_{2}i_{1})}.
$$

\vspace{8mm}

It is well known \cite{KlPl2} that under the standard conditions
the numerical schemes (\ref{al1})--(\ref{al5}) 
have strong orders of convergence 1.0, 1.5, 2.0, 2.5, and 3.0 
correspondingly.
Among these conditions we consider only the condition
for approximations of iterated Ito stochastic 
integrals from the numerical
schemes (\ref{al1})--(\ref{al5}) \cite{KlPl2}, \cite{20}-\cite{12aa-afterxxx}

$$
{\sf M}\Biggl\{\Biggl(I_{{l_{1}\ldots l_{k}}_{\tau_{p+1},\tau_p}}
^{(i_{1}\ldots i_{k})} 
-\hat I_{{l_{1}\ldots l_{k}}_{\tau_{p+1},\tau_p}}^{(i_{1}\ldots i_{k})}
\Biggr)^2\Biggr\}\le C\Delta^{r+1},
$$

\vspace{3mm}
\noindent
where $r/2$ are strong orders of convergence for the numerical schemes
(\ref{al1})--(\ref{al5}), i.e. $r/2=1.0, 1.5,$ $2.0, 2.5,$ and $3.0.$
Moreover, constant $C$ does not depends on $\Delta$.

As we mentioned above, the numerical schemes (\ref{al1})--(\ref{al5})
are unrealizable in practice without 
procedures for the numerical simulation 
of iterated Ito stochastic integrals
from (\ref{15.001}).
In Sect. 13,
we give a brief overview to the effective method
of the mean-square approximation of
iterated Ito and Stratonovich stochastic integrals
of arbitrary multiplicity $k$ ($k\in\mathbb{N}$).

\vspace{5mm}

\section{Application of First Form of the Unified 
Taylor--Stratonovich Expansion to the High-Order
Strong Numerical Methods for Ito SDEs}

\vspace{5mm}

Let us write (\ref{t100x}) for all $s, t\in[0,T]$ such that $s>t$
in the following from

\vspace{2mm}
$$
R({\bf x}_s,s)=R({\bf x}_t,t)+\sum_{q=1}^r\sum_{(k,j,l_1,\ldots,l_k)\in
{\rm D}_q}
\frac{(s-t)^j}{j!} \sum\limits_{i_1,\ldots,i_k=1}^m
\bar G_{l_1}^{(i_1)}\ldots \bar G_{l_k}^{(i_k)}
\bar L^j R({\bf x}_t,t)\
I^{*(i_1\ldots i_k)}_{{l_1\ldots l_k}_{s,t}}+
$$

\vspace{2mm}
\begin{equation}
\label{15.001x}
+{\bf 1}_{\{r=2d-1,
d\in N\}}\frac{(s-t)^{(r+1)/2}}{\left(
(r+1)/2\right)!}L^{(r+1)/2}R({\bf x}_t,t)
+\left(\bar H_{r+1}\right)_{s,t}\ \ \ \hbox{w.\ p.\ 1},
\end{equation}

\vspace{6mm}
\noindent
where

\vspace{-3mm}
$$
\left(\bar H_{r+1}\right)_{s,t}= \left(H_{r+1}\right)_{s,t}
-{\bf 1}_{\{r=2d-1,
d\in N\}}\frac{(s-t)^{(r+1)/2}}{\left(
(r+1)/2\right)!}L^{(r+1)/2}R({\bf x}_t,t).
$$

\vspace{7mm}

Consider the partition 
$\{\tau_p\}_{p=0}^N$ of the interval
$[0,T]$ such that

\vspace{-1mm}
$$
0=\tau_0<\tau_1<\ldots<\tau_N=T,\ \ \
\Delta_N=
\max\limits_{0\le j\le N-1}\left|\tau_{j+1}-\tau_j\right|.
$$

\vspace{3mm}

From (\ref{15.001x}) for $s=\tau_{p+1},$
$t=\tau_p$ we obtain the following representation of
explicit one-step strong numerical scheme for Ito SDE (\ref{1.5.2}),
which is based on the first form of the unified Taylor--Stratonovich
expansion

$$
{\bf y}_{p+1}={\bf y}_{p}+
\sum_{q=1}^r\sum_{(k,j,l_1,\ldots,l_k)\in{\rm D}_q}
\frac{(\tau_{p+1}-\tau_p)^j}{j!} \sum\limits_{i_1,\ldots,i_k=1}^m
\bar G_{l_1}^{(i_1)}\ldots \bar G_{l_k}^{(i_k)}
\bar L^j\hspace{0.4mm}{\bf y}_{p}\
{\hat I}^{*(i_1\ldots i_k)}_{{l_1\ldots l_k}_{\tau_{p+1},\tau_p}}+
$$

\vspace{2mm}
\begin{equation}
\label{15.002x}
+{\bf 1}_{\{r=2d-1,
d\in N\}}\frac{(\tau_{p+1}-\tau_p)^{(r+1)/2}}{\left(
(r+1)/2\right)!}L^{(r+1)/2}\hspace{0.4mm}{\bf y}_{p},
\end{equation}

\vspace{6mm}
\noindent
where
$\hat I^{*(i_1\ldots i_k)}_{{l_1\ldots l_k}_{\tau_{p+1},\tau_p}}$ 
is an approximation of the iterated Stratonovich
stochastic integral 
$I^{*(i_1\ldots i_k)}_{{l_1\ldots l_k}_{\tau_{p+1},\tau_p}}$
defined as
$$
I_{{l_1\ldots l_k}_{s,t}}^{*(i_1\ldots i_k)}=
{\int\limits_t^{*}}^s
(t-t_{k})^{l_{k}}\ldots 
{\int\limits_t^{*}}^{t_2}
(t-t_ {1}) ^ {l_ {1}} d
{\bf f} ^ {(i_ {1})} _ {t_ {1}} \ldots 
d {\bf f} _ {t_ {k}} ^ {(i_ {k})}.
$$

\vspace{2mm}

Note that we understand the equality (\ref{15.002x}) componentwise
with respect to the components ${\bf y}_p^{(i)}$ of the column
${\bf y}_p.$
Also for simplicity we put 
$\tau_p=p\Delta$,
$\Delta=T/N,$ $T=\tau_N,$ $p=0,1,\ldots,N.$

It is known \cite{KlPl2} that under the appropriate conditions
the numerical scheme (\ref{15.002x}) has strong order of convergence $r/2$
($r\in\mathbb{N}$).

Denote
$$
\bar{\bf a}({\bf x},t)={\bf a}({\bf x},t)-
\frac{1}{2}\sum\limits_{j=1}^m G_0^{(j)}B_j({\bf x},t),
$$

\vspace{2mm}
\noindent
where $B_j({\bf x},t)$ is the $j$-th column of the matrix
function $B({\bf x},t)$  (see (\ref{1.5.2})).

Below we consider particular cases of the numerical scheme
(\ref{15.002x}) for $r=2,3,4,5,$ and $6,$ i.e. 
explicit one-step strong numerical schemes for Ito SDE (\ref{1.5.2})
with orders $1.0, 1.5, 2.0, 2.5,$ and $3.0$ of convergence.
At that for simplicity 
we will write $\bar{\bf a},$ $\bar L\bar {\bf a},$ $L{\bf a},$
$B_i,$ $G_0^{(i)}B_{j}$ etc.
instead of $\bar{\bf a}({\bf y}_p,\tau_p),$ 
$\bar L \bar {\bf a}({\bf y}_p,\tau_p),$ 
$L{\bf a}({\bf y}_p,\tau_p),$  $B_i({\bf y}_p,\tau_p),$ 
$G_0^{(i)}B_{j}({\bf y}_p,\tau_p)$ etc. correspondingly.
Moreover, the operators $\bar L$ and $G_0^{(i)},$ $i=1,\ldots,m,$
are determined by the equalities
(\ref{2.3}), (\ref{2.4}), and (\ref{2.4a}) as before.

\vspace{10mm}

\centerline{\bf Scheme with strong order 1.0}

\vspace{2mm}
\begin{equation}
\label{al1x}
{\bf y}_{p+1}={\bf y}_{p}+\sum_{i_{1}=1}^{m}B_{i_{1}}
\hat I_{0_{\tau_{p+1},\tau_p}}^{*(i_{1})}+\Delta\bar{\bf a}
+\sum_{i_{1},i_{2}=1}^{m}G_0^{(i_{2})}
B_{i_{1}}\hat I_{00_{\tau_{p+1},\tau_p}}^{*(i_{2}i_{1})}.
\end{equation}

\vspace{10mm}

\centerline{\bf Scheme with strong order 1.5}

\vspace{4mm}
$$
{\bf y}_{p+1}={\bf y}_{p}+\sum_{i_{1}=1}^{m}B_{i_{1}}
\hat I_{0_{\tau_{p+1},\tau_p}}^{*(i_{1})}+\Delta\bar{\bf a}
+\sum_{i_{1},i_{2}=1}^{m}G_0^{(i_{2})}
B_{i_{1}}\hat I_{00_{\tau_{p+1},\tau_p}}^{*(i_{2}i_{1})}+
$$

\vspace{2mm}
$$
+
\sum_{i_{1}=1}^{m}\left[G_0^{(i_{1})}\bar{\bf a}\left(
\Delta \hat I_{0_{\tau_{p+1},\tau_p}}^{*(i_{1})}+
\hat I_{1_{\tau_{p+1},\tau_p}}^{*(i_{1})}\right)
-\bar LB_{i_{1}}\hat I_{1_{\tau_{p+1},\tau_p}}^{*(i_{1})}\right]+
$$

\vspace{2mm}
\begin{equation}
\label{al2x}
+\sum_{i_{1},i_{2},i_{3}=1}^{m} G_0^{(i_{3})}G_0^{(i_{2})}
B_{i_{1}}\hat I_{000_{\tau_{p+1},\tau_p}}^{*(i_{3}i_{2}i_{1})}+
\frac{\Delta^2}{2}L{\bf a}.
\end{equation}

\vspace{10mm}

\centerline{\bf Scheme with strong order 2.0}

\vspace{4mm}
$$
{\bf y}_{p+1}={\bf y}_{p}+\sum_{i_{1}=1}^{m}B_{i_{1}}
\hat I_{0_{\tau_{p+1},\tau_p}}^{*(i_{1})}+\Delta\bar{\bf a}
+\sum_{i_{1},i_{2}=1}^{m}G_0^{(i_{2})}
B_{i_{1}}\hat I_{00_{\tau_{p+1},\tau_p}}^{*(i_{2}i_{1})}+
$$

\vspace{2mm}
$$
+
\sum_{i_{1}=1}^{m}\left[G_0^{(i_{1})}\bar{\bf a}\left(
\Delta \hat I_{0_{\tau_{p+1},\tau_p}}^{*(i_{1})}+
\hat I_{1_{\tau_{p+1},\tau_p}}^{*(i_{1})}\right)
-\bar LB_{i_{1}}\hat I_{1_{\tau_{p+1},\tau_p}}^{*(i_{1})}\right]+
$$

\vspace{2mm}
$$
+\sum_{i_{1},i_{2},i_{3}=1}^{m} G_0^{(i_{3})}G_0^{(i_{2})}
B_{i_{1}}\hat I_{000_{\tau_{p+1},\tau_p}}^{*(i_{3}i_{2}i_{1})}+
\frac{\Delta^2}{2}\bar L\bar{\bf a}+
$$

\vspace{2mm}
$$
+\sum_{i_{1},i_{2}=1}^{m}
\left[G_0^{(i_{2})}\bar LB_{i_{1}}\left(
\hat I_{10_{\tau_{p+1},\tau_p}}^{*(i_{2}i_{1})}-
\hat I_{01_{\tau_{p+1},\tau_p}}^{*(i_{2}i_{1})}
\right)
-\bar LG_0^{(i_{2})}
B_{i_{1}}\hat I_{10_{\tau_{p+1},\tau_p}}^{*(i_{2}i_{1})}
+\right.
$$

\vspace{2mm}
$$
\left.+G_0^{(i_{2})}G_0^{(i_{1})}\bar{\bf a}\left(
\hat I_{01_{\tau_{p+1},\tau_p}}
^{*(i_{2}i_{1})}+\Delta \hat I_{00_{\tau_{p+1},\tau_p}}^{*(i_{2}i_{1})}
\right)\right]+
$$

\vspace{2mm}
\begin{equation}
\label{al3x}
+
\sum_{i_{1},i_{2},i_{3},i_{4}=1}^{m}G_0^{(i_{4})}G_0^{(i_{3})}G_0^{(i_{2})}
B_{i_{1}}\hat I_{0000_{\tau_{p+1},\tau_p}}^{*(i_{4}i_{3}i_{2}i_{1})}.
\end{equation}

\vspace{10mm}

\centerline{\bf Scheme with strong order 2.5}

\vspace{4mm}
$$
{\bf y}_{p+1}={\bf y}_{p}+\sum_{i_{1}=1}^{m}B_{i_{1}}
\hat I_{0_{\tau_{p+1},\tau_p}}^{*(i_{1})}+\Delta\bar{\bf a}
+\sum_{i_{1},i_{2}=1}^{m}G_0^{(i_{2})}
B_{i_{1}}\hat I_{00_{\tau_{p+1},\tau_p}}^{*(i_{2}i_{1})}+
$$

\vspace{2mm}
$$
+
\sum_{i_{1}=1}^{m}\left[G_0^{(i_{1})}\bar{\bf a}\left(
\Delta \hat I_{0_{\tau_{p+1},\tau_p}}^{*(i_{1})}+
\hat I_{1_{\tau_{p+1},\tau_p}}^{*(i_{1})}\right)
-\bar LB_{i_{1}}\hat I_{1_{\tau_{p+1},\tau_p}}^{*(i_{1})}\right]+
$$

\vspace{2mm}
$$
+\sum_{i_{1},i_{2},i_{3}=1}^{m} G_0^{(i_{3})}G_0^{(i_{2})}
B_{i_{1}}\hat I_{000_{\tau_{p+1},\tau_p}}^{*(i_{3}i_{2}i_{1})}+
\frac{\Delta^2}{2}\bar L\bar{\bf a}+
$$

\vspace{2mm}
$$
+\sum_{i_{1},i_{2}=1}^{m}
\left[G_0^{(i_{2})}\bar LB_{i_{1}}\left(
\hat I_{10_{\tau_{p+1},\tau_p}}^{*(i_{2}i_{1})}-
\hat I_{01_{\tau_{p+1},\tau_p}}^{*(i_{2}i_{1})}
\right)
-\bar LG_0^{(i_{2})}
B_{i_{1}}\hat I_{10_{\tau_{p+1},\tau_p}}^{*(i_{2}i_{1})}
+\right.
$$

\vspace{2mm}
$$
\left.+G_0^{(i_{2})}G_0^{(i_{1})}\bar{\bf a}\left(
\hat I_{01_{\tau_{p+1},\tau_p}}
^{*(i_{2}i_{1})}+\Delta \hat I_{00_{\tau_{p+1},\tau_p}}^{*(i_{2}i_{1})}
\right)\right]+
$$

$$
+
\sum_{i_{1},i_{2},i_{3},i_{4}=1}^{m}G_0^{(i_{4})}G_0^{(i_{3})}G_0^{(i_{2})}
B_{i_{1}}\hat I_{0000_{\tau_{p+1},\tau_p}}^{*(i_{4}i_{3}i_{2}i_{1})}+
$$

\vspace{2mm}
$$
+\sum_{i_{1}=1}^{m}\Biggl[G_0^{(i_{1})}\bar L\bar{\bf a}\left(\frac{1}{2}
\hat I_{2_{\tau_{p+1},\tau_p}}
^{*(i_{1})}+\Delta \hat I_{1_{\tau_{p+1},\tau_p}}^{*(i_{1})}+
\frac{\Delta^2}{2}\hat I_{0_{\tau_{p+1},\tau_p}}^{*(i_{1})}\right)\Biggr.+
$$

\vspace{2mm}
$$
+\frac{1}{2}\bar L\bar LB_{i_{1}}\hat I_{2_{\tau_{p+1},\tau_p}}^{*(i_{1})}-
\bar LG_0^{(i_{1})}\bar{\bf a}\Biggl.
\left(\hat I_{2_{\tau_{p+1},\tau_p}}^{*(i_{1})}+
\Delta \hat I_{1{\tau_{p+1},\tau_p}}^{*(i_{1})}\right)\Biggr]+
$$

\vspace{2mm}
$$
+
\sum_{i_{1},i_{2},i_{3}=1}^m\left[
G_0^{(i_{3})}\bar LG_0^{(i_{2})}B_{i_{1}}
\left(\hat I_{100_{\tau_{p+1},\tau_p}}
^{*(i_{3}i_{2}i_{1})}-\hat I_{010_{\tau_{p+1},\tau_p}}
^{*(i_{3}i_{2}i_{1})}\right)
\right.+
$$

\vspace{2mm}
$$
+G_0^{(i_{3})}G_0^{(i_{2})}\bar LB_{i_{1}}\left(
\hat I_{010_{\tau_{p+1},\tau_p}}^{*(i_{3}i_{2}i_{1})}-
\hat I_{001_{\tau_{p+1},\tau_p}}^{*(i_{3}i_{2}i_{1})}\right)+
$$

\vspace{2mm}
$$
+G_0^{(i_{3})}G_0^{(i_{2})}G_0^{(i_{1})}\bar {\bf a}
\left(\Delta \hat I_{000_{\tau_{p+1},\tau_p}}^{*(i_{3}i_{2}i_{1})}+
\hat I_{001_{\tau_{p+1},\tau_p}}^{*(i_{3}i_{2}i_{1})}\right)-
$$

\vspace{2mm}

$$
\left.-\bar LG_0^{(i_{3})}G_0^{(i_{2})}B_{i_{1}}
\hat I_{100_{\tau_{p+1},\tau_p}}^{*(i_{3}i_{2}i_{1})}\right]+
$$

\vspace{2mm}
$$
+\sum_{i_{1},i_{2},i_{3},i_{4},i_{5}=1}^m
G_0^{(i_{5})}G_0^{(i_{4})}G_0^{(i_{3})}G_0^{(i_{2})}B_{i_{1}}
\hat I_{00000_{\tau_{p+1},\tau_p}}^{*(i_{5}i_{4}i_{3}i_{2}i_{1})}+
$$

\vspace{1mm}
\begin{equation}
\label{al4x}
+
\frac{\Delta^3}{6}LL{\bf a}.
\end{equation}

\vspace{10mm}

\centerline{\bf Scheme with strong order 3.0}

\vspace{4mm}

$$
{\bf y}_{p+1}={\bf y}_{p}+\sum_{i_{1}=1}^{m}B_{i_{1}}
\hat I_{0_{\tau_{p+1},\tau_p}}^{*(i_{1})}+\Delta\bar{\bf a}
+\sum_{i_{1},i_{2}=1}^{m}G_0^{(i_{2})}
B_{i_{1}}\hat I_{00_{\tau_{p+1},\tau_p}}^{*(i_{2}i_{1})}+
$$

\vspace{2mm}
$$
+
\sum_{i_{1}=1}^{m}\left[G_0^{(i_{1})}\bar{\bf a}\left(
\Delta \hat I_{0_{\tau_{p+1},\tau_p}}^{*(i_{1})}+
\hat I_{1_{\tau_{p+1},\tau_p}}^{*(i_{1})}\right)
-\bar LB_{i_{1}}\hat I_{1_{\tau_{p+1},\tau_p}}^{*(i_{1})}\right]+
$$

\vspace{2mm}
$$
+\sum_{i_{1},i_{2},i_{3}=1}^{m} G_0^{(i_{3})}G_0^{(i_{2})}
B_{i_{1}}\hat I_{000_{\tau_{p+1},\tau_p}}^{*(i_{3}i_{2}i_{1})}+
\frac{\Delta^2}{2}\bar L\bar{\bf a}+
$$

\vspace{2mm}
$$
+\sum_{i_{1},i_{2}=1}^{m}
\left[G_0^{(i_{2})}\bar LB_{i_{1}}\left(
\hat I_{10_{\tau_{p+1},\tau_p}}^{*(i_{2}i_{1})}-
\hat I_{01_{\tau_{p+1},\tau_p}}^{*(i_{2}i_{1})}
\right)
-\bar LG_0^{(i_{2})}
B_{i_{1}}\hat I_{10_{\tau_{p+1},\tau_p}}^{*(i_{2}i_{1})}
+\right.
$$

\vspace{2mm}
$$
\left.+G_0^{(i_{2})}G_0^{(i_{1})}\bar{\bf a}\left(
\hat I_{01_{\tau_{p+1},\tau_p}}
^{*(i_{2}i_{1})}+\Delta \hat I_{00_{\tau_{p+1},\tau_p}}^{*(i_{2}i_{1})}
\right)\right]+
$$

\vspace{2mm}
\begin{equation}
\label{al5x}
+
\sum_{i_{1},i_{2},i_{3},i_{4}=1}^{m}G_0^{(i_{4})}G_0^{(i_{3})}G_0^{(i_{2})}
B_{i_{1}}\hat I_{0000_{\tau_{p+1},\tau_p}}^{*(i_{4}i_{3}i_{2}i_{1})}+
{\bf q}_{p+1,p}+{\bf r}_{p+1,p},
\end{equation}

\vspace{6mm}
\noindent
where

\vspace{1mm}

$$
{\bf q}_{p+1,p}=
\sum_{i_{1}=1}^{m}\Biggl[G_0^{(i_{1})}\bar L\bar{\bf a}\left(\frac{1}{2}
\hat I_{2_{\tau_{p+1},\tau_p}}
^{*(i_{1})}+\Delta \hat I_{1_{\tau_{p+1},\tau_p}}^{*(i_{1})}+
\frac{\Delta^2}{2}\hat I_{0_{\tau_{p+1},\tau_p}}^{*(i_{1})}\right)\Biggr.+
$$

\vspace{2mm}
$$
+\frac{1}{2}\bar L\bar LB_{i_{1}}\hat I_{2_{\tau_{p+1},\tau_p}}^{*(i_{1})}-
\bar LG_0^{(i_{1})}\bar{\bf a}\Biggl.
\left(\hat I_{2_{\tau_{p+1},\tau_p}}^{*(i_{1})}+
\Delta \hat I_{1{\tau_{p+1},\tau_p}}^{*(i_{1})}\right)\Biggr]+
$$

\vspace{2mm}
$$
+
\sum_{i_{1},i_{2},i_{3}=1}^m\left[
G_0^{(i_{3})}\bar LG_0^{(i_{2})}B_{i_{1}}
\left(\hat I_{100_{\tau_{p+1},\tau_p}}
^{*(i_{3}i_{2}i_{1})}-\hat I_{010_{\tau_{p+1},\tau_p}}
^{*(i_{3}i_{2}i_{1})}\right)
\right.+
$$

\vspace{2mm}
$$
+G_0^{(i_{3})}G_0^{(i_{2})}\bar LB_{i_{1}}\left(
\hat I_{010_{\tau_{p+1},\tau_p}}^{*(i_{3}i_{2}i_{1})}-
\hat I_{001_{\tau_{p+1},\tau_p}}^{*(i_{3}i_{2}i_{1})}\right)+
$$

\vspace{2mm}
$$
+G_0^{(i_{3})}G_0^{(i_{2})}G_0^{(i_{1})}\bar {\bf a}
\left(\Delta \hat I_{000_{\tau_{p+1},\tau_p}}^{*(i_{3}i_{2}i_{1})}+
\hat I_{001_{\tau_{p+1},\tau_p}}^{*(i_{3}i_{2}i_{1})}\right)-
$$

\vspace{2mm}

$$
\left.-\bar LG_0^{(i_{3})}G_0^{(i_{2})}B_{i_{1}}
\hat I_{100_{\tau_{p+1},\tau_p}}^{*(i_{3}i_{2}i_{1})}\right]+
$$

\vspace{2mm}
$$
+\sum_{i_{1},i_{2},i_{3},i_{4},i_{5}=1}^m
G_0^{(i_{5})}G_0^{(i_{4})}G_0^{(i_{3})}G_0^{(i_{2})}B_{i_{1}}
\hat I_{00000_{\tau_{p+1},\tau_p}}^{*(i_{5}i_{4}i_{3}i_{2}i_{1})}+
$$

\vspace{2mm}
$$
+
\frac{\Delta^3}{6}\bar L\bar L\bar {\bf a},
$$

\vspace{6mm}
\noindent
and

\vspace{2mm}

$$
{\bf r}_{p+1,p}=\sum_{i_{1},i_{2}=1}^{m}
\Biggl[G_0^{(i_{2})}G_0^{(i_{1})}\bar L\bar {\bf a}\Biggl(
\frac{1}{2}\hat I_{02_{\tau_{p+1},\tau_p}}^{*(i_{2}i_{1})}
+
\Delta \hat I_{01_{\tau_{p+1},\tau_p}}^{*(i_{2}i_{1})}
+
\frac{\Delta^2}{2}
\hat I_{00_{\tau_{p+1},\tau_p}}^{*(i_{2}i_{1})}\Biggr)+\Biggr.
$$

\vspace{2mm}
$$
+
\frac{1}{2}\bar L\bar LG_0^{(i_{2})}B_{i_{1}}
\hat I_{20_{\tau_{p+1},\tau_p}}^{*(i_{2}i_{1})}
$$

\vspace{2mm}
$$
+G_0^{(i_{2})}\bar LG_0^{(i_{1})}\bar {\bf a}\left(
\hat I_{11_{\tau_{p+1},\tau_p}}
^{*(i_{2}i_{1})}-\hat I_{02_{\tau_{p+1},\tau_p}}^{*(i_{2}i_{1})}+
\Delta\left(\hat I_{10_{\tau_{p+1},\tau_p}}
^{*(i_{2}i_{1})}-\hat I_{01_{\tau_{p+1},\tau_p}}^{*(i_{2}i_{1})}
\right)\right)+
$$

\vspace{2mm}
$$
+\bar LG_0^{(i_{2})}\bar LB_{i_1}\left(
\hat I_{11_{\tau_{p+1},\tau_p}}
^{*(i_{2}i_{1})}-\hat I_{20_{\tau_{p+1},\tau_p}}^{*(i_{2}i_{1})}\right)+
$$

\vspace{2mm}
$$
+G_0^{(i_{2})}\bar L\bar LB_{i_1}\Biggl(
\frac{1}{2}\hat I_{02_{\tau_{p+1},\tau_p}}^{*(i_{2}i_{1})}+
\frac{1}{2}\hat I_{20_{\tau_{p+1},\tau_p}}^{*(i_{2}i_{1})}-
\hat I_{11_{\tau_{p+1},\tau_p}}^{*(i_{2}i_{1})}\Biggr)-
$$

\vspace{2mm}
$$
\Biggl.-\bar LG_0^{(i_{2})}G_0^{(i_{1})}\bar{\bf a}\left(
\Delta \hat I_{10_{\tau_{p+1},\tau_p}}
^{*(i_{2}i_{1})}+\hat I_{11_{\tau_{p+1},\tau_p}}^{*(i_{2}i_{1})}\right)
\Biggr]+
$$

\vspace{2mm}
$$
+
\sum_{i_{1},i_2,i_3,i_{4}=1}^m\Biggl[
G_0^{(i_{4})}G_0^{(i_{3})}G_0^{(i_{2})}G_0^{(i_{1})}\bar{\bf a}
\left(\Delta \hat I_{0000_{\tau_{p+1},\tau_p}}
^{*(i_4i_{3}i_{2}i_{1})}+\hat I_{0001_{\tau_{p+1},\tau_p}}
^{*(i_4i_{3}i_{2}i_{1})}\right)
+\Biggr.
$$

\vspace{2mm}
$$
+G_0^{(i_{4})}G_0^{(i_{3})}\bar LG_0^{(i_{2})}B_{i_1}
\left(\hat I_{0100_{\tau_{p+1},\tau_p}}
^{*(i_4i_{3}i_{2}i_{1})}-\hat I_{0010_{\tau_{p+1},\tau_p}}
^{*(i_4i_{3}i_{2}i_{1})}\right)-
$$

\vspace{2mm}
$$
-\bar LG_0^{(i_{4})}G_0^{(i_{3})}G_0^{(i_{2})}B_{i_1}
\hat I_{1000_{\tau_{p+1},\tau_p}}
^{*(i_4i_{3}i_{2}i_{1})}+
$$

\vspace{2mm}
$$
+G_0^{(i_{4})}\bar LG_0^{(i_{3})}G_0^{(i_{2})}B_{i_1}
\left(\hat I_{1000_{\tau_{p+1},\tau_p}}
^{*(i_4i_{3}i_{2}i_{1})}-\hat I_{0100_{\tau_{p+1},\tau_p}}
^{*(i_4i_{3}i_{2}i_{1})}\right)+
$$

\vspace{2mm}
$$
\Biggl.+G_0^{(i_{4})}G_0^{(i_{3})}G_0^{(i_{2})}\bar LB_{i_1}
\left(\hat I_{0010_{\tau_{p+1},\tau_p}}
^{*(i_4i_{3}i_{2}i_{1})}-\hat I_{0001_{\tau_{p+1},\tau_p}}
^{*(i_4i_{3}i_{2}i_{1})}\right)\Biggr]+
$$

\vspace{2mm}
$$
+\sum_{i_{1},i_2,i_3,i_4,i_5,i_{6}=1}^m
G_0^{(i_{6})}G_0^{(i_{5})}
G_0^{(i_{4})}G_0^{(i_{3})}G_0^{(i_{2})}B_{i_{1}}
\hat I_{000000_{\tau_{p+1},\tau_p}}^{*(i_6i_{5}i_{4}i_{3}i_{2}i_{1})}.
$$

\vspace{8mm}

It is well known \cite{KlPl2} that under the standard conditions
the numerical schemes (\ref{al1x})--(\ref{al5x}) 
have strong orders of convergence 1.0, 1.5, 2.0, 2.5, and 3.0
correspondingly.
Among these conditions we consider only the condition
for approximations of iterated Stratonovich stochastic 
integrals from the numerical
schemes (\ref{al1x})--(\ref{al5x}) \cite{KlPl2}, \cite{20}-\cite{12aa-afterxxx} 

$$
{\sf M}\Biggl\{\Biggl(I_{{l_{1}\ldots l_{k}}_{\tau_{p+1},\tau_p}}
^{*(i_{1}\ldots i_{k})} 
-\hat I_{{l_{1}\ldots l_{k}}_{\tau_{p+1},\tau_p}}^{*(i_{1}\ldots i_{k})}
\Biggr)^2\Biggr\}\le C\Delta^{r+1},
$$

\vspace{3mm}
\noindent
where $r/2$ are strong orders of convergence for the numerical schemes
(\ref{al1x})--(\ref{al5x}), i.e. $r/2=1.0, 1.5,$ $2.0, 2.5,$ and $3.0.$
Moreover, constant $C$ does not depends on $\Delta$.

As we mentioned above, the numerical schemes (\ref{al1x})--(\ref{al5x})
are unrealizable in practice without 
procedures for the numerical simulation 
of iterated Stratonovich stochastic integrals
from (\ref{15.001x}).
In the next section,
we give a brief overview to the effective method
of the mean-square approximation of
iterated Ito and Stratonovich stochastic integrals
of arbitrary multiplicity $k$ ($k\in\mathbb{N}$).

\vspace{5mm}

\section{Method of the Mean-Square Approximation
of Iterated Ito and Stratonovich 
Stochastic Integrals Based on Generalized Multiple Fourier Series}

\vspace{5mm}

It should be noted that 
there is an approach to the mean-square approximation of
iterated stochastic integrals based on multiple integral sums
(see, for example, \cite{Mi2}). This method implies
the partitioning of the integration interval   
of the iterated stochastic integral under consideration; 
this interval is the integration step of the numerical methods 
used to solve Ito SDEs; therefore, it is already fairly small and does 
not need to be partitioned. 
Computational experiments \cite{7} show that the application of 
the method \cite{Mi2} to stochastic integrals with multiplicities
$k\ge 2$ leads to unacceptably high computational cost and 
accumulation of computation errors. 
Another well-known method 
is based on the Karhunen--Loeve expansion of the
Brownian bridge process \cite{Mi2}. This method 
has no the mentioned drawback
(also see \cite{KlPl2}, \cite{KPS})
but leads to iterated application of the 
operation of limit transition.
So, the mentioned method may not converge in 
the mean-square sense 
to appropriate iterated stochastic integrals
for some methods of series summation (see discussion in Sect.~14
for details).

The difficulties noted above can be overcome with a different 
and more effective method proposed and developed
by the author in \cite{7} (also see \cite{arxiv-1}-\cite{5-008}, 
\cite{8}-\cite{new-art-1xxy}).
The idea of this method is as follows: the iterated Ito stochastic 
integral $J[\psi^{(k)}]_{T,t}$ of the form (\ref{ito})
with multiplicity $k$ is represented as a multiple stochastic 
integral from the nonrandom discontinuous function 
$K(t_1,\ldots,t_k)$ of $k$ variables
(see (\ref{ppp}) below)
defined on the hypercube $[t, T]^k$, where here and further $[t, T]$ is an interval of 
integration of the iterated Ito stochastic integral. Then, 
the function $K(t_1,\ldots,t_k)$ 
is expanded in the hypercube $[t, T]^k$ into the generalized 
multiple Fourier series converging 
in the mean-square sense
in the space 
$L_2([t,T]^k)$. After a number of nontrivial transformations we come 
(see Theorems 7, 8 below) to the 
mean-square convergening expansion of the iterated Ito stochastic 
integral into the multiple 
series of products
of standard  Gaussian random 
variables. The coefficients of this 
series are the coefficients of 
generalized multiple Fourier series for the function $K(t_1,\ldots,t_k)$, 
which can be calculated using the explicit formula 
regardless of multiplicity $k$ of the iterated Ito stochastic integral.
Hereinafter, this method is referred to as the method of 
generalized multiple Fourier series.

Suppose that $\{\phi_j(x)\}_{j=0}^{\infty}$
is a complete orthonormal system of functions in 
the space $L_2([t, T])$. 
Define the following function on the hypercube $[t, T]^k$

\vspace{-3mm}
\begin{equation}
\label{ppp}
K(t_1,\ldots,t_k)=
\begin{cases}
\psi_1(t_1)\ldots \psi_k(t_k),\ &t_1<\ldots<t_k\\
~\\
~\\
0,\ &\hbox{\rm otherwise}
\end{cases}\ \ \ \ 
=\ \ \ \ 
\prod\limits_{l=1}^k
\psi_l(t_l)\ \prod\limits_{l=1}^{k-1}{\bf 1}_{\{t_l<t_{l+1}\}},\ 
\end{equation}

\vspace{3mm}
\noindent
where $t_1,\ldots,t_k\in [t, T]$ $(k\ge 2)$,  
$K(t_1)\equiv\psi_1(t_1)$ for $t_1\in[t, T],$
$\psi_1(\tau),\ldots,\psi_k(\tau): [t, T]\to\mathbb{R}$
are continuous nonrandom functions (the case 
$\psi_1(\tau),\ldots,\psi_k(\tau)\in L_2([t, T])$ will be considered below
in this section).
Here 
${\bf 1}_A$ denotes the indicator of the set $A$.

The function $K(t_1,\ldots,t_k)$ of the form (\ref{ppp})
is piecewise continuous in the 
hypercube $[t, T]^k.$
At this situation it is well known that the generalized 
multiple Fourier series 
of $K(t_1,\ldots,t_k)\in L_2([t, T]^k)$ is converging 
to $K(t_1,\ldots,t_k)$ in the hypercube $[t, T]^k$ in 
the mean-square sense, i.e.

\vspace{1mm}
\begin{equation}
\label{sos1z}
\lim\limits_{p_1,\ldots,p_k\to \infty}
\Biggl\Vert
K(t_1,\ldots,t_k)-
\sum_{j_1=0}^{p_1}\ldots \sum_{j_k=0}^{p_k}
C_{j_k\ldots j_1}\prod_{l=1}^{k} \phi_{j_l}(t_l)\Biggr\Vert_{L_2([t, T]^k)}
=0,
\end{equation}

\vspace{2mm}
\noindent
where

\vspace{-3mm}
\begin{equation}
\label{ppppa}
C_{j_k\ldots j_1}=\int\limits_{[t,T]^k}
K(t_1,\ldots,t_k)\prod_{l=1}^{k}\phi_{j_l}(t_l)dt_1\ldots dt_k
\end{equation}

\vspace{5mm}
\noindent
is the Fourier coefficient, and

$$
\left\Vert f\right\Vert_{L_2([t, T]^k)}=\left(\int\limits_{[t,T]^k}
f^2(t_1,\ldots,t_k)dt_1\ldots dt_k\right)^{1/2}.
$$

\vspace{5mm}

Consider the partition $\{\tau_j\}_{j=0}^N$ of $[t,T]$ such that

\begin{equation}
\label{1111}
t=\tau_0<\ldots <\tau_N=T,\ \ \
\Delta_N=
\hbox{\vtop{\offinterlineskip\halign{
\hfil#\hfil\cr
{\rm max}\cr
$\stackrel{}{{}_{0\le j\le N-1}}$\cr
}} }\Delta\tau_j\to 0\ \ \hbox{if}\ \ N\to \infty,\ \ \ 
\Delta\tau_j=\tau_{j+1}-\tau_j.
\end{equation}

\vspace{5mm}

{\bf Theorem 7} \cite{7} (2006) (also see
\cite{arxiv-1}-\cite{5-008}, \cite{8}-\cite{new-2023a}). {\it Suppose that
every $\psi_l(\tau)$ $(l=1,\ldots, k)$ is a continuous nonrandom function on 
$[t, T]$ and
$\{\phi_j(x)\}_{j=0}^{\infty}$ is a complete orthonormal system  
of continuous functions in the space $L_2([t,T]).$ 
Then

$$
J[\psi^{(k)}]_{T,t}\  =\ 
\hbox{\vtop{\offinterlineskip\halign{
\hfil#\hfil\cr
{\rm l.i.m.}\cr
$\stackrel{}{{}_{p_1,\ldots,p_k\to \infty}}$\cr
}} }\sum_{j_1=0}^{p_1}\ldots\sum_{j_k=0}^{p_k}
C_{j_k\ldots j_1}\Biggl(
\prod_{l=1}^k\zeta_{j_l}^{(i_l)}\ -
\Biggr.
$$

\vspace{2mm}
\begin{equation}
\label{tyyy}
-\ \Biggl.
\hbox{\vtop{\offinterlineskip\halign{
\hfil#\hfil\cr
{\rm l.i.m.}\cr
$\stackrel{}{{}_{N\to \infty}}$\cr
}} }\sum_{(l_1,\ldots,l_k)\in {\rm G}_k}
\phi_{j_{1}}(\tau_{l_1})
\Delta{\bf w}_{\tau_{l_1}}^{(i_1)}\ldots
\phi_{j_{k}}(\tau_{l_k})
\Delta{\bf w}_{\tau_{l_k}}^{(i_k)}\Biggr),
\end{equation}

\vspace{5mm}
\noindent
where

\vspace{-2mm}
$$
{\rm G}_k={\rm H}_k\backslash{\rm L}_k,\ \ \
{\rm H}_k=\biggl\{(l_1,\ldots,l_k):\ l_1,\ldots,l_k=0,\ 1,\ldots,N-1\biggr\},
$$

$$
{\rm L}_k=\biggl\{(l_1,\ldots,l_k):\ l_1,\ldots,l_k=0,\ 1,\ldots,N-1;\
l_g\ne l_r\ (g\ne r);\ g, r=1,\ldots,k\biggr\},
$$

\vspace{5mm}
\noindent
${\rm l.i.m.}$ is a limit in the mean-square sense,
$i_1,\ldots,i_k=0,1,\ldots,m,$ 

\begin{equation}
\label{rr23}
\zeta_{j}^{(i)}=
\int\limits_t^T \phi_{j}(s) d{\bf w}_s^{(i)}
\end{equation} 

\vspace{3mm}
\noindent
are independent standard Gaussian random variables
for various
$i$ or $j$ {\rm(}if $i\ne 0${\rm),}
$C_{j_k\ldots j_1}$ is the Fourier coefficient {\rm(\ref{ppppa}),}
$\Delta{\bf w}_{\tau_{j}}^{(i)}=
{\bf w}_{\tau_{j+1}}^{(i)}-{\bf w}_{\tau_{j}}^{(i)}$
$(i=0,\ 1,\ldots,m),$\
$\left\{\tau_{j}\right\}_{j=0}^{N}$ is a partition of
$[t,T],$ which satisfies the condition {\rm (\ref{1111})}.}

\vspace{2mm}

In order to evaluate the significance of Theorem 7 for practice we will
demonstrate its transformed particular cases for 
$k=1,\ldots,6$ \cite{7} (2006) (also see
\cite{arxiv-1}-\cite{5-008}, \cite{8}-\cite{new-2023a}) 

\vspace{1mm}
\begin{equation}
\label{za1}
J[\psi^{(1)}]_{T,t}
=\hbox{\vtop{\offinterlineskip\halign{
\hfil#\hfil\cr
{\rm l.i.m.}\cr
$\stackrel{}{{}_{p_1\to \infty}}$\cr
}} }\sum_{j_1=0}^{p_1}
C_{j_1}\zeta_{j_1}^{(i_1)},
\end{equation}

\vspace{3mm}
\begin{equation}
\label{za2}
J[\psi^{(2)}]_{T,t}
=\hbox{\vtop{\offinterlineskip\halign{
\hfil#\hfil\cr
{\rm l.i.m.}\cr
$\stackrel{}{{}_{p_1,p_2\to \infty}}$\cr
}} }\sum_{j_1=0}^{p_1}\sum_{j_2=0}^{p_2}
C_{j_2j_1}\Biggl(\zeta_{j_1}^{(i_1)}\zeta_{j_2}^{(i_2)}
-{\bf 1}_{\{i_1=i_2\ne 0\}}
{\bf 1}_{\{j_1=j_2\}}\Biggr),
\end{equation}

\vspace{5mm}
$$
J[\psi^{(3)}]_{T,t}=
\hbox{\vtop{\offinterlineskip\halign{
\hfil#\hfil\cr
{\rm l.i.m.}\cr
$\stackrel{}{{}_{p_1,\ldots,p_3\to \infty}}$\cr
}} }\sum_{j_1=0}^{p_1}\sum_{j_2=0}^{p_2}\sum_{j_3=0}^{p_3}
C_{j_3j_2j_1}\Biggl(
\zeta_{j_1}^{(i_1)}\zeta_{j_2}^{(i_2)}\zeta_{j_3}^{(i_3)}
-\Biggr.
$$
\begin{equation}
\label{za3}
-\Biggl.
{\bf 1}_{\{i_1=i_2\ne 0\}}
{\bf 1}_{\{j_1=j_2\}}
\zeta_{j_3}^{(i_3)}
-{\bf 1}_{\{i_2=i_3\ne 0\}}
{\bf 1}_{\{j_2=j_3\}}
\zeta_{j_1}^{(i_1)}-
{\bf 1}_{\{i_1=i_3\ne 0\}}
{\bf 1}_{\{j_1=j_3\}}
\zeta_{j_2}^{(i_2)}\Biggr),
\end{equation}

\vspace{5mm}
$$
J[\psi^{(4)}]_{T,t}
=
\hbox{\vtop{\offinterlineskip\halign{
\hfil#\hfil\cr
{\rm l.i.m.}\cr
$\stackrel{}{{}_{p_1,\ldots,p_4\to \infty}}$\cr
}} }\sum_{j_1=0}^{p_1}\ldots\sum_{j_4=0}^{p_4}
C_{j_4\ldots j_1}\Biggl(
\prod_{l=1}^4\zeta_{j_l}^{(i_l)}
\Biggr.
-
$$
$$
-
{\bf 1}_{\{i_1=i_2\ne 0\}}
{\bf 1}_{\{j_1=j_2\}}
\zeta_{j_3}^{(i_3)}
\zeta_{j_4}^{(i_4)}
-
{\bf 1}_{\{i_1=i_3\ne 0\}}
{\bf 1}_{\{j_1=j_3\}}
\zeta_{j_2}^{(i_2)}
\zeta_{j_4}^{(i_4)}-
$$
$$
-
{\bf 1}_{\{i_1=i_4\ne 0\}}
{\bf 1}_{\{j_1=j_4\}}
\zeta_{j_2}^{(i_2)}
\zeta_{j_3}^{(i_3)}
-
{\bf 1}_{\{i_2=i_3\ne 0\}}
{\bf 1}_{\{j_2=j_3\}}
\zeta_{j_1}^{(i_1)}
\zeta_{j_4}^{(i_4)}-
$$
$$
-
{\bf 1}_{\{i_2=i_4\ne 0\}}
{\bf 1}_{\{j_2=j_4\}}
\zeta_{j_1}^{(i_1)}
\zeta_{j_3}^{(i_3)}
-
{\bf 1}_{\{i_3=i_4\ne 0\}}
{\bf 1}_{\{j_3=j_4\}}
\zeta_{j_1}^{(i_1)}
\zeta_{j_2}^{(i_2)}+
$$
$$
+
{\bf 1}_{\{i_1=i_2\ne 0\}}
{\bf 1}_{\{j_1=j_2\}}
{\bf 1}_{\{i_3=i_4\ne 0\}}
{\bf 1}_{\{j_3=j_4\}}
+
$$
$$
+
{\bf 1}_{\{i_1=i_3\ne 0\}}
{\bf 1}_{\{j_1=j_3\}}
{\bf 1}_{\{i_2=i_4\ne 0\}}
{\bf 1}_{\{j_2=j_4\}}+
$$
\begin{equation}
\label{za4}
+\Biggl.
{\bf 1}_{\{i_1=i_4\ne 0\}}
{\bf 1}_{\{j_1=j_4\}}
{\bf 1}_{\{i_2=i_3\ne 0\}}
{\bf 1}_{\{j_2=j_3\}}\Biggr),
\end{equation}

\vspace{6mm}
$$
J[\psi^{(5)}]_{T,t}
=\hbox{\vtop{\offinterlineskip\halign{
\hfil#\hfil\cr
{\rm l.i.m.}\cr
$\stackrel{}{{}_{p_1,\ldots,p_5\to \infty}}$\cr
}} }\sum_{j_1=0}^{p_1}\ldots\sum_{j_5=0}^{p_5}
C_{j_5\ldots j_1}\Biggl(
\prod_{l=1}^5\zeta_{j_l}^{(i_l)}
-\Biggr.
$$
$$
-
{\bf 1}_{\{i_1=i_2\ne 0\}}
{\bf 1}_{\{j_1=j_2\}}
\zeta_{j_3}^{(i_3)}
\zeta_{j_4}^{(i_4)}
\zeta_{j_5}^{(i_5)}-
{\bf 1}_{\{i_1=i_3\ne 0\}}
{\bf 1}_{\{j_1=j_3\}}
\zeta_{j_2}^{(i_2)}
\zeta_{j_4}^{(i_4)}
\zeta_{j_5}^{(i_5)}-
$$
$$
-
{\bf 1}_{\{i_1=i_4\ne 0\}}
{\bf 1}_{\{j_1=j_4\}}
\zeta_{j_2}^{(i_2)}
\zeta_{j_3}^{(i_3)}
\zeta_{j_5}^{(i_5)}-
{\bf 1}_{\{i_1=i_5\ne 0\}}
{\bf 1}_{\{j_1=j_5\}}
\zeta_{j_2}^{(i_2)}
\zeta_{j_3}^{(i_3)}
\zeta_{j_4}^{(i_4)}-
$$
$$
-
{\bf 1}_{\{i_2=i_3\ne 0\}}
{\bf 1}_{\{j_2=j_3\}}
\zeta_{j_1}^{(i_1)}
\zeta_{j_4}^{(i_4)}
\zeta_{j_5}^{(i_5)}-
{\bf 1}_{\{i_2=i_4\ne 0\}}
{\bf 1}_{\{j_2=j_4\}}
\zeta_{j_1}^{(i_1)}
\zeta_{j_3}^{(i_3)}
\zeta_{j_5}^{(i_5)}-
$$
$$
-
{\bf 1}_{\{i_2=i_5\ne 0\}}
{\bf 1}_{\{j_2=j_5\}}
\zeta_{j_1}^{(i_1)}
\zeta_{j_3}^{(i_3)}
\zeta_{j_4}^{(i_4)}
-{\bf 1}_{\{i_3=i_4\ne 0\}}
{\bf 1}_{\{j_3=j_4\}}
\zeta_{j_1}^{(i_1)}
\zeta_{j_2}^{(i_2)}
\zeta_{j_5}^{(i_5)}-
$$
$$
-
{\bf 1}_{\{i_3=i_5\ne 0\}}
{\bf 1}_{\{j_3=j_5\}}
\zeta_{j_1}^{(i_1)}
\zeta_{j_2}^{(i_2)}
\zeta_{j_4}^{(i_4)}
-{\bf 1}_{\{i_4=i_5\ne 0\}}
{\bf 1}_{\{j_4=j_5\}}
\zeta_{j_1}^{(i_1)}
\zeta_{j_2}^{(i_2)}
\zeta_{j_3}^{(i_3)}+
$$
$$
+
{\bf 1}_{\{i_1=i_2\ne 0\}}
{\bf 1}_{\{j_1=j_2\}}
{\bf 1}_{\{i_3=i_4\ne 0\}}
{\bf 1}_{\{j_3=j_4\}}\zeta_{j_5}^{(i_5)}+
{\bf 1}_{\{i_1=i_2\ne 0\}}
{\bf 1}_{\{j_1=j_2\}}
{\bf 1}_{\{i_3=i_5\ne 0\}}
{\bf 1}_{\{j_3=j_5\}}\zeta_{j_4}^{(i_4)}+
$$
$$
+
{\bf 1}_{\{i_1=i_2\ne 0\}}
{\bf 1}_{\{j_1=j_2\}}
{\bf 1}_{\{i_4=i_5\ne 0\}}
{\bf 1}_{\{j_4=j_5\}}\zeta_{j_3}^{(i_3)}+
{\bf 1}_{\{i_1=i_3\ne 0\}}
{\bf 1}_{\{j_1=j_3\}}
{\bf 1}_{\{i_2=i_4\ne 0\}}
{\bf 1}_{\{j_2=j_4\}}\zeta_{j_5}^{(i_5)}+
$$
$$
+
{\bf 1}_{\{i_1=i_3\ne 0\}}
{\bf 1}_{\{j_1=j_3\}}
{\bf 1}_{\{i_2=i_5\ne 0\}}
{\bf 1}_{\{j_2=j_5\}}\zeta_{j_4}^{(i_4)}+
{\bf 1}_{\{i_1=i_3\ne 0\}}
{\bf 1}_{\{j_1=j_3\}}
{\bf 1}_{\{i_4=i_5\ne 0\}}
{\bf 1}_{\{j_4=j_5\}}\zeta_{j_2}^{(i_2)}+
$$
$$
+
{\bf 1}_{\{i_1=i_4\ne 0\}}
{\bf 1}_{\{j_1=j_4\}}
{\bf 1}_{\{i_2=i_3\ne 0\}}
{\bf 1}_{\{j_2=j_3\}}\zeta_{j_5}^{(i_5)}+
{\bf 1}_{\{i_1=i_4\ne 0\}}
{\bf 1}_{\{j_1=j_4\}}
{\bf 1}_{\{i_2=i_5\ne 0\}}
{\bf 1}_{\{j_2=j_5\}}\zeta_{j_3}^{(i_3)}+
$$
$$
+
{\bf 1}_{\{i_1=i_4\ne 0\}}
{\bf 1}_{\{j_1=j_4\}}
{\bf 1}_{\{i_3=i_5\ne 0\}}
{\bf 1}_{\{j_3=j_5\}}\zeta_{j_2}^{(i_2)}+
{\bf 1}_{\{i_1=i_5\ne 0\}}
{\bf 1}_{\{j_1=j_5\}}
{\bf 1}_{\{i_2=i_3\ne 0\}}
{\bf 1}_{\{j_2=j_3\}}\zeta_{j_4}^{(i_4)}+
$$
$$
+
{\bf 1}_{\{i_1=i_5\ne 0\}}
{\bf 1}_{\{j_1=j_5\}}
{\bf 1}_{\{i_2=i_4\ne 0\}}
{\bf 1}_{\{j_2=j_4\}}\zeta_{j_3}^{(i_3)}+
{\bf 1}_{\{i_1=i_5\ne 0\}}
{\bf 1}_{\{j_1=j_5\}}
{\bf 1}_{\{i_3=i_4\ne 0\}}
{\bf 1}_{\{j_3=j_4\}}\zeta_{j_2}^{(i_2)}+
$$
$$
+
{\bf 1}_{\{i_2=i_3\ne 0\}}
{\bf 1}_{\{j_2=j_3\}}
{\bf 1}_{\{i_4=i_5\ne 0\}}
{\bf 1}_{\{j_4=j_5\}}\zeta_{j_1}^{(i_1)}+
{\bf 1}_{\{i_2=i_4\ne 0\}}
{\bf 1}_{\{j_2=j_4\}}
{\bf 1}_{\{i_3=i_5\ne 0\}}
{\bf 1}_{\{j_3=j_5\}}\zeta_{j_1}^{(i_1)}+
$$
\begin{equation}
\label{za5}
+\Biggl.
{\bf 1}_{\{i_2=i_5\ne 0\}}
{\bf 1}_{\{j_2=j_5\}}
{\bf 1}_{\{i_3=i_4\ne 0\}}
{\bf 1}_{\{j_3=j_4\}}\zeta_{j_1}^{(i_1)}\Biggr),
\end{equation}

\vspace{9mm}

$$
J[\psi^{(6)}]_{T,t}
=\hbox{\vtop{\offinterlineskip\halign{
\hfil#\hfil\cr
{\rm l.i.m.}\cr
$\stackrel{}{{}_{p_1,\ldots,p_6\to \infty}}$\cr
}} }\sum_{j_1=0}^{p_1}\ldots\sum_{j_6=0}^{p_6}
C_{j_6\ldots j_1}\Biggl(
\prod_{l=1}^6
\zeta_{j_l}^{(i_l)}
-\Biggr.
$$
$$
-
{\bf 1}_{\{i_1=i_6\ne 0\}}
{\bf 1}_{\{j_1=j_6\}}
\zeta_{j_2}^{(i_2)}
\zeta_{j_3}^{(i_3)}
\zeta_{j_4}^{(i_4)}
\zeta_{j_5}^{(i_5)}-
{\bf 1}_{\{i_2=i_6\ne 0\}}
{\bf 1}_{\{j_2=j_6\}}
\zeta_{j_1}^{(i_1)}
\zeta_{j_3}^{(i_3)}
\zeta_{j_4}^{(i_4)}
\zeta_{j_5}^{(i_5)}-
$$
$$
-
{\bf 1}_{\{i_3=i_6\ne 0\}}
{\bf 1}_{\{j_3=j_6\}}
\zeta_{j_1}^{(i_1)}
\zeta_{j_2}^{(i_2)}
\zeta_{j_4}^{(i_4)}
\zeta_{j_5}^{(i_5)}-
{\bf 1}_{\{i_4=i_6\ne 0\}}
{\bf 1}_{\{j_4=j_6\}}
\zeta_{j_1}^{(i_1)}
\zeta_{j_2}^{(i_2)}
\zeta_{j_3}^{(i_3)}
\zeta_{j_5}^{(i_5)}-
$$
$$
-
{\bf 1}_{\{i_5=i_6\ne 0\}}
{\bf 1}_{\{j_5=j_6\}}
\zeta_{j_1}^{(i_1)}
\zeta_{j_2}^{(i_2)}
\zeta_{j_3}^{(i_3)}
\zeta_{j_4}^{(i_4)}-
{\bf 1}_{\{i_1=i_2\ne 0\}}
{\bf 1}_{\{j_1=j_2\}}
\zeta_{j_3}^{(i_3)}
\zeta_{j_4}^{(i_4)}
\zeta_{j_5}^{(i_5)}
\zeta_{j_6}^{(i_6)}-
$$
$$
-
{\bf 1}_{\{i_1=i_3\ne 0\}}
{\bf 1}_{\{j_1=j_3\}}
\zeta_{j_2}^{(i_2)}
\zeta_{j_4}^{(i_4)}
\zeta_{j_5}^{(i_5)}
\zeta_{j_6}^{(i_6)}-
{\bf 1}_{\{i_1=i_4\ne 0\}}
{\bf 1}_{\{j_1=j_4\}}
\zeta_{j_2}^{(i_2)}
\zeta_{j_3}^{(i_3)}
\zeta_{j_5}^{(i_5)}
\zeta_{j_6}^{(i_6)}-
$$
$$
-
{\bf 1}_{\{i_1=i_5\ne 0\}}
{\bf 1}_{\{j_1=j_5\}}
\zeta_{j_2}^{(i_2)}
\zeta_{j_3}^{(i_3)}
\zeta_{j_4}^{(i_4)}
\zeta_{j_6}^{(i_6)}-
{\bf 1}_{\{i_2=i_3\ne 0\}}
{\bf 1}_{\{j_2=j_3\}}
\zeta_{j_1}^{(i_1)}
\zeta_{j_4}^{(i_4)}
\zeta_{j_5}^{(i_5)}
\zeta_{j_6}^{(i_6)}-
$$
$$
-
{\bf 1}_{\{i_2=i_4\ne 0\}}
{\bf 1}_{\{j_2=j_4\}}
\zeta_{j_1}^{(i_1)}
\zeta_{j_3}^{(i_3)}
\zeta_{j_5}^{(i_5)}
\zeta_{j_6}^{(i_6)}-
{\bf 1}_{\{i_2=i_5\ne 0\}}
{\bf 1}_{\{j_2=j_5\}}
\zeta_{j_1}^{(i_1)}
\zeta_{j_3}^{(i_3)}
\zeta_{j_4}^{(i_4)}
\zeta_{j_6}^{(i_6)}-
$$
$$
-
{\bf 1}_{\{i_3=i_4\ne 0\}}
{\bf 1}_{\{j_3=j_4\}}
\zeta_{j_1}^{(i_1)}
\zeta_{j_2}^{(i_2)}
\zeta_{j_5}^{(i_5)}
\zeta_{j_6}^{(i_6)}-
{\bf 1}_{\{i_3=i_5\ne 0\}}
{\bf 1}_{\{j_3=j_5\}}
\zeta_{j_1}^{(i_1)}
\zeta_{j_2}^{(i_2)}
\zeta_{j_4}^{(i_4)}
\zeta_{j_6}^{(i_6)}-
$$
$$
-
{\bf 1}_{\{i_4=i_5\ne 0\}}
{\bf 1}_{\{j_4=j_5\}}
\zeta_{j_1}^{(i_1)}
\zeta_{j_2}^{(i_2)}
\zeta_{j_3}^{(i_3)}
\zeta_{j_6}^{(i_6)}+
$$
$$
+
{\bf 1}_{\{i_1=i_2\ne 0\}}
{\bf 1}_{\{j_1=j_2\}}
{\bf 1}_{\{i_3=i_4\ne 0\}}
{\bf 1}_{\{j_3=j_4\}}
\zeta_{j_5}^{(i_5)}
\zeta_{j_6}^{(i_6)}+
{\bf 1}_{\{i_1=i_2\ne 0\}}
{\bf 1}_{\{j_1=j_2\}}
{\bf 1}_{\{i_3=i_5\ne 0\}}
{\bf 1}_{\{j_3=j_5\}}
\zeta_{j_4}^{(i_4)}
\zeta_{j_6}^{(i_6)}+
$$
$$
+
{\bf 1}_{\{i_1=i_2\ne 0\}}
{\bf 1}_{\{j_1=j_2\}}
{\bf 1}_{\{i_4=i_5\ne 0\}}
{\bf 1}_{\{j_4=j_5\}}
\zeta_{j_3}^{(i_3)}
\zeta_{j_6}^{(i_6)}
+
{\bf 1}_{\{i_1=i_3\ne 0\}}
{\bf 1}_{\{j_1=j_3\}}
{\bf 1}_{\{i_2=i_4\ne 0\}}
{\bf 1}_{\{j_2=j_4\}}
\zeta_{j_5}^{(i_5)}
\zeta_{j_6}^{(i_6)}+
$$
$$
+
{\bf 1}_{\{i_1=i_3\ne 0\}}
{\bf 1}_{\{j_1=j_3\}}
{\bf 1}_{\{i_2=i_5\ne 0\}}
{\bf 1}_{\{j_2=j_5\}}
\zeta_{j_4}^{(i_4)}
\zeta_{j_6}^{(i_6)}
+{\bf 1}_{\{i_1=i_3\ne 0\}}
{\bf 1}_{\{j_1=j_3\}}
{\bf 1}_{\{i_4=i_5\ne 0\}}
{\bf 1}_{\{j_4=j_5\}}
\zeta_{j_2}^{(i_2)}
\zeta_{j_6}^{(i_6)}+
$$
$$
+
{\bf 1}_{\{i_1=i_4\ne 0\}}
{\bf 1}_{\{j_1=j_4\}}
{\bf 1}_{\{i_2=i_3\ne 0\}}
{\bf 1}_{\{j_2=j_3\}}
\zeta_{j_5}^{(i_5)}
\zeta_{j_6}^{(i_6)}
+
{\bf 1}_{\{i_1=i_4\ne 0\}}
{\bf 1}_{\{j_1=j_4\}}
{\bf 1}_{\{i_2=i_5\ne 0\}}
{\bf 1}_{\{j_2=j_5\}}
\zeta_{j_3}^{(i_3)}
\zeta_{j_6}^{(i_6)}+
$$
$$
+
{\bf 1}_{\{i_1=i_4\ne 0\}}
{\bf 1}_{\{j_1=j_4\}}
{\bf 1}_{\{i_3=i_5\ne 0\}}
{\bf 1}_{\{j_3=j_5\}}
\zeta_{j_2}^{(i_2)}
\zeta_{j_6}^{(i_6)}
+
{\bf 1}_{\{i_1=i_5\ne 0\}}
{\bf 1}_{\{j_1=j_5\}}
{\bf 1}_{\{i_2=i_3\ne 0\}}
{\bf 1}_{\{j_2=j_3\}}
\zeta_{j_4}^{(i_4)}
\zeta_{j_6}^{(i_6)}+
$$
$$
+
{\bf 1}_{\{i_1=i_5\ne 0\}}
{\bf 1}_{\{j_1=j_5\}}
{\bf 1}_{\{i_2=i_4\ne 0\}}
{\bf 1}_{\{j_2=j_4\}}
\zeta_{j_3}^{(i_3)}
\zeta_{j_6}^{(i_6)}
+
{\bf 1}_{\{i_1=i_5\ne 0\}}
{\bf 1}_{\{j_1=j_5\}}
{\bf 1}_{\{i_3=i_4\ne 0\}}
{\bf 1}_{\{j_3=j_4\}}
\zeta_{j_2}^{(i_2)}
\zeta_{j_6}^{(i_6)}+
$$
$$
+
{\bf 1}_{\{i_2=i_3\ne 0\}}
{\bf 1}_{\{j_2=j_3\}}
{\bf 1}_{\{i_4=i_5\ne 0\}}
{\bf 1}_{\{j_4=j_5\}}
\zeta_{j_1}^{(i_1)}
\zeta_{j_6}^{(i_6)}
+
{\bf 1}_{\{i_2=i_4\ne 0\}}
{\bf 1}_{\{j_2=j_4\}}
{\bf 1}_{\{i_3=i_5\ne 0\}}
{\bf 1}_{\{j_3=j_5\}}
\zeta_{j_1}^{(i_1)}
\zeta_{j_6}^{(i_6)}+
$$
$$
+
{\bf 1}_{\{i_2=i_5\ne 0\}}
{\bf 1}_{\{j_2=j_5\}}
{\bf 1}_{\{i_3=i_4\ne 0\}}
{\bf 1}_{\{j_3=j_4\}}
\zeta_{j_1}^{(i_1)}
\zeta_{j_6}^{(i_6)}
+
{\bf 1}_{\{i_6=i_1\ne 0\}}
{\bf 1}_{\{j_6=j_1\}}
{\bf 1}_{\{i_3=i_4\ne 0\}}
{\bf 1}_{\{j_3=j_4\}}
\zeta_{j_2}^{(i_2)}
\zeta_{j_5}^{(i_5)}+
$$
$$
+
{\bf 1}_{\{i_6=i_1\ne 0\}}
{\bf 1}_{\{j_6=j_1\}}
{\bf 1}_{\{i_3=i_5\ne 0\}}
{\bf 1}_{\{j_3=j_5\}}
\zeta_{j_2}^{(i_2)}
\zeta_{j_4}^{(i_4)}
+
{\bf 1}_{\{i_6=i_1\ne 0\}}
{\bf 1}_{\{j_6=j_1\}}
{\bf 1}_{\{i_2=i_5\ne 0\}}
{\bf 1}_{\{j_2=j_5\}}
\zeta_{j_3}^{(i_3)}
\zeta_{j_4}^{(i_4)}+
$$
$$
+
{\bf 1}_{\{i_6=i_1\ne 0\}}
{\bf 1}_{\{j_6=j_1\}}
{\bf 1}_{\{i_2=i_4\ne 0\}}
{\bf 1}_{\{j_2=j_4\}}
\zeta_{j_3}^{(i_3)}
\zeta_{j_5}^{(i_5)}
+
{\bf 1}_{\{i_6=i_1\ne 0\}}
{\bf 1}_{\{j_6=j_1\}}
{\bf 1}_{\{i_4=i_5\ne 0\}}
{\bf 1}_{\{j_4=j_5\}}
\zeta_{j_2}^{(i_2)}
\zeta_{j_3}^{(i_3)}+
$$
$$
+
{\bf 1}_{\{i_6=i_1\ne 0\}}
{\bf 1}_{\{j_6=j_1\}}
{\bf 1}_{\{i_2=i_3\ne 0\}}
{\bf 1}_{\{j_2=j_3\}}
\zeta_{j_4}^{(i_4)}
\zeta_{j_5}^{(i_5)}
+
{\bf 1}_{\{i_6=i_2\ne 0\}}
{\bf 1}_{\{j_6=j_2\}}
{\bf 1}_{\{i_3=i_5\ne 0\}}
{\bf 1}_{\{j_3=j_5\}}
\zeta_{j_1}^{(i_1)}
\zeta_{j_4}^{(i_4)}+
$$
$$
+
{\bf 1}_{\{i_6=i_2\ne 0\}}
{\bf 1}_{\{j_6=j_2\}}
{\bf 1}_{\{i_4=i_5\ne 0\}}
{\bf 1}_{\{j_4=j_5\}}
\zeta_{j_1}^{(i_1)}
\zeta_{j_3}^{(i_3)}
+
{\bf 1}_{\{i_6=i_2\ne 0\}}
{\bf 1}_{\{j_6=j_2\}}
{\bf 1}_{\{i_3=i_4\ne 0\}}
{\bf 1}_{\{j_3=j_4\}}
\zeta_{j_1}^{(i_1)}
\zeta_{j_5}^{(i_5)}+
$$
$$
+
{\bf 1}_{\{i_6=i_2\ne 0\}}
{\bf 1}_{\{j_6=j_2\}}
{\bf 1}_{\{i_1=i_5\ne 0\}}
{\bf 1}_{\{j_1=j_5\}}
\zeta_{j_3}^{(i_3)}
\zeta_{j_4}^{(i_4)}
+
{\bf 1}_{\{i_6=i_2\ne 0\}}
{\bf 1}_{\{j_6=j_2\}}
{\bf 1}_{\{i_1=i_4\ne 0\}}
{\bf 1}_{\{j_1=j_4\}}
\zeta_{j_3}^{(i_3)}
\zeta_{j_5}^{(i_5)}+
$$
$$
+
{\bf 1}_{\{i_6=i_2\ne 0\}}
{\bf 1}_{\{j_6=j_2\}}
{\bf 1}_{\{i_1=i_3\ne 0\}}
{\bf 1}_{\{j_1=j_3\}}
\zeta_{j_4}^{(i_4)}
\zeta_{j_5}^{(i_5)}
+
{\bf 1}_{\{i_6=i_3\ne 0\}}
{\bf 1}_{\{j_6=j_3\}}
{\bf 1}_{\{i_2=i_5\ne 0\}}
{\bf 1}_{\{j_2=j_5\}}
\zeta_{j_1}^{(i_1)}
\zeta_{j_4}^{(i_4)}+
$$
$$
+
{\bf 1}_{\{i_6=i_3\ne 0\}}
{\bf 1}_{\{j_6=j_3\}}
{\bf 1}_{\{i_4=i_5\ne 0\}}
{\bf 1}_{\{j_4=j_5\}}
\zeta_{j_1}^{(i_1)}
\zeta_{j_2}^{(i_2)}
+
{\bf 1}_{\{i_6=i_3\ne 0\}}
{\bf 1}_{\{j_6=j_3\}}
{\bf 1}_{\{i_2=i_4\ne 0\}}
{\bf 1}_{\{j_2=j_4\}}
\zeta_{j_1}^{(i_1)}
\zeta_{j_5}^{(i_5)}+
$$
$$
+
{\bf 1}_{\{i_6=i_3\ne 0\}}
{\bf 1}_{\{j_6=j_3\}}
{\bf 1}_{\{i_1=i_5\ne 0\}}
{\bf 1}_{\{j_1=j_5\}}
\zeta_{j_2}^{(i_2)}
\zeta_{j_4}^{(i_4)}
+
{\bf 1}_{\{i_6=i_3\ne 0\}}
{\bf 1}_{\{j_6=j_3\}}
{\bf 1}_{\{i_1=i_4\ne 0\}}
{\bf 1}_{\{j_1=j_4\}}
\zeta_{j_2}^{(i_2)}
\zeta_{j_5}^{(i_5)}+
$$
$$
+
{\bf 1}_{\{i_6=i_3\ne 0\}}
{\bf 1}_{\{j_6=j_3\}}
{\bf 1}_{\{i_1=i_2\ne 0\}}
{\bf 1}_{\{j_1=j_2\}}
\zeta_{j_4}^{(i_4)}
\zeta_{j_5}^{(i_5)}
+
{\bf 1}_{\{i_6=i_4\ne 0\}}
{\bf 1}_{\{j_6=j_4\}}
{\bf 1}_{\{i_3=i_5\ne 0\}}
{\bf 1}_{\{j_3=j_5\}}
\zeta_{j_1}^{(i_1)}
\zeta_{j_2}^{(i_2)}+
$$
$$
+
{\bf 1}_{\{i_6=i_4\ne 0\}}
{\bf 1}_{\{j_6=j_4\}}
{\bf 1}_{\{i_2=i_5\ne 0\}}
{\bf 1}_{\{j_2=j_5\}}
\zeta_{j_1}^{(i_1)}
\zeta_{j_3}^{(i_3)}
+
{\bf 1}_{\{i_6=i_4\ne 0\}}
{\bf 1}_{\{j_6=j_4\}}
{\bf 1}_{\{i_2=i_3\ne 0\}}
{\bf 1}_{\{j_2=j_3\}}
\zeta_{j_1}^{(i_1)}
\zeta_{j_5}^{(i_5)}+
$$
$$
+
{\bf 1}_{\{i_6=i_4\ne 0\}}
{\bf 1}_{\{j_6=j_4\}}
{\bf 1}_{\{i_1=i_5\ne 0\}}
{\bf 1}_{\{j_1=j_5\}}
\zeta_{j_2}^{(i_2)}
\zeta_{j_3}^{(i_3)}
+
{\bf 1}_{\{i_6=i_4\ne 0\}}
{\bf 1}_{\{j_6=j_4\}}
{\bf 1}_{\{i_1=i_3\ne 0\}}
{\bf 1}_{\{j_1=j_3\}}
\zeta_{j_2}^{(i_2)}
\zeta_{j_5}^{(i_5)}+
$$
$$
+
{\bf 1}_{\{i_6=i_4\ne 0\}}
{\bf 1}_{\{j_6=j_4\}}
{\bf 1}_{\{i_1=i_2\ne 0\}}
{\bf 1}_{\{j_1=j_2\}}
\zeta_{j_3}^{(i_3)}
\zeta_{j_5}^{(i_5)}
+
{\bf 1}_{\{i_6=i_5\ne 0\}}
{\bf 1}_{\{j_6=j_5\}}
{\bf 1}_{\{i_3=i_4\ne 0\}}
{\bf 1}_{\{j_3=j_4\}}
\zeta_{j_1}^{(i_1)}
\zeta_{j_2}^{(i_2)}+
$$
$$
+
{\bf 1}_{\{i_6=i_5\ne 0\}}
{\bf 1}_{\{j_6=j_5\}}
{\bf 1}_{\{i_2=i_4\ne 0\}}
{\bf 1}_{\{j_2=j_4\}}
\zeta_{j_1}^{(i_1)}
\zeta_{j_3}^{(i_3)}
+
{\bf 1}_{\{i_6=i_5\ne 0\}}
{\bf 1}_{\{j_6=j_5\}}
{\bf 1}_{\{i_2=i_3\ne 0\}}
{\bf 1}_{\{j_2=j_3\}}
\zeta_{j_1}^{(i_1)}
\zeta_{j_4}^{(i_4)}+
$$
$$
+
{\bf 1}_{\{i_6=i_5\ne 0\}}
{\bf 1}_{\{j_6=j_5\}}
{\bf 1}_{\{i_1=i_4\ne 0\}}
{\bf 1}_{\{j_1=j_4\}}
\zeta_{j_2}^{(i_2)}
\zeta_{j_3}^{(i_3)}
+
{\bf 1}_{\{i_6=i_5\ne 0\}}
{\bf 1}_{\{j_6=j_5\}}
{\bf 1}_{\{i_1=i_3\ne 0\}}
{\bf 1}_{\{j_1=j_3\}}
\zeta_{j_2}^{(i_2)}
\zeta_{j_4}^{(i_4)}+
$$
$$
+
{\bf 1}_{\{i_6=i_5\ne 0\}}
{\bf 1}_{\{j_6=j_5\}}
{\bf 1}_{\{i_1=i_2\ne 0\}}
{\bf 1}_{\{j_1=j_2\}}
\zeta_{j_3}^{(i_3)}
\zeta_{j_4}^{(i_4)}-
$$
$$
-
{\bf 1}_{\{i_6=i_1\ne 0\}}
{\bf 1}_{\{j_6=j_1\}}
{\bf 1}_{\{i_2=i_5\ne 0\}}
{\bf 1}_{\{j_2=j_5\}}
{\bf 1}_{\{i_3=i_4\ne 0\}}
{\bf 1}_{\{j_3=j_4\}}-
$$
$$
-
{\bf 1}_{\{i_6=i_1\ne 0\}}
{\bf 1}_{\{j_6=j_1\}}
{\bf 1}_{\{i_2=i_4\ne 0\}}
{\bf 1}_{\{j_2=j_4\}}
{\bf 1}_{\{i_3=i_5\ne 0\}}
{\bf 1}_{\{j_3=j_5\}}-
$$
$$
-
{\bf 1}_{\{i_6=i_1\ne 0\}}
{\bf 1}_{\{j_6=j_1\}}
{\bf 1}_{\{i_2=i_3\ne 0\}}
{\bf 1}_{\{j_2=j_3\}}
{\bf 1}_{\{i_4=i_5\ne 0\}}
{\bf 1}_{\{j_4=j_5\}}-
$$
$$
-
{\bf 1}_{\{i_6=i_2\ne 0\}}
{\bf 1}_{\{j_6=j_2\}}
{\bf 1}_{\{i_1=i_5\ne 0\}}
{\bf 1}_{\{j_1=j_5\}}
{\bf 1}_{\{i_3=i_4\ne 0\}}
{\bf 1}_{\{j_3=j_4\}}-
$$
$$
-
{\bf 1}_{\{i_6=i_2\ne 0\}}
{\bf 1}_{\{j_6=j_2\}}
{\bf 1}_{\{i_1=i_4\ne 0\}}
{\bf 1}_{\{j_1=j_4\}}
{\bf 1}_{\{i_3=i_5\ne 0\}}
{\bf 1}_{\{j_3=j_5\}}-
$$
$$
-
{\bf 1}_{\{i_6=i_2\ne 0\}}
{\bf 1}_{\{j_6=j_2\}}
{\bf 1}_{\{i_1=i_3\ne 0\}}
{\bf 1}_{\{j_1=j_3\}}
{\bf 1}_{\{i_4=i_5\ne 0\}}
{\bf 1}_{\{j_4=j_5\}}-
$$
$$
-
{\bf 1}_{\{i_6=i_3\ne 0\}}
{\bf 1}_{\{j_6=j_3\}}
{\bf 1}_{\{i_1=i_5\ne 0\}}
{\bf 1}_{\{j_1=j_5\}}
{\bf 1}_{\{i_2=i_4\ne 0\}}
{\bf 1}_{\{j_2=j_4\}}-
$$
$$
-
{\bf 1}_{\{i_6=i_3\ne 0\}}
{\bf 1}_{\{j_6=j_3\}}
{\bf 1}_{\{i_1=i_4\ne 0\}}
{\bf 1}_{\{j_1=j_4\}}
{\bf 1}_{\{i_2=i_5\ne 0\}}
{\bf 1}_{\{j_2=j_5\}}-
$$
$$
-
{\bf 1}_{\{i_3=i_6\ne 0\}}
{\bf 1}_{\{j_3=j_6\}}
{\bf 1}_{\{i_1=i_2\ne 0\}}
{\bf 1}_{\{j_1=j_2\}}
{\bf 1}_{\{i_4=i_5\ne 0\}}
{\bf 1}_{\{j_4=j_5\}}-
$$
$$
-
{\bf 1}_{\{i_6=i_4\ne 0\}}
{\bf 1}_{\{j_6=j_4\}}
{\bf 1}_{\{i_1=i_5\ne 0\}}
{\bf 1}_{\{j_1=j_5\}}
{\bf 1}_{\{i_2=i_3\ne 0\}}
{\bf 1}_{\{j_2=j_3\}}-
$$
$$
-
{\bf 1}_{\{i_6=i_4\ne 0\}}
{\bf 1}_{\{j_6=j_4\}}
{\bf 1}_{\{i_1=i_3\ne 0\}}
{\bf 1}_{\{j_1=j_3\}}
{\bf 1}_{\{i_2=i_5\ne 0\}}
{\bf 1}_{\{j_2=j_5\}}-
$$
$$
-
{\bf 1}_{\{i_6=i_4\ne 0\}}
{\bf 1}_{\{j_6=j_4\}}
{\bf 1}_{\{i_1=i_2\ne 0\}}
{\bf 1}_{\{j_1=j_2\}}
{\bf 1}_{\{i_3=i_5\ne 0\}}
{\bf 1}_{\{j_3=j_5\}}-
$$
$$
-
{\bf 1}_{\{i_6=i_5\ne 0\}}
{\bf 1}_{\{j_6=j_5\}}
{\bf 1}_{\{i_1=i_4\ne 0\}}
{\bf 1}_{\{j_1=j_4\}}
{\bf 1}_{\{i_2=i_3\ne 0\}}
{\bf 1}_{\{j_2=j_3\}}-
$$
$$
-
{\bf 1}_{\{i_6=i_5\ne 0\}}
{\bf 1}_{\{j_6=j_5\}}
{\bf 1}_{\{i_1=i_2\ne 0\}}
{\bf 1}_{\{j_1=j_2\}}
{\bf 1}_{\{i_3=i_4\ne 0\}}
{\bf 1}_{\{j_3=j_4\}}-
$$
\begin{equation}
\label{a6}
\Biggl.-
{\bf 1}_{\{i_6=i_5\ne 0\}}
{\bf 1}_{\{j_6=j_5\}}
{\bf 1}_{\{i_1=i_3\ne 0\}}
{\bf 1}_{\{j_1=j_3\}}
{\bf 1}_{\{i_2=i_4\ne 0\}}
{\bf 1}_{\{j_2=j_4\}}\Biggr),
\end{equation}

\vspace{5mm}
\noindent
where ${\bf 1}_A$ is the indicator of the set $A$.

It was shown in \cite{arxiv-1} (Sect.~6, 15, 16), \cite{20aa} (Sect.~1.1.9, 1.11, 1.12) that 
Theorem 7 is valid for convergence 
in the mean of degree $2n$ ($n\in \mathbb{N}$).
The convergence w.~p.~1 in Theorem~7 is proved in 
\cite{20aa} (Sect.~1.7.2)
for the complete orthonormal systems of Legendre polynomials 
and trigonometric functions in the space $L_2([t, T]).$
Moreover, the complete orthonormal systems of Haar and 
Rademacher--Walsh functions in $L_2([t,T])$ 
can also be applied in Theorem 7
\cite{arxiv-1}, \cite{9}-\cite{12aa-afterxxx}.
The generalization of Theorem 7 for 
complete orthonormal with weigth $r(t_1)\ldots r(t_k)\ge 0$ systems
of functions in the space $L_2([t,T]^k)$ can be found in 
\cite{arxiv-13}, \cite{20}-\cite{12aa-afterxxx}.
Another modification of Theorem 7 is connected with the mean-square
approximation of iterated stochastic integrals
with respect to the infinite-dimensional
$Q$-Wiener process \cite{arxiv-20}, \cite{arxiv-21}, 
\cite{5-008}. The latters play the key role
for implemetation of the high-order
strong numerical methods
for non-commutative semilinear
stochastic partial differential equations with 
multiplicative trace class noise \cite{xx1}-\cite{xx4}.

For further consideration, let us 
consider the generalization of formulas (\ref{za1})--(\ref{a6}) 
for the case of an arbitrary multiplicity $k$ $(k\in\mathbb{N})$ of 
the iterated Ito stochastic integral $J[\psi^{(k)}]_{T,t}$ defined by (\ref{ito}).
In order to do this, let us
introduce some notations. 
Let us consider the unordered
set $\{1, 2, \ldots, k\}$ 
and separate it into two parts:
the first part consists of $r$ unordered 
pairs (sequence order of these pairs is also unimportant) and the 
second one consists of the 
remaining $k-2r$ numbers.
So, we have

\begin{equation}
\label{leto5007}
(\{
\underbrace{\{g_1, g_2\}, \ldots, 
\{g_{2r-1}, g_{2r}\}}_{\small{\hbox{part 1}}}
\},
\{\underbrace{q_1, \ldots, q_{k-2r}}_{\small{\hbox{part 2}}}
\}),
\end{equation}

\vspace{3mm}
\noindent
where 
$$
\{g_1, g_2, \ldots, 
g_{2r-1}, g_{2r}, q_1, \ldots, q_{k-2r}\}=\{1, 2, \ldots, k\},
$$

\vspace{3mm}
\noindent
braces   
mean an unordered 
set, and pa\-ren\-the\-ses mean an ordered set.

We will say that (\ref{leto5007}) is a partition 
and consider the sum with respect to all possible
partitions

\begin{equation}
\label{leto5008}
\sum_{\stackrel{(\{\{g_1, g_2\}, \ldots, 
\{g_{2r-1}, g_{2r}\}\}, \{q_1, \ldots, q_{k-2r}\})}
{{}_{\{g_1, g_2, \ldots, 
g_{2r-1}, g_{2r}, q_1, \ldots, q_{k-2r}\}=\{1, 2, \ldots, k\}}}}
a_{g_1 g_2, \ldots, 
g_{2r-1} g_{2r}, q_1 \ldots q_{k-2r}},
\end{equation}

\vspace{4mm}
\noindent
where $a_{g_1 g_2, \ldots, 
g_{2r-1} g_{2r}, q_1 \ldots q_{k-2r}}\in \mathbb{R}^1.$

Below there are several examples of sums in the form (\ref{leto5008})

\vspace{2mm}
$$
\sum_{\stackrel{(\{g_1, g_2\})}{{}_{\{g_1, g_2\}=\{1, 2\}}}}
a_{g_1 g_2}=a_{12},
$$

\vspace{3mm}
$$
\sum_{\stackrel{(\{\{g_1, g_2\}, \{g_3, g_4\}\})}
{{}_{\{g_1, g_2, g_3, g_4\}=\{1, 2, 3, 4\}}}}
a_{g_1 g_2, g_3 g_4}=a_{12,34} + a_{13,24} + a_{23,14},
$$

\vspace{3mm}
$$
\sum_{\stackrel{(\{g_1, g_2\}, \{q_1, q_{2}\})}
{{}_{\{g_1, g_2, q_1, q_{2}\}=\{1, 2, 3, 4\}}}}
a_{g_1 g_2, q_1 q_{2}}=
$$

$$
=a_{12,34}+a_{13,24}+a_{14,23}
+a_{23,14}+a_{24,13}+a_{34,12},
$$

\vspace{3mm}
$$
\sum_{\stackrel{(\{g_1, g_2\}, \{q_1, q_{2}, q_3\})}
{{}_{\{g_1, g_2, q_1, q_{2}, q_3\}=\{1, 2, 3, 4, 5\}}}}
a_{g_1 g_2, q_1 q_{2}q_3}
=
$$

$$
=a_{12,345}+a_{13,245}+a_{14,235}
+a_{15,234}+a_{23,145}+a_{24,135}+
$$
$$
+a_{25,134}+a_{34,125}+a_{35,124}+a_{45,123},
$$

\vspace{2mm}
$$
\sum_{\stackrel{(\{\{g_1, g_2\}, \{g_3, g_{4}\}\}, \{q_1\})}
{{}_{\{g_1, g_2, g_3, g_{4}, q_1\}=\{1, 2, 3, 4, 5\}}}}
a_{g_1 g_2, g_3 g_{4},q_1}
=
$$

$$
=
a_{12,34,5}+a_{13,24,5}+a_{14,23,5}+
a_{12,35,4}+a_{13,25,4}+a_{15,23,4}+
$$
$$
+a_{12,54,3}+a_{15,24,3}+a_{14,25,3}+a_{15,34,2}+a_{13,54,2}+a_{14,53,2}+
$$
$$
+
a_{52,34,1}+a_{53,24,1}+a_{54,23,1}.
$$

\vspace{4mm}

Now we can write (\ref{tyyy}) as

\vspace{1mm}

$$
J[\psi^{(k)}]_{T,t}=
\hbox{\vtop{\offinterlineskip\halign{
\hfil#\hfil\cr
{\rm l.i.m.}\cr
$\stackrel{}{{}_{p_1,\ldots,p_k\to \infty}}$\cr
}} }
\sum\limits_{j_1=0}^{p_1}\ldots
\sum\limits_{j_k=0}^{p_k}
C_{j_k\ldots j_1}\Biggl(
\prod_{l=1}^k\zeta_{j_l}^{(i_l)}+\sum\limits_{r=1}^{[k/2]}
(-1)^r \times
\Biggr.
$$

\vspace{3mm}
\begin{equation}
\label{leto6000}
\times
\sum_{\stackrel{(\{\{g_1, g_2\}, \ldots, 
\{g_{2r-1}, g_{2r}\}\}, \{q_1, \ldots, q_{k-2r}\})}
{{}_{\{g_1, g_2, \ldots, 
g_{2r-1}, g_{2r}, q_1, \ldots, q_{k-2r}\}=\{1, 2, \ldots, k\}}}}
\prod\limits_{s=1}^r
{\bf 1}_{\{i_{g_{{}_{2s-1}}}=~i_{g_{{}_{2s}}}\ne 0\}}
\Biggl.{\bf 1}_{\{j_{g_{{}_{2s-1}}}=~j_{g_{{}_{2s}}}\}}
\prod_{l=1}^{k-2r}\zeta_{j_{q_l}}^{(i_{q_l})}\Biggr),
\end{equation}

\vspace{5mm}
\noindent
where $\prod\limits_{\emptyset}
\stackrel{\sf def}{=}1,$ $\sum\limits_{\emptyset}
\stackrel{\sf def}{=}0,$ $[x]$ is an integer part of a real number $x;$
another notations are the same as in Theorem~7.

In particular, from (\ref{leto6000}) for $k=5$ we obtain

\vspace{1mm}

$$
J[\psi^{(5)}]_{T,t}=
\hbox{\vtop{\offinterlineskip\halign{
\hfil#\hfil\cr
{\rm l.i.m.}\cr
$\stackrel{}{{}_{p_1,\ldots,p_5\to \infty}}$\cr
}} }\sum_{j_1=0}^{p_1}\ldots\sum_{j_5=0}^{p_5}
C_{j_5\ldots j_1}\Biggl(
\prod_{l=1}^5\zeta_{j_l}^{(i_l)}-\Biggr.
$$

$$
-
\sum\limits_{\stackrel{(\{g_1, g_2\}, \{q_1, q_{2}, q_3\})}
{{}_{\{g_1, g_2, q_{1}, q_{2}, q_3\}=\{1, 2, 3, 4, 5\}}}}
{\bf 1}_{\{i_{g_{{}_{1}}}=~i_{g_{{}_{2}}}\ne 0\}}
{\bf 1}_{\{j_{g_{{}_{1}}}=~j_{g_{{}_{2}}}\}}
\prod_{l=1}^{3}\zeta_{j_{q_l}}^{(i_{q_l})}+
$$

$$
+
\sum_{\stackrel{(\{\{g_1, g_2\}, 
\{g_{3}, g_{4}\}\}, \{q_1\})}
{{}_{\{g_1, g_2, g_{3}, g_{4}, q_1\}=\{1, 2, 3, 4, 5\}}}}
{\bf 1}_{\{i_{g_{{}_{1}}}=~i_{g_{{}_{2}}}\ne 0\}}
{\bf 1}_{\{j_{g_{{}_{1}}}=~j_{g_{{}_{2}}}\}}
\Biggl.{\bf 1}_{\{i_{g_{{}_{3}}}=~i_{g_{{}_{4}}}\ne 0\}}
{\bf 1}_{\{j_{g_{{}_{3}}}=~j_{g_{{}_{4}}}\}}
\zeta_{j_{q_1}}^{(i_{q_1})}\Biggr).
$$

\vspace{5mm}
\noindent
The last equality obviously agrees with
(\ref{za5}).

Let us consider the generalization of Theorem 7 for the case
of an arbitrary complete orthonormal systems  
of functions in the space $L_2([t,T])$ 
and $\psi_1(\tau),\ldots,\psi_k(\tau)\in L_2([t, T]).$

\vspace{2mm}  
         
{\bf Theorem 8}\ \cite{arxiv-1} (Sect.~15), \cite{20aa} (Sect.~1.11), \cite{new-2023a},
\cite{diffjournal}.\
{\it Suppose that
$\psi_1(\tau),\ldots,\psi_k(\tau)\in L_2([t, T])$ and
$\{\phi_j(x)\}_{j=0}^{\infty}$ is an arbitrary complete orthonormal system  
of functions in the space $L_2([t,T]).$
Then the following expansion

\vspace{1mm}

$$
J[\psi^{(k)}]_{T,t}=
\hbox{\vtop{\offinterlineskip\halign{
\hfil#\hfil\cr
{\rm l.i.m.}\cr
$\stackrel{}{{}_{p_1,\ldots,p_k\to \infty}}$\cr
}} }
\sum\limits_{j_1=0}^{p_1}\ldots
\sum\limits_{j_k=0}^{p_k}
C_{j_k\ldots j_1}\Biggl(
\prod_{l=1}^k\zeta_{j_l}^{(i_l)}+\sum\limits_{r=1}^{[k/2]}
(-1)^r \times
\Biggr.
$$

\vspace{3mm}
\begin{equation}
\label{chainxx90}
\times
\sum_{\stackrel{(\{\{g_1, g_2\}, \ldots, 
\{g_{2r-1}, g_{2r}\}\}, \{q_1, \ldots, q_{k-2r}\})}
{{}_{\{g_1, g_2, \ldots, 
g_{2r-1}, g_{2r}, q_1, \ldots, q_{k-2r}\}=\{1, 2, \ldots, k\}}}}
\prod\limits_{s=1}^r
{\bf 1}_{\{i_{g_{{}_{2s-1}}}=~i_{g_{{}_{2s}}}\ne 0\}}
\Biggl.{\bf 1}_{\{j_{g_{{}_{2s-1}}}=~j_{g_{{}_{2s}}}\}}
\prod_{l=1}^{k-2r}\zeta_{j_{q_l}}^{(i_{q_l})}\Biggr)
\end{equation}

\vspace{5mm}
\noindent
con\-verg\-ing in the mean-square sense is valid$,$
where $\prod\limits_{\emptyset}
\stackrel{\sf def}{=}1,$ $\sum\limits_{\emptyset}
\stackrel{\sf def}{=}0,$ $[x]$ is an integer part of a real number $x;$
another
notations are the same as in Theorem {\rm 7}.}

\vspace{2mm}

It should be noted that an analogue 
of Theorem 8 (based on 
explicit representation using Hermite polynomials) was considered 
in \cite{Rybakov1000}. 
Note that we use another notations in comparison with \cite{Rybakov1000}.
Moreover, the proof of an analogue of Theorem 8
from \cite{Rybakov1000} is different from the proof given in 
\cite{arxiv-1} (Sect.~15), \cite{20aa} (Sect.~1.11), \cite{new-2023a},
\cite{diffjournal}.

In a number of works of the author \cite{arxiv-5}-\cite{arxiv-11}, 
\cite{arxiv-19}, \cite{arxiv-23},
\cite{5-005}, \cite{5-007}, \cite{12}-\cite{12aa-afterxxx}, \cite{new-art-1xxy}
Theorems 7, 8 have been adapted for the iterated Stratonovich stochastic integrals
(\ref{str}). Let us collect some of these results (old results)
in the following statement.

\vspace{2mm}

{\bf Theorem 9}\ \cite{arxiv-5}-\cite{arxiv-11}, 
\cite{arxiv-19}, \cite{arxiv-23},
\cite{5-005}, \cite{5-007}, \cite{12}-\cite{12aa-afterxxx}.\
{\it Suppose that 
$\{\phi_j(x)\}_{j=0}^{\infty}$ is a complete or\-tho\-nor\-mal system of 
Legendre polynomials or trigonometric functions in the space $L_2([t, T]).$
At the same time $\psi_2(\tau)$ is a continuously differentiable 
function on $[t, T]$ and $\psi_1(\tau),\psi_3(\tau)$ are two times 
continuously differentiable functions on $[t, T]$. Then

\vspace{-2mm}
\begin{equation}
\label{a}
J^{*}[\psi^{(2)}]_{T,t}=
\hbox{\vtop{\offinterlineskip\halign{
\hfil#\hfil\cr
{\rm l.i.m.}\cr
$\stackrel{}{{}_{p_1,p_2\to \infty}}$\cr
}} }\sum_{j_1=0}^{p_1}\sum_{j_2=0}^{p_2}
C_{j_2j_1}\zeta_{j_1}^{(i_1)}\zeta_{j_2}^{(i_2)}\ \ \ (i_1,i_2=1,\ldots,m),
\end{equation}

\vspace{0.5mm}
\begin{equation}
\label{feto19000ab}
J^{*}[\psi^{(3)}]_{T,t}=
\hbox{\vtop{\offinterlineskip\halign{
\hfil#\hfil\cr
{\rm l.i.m.}\cr
$\stackrel{}{{}_{p_1,p_2,p_3\to \infty}}$\cr
}} }\sum_{j_1=0}^{p_1}\sum_{j_2=0}^{p_2}\sum_{j_3=0}^{p_3}
C_{j_3 j_2 j_1}\zeta_{j_1}^{(i_1)}\zeta_{j_2}^{(i_2)}\zeta_{j_3}^{(i_3)}\ \ \
(i_1,i_2,i_3=0, 1,\ldots,m),
\end{equation}

\vspace{0.5mm}
\begin{equation}
\label{feto19000a}
J^{*}[\psi^{(3)}]_{T,t}=
\hbox{\vtop{\offinterlineskip\halign{
\hfil#\hfil\cr
{\rm l.i.m.}\cr
$\stackrel{}{{}_{p\to \infty}}$\cr
}} }
\sum\limits_{j_1,j_2,j_3=0}^{p}
C_{j_3 j_2 j_1}\zeta_{j_1}^{(i_1)}\zeta_{j_2}^{(i_2)}\zeta_{j_3}^{(i_3)}\ \ \
(i_1,i_2,i_3=1,\ldots,m),
\end{equation}

\vspace{0.5mm}
\begin{equation}
\label{uu}
J^{*}[\psi^{(4)}]_{T,t}=
\hbox{\vtop{\offinterlineskip\halign{
\hfil#\hfil\cr
{\rm l.i.m.}\cr
$\stackrel{}{{}_{p\to \infty}}$\cr
}} }
\sum\limits_{j_1, \ldots, j_4=0}^{p}
C_{j_4 j_3 j_2 j_1}\zeta_{j_1}^{(i_1)}
\zeta_{j_2}^{(i_2)}\zeta_{j_3}^{(i_3)}\zeta_{j_4}^{(i_4)}\ \ \
(i_1,\ldots,i_4=0, 1,\ldots,m),
\end{equation}

\vspace{3mm}
\noindent
where $J^{*}[\psi^{(k)}]_{T,t}$ is defined by {\rm (\ref{str})} and
$\psi_l(\tau)\equiv 1$ $(l=1,\ldots,4)$ in {\rm (\ref{feto19000ab})}, 
{\rm (\ref{uu}),} ${\bf w}_{\tau}^{(0)}=\tau;$ 
another notations are the same as in Theorem {\rm 7.}}

\vspace{2mm}

Recently, a new approach to the expansion and mean-square 
approximation of iterated Stratonovich stochastic integrals has been obtained
\cite{20aa} (Sect.~2.10--2.24) (also see \cite{arxiv-5}, 
\cite{arxiv-10}, \cite{arxiv-11}, \cite{12aa-afterxxx},
\cite{new-art-1xxy}, \cite{new-art-1xxys}).
The mentioned approach covers the following cases
for iterated Stratonovich stochastic
integrals of multiplicities 1 to 6.

1. The case of continuously differentiable 
weight functions (multiplicities 1 to 5) and 
weight functions identically equal to one (multiplicity 6).
In this case, we use
a complete orthonormal system of Legendre polynomials or 
trigonometric functions in $L_2([t, T])$.

2. The case of continuous weight functions (multiplicities 1 and 2),
binomial weight functions (multiplicities 3 and 4)
and weight functions identically equal to one (multiplicities 5 and 6).
In this case, we use
an arbitrary complete orthonormal system of functions in $L_2([t, T])$.

Recently (in 2024), the above approach has been generalized for 
iterated Stratonovich stochastic
integrals of multiplicity $k,$ $k\in \mathbb{N}$ \cite{20aa} (Chapter 2, Theorems~2.59, 2.61) 
but under one additional condition
(the case of continuous weight functions 
and an arbitrary complete orthonormal system of functions in $L_2([t, T])$).

Let us consider a number of theorems
that relate to Cases 1 and 2.

\vspace{2mm}

{\bf Theorem 10}\ \cite{arxiv-5}, 
\cite{arxiv-10}, \cite{arxiv-11}, \cite{20aa}, \cite{new-art-1xxy}.\
{\it Suppose 
that $\{\phi_j(x)\}_{j=0}^{\infty}$ is a complete orthonormal system of 
Legendre polynomials or trigonometric functions in the space $L_2([t, T]).$
Furthermore, let $\psi_1(\tau), \psi_2(\tau),$ $\psi_3(\tau)$ are continuously dif\-ferentiable 
nonrandom functions on $[t, T].$ 
Then, for the 
iterated Stra\-to\-no\-vich stochastic integral of third multiplicity

\vspace{-1mm}
$$
J^{*}[\psi^{(3)}]_{T,t}={\int\limits_t^{*}}^T\psi_3(t_3)
{\int\limits_t^{*}}^{t_3}\psi_2(t_2)
{\int\limits_t^{*}}^{t_2}\psi_1(t_1)
d{\bf w}_{t_1}^{(i_1)}
d{\bf w}_{t_2}^{(i_2)}d{\bf w}_{t_3}^{(i_3)}\ \ \ (i_1,i_2,i_3=0,1,\ldots,m)
$$

\vspace{3mm}
\noindent
defined by {\rm (\ref{str})} the following 
relations

\begin{equation}
\label{fin1}
J^{*}[\psi^{(3)}]_{T,t}
=\hbox{\vtop{\offinterlineskip\halign{
\hfil#\hfil\cr
{\rm l.i.m.}\cr
$\stackrel{}{{}_{p\to \infty}}$\cr
}} }
\sum\limits_{j_1, j_2, j_3=0}^{p}
C_{j_3 j_2 j_1}\zeta_{j_1}^{(i_1)}\zeta_{j_2}^{(i_2)}\zeta_{j_3}^{(i_3)},
\end{equation}

\vspace{2mm}
\begin{equation}
\label{fin2}
{\sf M}\left\{\left(
J^{*}[\psi^{(3)}]_{T,t}-
\sum\limits_{j_1, j_2, j_3=0}^{p}
C_{j_3 j_2 j_1}\zeta_{j_1}^{(i_1)}\zeta_{j_2}^{(i_2)}\zeta_{j_3}^{(i_3)}\right)^2\right\}
\le \frac{C}{p}
\end{equation}

\vspace{5mm}
\noindent
are fulfilled, where $i_1, i_2, i_3=0,1,\ldots,m$ in {\rm (\ref{fin1})} and 
$i_1, i_2, i_3=1,\ldots,m$ in {\rm (\ref{fin2})},
constant $C$ is independent of $p,$

$$
C_{j_3 j_2 j_1}=\int\limits_t^T\psi_3(t_3)\phi_{j_3}(t_3)
\int\limits_t^{t_3}\psi_2(t_2)\phi_{j_2}(t_2)
\int\limits_t^{t_2}\psi_1(t_1)\phi_{j_1}(t_1)dt_1dt_2dt_3
$$

\vspace{3mm}
\noindent
and
$$
\zeta_{j}^{(i)}=
\int\limits_t^T \phi_{j}(\tau) d{\bf f}_{\tau}^{(i)}
$$ 

\vspace{2mm}
\noindent
are independent standard Gaussian random variables for various 
$i$ or $j$ {\rm (}in the case when $i\ne 0${\rm ),} 
${\bf w}_{\tau}^{(0)}=\tau;$ 
another notations are the same as in Theorem~{\rm 7}.}

\vspace{2mm}

{\bf Theorem 11}\ \cite{arxiv-5}, 
\cite{arxiv-10}, \cite{arxiv-11}, \cite{20aa}, \cite{new-art-1xxy}.\ {\it Let
$\{\phi_j(x)\}_{j=0}^{\infty}$ be a complete orthonormal system of 
Legendre polynomials or trigonometric functions in the space $L_2([t, T]).$
Furthermore, let $\psi_1(\tau), \ldots,$ $\psi_4(\tau)$ be continuously dif\-ferentiable 
nonrandom functions on $[t, T].$ 
Then, for the 
iterated Stra\-to\-no\-vich stochastic integral of fourth multiplicity

\vspace{-1mm}
\begin{equation}
\label{fin0}
J^{*}[\psi^{(4)}]_{T,t}={\int\limits_t^{*}}^T\psi_4(t_4)
{\int\limits_t^{*}}^{t_4}\psi_3(t_3)
{\int\limits_t^{*}}^{t_3}\psi_2(t_2)
{\int\limits_t^{*}}^{t_2}\psi_1(t_1)
d{\bf w}_{t_1}^{(i_1)}
d{\bf w}_{t_2}^{(i_2)}d{\bf w}_{t_3}^{(i_3)}d{\bf w}_{t_4}^{(i_4)}
\end{equation}

\vspace{3mm}
\noindent
defined by {\rm (\ref{str})} the following 
relations

\begin{equation}
\label{fin3}
J^{*}[\psi^{(4)}]_{T,t}
=\hbox{\vtop{\offinterlineskip\halign{
\hfil#\hfil\cr
{\rm l.i.m.}\cr
$\stackrel{}{{}_{p\to \infty}}$\cr
}} }
\sum\limits_{j_1, j_2, j_3,j_4=0}^{p}
C_{j_4j_3 j_2 j_1}\zeta_{j_1}^{(i_1)}\zeta_{j_2}^{(i_2)}\zeta_{j_3}^{(i_3)}\zeta_{j_4}^{(i_4)},
\end{equation}

\vspace{2mm}

\begin{equation}
\label{fin4}
{\sf M}\left\{\left(
J^{*}[\psi^{(4)}]_{T,t}-
\sum\limits_{j_1, j_2, j_3, j_4=0}^{p}
C_{j_4 j_3 j_2 j_1}\zeta_{j_1}^{(i_1)}\zeta_{j_2}^{(i_2)}\zeta_{j_3}^{(i_3)}
\zeta_{j_4}^{(i_4)}
\right)^2\right\}
\le \frac{C}{p^{1-\varepsilon}}
\end{equation}

\vspace{4mm}
\noindent
are fulfilled, where $i_1, \ldots , i_4=0,1,\ldots,m$ in {\rm (\ref{fin0}),} {\rm (\ref{fin3})} 
and $i_1, \ldots, i_4=1,\ldots,m$ in {\rm (\ref{fin4})},
constant $C$ does not depend on $p,$
$\varepsilon$ is an arbitrary
small positive real number 
for the case of complete orthonormal system of 
Legendre polynomials in the space $L_2([t, T])$
and $\varepsilon=0$ for the case of
complete orthonormal system of 
trigonometric functions in the space $L_2([t, T]),$

$$
C_{j_4 j_3 j_2 j_1}=
\int\limits_t^T\psi_4(t_4)\phi_{j_4}(t_4)
\int\limits_t^{t_4}\psi_3(t_3)\phi_{j_3}(t_3)
\int\limits_t^{t_3}\psi_2(t_2)\phi_{j_2}(t_2)
\int\limits_t^{t_2}\psi_1(t_1)\phi_{j_1}(t_1)dt_1dt_2dt_3dt_4,
$$

\vspace{4mm}
\noindent
${\bf w}_{\tau}^{(0)}=\tau;$ another notations are the same as in Theorem~{\rm 10}.}

\vspace{2mm}

{\bf Theorem 12}\ \cite{arxiv-5}, 
\cite{arxiv-10}, \cite{arxiv-11}, \cite{20aa}, \cite{new-art-1xxy}.\
{\it Assume 
that $\{\phi_j(x)\}_{j=0}^{\infty}$ is a complete orthonormal system of 
Legendre polynomials or trigonometric functions in the space $L_2([t, T])$
and $\psi_1(\tau), \ldots,$ $\psi_5(\tau)$ are continuously dif\-ferentiable 
nonrandom functions on $[t, T].$ 
Then, for the 
iterated Stra\-to\-no\-vich stochastic integral of fifth multiplicity

\vspace{-1mm}
\begin{equation}
\label{fin7}
J^{*}[\psi^{(5)}]_{T,t}={\int\limits_t^{*}}^T\psi_5(t_5)
\ldots
{\int\limits_t^{*}}^{t_2}\psi_1(t_1)
d{\bf w}_{t_1}^{(i_1)}
\ldots d{\bf w}_{t_5}^{(i_5)}
\end{equation}

\vspace{3mm}
\noindent
the following 
relations

\begin{equation}
\label{fin8}
J^{*}[\psi^{(5)}]_{T,t}
=\hbox{\vtop{\offinterlineskip\halign{
\hfil#\hfil\cr
{\rm l.i.m.}\cr
$\stackrel{}{{}_{p\to \infty}}$\cr
}} }
\sum\limits_{j_1,\ldots,j_5=0}^{p}
C_{j_5 \ldots j_1}\zeta_{j_1}^{(i_1)}\ldots \zeta_{j_5}^{(i_5)},
\end{equation}

\vspace{2mm}

\begin{equation}
\label{fin9}
{\sf M}\left\{\left(
J^{*}[\psi^{(5)}]_{T,t}-
\sum\limits_{j_1, \ldots, j_5=0}^{p}
C_{j_5 \ldots j_1}\zeta_{j_1}^{(i_1)}\ldots
\zeta_{j_5}^{(i_5)}
\right)^2\right\}
\le \frac{C}{p^{1-\varepsilon}}
\end{equation}

\vspace{5mm}
\noindent
are fulfilled, where $i_1, \ldots , i_5=0,1,\ldots,m$ in {\rm (\ref{fin7}),} {\rm (\ref{fin8})} 
and $i_1, \ldots, i_5=1,\ldots,m$ in {\rm (\ref{fin9})},
constant $C$ is independent of $p,$
$\varepsilon$ is an arbitrary
small positive real number 
for the case of complete orthonormal system of 
Legendre polynomials in the space $L_2([t, T])$
and $\varepsilon=0$ for the case of
complete orthonormal system of 
trigonometric functions in the space $L_2([t, T]),$

$$
C_{j_5 \ldots j_1}=
\int\limits_t^T\psi_5(t_5)\phi_{j_5}(t_5)\ldots
\int\limits_t^{t_2}\psi_1(t_1)\phi_{j_1}(t_1)dt_1\ldots dt_5,
$$

\vspace{3mm}
\noindent
${\bf w}_{\tau}^{(0)}=\tau;$ 
another notations are the same as in Theorem~{\rm 10, 11}.}

\vspace{2mm}

{\bf Theorem 13}\ \cite{arxiv-5}, 
\cite{arxiv-10}, \cite{arxiv-11}, \cite{20aa}, \cite{new-art-1xxys}.\
{\it Suppose that 
$\{\phi_j(x)\}_{j=0}^{\infty}$ is a complete orthonormal system of 
Legendre polynomials or trigonometric functions in the space $L_2([t, T]).$
Then, for the 
iterated Stratonovich stochastic integral of sixth multiplicity

\begin{equation}
\label{after10001qu1}
J^{*}[\psi^{(6)}]_{T,t}={\int\limits_t^{*}}^T
\ldots
{\int\limits_t^{*}}^{t_2}
d{\bf w}_{t_1}^{(i_1)}
\ldots d{\bf w}_{t_6}^{(i_6)}
\end{equation}

\vspace{3mm}
\noindent
the following 
expansion 

\vspace{-1mm}
$$
J^{*}[\psi^{(6)}]_{T,t}
=\hbox{\vtop{\offinterlineskip\halign{
\hfil#\hfil\cr
{\rm l.i.m.}\cr
$\stackrel{}{{}_{p\to \infty}}$\cr
}} }
\sum\limits_{j_1, \ldots, j_6=0}^{p}
C_{j_6 \ldots j_1}\zeta_{j_1}^{(i_1)}\ldots
\zeta_{j_6}^{(i_6)}
$$

\vspace{4mm}
\noindent
that converges in the mean-square sense is valid, where
$i_1, \ldots, i_6=0, 1,\ldots,m,$

$$
C_{j_6 \ldots j_1}=
\int\limits_t^T\phi_{j_6}(t_6)\ldots
\int\limits_t^{t_2}\phi_{j_1}(t_1)dt_1\ldots dt_6.
$$

\vspace{2mm}
\noindent
${\bf w}_{\tau}^{(0)}=\tau;$ 
another notations are the same as in Theorems~{\rm 10--12}.}

\vspace{2mm}

{\bf Theorem 14}\ \cite{arxiv-5}, 
\cite{arxiv-10}, \cite{arxiv-11}, \cite{20aa}.\  {\it Suppose that 
$\{\phi_j(x)\}_{j=0}^{\infty}$ is an arbitrary complete orthonormal system of 
functions in the space $L_2([t, T]).$
Moreover$,$ $\psi_1(\tau), \psi_2(\tau)$ are continuous 
functions on $[t, T].$ 
Then$,$ 
for the iterated Stra\-to\-novich stochastic integral

$$
J^{*}[\psi^{(2)}]_{T,t}={\int\limits_t^{*}}^T\psi_2(t_2)
{\int\limits_t^{*}}^{t_2}\psi_1(t_1)d{\bf f}_{t_1}^{(i_1)}
d{\bf f}_{t_2}^{(i_2)}\ \ \ (i_1, i_2=1,\ldots,m)
$$

\vspace{3mm}
\noindent
defined by {\rm (\ref{str})} the following expansion  

$$
J^{*}[\psi^{(2)}]_{T,t}=\hbox{\vtop{\offinterlineskip\halign{
\hfil#\hfil\cr
{\rm l.i.m.}\cr
$\stackrel{}{{}_{p_1,p_2\to \infty}}$\cr
}} }\sum_{j_1=0}^{p_1}\sum_{j_2=0}^{p_2}
C_{j_2j_1}\zeta_{j_1}^{(i_1)}\zeta_{j_2}^{(i_2)}
$$

\vspace{3mm}
\noindent
that converges in the mean-square
sence is valid$,$ where the notations are the same as in Theorems {\rm 7, 8.}
}

\vspace{2mm}

{\bf Theorem~15}\ \cite{arxiv-5}, 
\cite{arxiv-10}, \cite{arxiv-11}, \cite{20aa}, \cite{2024xxx}.\ {\it Suppose that
$\{\phi_j(x)\}_{j=0}^{\infty}$ is an arbitrary complete orthonormal system of 
functions in the space $L_2([t,T]).$
Then$,$ for the iterated Stra\-to\-no\-vich stochastic integrals

$$
I_{{l_1l_2l_3}_{T,t}}^{*(i_1i_2i_3)}=
{\int\limits_t^{*}}^T (t_3-t)^{l_3}
{\int\limits_t^{*}}^{t_3}(t_2-t)^{l_2}
{\int\limits_t^{*}}^{t_2}(t_1-t)^{l_1}
d{\bf w}_{t_1}^{(i_1)}
d{\bf w}_{t_2}^{(i_2)}d{\bf w}_{t_3}^{(i_3)},
$$

$$
I_{{l_1l_2l_3 l_4}_{T,t}}^{*(i_1i_2i_3 i_4)}=
{\int\limits_t^{*}}^T (t_4-t)^{l_4}{\int\limits_t^{*}}^{t_4} (t_3-t)^{l_3}
{\int\limits_t^{*}}^{t_3}(t_2-t)^{l_2}
{\int\limits_t^{*}}^{t_2}(t_1-t)^{l_1}
d{\bf w}_{t_1}^{(i_1)}
d{\bf w}_{t_2}^{(i_2)}d{\bf w}_{t_3}^{(i_3)}d{\bf w}_{t_4}^{(i_4)}
$$

\vspace{3mm}
\noindent
the following expansions 

$$
I_{{l_1l_2l_3}_{T,t}}^{*(i_1i_2i_3)}=
\hbox{\vtop{\offinterlineskip\halign{
\hfil#\hfil\cr
{\rm l.i.m.}\cr
$\stackrel{}{{}_{p\to \infty}}$\cr
}} }\sum_{j_1,j_2,j_3=0}^{p}
C_{j_3 j_2 j_1}\zeta_{j_1}^{(i_1)}\zeta_{j_2}^{(i_2)}\zeta_{j_3}^{(i_3)},
$$

\vspace{1mm}
$$
I_{{l_1l_2l_3 l_4}_{T,t}}^{*(i_1i_2i_3i_4)}=
\hbox{\vtop{\offinterlineskip\halign{
\hfil#\hfil\cr
{\rm l.i.m.}\cr
$\stackrel{}{{}_{p\to \infty}}$\cr
}} }\sum_{j_1,j_2,j_3,j_4=0}^{p}
C_{j_4 j_3 j_2 j_1}\zeta_{j_1}^{(i_1)}\zeta_{j_2}^{(i_2)}\zeta_{j_3}^{(i_3)}\zeta_{j_4}^{(i_4)}
$$

\vspace{5mm}
\noindent
that converge in the mean-square sense are valid, where 
$i_1,i_2,i_3,i_4=0,1,\ldots,m;$ $l_1,l_2,l_3,l_4=0,1,2,\ldots,$

$$
C_{j_3 j_2 j_1}=\int\limits_t^T
(t_3-t)^{l_3}\phi_{j_3}(t_3)\int\limits_t^{t_3}
(t_2-t)^{l_2}
\phi_{j_2}(t_2)
\int\limits_t^{t_2}
(t_1-t)^{l_1}\phi_{j_1}(t_1)
dt_1dt_2dt_3,
$$

$$
C_{j_4 j_3 j_2 j_1}=\int\limits_t^T
(t_4-t)^{l_4}\phi_{j_4}(t_4)\int\limits_t^{t_4}
(t_3-t)^{l_3}\phi_{j_3}(t_3)\int\limits_t^{t_3}
(t_2-t)^{l_2}
\phi_{j_2}(t_2)
\int\limits_t^{t_2}
(t_1-t)^{l_1}\phi_{j_1}(t_1)
dt_1dt_2dt_3dt_4,
$$

\vspace{3mm}
\noindent
${\bf w}_{\tau}^{(0)}=\tau;$ another notations are the same as in Theorems {\rm 7, 8.}
}

\vspace{2mm}

{\bf Theorem~16}\ \cite{arxiv-5}, 
\cite{arxiv-10}, \cite{arxiv-11}, \cite{20aa}, \cite{2024xxx}.\ {\it Suppose that
$\{\phi_j(x)\}_{j=0}^{\infty}$ is an arbitrary complete orthonormal system of 
func\-ti\-ons in the space $L_2([t,T]).$
Then$,$ for the iterated Stra\-to\-no\-vich stochastic integrals

$$
J^{*}[\psi^{(5)}]_{T,t}=
{\int\limits_t^{*}}^T
\ldots
{\int\limits_t^{*}}^{t_2}
d{\bf w}_{t_1}^{(i_1)}
\ldots d{\bf w}_{t_5}^{(i_5)},\ \ \
J^{*}[\psi^{(6)}]_{T,t}=
{\int\limits_t^{*}}^T
\ldots
{\int\limits_t^{*}}^{t_2}
d{\bf w}_{t_1}^{(i_1)}
\ldots d{\bf w}_{t_6}^{(i_6)}
$$

\vspace{3mm}
\noindent
the following 
expansions 

$$
J^{*}[\psi^{(5)}]_{T,t}=
\hbox{\vtop{\offinterlineskip\halign{
\hfil#\hfil\cr
{\rm l.i.m.}\cr
$\stackrel{}{{}_{p\to \infty}}$\cr
}} }
\sum\limits_{j_1,\ldots, j_5=0}^{p}
C_{j_5 \ldots j_1}\zeta_{j_1}^{(i_1)}\ldots \zeta_{j_5}^{(i_5)},\ \ \ 
J^{*}[\psi^{(6)}]_{T,t}=
\hbox{\vtop{\offinterlineskip\halign{
\hfil#\hfil\cr
{\rm l.i.m.}\cr
$\stackrel{}{{}_{p\to \infty}}$\cr
}} }
\sum\limits_{j_1,\ldots, j_6=0}^{p}
C_{j_6 \ldots j_1}\zeta_{j_1}^{(i_1)}\ldots \zeta_{j_6}^{(i_6)}
$$

\vspace{5mm}
\noindent
that converge in the mean-square sense are valid, where 
$i_1,\ldots,i_6=0, 1,\ldots,m,$

$$
C_{j_5\ldots j_1}=\int\limits_t^T
\phi_{j_5}(t_5)
\ldots
\int\limits_t^{t_2}
\phi_{j_1}(t_1)dt_1\ldots dt_5,\ \ \
C_{j_6\ldots j_1}=\int\limits_t^T
\phi_{j_6}(t_6)
\ldots
\int\limits_t^{t_2}
\phi_{j_1}(t_1)dt_1\ldots dt_6,
$$

\vspace{3mm}
\noindent
${\bf w}_{\tau}^{(0)}=\tau;$
another notations are the same as in Theorems {\rm 7, 8.}
}

\vspace{2mm}

On the base of Theorems 14--16 in
\cite{20aa}
the following hypothesis was formulated.

\vspace{2mm}

{\bf Hypothesis 1}\  \cite{20aa}.
{\it Assume that
$\{\phi_j(x)\}_{j=0}^{\infty}$ is an arbitrary complete orthonormal system of 
functions in $L_2([t,T])$
and $\psi_1(\tau),\ldots,\psi_k(\tau)$ are continuous functions on $[t, T].$
Then, for the iterated Stra\-to\-no\-vich stochastic 
integral of $k$th multiplicity

$$
J^{*}[\psi^{(k)}]_{T,t}=
{\int\limits_t^{*}}^T \psi_k(t_k) \ldots 
{\int\limits_t^{*}}^{t_2}
\psi_1(t_1) d{\bf w}_{t_1}^{(i_1)}\ldots
d{\bf w}_{t_k}^{(i_k)}
\ \ \ (i_1,\ldots, i_k=0, 1,\ldots,m)
$$

\vspace{3mm}
\noindent
defined by {\rm (\ref{str})} the following converging in the mean-square sense 
expansion 
         
\vspace{-2mm}
$$
J^{*}[\psi^{(k)}]_{T,t}=
\hbox{\vtop{\offinterlineskip\halign{
\hfil#\hfil\cr
{\rm l.i.m.}\cr
$\stackrel{}{{}_{p\to \infty}}$\cr
}} }
\sum\limits_{j_1, \ldots j_k=0}^{p}
C_{j_k \ldots j_1}\zeta_{j_1}^{(i_1)}
\ldots
\zeta_{j_k}^{(i_k)}
$$

\vspace{3mm}
\noindent
holds, where the Fourier 
coefficient $C_{j_k \ldots j_1}$ has the form

$$
C_{j_k \ldots j_1}=\int\limits_t^T\psi_k(t_k)\phi_{j_k}(t_k)\ldots
\int\limits_t^{t_2}
\psi_1(t_1)\phi_{j_1}(t_1)
dt_1\ldots dt_k,
$$

\vspace{3mm}
\noindent
${\rm l.i.m.}$ is a limit in the mean-square sense,

\vspace{-2mm}
$$
\zeta_{j}^{(i)}=
\int\limits_t^T \phi_{j}(\tau) d{\bf w}_{\tau}^{(i)}
$$ 

\vspace{2mm}
\noindent
are independent standard Gaussian random variables for various 
$i$ or $j$\ {\rm (}if $i\ne 0${\rm ),}
${\bf w}_{\tau}^{(i)}={\bf f}_{\tau}^{(i)}$ are independent 
standard Wiener processes
$(i=1,\ldots,m)$ and 
${\bf w}_{\tau}^{(0)}=\tau.$}

\vspace{2mm}

The hypothesis 1 allows to approximate the iterated
Stratonovich stochastic 
integral $J^{*}[\psi^{(k)}]_{T,t}$
by the sum

\vspace{-3mm}
$$
J^{*}[\psi^{(k)}]_{T,t}^{p}=
\sum\limits_{j_1, \ldots j_k=0}^{p}
C_{j_k \ldots j_1}\zeta_{j_1}^{(i_1)}
\ldots
\zeta_{j_k}^{(i_k)},
$$

\vspace{3mm}
\noindent
where
$$
\lim_{p\to\infty}{\sf M}\left\{\left(
J^{*}[\psi^{(k)}]_{T,t}-
J^{*}[\psi^{(k)}]_{T,t}^{p}\right)^2\right\}=0.
$$

\vspace{4mm}

Hypothesis~1 has been proved in \cite{20aa} (Chapter 2, Theorems~2.59, 2.61) 
but under one additional condition.

The following theorem shows how to calculate exactly
the mean-square approximation error for iterated Ito
stochastic integrals in Theorems 7, 8.

\vspace{2mm}

{\bf Theorem 17} \cite{20aa} (also see \cite{arxiv-3}).\
{\it Suppose that
$\psi_1(\tau),\ldots,\psi_k(\tau)\in L_2([t, T])$ 
and
$\{\phi_j(x)\}_{j=0}^{\infty}$ is an arbitrary complete orthonormal system  
of functions in the space $L_2([t,T]).$ 
Then

$$
{\sf M}\left\{\left(J[\psi^{(k)}]_{T,t}-
J[\psi^{(k)}]_{T,t}^p\right)^2\right\}
= \int\limits_{[t,T]^k} K^2(t_1,\ldots,t_k)
dt_1\ldots dt_k - 
$$

\begin{equation}
\label{tttr11}
- \sum_{j_1=0}^{p}\ldots\sum_{j_k=0}^{p}
C_{j_k\ldots j_1}
{\sf M}\left\{J[\psi^{(k)}]_{T,t}
\sum\limits_{(j_1,\ldots,j_k)}
\int\limits_t^T \phi_{j_k}(t_k)
\ldots
\int\limits_t^{t_{2}}\phi_{j_{1}}(t_{1})
d{\bf f}_{t_1}^{(i_1)}\ldots
d{\bf f}_{t_k}^{(i_k)}\right\},
\end{equation}

\vspace{3mm}
\noindent
where $J[\psi^{(k)}]_{T,t}$ is the iterated Ito stochastic integral {\rm (\ref{ito}),}

$$
J[\psi^{(k)}]_{T,t}^{p}=
\sum\limits_{j_1=0}^{p}\ldots
\sum\limits_{j_k=0}^{p}
C_{j_k\ldots j_1}\Biggl(
\prod_{l=1}^k\zeta_{j_l}^{(i_l)}+\sum\limits_{r=1}^{[k/2]}
(-1)^r \times
\Biggr.
$$

\vspace{2mm}
\begin{equation}
\label{yeee2}
\times
\sum_{\stackrel{(\{\{g_1, g_2\}, \ldots, 
\{g_{2r-1}, g_{2r}\}\}, \{q_1, \ldots, q_{k-2r}\})}
{{}_{\{g_1, g_2, \ldots, 
g_{2r-1}, g_{2r}, q_1, \ldots, q_{k-2r}\}=\{1, 2, \ldots, k\}}}}
\prod\limits_{s=1}^r
{\bf 1}_{\{i_{g_{{}_{2s-1}}}=~i_{g_{{}_{2s}}}\ne 0\}}
\Biggl.{\bf 1}_{\{j_{g_{{}_{2s-1}}}=~j_{g_{{}_{2s}}}\}}
\prod_{l=1}^{k-2r}\zeta_{j_{q_l}}^{(i_{q_l})}\Biggr),
\end{equation}

\vspace{4mm}
\noindent
$i_1,\ldots,i_k=1,\ldots,m,$ the expression

\vspace{-2mm}
$$
\sum\limits_{(j_1,\ldots,j_k)}
$$ 

\vspace{2mm}
\noindent
means the sum with respect to all
possible permutations 
$(j_1,\ldots,j_k)$. At the same time if 
$j_r$ swapped with $j_q$ in the permutation $(j_1,\ldots,j_k)$,
then $i_r$ swapped with $i_q$ in the permutation
$(i_1,\ldots,i_k);$
another notations are the same as in Theorems {\rm 7, 8.}}

\vspace{2mm}

Denote
$$
E_k^p \stackrel{\sf def}{=}{\sf M}\left\{\left(J[\psi^{(k)}]_{T,t}-
J[\psi^{(k)}]_{T,t}^p\right)^2\right\}.
$$

\vspace{3mm}

Note that
$$
{\sf M}\left\{J[\psi^{(k)}]_{T,t}
\int\limits_t^T \phi_{j_k}(t_k)
\ldots
\int\limits_t^{t_{2}}\phi_{j_{1}}(t_{1})
d{\bf f}_{t_1}^{(i_1)}\ldots
d{\bf f}_{t_k}^{(i_k)}\right\}=
$$

$$
=
\int\limits_t^T\psi_k(t_k) \phi_{j_k}(t_k)\ldots \int\limits_t^{t_{2}}
\psi_1(t_1)\phi_{j_1}(t_1) dt_1\ldots dt_k=
C_{j_k\ldots j_1}.
$$

\vspace{4mm}

Then from Theorem 17 for pairwise different $i_1,\ldots,i_k$ 
and for $i_1=\ldots=i_k$
we obtain

\begin{equation}
\label{zzz4}
E_k^p= I_k- \sum_{j_1,\ldots,j_k=0}^{p}
C_{j_k\ldots j_1}^2,
\end{equation}

$$
E_k^p= I_k - \sum_{j_1,\ldots,j_k=0}^{p}
C_{j_k\ldots j_1}\left(\sum\limits_{(j_1,\ldots,j_k)}
C_{j_k\ldots j_1}\right)
$$

\vspace{4mm}
\noindent
correspondingly.

Consider some examples of application of Theorem 17
$(i_1, \ldots ,i_5=1,\ldots,m)$

\vspace{-1mm}
\begin{equation}
\label{zzz4x}
E_2^p
=I_2
-\sum_{j_1,j_2=0}^p
C_{j_2j_1}^2-
\sum_{j_1,j_2=0}^p
C_{j_2j_1}C_{j_1j_2}\ \ \ (i_1=i_2),
\end{equation}

\vspace{1mm}
\begin{equation}
\label{zzz5}
E_3^p=I_3
-\sum_{j_3,j_2,j_1=0}^p C_{j_3j_2j_1}^2-
\sum_{j_3,j_2,j_1=0}^p C_{j_3j_1j_2}C_{j_3j_2j_1}\ \ \ (i_1=i_2\ne i_3),
\end{equation}

\vspace{1mm}
\begin{equation}
\label{zzz6}
E_3^p=I_3-
\sum_{j_3,j_2,j_1=0}^p C_{j_3j_2j_1}^2-
\sum_{j_3,j_2,j_1=0}^p C_{j_2j_3j_1}C_{j_3j_2j_1}\ \ \ (i_1\ne i_2=i_3),
\end{equation}

\vspace{1mm}
$$
E_3^p=I_3
-\sum_{j_3,j_2,j_1=0}^p C_{j_3j_2j_1}^2-
\sum_{j_3,j_2,j_1=0}^p C_{j_3j_2j_1}C_{j_1j_2j_3}\ \ \ (i_1=i_3\ne i_2),
$$

\vspace{2mm}
$$
E^p_4 = I_4 - \sum_{j_1,\ldots,j_4=0}^{p}
C_{j_4\ldots j_1}\Biggl(\sum\limits_{(j_1,j_2)}
C_{j_4\ldots j_1}\Biggr)\ \ \ (i_1=i_2\ne i_3, i_4;\ i_3\ne i_4),
$$

\vspace{2mm}
$$
E^p_4 = I_4 - \sum_{j_1,\ldots,j_4=0}^{p}
C_{j_4\ldots j_1}\Biggl(\sum\limits_{(j_1,j_3)}
C_{j_4\ldots j_1}\Biggr)\ \ \ (i_1=i_3\ne i_2, i_4;\ i_2\ne i_4),
$$

\vspace{2mm}
$$
E^p_4 = I_4 - \sum_{j_1,\ldots,j_4=0}^{p}
C_{j_4\ldots j_1}\Biggl(\sum\limits_{(j_1,j_4)}
C_{j_4\ldots j_1}\Biggr)\ \ \ (i_1=i_4\ne i_2, i_3;\ i_2\ne i_3),
$$

\vspace{2mm}
$$
E^p_4 = I_4 - \sum_{j_1,\ldots,j_4=0}^{p}
C_{j_4\ldots j_1}\Biggl(\sum\limits_{(j_2,j_3)}
C_{j_4\ldots j_1}\Biggr)\ \ \ (i_2=i_3\ne i_1, i_4;\ i_1\ne i_4),
$$

\vspace{2mm}
$$
E^p_4 = I_4 - \sum_{j_1,\ldots,j_4=0}^{p}
C_{j_4\ldots j_1}\Biggl(\sum\limits_{(j_2,j_4)}
C_{j_4\ldots j_1}\Biggr)\ \ \ (i_2=i_4\ne i_1, i_3;\ i_1\ne i_3),
$$

\vspace{2mm}
$$
E^p_4 = I_4 - \sum_{j_1,\ldots,j_4=0}^{p}
C_{j_4\ldots j_1}\Biggl(\sum\limits_{(j_3,j_4)}
C_{j_4\ldots j_1}\Biggr)\ \ \ (i_3=i_4\ne i_1, i_2;\ i_1\ne i_2),
$$

\vspace{2mm}
$$
E_4^p = I_4 -
 \sum_{j_1,\ldots,j_4=0}^{p}
C_{j_4\ldots j_1}\Biggl(\sum\limits_{(j_1,j_3,j_4)}
C_{j_4\ldots j_1}\Biggr)\ \ \ (i_1=i_3=i_4\ne i_2),
$$

\vspace{3mm}
$$
E_5^p = I_5 - \sum_{j_1,\ldots,j_5=0}^{p}
C_{j_5\ldots j_1}\Biggl(\sum\limits_{(j_2,j_4)}\Biggl(
\sum\limits_{(j_3,j_5)}
C_{j_5\ldots j_1}\Biggr)\Biggr)\ \ \ (i_1\ne i_2=i_4\ne i_3=i_5\ne i_1),
$$

\vspace{3mm}
$$
E^p_5 = I_5 - \sum_{j_1,\ldots,j_5=0}^{p}
C_{j_5\ldots j_1}\Biggl(\sum\limits_{(j_4,j_5)}\Biggl(
\sum\limits_{(j_1,j_2,j_3)}
C_{j_5\ldots j_1}\Biggr)\Biggr)\ \ \ (i_1=i_2=i_3\ne i_4=i_5),
$$

\vspace{3mm}
$$
E^p_5 = I_5 - \sum_{j_1,\ldots,j_5=0}^{p}
C_{j_5\ldots j_1}\Biggl(\sum\limits_{(j_1,j_3,j_4,j_5)}
C_{j_5\ldots j_1}\Biggr)\ \ \ (i_1=i_3=i_4=i_5\ne i_2).
$$

\vspace{5mm}

Let $J[\psi^{(k)}]_{T,t}^{p_1,\ldots,p_k}$ be the 
expression on the right-hand side of (\ref{chainxx90}) before passing
to the limit. Denote 
$$
I_k\stackrel{{\rm def}}{=}\int\limits_{[t,T]^k}
K^2(t_1,\ldots,t_k)dt_1\ldots dt_k.
$$

\vspace{4mm}

In \cite{arxiv-1}, \cite{arxiv-3}, \cite{20}-\cite{12aa-afterxxx} it was shown that 

\vspace{-1mm}
\begin{equation}
\label{zzz0}
{\sf M}\left\{\left(J[\psi^{(k)}]_{T,t}-
J[\psi^{(k)}]_{T,t}^{p_1,\ldots,p_k}\right)^2\right\}\le k!\left(I_k-\sum_{j_1=0}^{p_1}\ldots
\sum_{j_k=0}^{p_k}C^2_{j_k\ldots j_1}\right),
\end{equation}

\vspace{3mm}
\noindent
where $i_1,\ldots,i_k=1,\ldots,m$ and $0<T-t<\infty$ or 
$i_1,\ldots,i_k=0, 1,\ldots,m$ and $0<T-t<1.$

Moreover  \cite{arxiv-1} (Sect.~6, 15, 16), \cite{20aa} (Sect.~1.1.9, 1.11, 1.12),

\vspace{-1mm}
$$
{\sf M}\left\{\left(J[\psi^{(k)}]_{T,t}-
J[\psi^{(k)}]_{T,t}^{p_1,\ldots,p_k}\right)^{2n}\right\}\le C_{n,k}
\left(I_k-\sum_{j_1=0}^{p_1}\ldots
\sum_{j_k=0}^{p_k}C^2_{j_k\ldots j_1}\right)^n,
$$

\vspace{3mm}
\noindent
where 
$C_{n,k}=(k!)^{n} (2n-1)^{nk},$\ \ $n\in \mathbb{N}.$

Below 
we provide practical material concerning expansions and approximations of 
iterated Ito and 
Strato\-no\-vich stochastic integrals of multiplicities 1 to 6
(the case of Legendre
polynomial system).
The question about what kind of functions (polynomial
or trigonometric) is more convenient for the mean-square approximation
of iterated stochastic integrals is also considered.

Let us introduce more convenient (for further) notations
for the iterated Ito and Stratonovich stochastic integrals
\begin{equation}
\label{k1000}
I_{(l_1\ldots l_k)T,t}^{(i_1\ldots i_k)}
=\int\limits_t^T(t-t_k)^{l_k} \ldots \int\limits_t^{t_{2}}
(t-t_1)^{l_1} d{\bf f}_{t_1}^{(i_1)}\ldots
d{\bf f}_{t_k}^{(i_k)},
\end{equation}

\begin{equation}
\label{k1001}
I_{(l_1\ldots l_k)T,t}^{*(i_1\ldots i_k)}
={\int\limits_t^{*}}^T (t-t_k)^{l_k} \ldots {\int\limits_t^{*}}^{t_2}
(t-t_1)^{l_1} d{\bf f}_{t_1}^{(i_1)}\ldots
d{\bf f}_{t_k}^{(i_k)},
\end{equation}

\vspace{3mm}
\noindent
where $i_1,\ldots, i_k=1,\dots,m,$\ \  $l_1,\ldots,l_k=0, 1,\ldots$

The complete orthonormal system of Legendre polynomials in the 
space $L_2([t,T])$ looks as follows

\begin{equation}
\label{4009}
\phi_j(x)=\sqrt{\frac{2j+1}{T-t}}P_j\left(\left(
x-\frac{T+t}{2}\right)\frac{2}{T-t}\right),\ \ \ j=0, 1, 2,\ldots,
\end{equation}

\vspace{3mm}
\noindent
where $P_j(x)$ is the Legendre polynomial. 

Using the
system of 
functions (\ref{4009}) and Theorems 7--16,
we obtain the following expansions of iterated 
Ito and Stratonovich stochastic integrals (\ref{k1000}),
(\ref{k1001}) \cite{arxiv-2}, \cite{arxiv-4},
\cite{5-001}-\cite{5-004}, \cite{7}-\cite{new-art-1xxy}
(also see early publications
\cite{5-009} (2000), \cite{5-010} (2001),
\cite{322} (1997), \cite{444} (1998))

\vspace{2mm}
$$
I_{(0)T,t}^{(i_1)}=\sqrt{T-t}\zeta_0^{(i_1)},
$$

\vspace{2mm}
\begin{equation}
\label{4002}
I_{(1)T,t}^{(i_1)}=-\frac{(T-t)^{3/2}}{2}\left(\zeta_0^{(i_1)}+
\frac{1}{\sqrt{3}}\zeta_1^{(i_1)}\right),
\end{equation}

\vspace{3mm}
\begin{equation}
\label{4003}
I_{(2)T,t}^{(i_1)}=\frac{(T-t)^{5/2}}{3}\left(\zeta_0^{(i_1)}+
\frac{\sqrt{3}}{2}\zeta_1^{(i_1)}+
\frac{1}{2\sqrt{5}}\zeta_2^{(i_1)}\right),
\end{equation}

\vspace{3mm}
\begin{equation}
I_{(00)T,t}^{*(i_1 i_2)}=
\frac{T-t}{2}\left(\zeta_0^{(i_1)}\zeta_0^{(i_2)}+\sum_{i=1}^{\infty}
\frac{1}{\sqrt{4i^2-1}}\left(
\zeta_{i-1}^{(i_1)}\zeta_{i}^{(i_2)}-
\zeta_i^{(i_1)}\zeta_{i-1}^{(i_2)}\right)\right),
\label{4004}
\end{equation}

\vspace{6mm}

$$
I_{(00)T,t}^{(i_1 i_2)}=
\frac{T-t}{2}\left(\zeta_0^{(i_1)}\zeta_0^{(i_2)}+\sum_{i=1}^{\infty}
\frac{1}{\sqrt{4i^2-1}}\left(
\zeta_{i-1}^{(i_1)}\zeta_{i}^{(i_2)}-
\zeta_i^{(i_1)}\zeta_{i-1}^{(i_2)}\right) - {\bf 1}_{\{i_1=i_2\}}\right),
$$

\vspace{7mm}

$$
I_{(01)T,t}^{*(i_1 i_2)}=-\frac{T-t}{2}I_{(00)T,t}^{*(i_1 i_2)}
-\frac{(T-t)^2}{4}\Biggl(
\frac{1}{\sqrt{3}}\zeta_0^{(i_1)}\zeta_1^{(i_2)}+\Biggr.
$$

\vspace{2mm}
$$
+\Biggl.\sum_{i=0}^{\infty}\Biggl(
\frac{(i+2)\zeta_i^{(i_1)}\zeta_{i+2}^{(i_2)}
-(i+1)\zeta_{i+2}^{(i_1)}\zeta_{i}^{(i_2)}}
{\sqrt{(2i+1)(2i+5)}(2i+3)}-
\frac{\zeta_i^{(i_1)}\zeta_{i}^{(i_2)}}{(2i-1)(2i+3)}\Biggr)\Biggr),
$$

\vspace{7mm}

$$
I_{(10)T,t}^{*(i_1 i_2)}=-\frac{T-t}{2}I_{(00)T,t}^{*(i_1 i_2)}
-\frac{(T-t)^2}{4}\Biggl(
\frac{1}{\sqrt{3}}\zeta_0^{(i_2)}\zeta_1^{(i_1)}+\Biggr.
$$

\vspace{2mm}
\begin{equation}
\label{4006}
+\Biggl.\sum_{i=0}^{\infty}\Biggl(
\frac{(i+1)\zeta_{i+2}^{(i_2)}\zeta_{i}^{(i_1)}
-(i+2)\zeta_{i}^{(i_2)}\zeta_{i+2}^{(i_1)}}
{\sqrt{(2i+1)(2i+5)}(2i+3)}+
\frac{\zeta_i^{(i_1)}\zeta_{i}^{(i_2)}}{(2i-1)(2i+3)}\Biggr)\Biggr),
\end{equation}

\vspace{6mm}
\noindent 
or
$$
I_{(01)T,t}^{*(i_1i_2)}
=\hbox{\vtop{\offinterlineskip\halign{
\hfil#\hfil\cr
{\rm l.i.m.}\cr
$\stackrel{}{{}_{p\to \infty}}$\cr
}} }
\sum_{j_1,j_2=0}^{p}
C_{j_2j_1}^{01}
\zeta_{j_1}^{(i_1)}\zeta_{j_2}^{(i_2)},
$$

\vspace{3mm}
$$
I_{(10)T,t}^{*(i_1i_2)}
=\hbox{\vtop{\offinterlineskip\halign{
\hfil#\hfil\cr
{\rm l.i.m.}\cr
$\stackrel{}{{}_{p\to \infty}}$\cr
}} }
\sum_{j_1,j_2=0}^{p}
C_{j_2j_1}^{10}
\zeta_{j_1}^{(i_1)}\zeta_{j_2}^{(i_2)},
$$

\vspace{4mm}
\noindent
where

\vspace{-2mm}
$$
C_{j_2j_1}^{01}
=\frac{\sqrt{(2j_1+1)(2j_2+1)}}{8}(T-t)^{2}\bar
C_{j_2j_1}^{01},
$$

\vspace{3mm}
$$
C_{j_2j_1}^{10}
=\frac{\sqrt{(2j_1+1)(2j_2+1)}}{8}(T-t)^{2}\bar
C_{j_2j_1}^{10},
$$

\vspace{3mm}
$$
\bar C_{j_2j_1}^{01}=-\int\limits_{-1}^{1}(1+y)P_{j_2}(y)
\int\limits_{-1}^{y}
P_{j_1}(x)dx dy,
$$

\vspace{3mm}
$$
\bar C_{j_2j_1}^{10}=-\int\limits_{-1}^{1}P_{j_2}(y)
\int\limits_{-1}^{y}
(1+x)P_{j_1}(x)dx dy;
$$

\vspace{6mm}

$$
I_{(10)T,t}^{(i_1 i_2)}=
I_{(10)T,t}^{*(i_1 i_2)}+
\frac{1}{4}{\bf 1}_{\{i_1=i_2\}}(T-t)^2\ \ \ {\rm w.\ p.\ 1},
$$

\vspace{5mm}

$$
I_{(01)T,t}^{(i_1 i_2)}=
I_{(01)T,t}^{*(i_1 i_2)}+
\frac{1}{4}{\bf 1}_{\{i_1=i_2\}}(T-t)^2\ \ \ {\rm w.\ p.\ 1},
$$

\vspace{8mm}

$$
I_{(01)T,t}^{(i_1 i_2)}=
-\frac{T-t}{2}
I_{(00)T,t}^{(i_1 i_2)}
-\frac{(T-t)^2}{4}\Biggl(
\frac{1}{\sqrt{3}}\zeta_0^{(i_1)}\zeta_1^{(i_2)}+\Biggr.
$$

\vspace{2mm}
$$
+\Biggl.\sum_{i=0}^{\infty}\Biggl(
\frac{(i+2)\zeta_i^{(i_1)}\zeta_{i+2}^{(i_2)}
-(i+1)\zeta_{i+2}^{(i_1)}\zeta_{i}^{(i_2)}}
{\sqrt{(2i+1)(2i+5)}(2i+3)}-
\frac{\zeta_i^{(i_1)}\zeta_{i}^{(i_2)}}{(2i-1)(2i+3)}\Biggr)\Biggr),
$$

\vspace{7mm}

$$
I_{(10)T,t}^{(i_1 i_2)}=
-\frac{T-t}{2}I_{(00)T,t}^{(i_1 i_2)}
-\frac{(T-t)^2}{4}\Biggl(
\frac{1}{\sqrt{3}}\zeta_0^{(i_2)}\zeta_1^{(i_1)}+\Biggr.
$$

\vspace{2mm}
$$
+\Biggl.\sum_{i=0}^{\infty}\Biggl(
\frac{(i+1)\zeta_{i+2}^{(i_2)}\zeta_{i}^{(i_1)}
-(i+2)\zeta_{i}^{(i_2)}\zeta_{i+2}^{(i_1)}}
{\sqrt{(2i+1)(2i+5)}(2i+3)}+
\frac{\zeta_i^{(i_1)}\zeta_{i}^{(i_2)}}{(2i-1)(2i+3)}\Biggr)\Biggr),
$$

\vspace{7mm}
\noindent
or
$$
I_{(01)T,t}^{(i_1 i_2)}=
\hbox{\vtop{\offinterlineskip\halign{
\hfil#\hfil\cr
{\rm l.i.m.}\cr
$\stackrel{}{{}_{p\to \infty}}$\cr
}} }
\sum_{j_1,j_2=0}^{p}
C_{j_2j_1}^{01}\Biggl(\zeta_{j_1}^{(i_1)}\zeta_{j_2}^{(i_2)}
-{\bf 1}_{\{i_1=i_2\}}
{\bf 1}_{\{j_1=j_2\}}\Biggr),
$$

\vspace{5mm}
$$
I_{(10)T,t}^{(i_1 i_2)}=
\hbox{\vtop{\offinterlineskip\halign{
\hfil#\hfil\cr
{\rm l.i.m.}\cr
$\stackrel{}{{}_{p\to \infty}}$\cr
}} }
\sum_{j_1,j_2=0}^{p}
C_{j_2j_1}^{10}\Biggl(\zeta_{j_1}^{(i_1)}\zeta_{j_2}^{(i_2)}
-{\bf 1}_{\{i_1=i_2\}}
{\bf 1}_{\{j_1=j_2\}}\Biggr),
$$

\vspace{8mm}

$$
I_{(000)T,t}^{*(i_1 i_2 i_3)}=
\hbox{\vtop{\offinterlineskip\halign{
\hfil#\hfil\cr
{\rm l.i.m.}\cr
$\stackrel{}{{}_{p\to \infty}}$\cr
}} }
\sum\limits_{j_1, j_2, j_3=0}^{p}
C_{j_3 j_2 j_1}\zeta_{j_1}^{(i_1)}\zeta_{j_2}^{(i_2)}\zeta_{j_3}^{(i_3)},
$$

\vspace{7mm}

$$
I_{(000)T,t}^{(i_1i_2i_3)}
=\hbox{\vtop{\offinterlineskip\halign{
\hfil#\hfil\cr
{\rm l.i.m.}\cr
$\stackrel{}{{}_{p\to \infty}}$\cr
}} }
\sum_{j_1,j_2,j_3=0}^{p}
C_{j_3j_2j_1}
\Biggl(
\zeta_{j_1}^{(i_1)}\zeta_{j_2}^{(i_2)}\zeta_{j_3}^{(i_3)}
-{\bf 1}_{\{i_1=i_2\}}
{\bf 1}_{\{j_1=j_2\}}
\zeta_{j_3}^{(i_3)}-
\Biggr.
$$

\vspace{1mm}
\begin{equation}
\label{zzz1}
\Biggl.
-{\bf 1}_{\{i_2=i_3\}}
{\bf 1}_{\{j_2=j_3\}}
\zeta_{j_1}^{(i_1)}-
{\bf 1}_{\{i_1=i_3\}}
{\bf 1}_{\{j_1=j_3\}}
\zeta_{j_2}^{(i_2)}\Biggr),
\end{equation}

\vspace{7mm}

$$
I_{(000)T,t}^{(i_1 i_1 i_1)}
=\frac{1}{6}(T-t)^{3/2}\biggl(
\left(\zeta_0^{(i_1)}\right)^3-3
\zeta_0^{(i_1)}\biggr)\ \ \ \hbox{w.\ p.\ 1},
$$

\vspace{6mm}

$$
I_{(000)T,t}^{*(i_1 i_1 i_1)}
=\frac{1}{6}(T-t)^{3/2}
\left(\zeta_0^{(i_1)}\right)^3\ \ \ \hbox{w.\ p.\ 1},
$$

\vspace{5mm}
\noindent
where

\begin{equation}
\label{zzz2}
C_{j_3j_2j_1}
=\frac{\sqrt{(2j_1+1)(2j_2+1)(2j_3+1)}}{8}(T-t)^{3/2}\bar
C_{j_3j_2j_1},
\end{equation}

\vspace{2mm}
\begin{equation}
\label{zzz3}
\bar C_{j_3j_2j_1}=\int\limits_{-1}^{1}P_{j_3}(z)
\int\limits_{-1}^{z}P_{j_2}(y)
\int\limits_{-1}^{y}
P_{j_1}(x)dx dy dz;
\end{equation}

\vspace{7mm}

$$
I_{(000)T,t}^{(i_1 i_2 i_3)}=I_{(000)T,t}^{*(i_1 i_2 i_3)}
+{\bf 1}_{\{i_1=i_2\ne 0\}}\frac{1}{2}I_{(1)T,t}^{(i_3)}-
$$

\vspace{3mm}

$$
-
{\bf 1}_{\{i_2=i_3\ne 0\}}\frac{1}{2}\left((T-t)
I_{(0)T,t}^{(i_1)}+I_{(1)T,t}^{(i_1)}\right)\ \ \ \hbox{w.\ p.\ 1},
$$

\vspace{8mm}

$$
I_{(02)T,t}^{*(i_1 i_2)}=-\frac{(T-t)^2}{4}I_{(00)T,t}^{*(i_1 i_2)}
-(T-t)I_{(01)T,t}^{*(i_1 i_2)}+
\frac{(T-t)^3}{8}\Biggl[
\frac{2}{3\sqrt{5}}\zeta_2^{(i_2)}\zeta_0^{(i_1)}+\Biggr.
$$

\vspace{5mm}
$$
+\frac{1}{3}\zeta_0^{(i_1)}\zeta_0^{(i_2)}+
\sum_{i=0}^{\infty}\Biggl(
\frac{(i+2)(i+3)\zeta_{i+3}^{(i_2)}\zeta_{i}^{(i_1)}
-(i+1)(i+2)\zeta_{i}^{(i_2)}\zeta_{i+3}^{(i_1)}}
{\sqrt{(2i+1)(2i+7)}(2i+3)(2i+5)}+
\Biggr.
$$

\vspace{5mm}
\begin{equation}
\label{leto1}
\Biggl.\Biggl.+
\frac{(i^2+i-3)\zeta_{i+1}^{(i_2)}\zeta_{i}^{(i_1)}
-(i^2+3i-1)\zeta_{i}^{(i_2)}\zeta_{i+1}^{(i_1)}}
{\sqrt{(2i+1)(2i+3)}(2i-1)(2i+5)}\Biggr)\Biggr],
\end{equation}

\vspace{12mm}

$$
I_{(20)T,t}^{*(i_1 i_2)}=-\frac{(T-t)^2}{4}I_{(00)T,t}^{*(i_1 i_2)}
-(T-t)I_{(10)T,t}^{*(i_1 i_2)}+
\frac{(T-t)^3}{8}\Biggl[
\frac{2}{3\sqrt{5}}\zeta_0^{(i_2)}\zeta_2^{(i_1)}+\Biggr.
$$

\vspace{5mm}
$$
+\frac{1}{3}\zeta_0^{(i_1)}\zeta_0^{(i_2)}+
\sum_{i=0}^{\infty}\Biggl(
\frac{(i+1)(i+2)\zeta_{i+3}^{(i_2)}\zeta_{i}^{(i_1)}
-(i+2)(i+3)\zeta_{i}^{(i_2)}\zeta_{i+3}^{(i_1)}}
{\sqrt{(2i+1)(2i+7)}(2i+3)(2i+5)}+
\Biggr.
$$

\vspace{5mm}

\begin{equation}
\label{leto2}
\Biggl.\Biggl.+
\frac{(i^2+3i-1)\zeta_{i+1}^{(i_2)}\zeta_{i}^{(i_1)}
-(i^2+i-3)\zeta_{i}^{(i_2)}\zeta_{i+1}^{(i_1)}}
{\sqrt{(2i+1)(2i+3)}(2i-1)(2i+5)}\Biggr)\Biggr],
\end{equation}

\vspace{14mm}

$$
I_{(11)T,t}^{*(i_1 i_2)}=-\frac{(T-t)^2}{4}I_{(00)T,t}^{*(i_1 i_2)}-
\frac{(T-t)}{2}\left(I_{(10)T,t}^{*(i_1 i_2)}+
I_{(01)T,t}^{*(i_1 i_2)}\right)+
$$

\vspace{4mm}
$$
+\frac{(T-t)^3}{8}\Biggl[
\frac{1}{3}\zeta_1^{(i_1)}\zeta_1^{(i_2)}+
\sum_{i=0}^{\infty}\Biggl(
\frac{(i+1)(i+3)\left(\zeta_{i+3}^{(i_2)}\zeta_{i}^{(i_1)}
-\zeta_{i}^{(i_2)}\zeta_{i+3}^{(i_1)}\right)}
{\sqrt{(2i+1)(2i+7)}(2i+3)(2i+5)}+
\Biggr.
$$

\vspace{4mm}
\begin{equation}
\label{leto3}
\Biggl.\Biggl.+
\frac{(i+1)^2\left(\zeta_{i+1}^{(i_2)}\zeta_{i}^{(i_1)}
-\zeta_{i}^{(i_2)}\zeta_{i+1}^{(i_1)}\right)}
{\sqrt{(2i+1)(2i+3)}(2i-1)(2i+5)}\Biggr)\Biggr],
\end{equation}

\vspace{13mm}
or
$$
I_{(02)T,t}^{*(i_1 i_2)}=
\hbox{\vtop{\offinterlineskip\halign{
\hfil#\hfil\cr
{\rm l.i.m.}\cr
$\stackrel{}{{}_{p\to \infty}}$\cr
}} }
\sum_{j_1,j_2=0}^{p}
C_{j_2j_1}^{02}\zeta_{j_1}^{(i_1)}\zeta_{j_2}^{(i_2)},
$$

\vspace{3mm}

$$
I_{(20)T,t}^{*(i_1 i_2)}=
\hbox{\vtop{\offinterlineskip\halign{
\hfil#\hfil\cr
{\rm l.i.m.}\cr
$\stackrel{}{{}_{p\to \infty}}$\cr
}} }
\sum_{j_1,j_2=0}^{p}
C_{j_2j_1}^{20}\zeta_{j_1}^{(i_1)}\zeta_{j_2}^{(i_2)},
$$

\vspace{3mm}

$$
I_{(11)T,t}^{*(i_1 i_2)}=
\hbox{\vtop{\offinterlineskip\halign{
\hfil#\hfil\cr
{\rm l.i.m.}\cr
$\stackrel{}{{}_{p\to \infty}}$\cr
}} }
\sum_{j_1,j_2=0}^{p}
C_{j_2j_1}^{11}\zeta_{j_1}^{(i_1)}\zeta_{j_2}^{(i_2)},
$$

\vspace{5mm}
\noindent
where

\vspace{-2mm}
$$
C_{j_2j_1}^{02}=
\frac{\sqrt{(2j_1+1)(2j_2+1)}}{16}(T-t)^{3}\bar
C_{j_2j_1}^{02},
$$

\vspace{3mm}
$$
C_{j_2j_1}^{20}=
\frac{\sqrt{(2j_1+1)(2j_2+1)}}{16}(T-t)^{3}\bar
C_{j_2j_1}^{20},
$$

\vspace{3mm}

$$
C_{j_2j_1}^{11}=
\frac{\sqrt{(2j_1+1)(2j_2+1)}}{16}(T-t)^{3}\bar
C_{j_2j_1}^{11}, 
$$

\vspace{3mm}
$$
\bar C_{j_2j_1}^{02}=
\int\limits_{-1}^{1}P_{j_2}(y)(y+1)^2
\int\limits_{-1}^{y}
P_{j_1}(x)dx dy,
$$

\vspace{3mm}
$$
\bar C_{j_2j_1}^{20}=
\int\limits_{-1}^{1}P_{j_2}(y)
\int\limits_{-1}^{y}
P_{j_1}(x)(x+1)^2 dx dy,
$$

\vspace{3mm}
$$
\bar C_{j_2j_1}^{11}=
\int\limits_{-1}^{1}P_{j_2}(y)(y+1)
\int\limits_{-1}^{y}
P_{j_1}(x)(x+1)dx dy,
$$

\vspace{6mm}

$$
I_{(11)T,t}^{*(i_1 i_1)}=\frac{1}{2}\left(I_{(1)T,t}^{(i_1)}
\right)^2\ \ \ \hbox{w.\ p.\ 1,}
$$

\vspace{3mm}

$$
I_{(02)T,t}^{(i_1 i_2)}=
I_{(02)T,t}^{*(i_1 i_2)}-
\frac{1}{6}{\bf 1}_{\{i_1=i_2\}}(T-t)^3\ \ \ \hbox{w.\ p.\ 1},
$$

\vspace{3mm}

$$
I_{(20)T,t}^{(i_1 i_2)}=
I_{(20)T,t}^{*(i_1 i_2)}-
\frac{1}{6}{\bf 1}_{\{i_1=i_2\}}(T-t)^3\ \ \ \hbox{w.\ p.\ 1},
$$

\vspace{3mm}

$$
I_{(11)T,t}^{(i_1 i_2)}=
I_{(11)T,t}^{*(i_1 i_2)}-
\frac{1}{6}{\bf 1}_{\{i_1=i_2\}}(T-t)^3\ \ \ \hbox{w.\ p.\ 1},
$$

\vspace{12mm}

$$
I_{(02)T,t}^{(i_1 i_2)}
=-\frac{(T-t)^2}{4}I_{(00)T,t}^{(i_1 i_2)}
-(T-t) I_{01_{T,t}}^{(i_1 i_2)}+
\frac{(T-t)^3}{8}\Biggl[
\frac{2}{3\sqrt{5}}\zeta_2^{(i_2)}\zeta_0^{(i_1)}+\Biggr.
$$

\vspace{6mm}
$$
+\frac{1}{3}\zeta_0^{(i_1)}\zeta_0^{(i_2)}+
\sum_{i=0}^{\infty}\Biggl(
\frac{(i+2)(i+3)\zeta_{i+3}^{(i_2)}\zeta_{i}^{(i_1)}
-(i+1)(i+2)\zeta_{i}^{(i_2)}\zeta_{i+3}^{(i_1)}}
{\sqrt{(2i+1)(2i+7)}(2i+3)(2i+5)}+
\Biggr.
$$

\vspace{6mm}
$$
\Biggl.\Biggl.+
\frac{(i^2+i-3)\zeta_{i+1}^{(i_2)}\zeta_{i}^{(i_1)}
-(i^2+3i-1)\zeta_{i}^{(i_2)}\zeta_{i+1}^{(i_1)}}
{\sqrt{(2i+1)(2i+3)}(2i-1)(2i+5)}\Biggr)\Biggr] - 
$$

\vspace{6mm}
$$
-\frac{1}{24}{\bf 1}_{\{i_1=i_2\}}{(T-t)^3},
$$

\vspace{12mm}
$$
I_{(20)T,t}^{(i_1 i_2)}=-\frac{(T-t)^2}{4}
I_{(00)T,t}^{(i_1 i_2)}
-(T-t) I_{(10)T,t}^{(i_1 i_2)}+
\frac{(T-t)^3}{8}\Biggl[
\frac{2}{3\sqrt{5}}\zeta_0^{(i_2)}\zeta_2^{(i_1)}+\Biggr.
$$

\vspace{5mm}
$$
+\frac{1}{3}\zeta_0^{(i_1)}\zeta_0^{(i_2)}+
\sum_{i=0}^{\infty}\Biggl(
\frac{(i+1)(i+2)\zeta_{i+3}^{(i_2)}\zeta_{i}^{(i_1)}
-(i+2)(i+3)\zeta_{i}^{(i_2)}\zeta_{i+3}^{(i_1)}}
{\sqrt{(2i+1)(2i+7)}(2i+3)(2i+5)}+
\Biggr.
$$

\vspace{5mm}
$$
\Biggl.\Biggl.+
\frac{(i^2+3i-1)\zeta_{i+1}^{(i_2)}\zeta_{i}^{(i_1)}
-(i^2+i-3)\zeta_{i}^{(i_2)}\zeta_{i+1}^{(i_1)}}
{\sqrt{(2i+1)(2i+3)}(2i-1)(2i+5)}\Biggr)\Biggr] - 
$$

\vspace{4mm}
$$
-
\frac{1}{24}{\bf 1}_{\{i_1=i_2\}}{(T-t)^3},
$$

\vspace{11mm}

$$
I_{(11)T,t}^{(i_1 i_2)}
=-\frac{(T-t)^2}{4}I_{(00)T,t}^{(i_1 i_2)}
-\frac{T-t}{2}\left(
I_{(10)T,t}^{(i_1 i_2)}+
I_{(01)T,t}^{(i_1 i_2)}\right)+
\frac{(T-t)^3}{8}\Biggl[
\frac{1}{3}\zeta_1^{(i_1)}\zeta_1^{(i_2)}+\Biggr.
$$

\vspace{4mm}
$$
+
\sum_{i=0}^{\infty}\Biggl(
\frac{(i+1)(i+3)\left(\zeta_{i+3}^{(i_2)}\zeta_{i}^{(i_1)}
-\zeta_{i}^{(i_2)}\zeta_{i+3}^{(i_1)}\right)}
{\sqrt{(2i+1)(2i+7)}(2i+3)(2i+5)}+
\Biggr.
$$

\vspace{4mm}
$$
\Biggl.\Biggl.
+\frac{(i+1)^2\left(\zeta_{i+1}^{(i_2)}\zeta_{i}^{(i_1)}
-\zeta_{i}^{(i_2)}\zeta_{i+1}^{(i_1)}\right)}
{\sqrt{(2i+1)(2i+3)}(2i-1)(2i+5)}\Biggr)\Biggr] - 
$$

\vspace{4mm}
$$
-
\frac{1}{24}{\bf 1}_{\{i_1=i_2\}}{(T-t)^3},
$$

\vspace{12mm}
or
$$
I_{(02)T,t}^{(i_1 i_2)}=
\hbox{\vtop{\offinterlineskip\halign{
\hfil#\hfil\cr
{\rm l.i.m.}\cr
$\stackrel{}{{}_{p\to \infty}}$\cr
}} }
\sum_{j_1,j_2=0}^p
C_{j_2j_1}^{02}\Biggl(\zeta_{j_1}^{(i_1)}\zeta_{j_2}^{(i_2)}
-{\bf 1}_{\{i_1=i_2\}}
{\bf 1}_{\{j_1=j_2\}}\Biggr),
$$

\vspace{4mm}
$$
I_{(20)T,t}^{(i_1 i_2)}=
\hbox{\vtop{\offinterlineskip\halign{
\hfil#\hfil\cr
{\rm l.i.m.}\cr
$\stackrel{}{{}_{p\to \infty}}$\cr
}} }
\sum_{j_1,j_2=0}^{p}
C_{j_2j_1}^{20}\Biggl(\zeta_{j_1}^{(i_1)}\zeta_{j_2}^{(i_2)}
-{\bf 1}_{\{i_1=i_2\}}
{\bf 1}_{\{j_1=j_2\}}\Biggr),
$$

\vspace{4mm}
$$
I_{(11)T,t}^{(i_1 i_2)}=
\hbox{\vtop{\offinterlineskip\halign{
\hfil#\hfil\cr
{\rm l.i.m.}\cr
$\stackrel{}{{}_{p\to \infty}}$\cr
}} }
\sum_{j_1,j_2=0}^{p}
C_{j_2j_1}^{11}\Biggl(\zeta_{j_1}^{(i_1)}\zeta_{j_2}^{(i_2)}
-{\bf 1}_{\{i_1=i_2\}}
{\bf 1}_{\{j_1=j_2\}}\Biggr),
$$

\vspace{12mm}

$$
I_{(3)T,t}^{(i_1)}=-\frac{(T-t)^{7/2}}{4}\left(\zeta_0^{(i_1)}+
\frac{3\sqrt{3}}{5}\zeta_1^{(i_1)}+
\frac{1}{\sqrt{5}}\zeta_2^{(i_1)}+
\frac{1}{5\sqrt{7}}\zeta_3^{(i_1)}\right),
$$

\vspace{7mm}

$$
I_{(0000)T,t}^{*(i_1 i_2 i_3 i_4)}=
\hbox{\vtop{\offinterlineskip\halign{
\hfil#\hfil\cr
{\rm l.i.m.}\cr
$\stackrel{}{{}_{p\to \infty}}$\cr
}} }
\sum\limits_{j_1, j_2, j_3, j_4=0}^{p}
C_{j_4 j_3 j_2 j_1}\zeta_{j_1}^{(i_1)}\zeta_{j_2}^{(i_2)}\zeta_{j_3}^{(i_3)}
\zeta_{j_4}^{(i_4)},
$$

\vspace{10mm}

$$
I_{(0000)T,t}^{(i_1 i_2 i_3 i_4)}
=
\hbox{\vtop{\offinterlineskip\halign{
\hfil#\hfil\cr
{\rm l.i.m.}\cr
$\stackrel{}{{}_{p\to \infty}}$\cr
}} }
\sum_{j_1,j_2,j_3,j_4=0}^{p}
C_{j_4 j_3 j_2 j_1}\Biggl(
\zeta_{j_1}^{(i_1)}\zeta_{j_2}^{(i_2)}\zeta_{j_3}^{(i_3)}\zeta_{j_4}^{(i_4)}
-\Biggr.
$$
$$
-
{\bf 1}_{\{i_1=i_2\}}
{\bf 1}_{\{j_1=j_2\}}
\zeta_{j_3}^{(i_3)}
\zeta_{j_4}^{(i_4)}
-
{\bf 1}_{\{i_1=i_3\}}
{\bf 1}_{\{j_1=j_3\}}
\zeta_{j_2}^{(i_2)}
\zeta_{j_4}^{(i_4)}-
$$

\vspace{-2mm}
$$
-
{\bf 1}_{\{i_1=i_4\}}
{\bf 1}_{\{j_1=j_4\}}
\zeta_{j_2}^{(i_2)}
\zeta_{j_3}^{(i_3)}
-
{\bf 1}_{\{i_2=i_3\}}
{\bf 1}_{\{j_2=j_3\}}
\zeta_{j_1}^{(i_1)}
\zeta_{j_4}^{(i_4)}-
$$

\vspace{-2mm}
$$
-
{\bf 1}_{\{i_2=i_4\}}
{\bf 1}_{\{j_2=j_4\}}
\zeta_{j_1}^{(i_1)}
\zeta_{j_3}^{(i_3)}
-
{\bf 1}_{\{i_3=i_4\}}
{\bf 1}_{\{j_3=j_4\}}
\zeta_{j_1}^{(i_1)}
\zeta_{j_2}^{(i_2)}+
$$

\vspace{-2mm}
$$
+
{\bf 1}_{\{i_1=i_2\}}
{\bf 1}_{\{j_1=j_2\}}
{\bf 1}_{\{i_3=i_4\}}
{\bf 1}_{\{j_3=j_4\}}+
$$

\vspace{-2mm}
$$
+
{\bf 1}_{\{i_1=i_3\}}
{\bf 1}_{\{j_1=j_3\}}
{\bf 1}_{\{i_2=i_4\}}
{\bf 1}_{\{j_2=j_4\}}+
$$
\begin{equation}
\label{zzz10}
+\Biggl.
{\bf 1}_{\{i_1=i_4\}}
{\bf 1}_{\{j_1=j_4\}}
{\bf 1}_{\{i_2=i_3\}}
{\bf 1}_{\{j_2=j_3\}}\Biggr),
\end{equation}

\vspace{9mm}

$$
I_{(0000)T,t}^{(i_1i_1i_1i_1)}=
\frac{1}{24}(T-t)^2
\left(\left(\zeta_0^{(i_1)}\right)^4-
6\left(\zeta_0^{(i_1)}\right)^2+3\right)\ \ \ \hbox{w.\ p.\ 1},
$$

\vspace{5mm}

$$
I_{(0000)T,t}^{*(i_1i_1i_1i_1)}=
\frac{1}{24}(T-t)^2
\left(\zeta_0^{(i_1)}\right)^4\ \ \ \hbox{w.\ p.\ 1},
$$

\vspace{5mm}
\noindent
where

\begin{equation}
\label{zzz11}
C_{j_4j_3j_2j_1}=
\frac{\sqrt{(2j_1+1)(2j_2+1)(2j_3+1)(2j_4+1)}}{16}(T-t)^{2}\bar
C_{j_4j_3j_2j_1},
\end{equation}

\vspace{3mm}

\begin{equation}
\label{zzz12}
\bar C_{j_4j_3j_2j_1}=\int\limits_{-1}^{1}P_{j_4}(u)
\int\limits_{-1}^{u}P_{j_3}(z)
\int\limits_{-1}^{z}P_{j_2}(y)
\int\limits_{-1}^{y}
P_{j_1}(x)dx dy dz du;
\end{equation}

\vspace{7mm}

$$
I_{(001)T,t}^{*(i_1i_2i_3)}
=\hbox{\vtop{\offinterlineskip\halign{
\hfil#\hfil\cr
{\rm l.i.m.}\cr
$\stackrel{}{{}_{p\to \infty}}$\cr
}} }
\sum_{j_1,j_2,j_3=0}^{p}
C_{j_3 j_2 j_1}^{001}
\zeta_{j_1}^{(i_1)}\zeta_{j_2}^{(i_2)}\zeta_{j_3}^{(i_3)},
$$

\vspace{5mm}

$$
I_{(010)T,t}^{*(i_1i_2i_3)}
=\hbox{\vtop{\offinterlineskip\halign{
\hfil#\hfil\cr
{\rm l.i.m.}\cr
$\stackrel{}{{}_{p\to \infty}}$\cr
}} }
\sum_{j_1,j_2,j_3=0}^{p}
C_{j_3 j_2 j_1}^{010}
\zeta_{j_1}^{(i_1)}\zeta_{j_2}^{(i_2)}\zeta_{j_3}^{(i_3)},
$$

\vspace{5mm}

$$
I_{(100)T,t}^{*(i_1i_2i_3)}
=\hbox{\vtop{\offinterlineskip\halign{
\hfil#\hfil\cr
{\rm l.i.m.}\cr
$\stackrel{}{{}_{p\to \infty}}$\cr
}} }
\sum_{j_1,j_2,j_3=0}^{p}
C_{j_3 j_2 j_1}^{100}
\zeta_{j_1}^{(i_1)}\zeta_{j_2}^{(i_2)}\zeta_{j_3}^{(i_3)},
$$

\vspace{10mm}

$$
I_{(001)T,t}^{(i_1i_2i_3)}
=\hbox{\vtop{\offinterlineskip\halign{
\hfil#\hfil\cr
{\rm l.i.m.}\cr
$\stackrel{}{{}_{p\to \infty}}$\cr
}} }
\sum_{j_1,j_2,j_3=0}^{p}
C_{j_3j_2j_1}^{001}\Biggl(
\zeta_{j_1}^{(i_1)}\zeta_{j_2}^{(i_2)}\zeta_{j_3}^{(i_3)}
-{\bf 1}_{\{i_1=i_2\}}
{\bf 1}_{\{j_1=j_2\}}
\zeta_{j_3}^{(i_3)}-
\Biggr.
$$

\vspace{1mm}
\begin{equation}
\label{sss1}
\Biggl.
-{\bf 1}_{\{i_2=i_3\}}
{\bf 1}_{\{j_2=j_3\}}
\zeta_{j_1}^{(i_1)}-
{\bf 1}_{\{i_1=i_3\}}
{\bf 1}_{\{j_1=j_3\}}
\zeta_{j_2}^{(i_2)}\Biggr),
\end{equation}

\vspace{8mm}

$$
I_{(010)T,t}^{(i_1i_2i_3)}
=\hbox{\vtop{\offinterlineskip\halign{
\hfil#\hfil\cr
{\rm l.i.m.}\cr
$\stackrel{}{{}_{p\to \infty}}$\cr
}} }
\sum_{j_1,j_2,j_3=0}^{p}
C_{j_3j_2j_1}^{010}\Biggl(
\zeta_{j_1}^{(i_1)}\zeta_{j_2}^{(i_2)}\zeta_{j_3}^{(i_3)}
-{\bf 1}_{\{i_1=i_2\}}
{\bf 1}_{\{j_1=j_2\}}
\zeta_{j_3}^{(i_3)}-
\Biggr.
$$

\vspace{1mm}
\begin{equation}
\label{sss2}
\Biggl.
-{\bf 1}_{\{i_2=i_3\}}
{\bf 1}_{\{j_2=j_3\}}
\zeta_{j_1}^{(i_1)}-
{\bf 1}_{\{i_1=i_3\}}
{\bf 1}_{\{j_1=j_3\}}
\zeta_{j_2}^{(i_2)}\Biggr),
\end{equation}

\vspace{8mm}

$$
I_{(100)T,t}^{(i_1i_2i_3)}
=\hbox{\vtop{\offinterlineskip\halign{
\hfil#\hfil\cr
{\rm l.i.m.}\cr
$\stackrel{}{{}_{p\to \infty}}$\cr
}} }
\sum_{j_1,j_2,j_3=0}^{p}
C_{j_3j_2j_1}^{100}\Biggl(
\zeta_{j_1}^{(i_1)}\zeta_{j_2}^{(i_2)}\zeta_{j_3}^{(i_3)}
-{\bf 1}_{\{i_1=i_2\}}
{\bf 1}_{\{j_1=j_2\}}
\zeta_{j_3}^{(i_3)}-
\Biggr.
$$

\vspace{1mm}
\begin{equation}
\label{sss3}
\Biggl.
-{\bf 1}_{\{i_2=i_3\}}
{\bf 1}_{\{j_2=j_3\}}
\zeta_{j_1}^{(i_1)}-
{\bf 1}_{\{i_1=i_3\}}
{\bf 1}_{\{j_1=j_3\}}
\zeta_{j_2}^{(i_2)}\Biggr),
\end{equation}

\vspace{6mm}
\noindent
where

\vspace{-2mm}
$$
C_{j_3j_2j_1}^{001}
=\frac{\sqrt{(2j_1+1)(2j_2+1)(2j_3+1)}}{16}(T-t)^{5/2}\bar
C_{j_3j_2j_1}^{001},
$$

\vspace{4mm}

$$
C_{j_3j_2j_1}^{010}
=\frac{\sqrt{(2j_1+1)(2j_2+1)(2j_3+1)}}{16}(T-t)^{5/2}\bar
C_{j_3j_2j_1}^{010},
$$

\vspace{4mm}

$$
C_{j_3j_2j_1}^{100}
=\frac{\sqrt{(2j_1+1)(2j_2+1)(2j_3+1)}}{16}(T-t)^{5/2}\bar
C_{j_3j_2j_1}^{100},
$$

\vspace{6mm}

$$
\bar C_{j_3j_2j_1}^{100}=-
\int\limits_{-1}^{1}P_{j_3}(z)
\int\limits_{-1}^{z}P_{j_2}(y)
\int\limits_{-1}^{y}
P_{j_1}(x)(x+1)dx dy dz,
$$

\vspace{4mm}

$$
\bar C_{j_3j_2j_1}^{010}=-
\int\limits_{-1}^{1}P_{j_3}(z)
\int\limits_{-1}^{z}P_{j_2}(y)(y+1)
\int\limits_{-1}^{y}
P_{j_1}(x)dx dy dz,
$$

\vspace{4mm}

$$
\bar C_{j_3j_2j_1}^{001}=-
\int\limits_{-1}^{1}P_{j_3}(z)(z+1)
\int\limits_{-1}^{z}P_{j_2}(y)
\int\limits_{-1}^{y}
P_{j_1}(x)dx dy dz;
$$

\vspace{9mm}

$$
I_{(lll)T,t}^{(i_1i_1i_1)}=
\frac{1}{6}\left(\left(I_{(l)T,t}^{(i_1)}\right)^3-
3I_{(l)T,t}^{(i_1)}\Delta_{l(T,t)}\right)\ \ \ \hbox{w.\ p.\ 1},
$$

\vspace{4mm}
$$
I_{(lll)T,t}^{*(i_1i_1i_1)}=
\frac{1}{6}\left(I_{(l)T,t}^{(i_1)}\right)^3\ \ \ \hbox{w.\ p.\ 1},
$$

\vspace{4mm}

$$
I_{(llll)T,t}^{(i_1i_1i_1i_1)}=
\frac{1}{24}\left(\left(I_{(l)T,t}^{(i_1)}\right)^4-
6\left(I_{(l)T,t}^{(i_1)}\right)^2\Delta_{(l)T,t}+3
\left(\Delta_{(l)T,t}\right)^2\right)\ \ \ \hbox{w.\ p.\ 1},
$$

\vspace{5mm}

$$
I_{(llll)T,t}^{*(i_1i_1i_1i_1)}=
\frac{1}{24}\left(I_{(l)T,t}^{(i_1)}\right)^4\ \ \ \hbox{w.\ p.\ 1},
$$

\vspace{5mm}
\noindent
where

\vspace{-2mm}
$$
I_{(l)T,t}^{(i_1)}=\sum_{j=0}^l C_j^l \zeta_j^{(i_1)}\ \ \ \hbox{w.\ p.\ 1},
$$       

\vspace{5mm}

$$
\Delta_{l(T,t)}=\int\limits_t^T(t-s)^{2l}ds,\ \ \
C_j^l=\int\limits_t^T(t-s)^l\phi_j(s)ds;
$$

\vspace{6mm}

$$
I_{(00000)T,t}^{*(i_1 i_2 i_3 i_4 i_5)}=
\hbox{\vtop{\offinterlineskip\halign{
\hfil#\hfil\cr
{\rm l.i.m.}\cr
$\stackrel{}{{}_{p\to \infty}}$\cr
}} }
\sum\limits_{j_1, j_2, j_3, j_4, j_5=0}^{p}
C_{j_5j_4 j_3 j_2 j_1}
\zeta_{j_1}^{(i_1)}\zeta_{j_2}^{(i_2)}\zeta_{j_3}^{(i_3)}
\zeta_{j_4}^{(i_4)}\zeta_{j_5}^{(i_5)},
$$

\vspace{8mm}

$$
I_{(00000)T,t}^{(i_1 i_2 i_3 i_4 i_5)}
=
\hbox{\vtop{\offinterlineskip\halign{
\hfil#\hfil\cr
{\rm l.i.m.}\cr
$\stackrel{}{{}_{p\to \infty}}$\cr
}} }
\sum_{j_1,j_2,j_3,j_4,j_5=0}^p
C_{j_5 j_4 j_3 j_2 j_1}\Biggl(
\prod_{l=1}^5\zeta_{j_l}^{(i_l)}
-\Biggr.
$$
$$
-
{\bf 1}_{\{i_1=i_2\}}
{\bf 1}_{\{j_1=j_2\}}
\zeta_{j_3}^{(i_3)}
\zeta_{j_4}^{(i_4)}
\zeta_{j_5}^{(i_5)}-
{\bf 1}_{\{i_1=i_3\}}
{\bf 1}_{\{j_1=j_3\}}
\zeta_{j_2}^{(i_2)}
\zeta_{j_4}^{(i_4)}
\zeta_{j_5}^{(i_5)}-
$$
$$
-
{\bf 1}_{\{i_1=i_4\}}
{\bf 1}_{\{j_1=j_4\}}
\zeta_{j_2}^{(i_2)}
\zeta_{j_3}^{(i_3)}
\zeta_{j_5}^{(i_5)}-
{\bf 1}_{\{i_1=i_5\}}
{\bf 1}_{\{j_1=j_5\}}
\zeta_{j_2}^{(i_2)}
\zeta_{j_3}^{(i_3)}
\zeta_{j_4}^{(i_4)}-
$$
$$
-
{\bf 1}_{\{i_2=i_3\}}
{\bf 1}_{\{j_2=j_3\}}
\zeta_{j_1}^{(i_1)}
\zeta_{j_4}^{(i_4)}
\zeta_{j_5}^{(i_5)}-
{\bf 1}_{\{i_2=i_4\}}
{\bf 1}_{\{j_2=j_4\}}
\zeta_{j_1}^{(i_1)}
\zeta_{j_3}^{(i_3)}
\zeta_{j_5}^{(i_5)}-
$$
$$
-
{\bf 1}_{\{i_2=i_5\}}
{\bf 1}_{\{j_2=j_5\}}
\zeta_{j_1}^{(i_1)}
\zeta_{j_3}^{(i_3)}
\zeta_{j_4}^{(i_4)}
-{\bf 1}_{\{i_3=i_4\}}
{\bf 1}_{\{j_3=j_4\}}
\zeta_{j_1}^{(i_1)}
\zeta_{j_2}^{(i_2)}
\zeta_{j_5}^{(i_5)}-
$$
$$
-
{\bf 1}_{\{i_3=i_5\}}
{\bf 1}_{\{j_3=j_5\}}
\zeta_{j_1}^{(i_1)}
\zeta_{j_2}^{(i_2)}
\zeta_{j_4}^{(i_4)}
-{\bf 1}_{\{i_4=i_5\}}
{\bf 1}_{\{j_4=j_5\}}
\zeta_{j_1}^{(i_1)}
\zeta_{j_2}^{(i_2)}
\zeta_{j_3}^{(i_3)}+
$$
$$
+
{\bf 1}_{\{i_1=i_2\}}
{\bf 1}_{\{j_1=j_2\}}
{\bf 1}_{\{i_3=i_4\}}
{\bf 1}_{\{j_3=j_4\}}\zeta_{j_5}^{(i_5)}+
{\bf 1}_{\{i_1=i_2\}}
{\bf 1}_{\{j_1=j_2\}}
{\bf 1}_{\{i_3=i_5\}}
{\bf 1}_{\{j_3=j_5\}}\zeta_{j_4}^{(i_4)}+
$$
$$
+
{\bf 1}_{\{i_1=i_2\}}
{\bf 1}_{\{j_1=j_2\}}
{\bf 1}_{\{i_4=i_5\}}
{\bf 1}_{\{j_4=j_5\}}\zeta_{j_3}^{(i_3)}+
{\bf 1}_{\{i_1=i_3\}}
{\bf 1}_{\{j_1=j_3\}}
{\bf 1}_{\{i_2=i_4\}}
{\bf 1}_{\{j_2=j_4\}}\zeta_{j_5}^{(i_5)}+
$$
$$
+
{\bf 1}_{\{i_1=i_3\}}
{\bf 1}_{\{j_1=j_3\}}
{\bf 1}_{\{i_2=i_5\}}
{\bf 1}_{\{j_2=j_5\}}\zeta_{j_4}^{(i_4)}+
{\bf 1}_{\{i_1=i_3\}}
{\bf 1}_{\{j_1=j_3\}}
{\bf 1}_{\{i_4=i_5\}}
{\bf 1}_{\{j_4=j_5\}}\zeta_{j_2}^{(i_2)}+
$$
$$
+
{\bf 1}_{\{i_1=i_4\}}
{\bf 1}_{\{j_1=j_4\}}
{\bf 1}_{\{i_2=i_3\}}
{\bf 1}_{\{j_2=j_3\}}\zeta_{j_5}^{(i_5)}+
{\bf 1}_{\{i_1=i_4\}}
{\bf 1}_{\{j_1=j_4\}}
{\bf 1}_{\{i_2=i_5\}}
{\bf 1}_{\{j_2=j_5\}}\zeta_{j_3}^{(i_3)}+
$$
$$
+
{\bf 1}_{\{i_1=i_4\}}
{\bf 1}_{\{j_1=j_4\}}
{\bf 1}_{\{i_3=i_5\}}
{\bf 1}_{\{j_3=j_5\}}\zeta_{j_2}^{(i_2)}+
{\bf 1}_{\{i_1=i_5\}}
{\bf 1}_{\{j_1=j_5\}}
{\bf 1}_{\{i_2=i_3\}}
{\bf 1}_{\{j_2=j_3\}}\zeta_{j_4}^{(i_4)}+
$$
$$
+
{\bf 1}_{\{i_1=i_5\}}
{\bf 1}_{\{j_1=j_5\}}
{\bf 1}_{\{i_2=i_4\}}
{\bf 1}_{\{j_2=j_4\}}\zeta_{j_3}^{(i_3)}+
{\bf 1}_{\{i_1=i_5\}}
{\bf 1}_{\{j_1=j_5\}}
{\bf 1}_{\{i_3=i_4\}}
{\bf 1}_{\{j_3=j_4\}}\zeta_{j_2}^{(i_2)}+
$$
$$
+
{\bf 1}_{\{i_2=i_3\}}
{\bf 1}_{\{j_2=j_3\}}
{\bf 1}_{\{i_4=i_5\}}
{\bf 1}_{\{j_4=j_5\}}\zeta_{j_1}^{(i_1)}+
{\bf 1}_{\{i_2=i_4\}}
{\bf 1}_{\{j_2=j_4\}}
{\bf 1}_{\{i_3=i_5\}}
{\bf 1}_{\{j_3=j_5\}}\zeta_{j_1}^{(i_1)}+
$$
\begin{equation}
\label{sss4}
+\Biggl.
{\bf 1}_{\{i_2=i_5\}}
{\bf 1}_{\{j_2=j_5\}}
{\bf 1}_{\{i_3=i_4\}}
{\bf 1}_{\{j_3=j_4\}}\zeta_{j_1}^{(i_1)}\Biggr),
\end{equation}

\vspace{8mm}

$$         
I_{(00000)T,t}^{(i_1i_1i_1i_1i_1)}=
\frac{1}{120}(T-t)^{5/2}
\left(\left(\zeta_0^{(i_1)}\right)^5-
10\left(\zeta_0^{(i_1)}\right)^3+15\zeta_0^{(i_1)}\right)\ \ \ 
\hbox{w.\ p.\ 1},
$$

\vspace{2mm}

$$
I_{(00000)T,t}^{*(i_1i_1i_1i_1i_1)}=
\frac{1}{120}(T-t)^{5/2}\left(\zeta_0^{(i_1)}\right)^5\ \ \ \hbox{w.\ p.\ 1},
$$

\vspace{5mm}
\noindent
where

\vspace{4mm}
$$
C_{j_5j_4 j_3 j_2 j_1}=
\frac{\sqrt{(2j_1+1)(2j_2+1)(2j_3+1)(2j_4+1)(2j_5+1)}}{32}(T-t)^{5/2}\bar
C_{j_5j_4 j_3 j_2 j_1},
$$

\vspace{4mm}

$$
\bar C_{j_5j_4 j_3 j_2 j_1}=
\int\limits_{-1}^{1}P_{j_5}(v)
\int\limits_{-1}^{v}P_{j_4}(u)
\int\limits_{-1}^{u}P_{j_3}(z)
\int\limits_{-1}^{z}P_{j_2}(y)
\int\limits_{-1}^{y}
P_{j_1}(x)dx dy dz du dv;
$$

\vspace{9mm}

$$
I_{(0001)T,t}^{*(i_1i_2i_3)}
=\hbox{\vtop{\offinterlineskip\halign{
\hfil#\hfil\cr
{\rm l.i.m.}\cr
$\stackrel{}{{}_{p\to \infty}}$\cr
}} }
\sum_{j_1,j_2,j_3,j_4=0}^{p}
C_{j_4j_3 j_2 j_1}^{0001}
\zeta_{j_1}^{(i_1)}\zeta_{j_2}^{(i_2)}\zeta_{j_3}^{(i_3)}\zeta_{j_4}^{(i_4)},
$$

\vspace{4mm}

$$
I_{(0010)T,t}^{*(i_1i_2i_3)}
=\hbox{\vtop{\offinterlineskip\halign{
\hfil#\hfil\cr
{\rm l.i.m.}\cr
$\stackrel{}{{}_{p\to \infty}}$\cr
}} }
\sum_{j_1,j_2,j_3,j_4=0}^{p}
C_{j_4j_3 j_2 j_1}^{0010}
\zeta_{j_1}^{(i_1)}\zeta_{j_2}^{(i_2)}\zeta_{j_3}^{(i_3)}\zeta_{j_4}^{(i_4)},
$$

\vspace{4mm}

$$
I_{(0100)T,t}^{*(i_1i_2i_3)}
=\hbox{\vtop{\offinterlineskip\halign{
\hfil#\hfil\cr
{\rm l.i.m.}\cr
$\stackrel{}{{}_{p\to \infty}}$\cr
}} }
\sum_{j_1,j_2,j_3,j_4=0}^{p}
C_{j_4j_3 j_2 j_1}^{0100}
\zeta_{j_1}^{(i_1)}\zeta_{j_2}^{(i_2)}\zeta_{j_3}^{(i_3)}\zeta_{j_4}^{(i_4)},
$$

\vspace{4mm}

$$
I_{(1000)T,t}^{*(i_1i_2i_3)}
=\hbox{\vtop{\offinterlineskip\halign{
\hfil#\hfil\cr
{\rm l.i.m.}\cr
$\stackrel{}{{}_{p\to \infty}}$\cr
}} }
\sum_{j_1,j_2,j_3,j_4=0}^{p}
C_{j_4j_3 j_2 j_1}^{1000}
\zeta_{j_1}^{(i_1)}\zeta_{j_2}^{(i_2)}\zeta_{j_3}^{(i_3)}\zeta_{j_4}^{(i_4)},
$$

\vspace{9mm}

$$
I_{(0001)T,t}^{(i_1 i_2 i_3 i_4)}
=\hbox{\vtop{\offinterlineskip\halign{
\hfil#\hfil\cr
{\rm l.i.m.}\cr
$\stackrel{}{{}_{p\to \infty}}$\cr
}} }
\sum_{j_1,j_2,j_3,j_4=0}^{p}
C_{j_4 j_3 j_2 j_1}^{0001}\Biggl(
\zeta_{j_1}^{(i_1)}\zeta_{j_2}^{(i_2)}\zeta_{j_3}^{(i_3)}\zeta_{j_4}^{(i_4)}
-\Biggr.
$$
$$
-
{\bf 1}_{\{i_1=i_2\}}
{\bf 1}_{\{j_1=j_2\}}
\zeta_{j_3}^{(i_3)}
\zeta_{j_4}^{(i_4)}
-
{\bf 1}_{\{i_1=i_3\}}
{\bf 1}_{\{j_1=j_3\}}
\zeta_{j_2}^{(i_2)}
\zeta_{j_4}^{(i_4)}-
$$
$$
-
{\bf 1}_{\{i_1=i_4\}}
{\bf 1}_{\{j_1=j_4\}}
\zeta_{j_2}^{(i_2)}
\zeta_{j_3}^{(i_3)}
-
{\bf 1}_{\{i_2=i_3\}}
{\bf 1}_{\{j_2=j_3\}}
\zeta_{j_1}^{(i_1)}
\zeta_{j_4}^{(i_4)}-
$$
$$
-
{\bf 1}_{\{i_2=i_4\}}
{\bf 1}_{\{j_2=j_4\}}
\zeta_{j_1}^{(i_1)}
\zeta_{j_3}^{(i_3)}
-
{\bf 1}_{\{i_3=i_4\}}
{\bf 1}_{\{j_3=j_4\}}
\zeta_{j_1}^{(i_1)}
\zeta_{j_2}^{(i_2)}+
$$

\vspace{-3mm}
$$
+
{\bf 1}_{\{i_1=i_2\}}
{\bf 1}_{\{j_1=j_2\}}
{\bf 1}_{\{i_3=i_4\}}
{\bf 1}_{\{j_3=j_4\}}+
$$

\vspace{-3mm}
$$
+
{\bf 1}_{\{i_1=i_3\}}
{\bf 1}_{\{j_1=j_3\}}
{\bf 1}_{\{i_2=i_4\}}
{\bf 1}_{\{j_2=j_4\}}+
$$
$$
+\Biggl.
{\bf 1}_{\{i_1=i_4\}}
{\bf 1}_{\{j_1=j_4\}}
{\bf 1}_{\{i_2=i_3\}}
{\bf 1}_{\{j_2=j_3\}}\Biggr),
$$

\vspace{7mm}

$$
I_{(0010)T,t}^{(i_1 i_2 i_3 i_4)}
=\hbox{\vtop{\offinterlineskip\halign{
\hfil#\hfil\cr
{\rm l.i.m.}\cr
$\stackrel{}{{}_{p\to \infty}}$\cr
}} }
\sum_{j_1,j_2,j_3,j_4=0}^{p}
C_{j_4 j_3 j_2 j_1}^{0010}\Biggl(
\zeta_{j_1}^{(i_1)}\zeta_{j_2}^{(i_2)}\zeta_{j_3}^{(i_3)}\zeta_{j_4}^{(i_4)}
-\Biggr.
$$
$$
-
{\bf 1}_{\{i_1=i_2\}}
{\bf 1}_{\{j_1=j_2\}}
\zeta_{j_3}^{(i_3)}
\zeta_{j_4}^{(i_4)}
-
{\bf 1}_{\{i_1=i_3\}}
{\bf 1}_{\{j_1=j_3\}}
\zeta_{j_2}^{(i_2)}
\zeta_{j_4}^{(i_4)}-
$$
$$
-
{\bf 1}_{\{i_1=i_4\}}
{\bf 1}_{\{j_1=j_4\}}
\zeta_{j_2}^{(i_2)}
\zeta_{j_3}^{(i_3)}
-
{\bf 1}_{\{i_2=i_3\}}
{\bf 1}_{\{j_2=j_3\}}
\zeta_{j_1}^{(i_1)}
\zeta_{j_4}^{(i_4)}-
$$
$$
-
{\bf 1}_{\{i_2=i_4\}}
{\bf 1}_{\{j_2=j_4\}}
\zeta_{j_1}^{(i_1)}
\zeta_{j_3}^{(i_3)}
-
{\bf 1}_{\{i_3=i_4\}}
{\bf 1}_{\{j_3=j_4\}}
\zeta_{j_1}^{(i_1)}
\zeta_{j_2}^{(i_2)}+
$$

\vspace{-3mm}
$$
+
{\bf 1}_{\{i_1=i_2\}}
{\bf 1}_{\{j_1=j_2\}}
{\bf 1}_{\{i_3=i_4\}}
{\bf 1}_{\{j_3=j_4\}}+
$$

\vspace{-3mm}
$$
+
{\bf 1}_{\{i_1=i_3\}}
{\bf 1}_{\{j_1=j_3\}}
{\bf 1}_{\{i_2=i_4\}}
{\bf 1}_{\{j_2=j_4\}}+
$$
$$
+\Biggl.
{\bf 1}_{\{i_1=i_4\}}
{\bf 1}_{\{j_1=j_4\}}
{\bf 1}_{\{i_2=i_3\}}
{\bf 1}_{\{j_2=j_3\}}\Biggr),
$$

\vspace{7mm}

$$
I_{(0100)T,t}^{(i_1 i_2 i_3 i_4)}
=\hbox{\vtop{\offinterlineskip\halign{
\hfil#\hfil\cr
{\rm l.i.m.}\cr
$\stackrel{}{{}_{p\to \infty}}$\cr
}} }
\sum_{j_1,j_2,j_3,j_4=0}^{p}
C_{j_4 j_3 j_2 j_1}^{0100}\Biggl(
\zeta_{j_1}^{(i_1)}\zeta_{j_2}^{(i_2)}\zeta_{j_3}^{(i_3)}\zeta_{j_4}^{(i_4)}
-\Biggr.
$$
$$
-
{\bf 1}_{\{i_1=i_2\}}
{\bf 1}_{\{j_1=j_2\}}
\zeta_{j_3}^{(i_3)}
\zeta_{j_4}^{(i_4)}
-
{\bf 1}_{\{i_1=i_3\}}
{\bf 1}_{\{j_1=j_3\}}
\zeta_{j_2}^{(i_2)}
\zeta_{j_4}^{(i_4)}-
$$
$$
-
{\bf 1}_{\{i_1=i_4\}}
{\bf 1}_{\{j_1=j_4\}}
\zeta_{j_2}^{(i_2)}
\zeta_{j_3}^{(i_3)}
-
{\bf 1}_{\{i_2=i_3\}}
{\bf 1}_{\{j_2=j_3\}}
\zeta_{j_1}^{(i_1)}
\zeta_{j_4}^{(i_4)}-
$$
$$
-
{\bf 1}_{\{i_2=i_4\}}
{\bf 1}_{\{j_2=j_4\}}
\zeta_{j_1}^{(i_1)}
\zeta_{j_3}^{(i_3)}
-
{\bf 1}_{\{i_3=i_4\}}
{\bf 1}_{\{j_3=j_4\}}
\zeta_{j_1}^{(i_1)}
\zeta_{j_2}^{(i_2)}+
$$

\vspace{-3mm}
$$
+
{\bf 1}_{\{i_1=i_2\}}
{\bf 1}_{\{j_1=j_2\}}
{\bf 1}_{\{i_3=i_4\}}
{\bf 1}_{\{j_3=j_4\}}+
$$

\vspace{-3mm}
$$
+
{\bf 1}_{\{i_1=i_3\}}
{\bf 1}_{\{j_1=j_3\}}
{\bf 1}_{\{i_2=i_4\}}
{\bf 1}_{\{j_2=j_4\}}+
$$
$$
+\Biggl.
{\bf 1}_{\{i_1=i_4\}}
{\bf 1}_{\{j_1=j_4\}}
{\bf 1}_{\{i_2=i_3\}}
{\bf 1}_{\{j_2=j_3\}}\Biggr),
$$

\vspace{7mm}

$$
I_{(1000)T,t}^{(i_1 i_2 i_3 i_4)}
=\hbox{\vtop{\offinterlineskip\halign{
\hfil#\hfil\cr
{\rm l.i.m.}\cr
$\stackrel{}{{}_{p\to \infty}}$\cr
}} }
\sum_{j_1,j_2,j_3,j_4=0}^{p}
C_{j_4 j_3 j_2 j_1}^{1000}\Biggl(
\zeta_{j_1}^{(i_1)}\zeta_{j_2}^{(i_2)}\zeta_{j_3}^{(i_3)}\zeta_{j_4}^{(i_4)}
-\Biggr.
$$
$$
-
{\bf 1}_{\{i_1=i_2\}}
{\bf 1}_{\{j_1=j_2\}}
\zeta_{j_3}^{(i_3)}
\zeta_{j_4}^{(i_4)}
-
{\bf 1}_{\{i_1=i_3\}}
{\bf 1}_{\{j_1=j_3\}}
\zeta_{j_2}^{(i_2)}
\zeta_{j_4}^{(i_4)}-
$$
$$
-
{\bf 1}_{\{i_1=i_4\}}
{\bf 1}_{\{j_1=j_4\}}
\zeta_{j_2}^{(i_2)}
\zeta_{j_3}^{(i_3)}
-
{\bf 1}_{\{i_2=i_3\}}
{\bf 1}_{\{j_2=j_3\}}
\zeta_{j_1}^{(i_1)}
\zeta_{j_4}^{(i_4)}-
$$
$$
-
{\bf 1}_{\{i_2=i_4\}}
{\bf 1}_{\{j_2=j_4\}}
\zeta_{j_1}^{(i_1)}
\zeta_{j_3}^{(i_3)}
-
{\bf 1}_{\{i_3=i_4\}}
{\bf 1}_{\{j_3=j_4\}}
\zeta_{j_1}^{(i_1)}
\zeta_{j_2}^{(i_2)}+
$$

\vspace{-3mm}
$$
+
{\bf 1}_{\{i_1=i_2\}}
{\bf 1}_{\{j_1=j_2\}}
{\bf 1}_{\{i_3=i_4\}}
{\bf 1}_{\{j_3=j_4\}}+
$$

\vspace{-3mm}
$$
+{\bf 1}_{\{i_1=i_3\}}
{\bf 1}_{\{j_1=j_3\}}
{\bf 1}_{\{i_2=i_4\}}
{\bf 1}_{\{j_2=j_4\}}+
$$
$$
+\Biggl.
{\bf 1}_{\{i_1=i_4\}}
{\bf 1}_{\{j_1=j_4\}}
{\bf 1}_{\{i_2=i_3\}}
{\bf 1}_{\{j_2=j_3\}}\Biggr),
$$

\vspace{7mm}
\noindent
where

\vspace{4mm}

$$
C_{j_4j_3j_2j_1}^{0001}
=\frac{\sqrt{(2j_1+1)(2j_2+1)(2j_3+1)(2j_4+1)}}{32}(T-t)^{3}\bar
C_{j_4j_3j_2j_1}^{0001},
$$

\vspace{4mm}

$$
C_{j_3j_2j_1}^{0010}
=\frac{\sqrt{(2j_1+1)(2j_2+1)(2j_3+1)(2j_4+1)}}{32}(T-t)^{3}\bar
C_{j_4j_3j_2j_1}^{0010},
$$

\vspace{4mm}

$$
C_{j_4j_3j_2j_1}^{0100}=
\frac{\sqrt{(2j_1+1)(2j_2+1)(2j_3+1)(2j_4+1)}}{32}(T-t)^{3}\bar
C_{j_3j_2j_1}^{0100},
$$

\vspace{4mm}

$$
C_{j_4j_3j_2j_1}^{1000}
=\frac{\sqrt{(2j_1+1)(2j_2+1)(2j_3+1)(2j_4+1)}}{32}(T-t)^{3}\bar
C_{j_4j_3j_2j_1}^{1000},
$$

\vspace{6mm}

$$
\bar C_{j_4j_3j_2j_1}^{1000}=-
\int\limits_{-1}^{1}P_{j_4}(u)
\int\limits_{-1}^{u}P_{j_3}(z)
\int\limits_{-1}^{z}P_{j_2}(y)
\int\limits_{-1}^{y}
P_{j_1}(x)(x+1)dx dy dz du,
$$

\vspace{4mm}

$$
\bar C_{j_4j_3j_2j_1}^{0100}=-
\int\limits_{-1}^{1}P_{j_4}(u)
\int\limits_{-1}^{u}P_{j_3}(z)
\int\limits_{-1}^{z}P_{j_2}(y)(y+1)
\int\limits_{-1}^{y}
P_{j_1}(x)dx dy dz du,
$$

\vspace{4mm}

$$
\bar C_{j_4j_3j_2j_1}^{0010}=-
\int\limits_{-1}^{1}P_{j_4}(u)
\int\limits_{-1}^{u}P_{j_3}(z)(z+1)
\int\limits_{-1}^{z}P_{j_2}(y)
\int\limits_{-1}^{y}
P_{j_1}(x)dx dy dz du,
$$

\vspace{4mm}

$$
\bar C_{j_4j_3j_2j_1}^{0001}=-
\int\limits_{-1}^{1}P_{j_4}(u)(u+1)
\int\limits_{-1}^{u}P_{j_3}(z)
\int\limits_{-1}^{z}P_{j_2}(y)
\int\limits_{-1}^{y}
P_{j_1}(x)dx dy dz du;
$$

\vspace{9mm}

$$
I_{(000000)T,t}^{*(i_1 i_2 i_3 i_4 i_5 i_6)}=
\hbox{\vtop{\offinterlineskip\halign{
\hfil#\hfil\cr
{\rm l.i.m.}\cr
$\stackrel{}{{}_{p\to \infty}}$\cr
}} }
\sum\limits_{j_1, j_2, j_3, j_4, j_5, j_6=0}^{p}
C_{j_6j_5j_4 j_3 j_2 j_1}
\zeta_{j_1}^{(i_1)}\zeta_{j_2}^{(i_2)}\zeta_{j_3}^{(i_3)}
\zeta_{j_4}^{(i_4)}\zeta_{j_5}^{(i_5)}\zeta_{j_6}^{(i_6)},
$$

\vspace{9mm}

$$
I_{(000000)T,t}^{(i_1 i_2 i_3 i_4 i_5 i_6)}
=\hbox{\vtop{\offinterlineskip\halign{
\hfil#\hfil\cr
{\rm l.i.m.}\cr
$\stackrel{}{{}_{p\to \infty}}$\cr
}} }\sum_{j_1,j_2,j_3,j_4,j_5,j_6=0}^{p}
C_{j_6 j_5 j_4 j_3 j_2 j_1}\Biggl(
\prod_{l=1}^6
\zeta_{j_l}^{(i_l)}
-\Biggr.
$$
$$
-
{\bf 1}_{\{j_1=j_6\}}
{\bf 1}_{\{i_1=i_6\}}
\zeta_{j_2}^{(i_2)}
\zeta_{j_3}^{(i_3)}
\zeta_{j_4}^{(i_4)}
\zeta_{j_5}^{(i_5)}-
{\bf 1}_{\{j_2=j_6\}}
{\bf 1}_{\{i_2=i_6\}}
\zeta_{j_1}^{(i_1)}
\zeta_{j_3}^{(i_3)}
\zeta_{j_4}^{(i_4)}
\zeta_{j_5}^{(i_5)}-
$$
$$
-
{\bf 1}_{\{j_3=j_6\}}
{\bf 1}_{\{i_3=i_6\}}
\zeta_{j_1}^{(i_1)}
\zeta_{j_2}^{(i_2)}
\zeta_{j_4}^{(i_4)}
\zeta_{j_5}^{(i_5)}-
{\bf 1}_{\{j_4=j_6\}}
{\bf 1}_{\{i_4=i_6\}}
\zeta_{j_1}^{(i_1)}
\zeta_{j_2}^{(i_2)}
\zeta_{j_3}^{(i_3)}
\zeta_{j_5}^{(i_5)}-
$$
$$
-
{\bf 1}_{\{j_5=j_6\}}
{\bf 1}_{\{i_5=i_6\}}
\zeta_{j_1}^{(i_1)}
\zeta_{j_2}^{(i_2)}
\zeta_{j_3}^{(i_3)}
\zeta_{j_4}^{(i_4)}-
{\bf 1}_{\{j_1=j_2\}}
{\bf 1}_{\{i_1=i_2\}}
\zeta_{j_3}^{(i_3)}
\zeta_{j_4}^{(i_4)}
\zeta_{j_5}^{(i_5)}
\zeta_{j_6}^{(i_6)}-
$$
$$
-
{\bf 1}_{\{j_1=j_3\}}
{\bf 1}_{\{i_1=i_3\}}
\zeta_{j_2}^{(i_2)}
\zeta_{j_4}^{(i_4)}
\zeta_{j_5}^{(i_5)}
\zeta_{j_6}^{(i_6)}-
{\bf 1}_{\{j_1=j_4\}}
{\bf 1}_{\{i_1=i_4\}}
\zeta_{j_2}^{(i_2)}
\zeta_{j_3}^{(i_3)}
\zeta_{j_5}^{(i_5)}
\zeta_{j_6}^{(i_6)}-
$$
$$
-
{\bf 1}_{\{j_1=j_5\}}
{\bf 1}_{\{i_1=i_5\}}
\zeta_{j_2}^{(i_2)}
\zeta_{j_3}^{(i_3)}
\zeta_{j_4}^{(i_4)}
\zeta_{j_6}^{(i_6)}-
{\bf 1}_{\{j_2=j_3\}}
{\bf 1}_{\{i_2=i_3\}}
\zeta_{j_1}^{(i_1)}
\zeta_{j_4}^{(i_4)}
\zeta_{j_5}^{(i_5)}
\zeta_{j_6}^{(i_6)}-
$$
$$
-
{\bf 1}_{\{j_2=j_4\}}
{\bf 1}_{\{i_2=i_4\}}
\zeta_{j_1}^{(i_1)}
\zeta_{j_3}^{(i_3)}
\zeta_{j_5}^{(i_5)}
\zeta_{j_6}^{(i_6)}-
{\bf 1}_{\{j_2=j_5\}}
{\bf 1}_{\{i_2=i_5\}}
\zeta_{j_1}^{(i_1)}
\zeta_{j_3}^{(i_3)}
\zeta_{j_4}^{(i_4)}
\zeta_{j_6}^{(i_6)}-
$$
$$
-
{\bf 1}_{\{j_3=j_4\}}
{\bf 1}_{\{i_3=i_4\}}
\zeta_{j_1}^{(i_1)}
\zeta_{j_2}^{(i_2)}
\zeta_{j_5}^{(i_5)}
\zeta_{j_6}^{(i_6)}-
{\bf 1}_{\{j_3=j_5\}}
{\bf 1}_{\{i_3=i_5\}}
\zeta_{j_1}^{(i_1)}
\zeta_{j_2}^{(i_2)}
\zeta_{j_4}^{(i_4)}
\zeta_{j_6}^{(i_6)}-
$$
$$
-
{\bf 1}_{\{j_4=j_5\}}
{\bf 1}_{\{i_4=i_5\}}
\zeta_{j_1}^{(i_1)}
\zeta_{j_2}^{(i_2)}
\zeta_{j_3}^{(i_3)}
\zeta_{j_6}^{(i_6)}+
$$
$$
+
{\bf 1}_{\{j_1=j_2\}}
{\bf 1}_{\{i_1=i_2\}}
{\bf 1}_{\{j_3=j_4\}}
{\bf 1}_{\{i_3=i_4\}}
\zeta_{j_5}^{(i_5)}
\zeta_{j_6}^{(i_6)}
+
{\bf 1}_{\{j_1=j_2\}}
{\bf 1}_{\{i_1=i_2\}}
{\bf 1}_{\{j_3=j_5\}}
{\bf 1}_{\{i_3=i_5\}}
\zeta_{j_4}^{(i_4)}
\zeta_{j_6}^{(i_6)}+
$$
$$
+
{\bf 1}_{\{j_1=j_2\}}
{\bf 1}_{\{i_1=i_2\}}
{\bf 1}_{\{j_4=j_5\}}
{\bf 1}_{\{i_4=i_5\}}
\zeta_{j_3}^{(i_3)}
\zeta_{j_6}^{(i_6)}
+
{\bf 1}_{\{j_1=j_3\}}
{\bf 1}_{\{i_1=i_3\}}
{\bf 1}_{\{j_2=j_4\}}
{\bf 1}_{\{i_2=i_4\}}
\zeta_{j_5}^{(i_5)}
\zeta_{j_6}^{(i_6)}+
$$
$$
+
{\bf 1}_{\{j_1=j_3\}}
{\bf 1}_{\{i_1=i_3\}}
{\bf 1}_{\{j_2=j_5\}}
{\bf 1}_{\{i_2=i_5\}}
\zeta_{j_4}^{(i_4)}
\zeta_{j_6}^{(i_6)}
+
{\bf 1}_{\{j_1=j_3\}}
{\bf 1}_{\{i_1=i_3\}}
{\bf 1}_{\{j_4=j_5\}}
{\bf 1}_{\{i_4=i_5\}}
\zeta_{j_2}^{(i_2)}
\zeta_{j_6}^{(i_6)}+
$$
$$
+
{\bf 1}_{\{j_1=j_4\}}
{\bf 1}_{\{i_1=i_4\}}
{\bf 1}_{\{j_2=j_3\}}
{\bf 1}_{\{i_2=i_3\}}
\zeta_{j_5}^{(i_5)}
\zeta_{j_6}^{(i_6)}
+
{\bf 1}_{\{j_1=j_4\}}
{\bf 1}_{\{i_1=i_4\}}
{\bf 1}_{\{j_2=j_5\}}
{\bf 1}_{\{i_2=i_5\}}
\zeta_{j_3}^{(i_3)}
\zeta_{j_6}^{(i_6)}+
$$
$$
+
{\bf 1}_{\{j_1=j_4\}}
{\bf 1}_{\{i_1=i_4\}}
{\bf 1}_{\{j_3=j_5\}}
{\bf 1}_{\{i_3=i_5\}}
\zeta_{j_2}^{(i_2)}
\zeta_{j_6}^{(i_6)}
+
{\bf 1}_{\{j_1=j_5\}}
{\bf 1}_{\{i_1=i_5\}}
{\bf 1}_{\{j_2=j_3\}}
{\bf 1}_{\{i_2=i_3\}}
\zeta_{j_4}^{(i_4)}
\zeta_{j_6}^{(i_6)}+
$$
$$
+
{\bf 1}_{\{j_1=j_5\}}
{\bf 1}_{\{i_1=i_5\}}
{\bf 1}_{\{j_2=j_4\}}
{\bf 1}_{\{i_2=i_4\}}
\zeta_{j_3}^{(i_3)}
\zeta_{j_6}^{(i_6)}
+
{\bf 1}_{\{j_1=j_5\}}
{\bf 1}_{\{i_1=i_5\}}
{\bf 1}_{\{j_3=j_4\}}
{\bf 1}_{\{i_3=i_4\}}
\zeta_{j_2}^{(i_2)}
\zeta_{j_6}^{(i_6)}+
$$
$$
+
{\bf 1}_{\{j_2=j_3\}}
{\bf 1}_{\{i_2=i_3\}}
{\bf 1}_{\{j_4=j_5\}}
{\bf 1}_{\{i_4=i_5\}}
\zeta_{j_1}^{(i_1)}
\zeta_{j_6}^{(i_6)}
+
{\bf 1}_{\{j_2=j_4\}}
{\bf 1}_{\{i_2=i_4\}}
{\bf 1}_{\{j_3=j_5\}}
{\bf 1}_{\{i_3=i_5\}}
\zeta_{j_1}^{(i_1)}
\zeta_{j_6}^{(i_6)}+
$$
$$
+
{\bf 1}_{\{j_2=j_5\}}
{\bf 1}_{\{i_2=i_5\}}
{\bf 1}_{\{j_3=j_4\}}
{\bf 1}_{\{i_3=i_4\}}
\zeta_{j_1}^{(i_1)}
\zeta_{j_6}^{(i_6)}
+
{\bf 1}_{\{j_6=j_1\}}
{\bf 1}_{\{i_6=i_1\}}
{\bf 1}_{\{j_3=j_4\}}
{\bf 1}_{\{i_3=i_4\}}
\zeta_{j_2}^{(i_2)}
\zeta_{j_5}^{(i_5)}+
$$
$$
+
{\bf 1}_{\{j_6=j_1\}}
{\bf 1}_{\{i_6=i_1\}}
{\bf 1}_{\{j_3=j_5\}}
{\bf 1}_{\{i_3=i_5\}}
\zeta_{j_2}^{(i_2)}
\zeta_{j_4}^{(i_4)}
+
{\bf 1}_{\{j_6=j_1\}}
{\bf 1}_{\{i_6=i_1\}}
{\bf 1}_{\{j_2=j_5\}}
{\bf 1}_{\{i_2=i_5\}}
\zeta_{j_3}^{(i_3)}
\zeta_{j_4}^{(i_4)}+
$$
$$
+
{\bf 1}_{\{j_6=j_1\}}
{\bf 1}_{\{i_6=i_1\}}
{\bf 1}_{\{j_2=j_4\}}
{\bf 1}_{\{i_2=i_4\}}
\zeta_{j_3}^{(i_3)}
\zeta_{j_5}^{(i_5)}
+
{\bf 1}_{\{j_6=j_1\}}
{\bf 1}_{\{i_6=i_1\}}
{\bf 1}_{\{j_4=j_5\}}
{\bf 1}_{\{i_4=i_5\}}
\zeta_{j_2}^{(i_2)}
\zeta_{j_3}^{(i_3)}+
$$
$$
+
{\bf 1}_{\{j_6=j_1\}}
{\bf 1}_{\{i_6=i_1\}}
{\bf 1}_{\{j_2=j_3\}}
{\bf 1}_{\{i_2=i_3\}}
\zeta_{j_4}^{(i_4)}
\zeta_{j_5}^{(i_5)}
+
{\bf 1}_{\{j_6=j_2\}}
{\bf 1}_{\{i_6=i_2\}}
{\bf 1}_{\{j_3=j_5\}}
{\bf 1}_{\{i_3=i_5\}}
\zeta_{j_1}^{(i_1)}
\zeta_{j_4}^{(i_4)}+
$$
$$
+
{\bf 1}_{\{j_6=j_2\}}
{\bf 1}_{\{i_6=i_2\}}
{\bf 1}_{\{j_4=j_5\}}
{\bf 1}_{\{i_4=i_5\}}
\zeta_{j_1}^{(i_1)}
\zeta_{j_3}^{(i_3)}
+
{\bf 1}_{\{j_6=j_2\}}
{\bf 1}_{\{i_6=i_2\}}
{\bf 1}_{\{j_3=j_4\}}
{\bf 1}_{\{i_3=i_4\}}
\zeta_{j_1}^{(i_1)}
\zeta_{j_5}^{(i_5)}+
$$
$$
+
{\bf 1}_{\{j_6=j_2\}}
{\bf 1}_{\{i_6=i_2\}}
{\bf 1}_{\{j_1=j_5\}}
{\bf 1}_{\{i_1=i_5\}}
\zeta_{j_3}^{(i_3)}
\zeta_{j_4}^{(i_4)}
+
{\bf 1}_{\{j_6=j_2\}}
{\bf 1}_{\{i_6=i_2\}}
{\bf 1}_{\{j_1=j_4\}}
{\bf 1}_{\{i_1=i_4\}}
\zeta_{j_3}^{(i_3)}
\zeta_{j_5}^{(i_5)}+
$$
$$
+
{\bf 1}_{\{j_6=j_2\}}
{\bf 1}_{\{i_6=i_2\}}
{\bf 1}_{\{j_1=j_3\}}
{\bf 1}_{\{i_1=i_3\}}
\zeta_{j_4}^{(i_4)}
\zeta_{j_5}^{(i_5)}
+
{\bf 1}_{\{j_6=j_3\}}
{\bf 1}_{\{i_6=i_3\}}
{\bf 1}_{\{j_2=j_5\}}
{\bf 1}_{\{i_2=i_5\}}
\zeta_{j_1}^{(i_1)}
\zeta_{j_4}^{(i_4)}+
$$
$$
+
{\bf 1}_{\{j_6=j_3\}}
{\bf 1}_{\{i_6=i_3\}}
{\bf 1}_{\{j_4=j_5\}}
{\bf 1}_{\{i_4=i_5\}}
\zeta_{j_1}^{(i_1)}
\zeta_{j_2}^{(i_2)}
+
{\bf 1}_{\{j_6=j_3\}}
{\bf 1}_{\{i_6=i_3\}}
{\bf 1}_{\{j_2=j_4\}}
{\bf 1}_{\{i_2=i_4\}}
\zeta_{j_1}^{(i_1)}
\zeta_{j_5}^{(i_5)}+
$$
$$
+
{\bf 1}_{\{j_6=j_3\}}
{\bf 1}_{\{i_6=i_3\}}
{\bf 1}_{\{j_1=j_5\}}
{\bf 1}_{\{i_1=i_5\}}
\zeta_{j_2}^{(i_2)}
\zeta_{j_4}^{(i_4)}
+
{\bf 1}_{\{j_6=j_3\}}
{\bf 1}_{\{i_6=i_3\}}
{\bf 1}_{\{j_1=j_4\}}
{\bf 1}_{\{i_1=i_4\}}
\zeta_{j_2}^{(i_2)}
\zeta_{j_5}^{(i_5)}+
$$
$$
+
{\bf 1}_{\{j_6=j_3\}}
{\bf 1}_{\{i_6=i_3\}}
{\bf 1}_{\{j_1=j_2\}}
{\bf 1}_{\{i_1=i_2\}}
\zeta_{j_4}^{(i_4)}
\zeta_{j_5}^{(i_5)}
+
{\bf 1}_{\{j_6=j_4\}}
{\bf 1}_{\{i_6=i_4\}}
{\bf 1}_{\{j_3=j_5\}}
{\bf 1}_{\{i_3=i_5\}}
\zeta_{j_1}^{(i_1)}
\zeta_{j_2}^{(i_2)}+
$$
$$
+
{\bf 1}_{\{j_6=j_4\}}
{\bf 1}_{\{i_6=i_4\}}
{\bf 1}_{\{j_2=j_5\}}
{\bf 1}_{\{i_2=i_5\}}
\zeta_{j_1}^{(i_1)}
\zeta_{j_3}^{(i_3)}
+
{\bf 1}_{\{j_6=j_4\}}
{\bf 1}_{\{i_6=i_4\}}
{\bf 1}_{\{j_2=j_3\}}
{\bf 1}_{\{i_2=i_3\}}
\zeta_{j_1}^{(i_1)}
\zeta_{j_5}^{(i_5)}+
$$
$$
+
{\bf 1}_{\{j_6=j_4\}}
{\bf 1}_{\{i_6=i_4\}}
{\bf 1}_{\{j_1=j_5\}}
{\bf 1}_{\{i_1=i_5\}}
\zeta_{j_2}^{(i_2)}
\zeta_{j_3}^{(i_3)}
+
{\bf 1}_{\{j_6=j_4\}}
{\bf 1}_{\{i_6=i_4\}}
{\bf 1}_{\{j_1=j_3\}}
{\bf 1}_{\{i_1=i_3\}}
\zeta_{j_2}^{(i_2)}
\zeta_{j_5}^{(i_5)}+
$$
$$
+
{\bf 1}_{\{j_6=j_4\}}
{\bf 1}_{\{i_6=i_4\}}
{\bf 1}_{\{j_1=j_2\}}
{\bf 1}_{\{i_1=i_2\}}
\zeta_{j_3}^{(i_3)}
\zeta_{j_5}^{(i_5)}
+
{\bf 1}_{\{j_6=j_5\}}
{\bf 1}_{\{i_6=i_5\}}
{\bf 1}_{\{j_3=j_4\}}
{\bf 1}_{\{i_3=i_4\}}
\zeta_{j_1}^{(i_1)}
\zeta_{j_2}^{(i_2)}+
$$
$$
+
{\bf 1}_{\{j_6=j_5\}}
{\bf 1}_{\{i_6=i_5\}}
{\bf 1}_{\{j_2=j_4\}}
{\bf 1}_{\{i_2=i_4\}}
\zeta_{j_1}^{(i_1)}
\zeta_{j_3}^{(i_3)}
+
{\bf 1}_{\{j_6=j_5\}}
{\bf 1}_{\{i_6=i_5\}}
{\bf 1}_{\{j_2=j_3\}}
{\bf 1}_{\{i_2=i_3\}}
\zeta_{j_1}^{(i_1)}
\zeta_{j_4}^{(i_4)}+
$$
$$
+
{\bf 1}_{\{j_6=j_5\}}
{\bf 1}_{\{i_6=i_5\}}
{\bf 1}_{\{j_1=j_4\}}
{\bf 1}_{\{i_1=i_4\}}
\zeta_{j_2}^{(i_2)}
\zeta_{j_3}^{(i_3)}
+
{\bf 1}_{\{j_6=j_5\}}
{\bf 1}_{\{i_6=i_5\}}
{\bf 1}_{\{j_1=j_3\}}
{\bf 1}_{\{i_1=i_3\}}
\zeta_{j_2}^{(i_2)}
\zeta_{j_4}^{(i_4)}+
$$
$$
+
{\bf 1}_{\{j_6=j_5\}}
{\bf 1}_{\{i_6=i_5\}}
{\bf 1}_{\{j_1=j_2\}}
{\bf 1}_{\{i_1=i_2\}}
\zeta_{j_3}^{(i_3)}
\zeta_{j_4}^{(i_4)}-
$$
$$
-
{\bf 1}_{\{j_6=j_1\}}
{\bf 1}_{\{i_6=i_1\}}
{\bf 1}_{\{j_2=j_5\}}
{\bf 1}_{\{i_2=i_5\}}
{\bf 1}_{\{j_3=j_4\}}
{\bf 1}_{\{i_3=i_4\}}-
$$
$$
-
{\bf 1}_{\{j_6=j_1\}}
{\bf 1}_{\{i_6=i_1\}}
{\bf 1}_{\{j_2=j_4\}}
{\bf 1}_{\{i_2=i_4\}}
{\bf 1}_{\{j_3=j_5\}}
{\bf 1}_{\{i_3=i_5\}}-
$$
$$
-
{\bf 1}_{\{j_6=j_1\}}
{\bf 1}_{\{i_6=i_1\}}
{\bf 1}_{\{j_2=j_3\}}
{\bf 1}_{\{i_2=i_3\}}
{\bf 1}_{\{j_4=j_5\}}
{\bf 1}_{\{i_4=i_5\}}-
$$
$$
-               
{\bf 1}_{\{j_6=j_2\}}
{\bf 1}_{\{i_6=i_2\}}
{\bf 1}_{\{j_1=j_5\}}
{\bf 1}_{\{i_1=i_5\}}
{\bf 1}_{\{j_3=j_4\}}
{\bf 1}_{\{i_3=i_4\}}-
$$
$$
-
{\bf 1}_{\{j_6=j_2\}}
{\bf 1}_{\{i_6=i_2\}}
{\bf 1}_{\{j_1=j_4\}}
{\bf 1}_{\{i_1=i_4\}}
{\bf 1}_{\{j_3=j_5\}}
{\bf 1}_{\{i_3=i_5\}}-
$$
$$
-
{\bf 1}_{\{j_6=j_2\}}
{\bf 1}_{\{i_6=i_2\}}
{\bf 1}_{\{j_1=j_3\}}
{\bf 1}_{\{i_1=i_3\}}
{\bf 1}_{\{j_4=j_5\}}
{\bf 1}_{\{i_4=i_5\}}-
$$
$$
-
{\bf 1}_{\{j_6=j_3\}}
{\bf 1}_{\{i_6=i_3\}}
{\bf 1}_{\{j_1=j_5\}}
{\bf 1}_{\{i_1=i_5\}}
{\bf 1}_{\{j_2=j_4\}}
{\bf 1}_{\{i_2=i_4\}}-
$$
$$
-
{\bf 1}_{\{j_6=j_3\}}
{\bf 1}_{\{i_6=i_3\}}
{\bf 1}_{\{j_1=j_4\}}
{\bf 1}_{\{i_1=i_4\}}
{\bf 1}_{\{j_2=j_5\}}
{\bf 1}_{\{i_2=i_5\}}-
$$
$$
-
{\bf 1}_{\{j_3=j_6\}}
{\bf 1}_{\{i_3=i_6\}}
{\bf 1}_{\{j_1=j_2\}}
{\bf 1}_{\{i_1=i_2\}}
{\bf 1}_{\{j_4=j_5\}}
{\bf 1}_{\{i_4=i_5\}}-
$$
$$
-
{\bf 1}_{\{j_6=j_4\}}
{\bf 1}_{\{i_6=i_4\}}
{\bf 1}_{\{j_1=j_5\}}
{\bf 1}_{\{i_1=i_5\}}
{\bf 1}_{\{j_2=j_3\}}
{\bf 1}_{\{i_2=i_3\}}-
$$
$$
-
{\bf 1}_{\{j_6=j_4\}}
{\bf 1}_{\{i_6=i_4\}}
{\bf 1}_{\{j_1=j_3\}}
{\bf 1}_{\{i_1=i_3\}}
{\bf 1}_{\{j_2=j_5\}}
{\bf 1}_{\{i_2=i_5\}}-
$$
$$
-
{\bf 1}_{\{j_6=j_4\}}
{\bf 1}_{\{i_6=i_4\}}
{\bf 1}_{\{j_1=j_2\}}
{\bf 1}_{\{i_1=i_2\}}
{\bf 1}_{\{j_3=j_5\}}
{\bf 1}_{\{i_3=i_5\}}-
$$
$$
-
{\bf 1}_{\{j_6=j_5\}}
{\bf 1}_{\{i_6=i_5\}}
{\bf 1}_{\{j_1=j_4\}}
{\bf 1}_{\{i_1=i_4\}}
{\bf 1}_{\{j_2=j_3\}}
{\bf 1}_{\{i_2=i_3\}}-
$$
$$
-
{\bf 1}_{\{j_6=j_5\}}
{\bf 1}_{\{i_6=i_5\}}
{\bf 1}_{\{j_1=j_2\}}
{\bf 1}_{\{i_1=i_2\}}
{\bf 1}_{\{j_3=j_4\}}
{\bf 1}_{\{i_3=i_4\}}-
$$
$$
\Biggl.-
{\bf 1}_{\{j_6=j_5\}}
{\bf 1}_{\{i_6=i_5\}}
{\bf 1}_{\{j_1=j_3\}}
{\bf 1}_{\{i_1=i_3\}}
{\bf 1}_{\{j_2=j_4\}}
{\bf 1}_{\{i_2=i_4\}}\Biggr),
$$

\vspace{8mm}

$$         
I_{(000000)T,t}^{(i_1i_1i_1i_1i_1i_1)}=
\frac{1}{720}(T-t)^{3}
\left(\left(\zeta_0^{(i_1)}\right)^6-
15\left(\zeta_0^{(i_1)}\right)^4+45\left(\zeta_0^{(i_1)}\right)^2-
15\right)\ \ \ 
\hbox{w.\ p.\ 1},
$$

\vspace{4mm}

$$
I_{(000000)T,t}^{*(i_1i_1i_1i_1i_1i_1)}=
\frac{1}{720}(T-t)^{3}\left(\zeta_0^{(i_1)}\right)^6\ \ \ \hbox{w.\ p.\ 1},
$$

\vspace{5mm}
\noindent
where

\vspace{4mm}

$$
C_{j_6j_5j_4 j_3 j_2 j_1}
=\frac{\sqrt{(2j_1+1)(2j_2+1)(2j_3+1)
(2j_4+1)(2j_5+1)(2j_6+1)}}{64}(T-t)^{3}\bar
C_{j_6j_5j_4 j_3 j_2 j_1},
$$

\vspace{2mm}

$$
\bar C_{j_6j_5j_4 j_3 j_2 j_1}=
\int\limits_{-1}^{1}P_{j_6}(w)
\int\limits_{-1}^{w}P_{j_5}(v)
\int\limits_{-1}^{v}P_{j_4}(u)
\int\limits_{-1}^{u}P_{j_3}(z)
\int\limits_{-1}^{z}P_{j_2}(y)
\int\limits_{-1}^{y}
P_{j_1}(x)dx dy dz du dv dw.
$$

\vspace{10mm}

Let us analyze the approximation $I_{(00)T,t}^{*(i_1 i_2)q}$
of the iterated stochastic 
integral $I_{(00)T,t}^{*(i_1 i_2)}$ obtained from 
(\ref{4004}) by replacing $\infty$ on $q$.

It is easy to prove that

\vspace{3mm}
\begin{equation}
\label{fff09}
{\sf M}\biggl\{\left(I_{(00)T,t}^{*(i_1 i_2)}-
I_{(00)T,t}^{*(i_1 i_2)q}
\right)^2\biggr\}
=\frac{(T-t)^2}{2}\left(\frac{1}{2}-\sum_{i=1}^q
\frac{1}{4i^2-1}\right)\ \ \ (i_1\ne i_2).
\end{equation}

\vspace{6mm}

Moreover, using Theorem 17, we obtain for $i_1\ne i_2$

\vspace{3mm}
$$
{\sf M}\biggl\{\left(I_{(10)T,t}^{*(i_1 i_2)}-I_{(10)T,t}^{*(i_1 i_2)q}
\right)^2\biggr\}=
{\sf M}\biggl\{\left(I_{(01)T,t}^{*(i_1 i_2)}-
I_{(01)T,t}^{*(i_1 i_2)q}\right)^2\biggr\}=
$$

\vspace{3mm}
$$
=\frac{(T-t)^4}{16}\left(\frac{5}{9}-
2\sum_{i=2}^q\frac{1}{4i^2-1}-
\sum_{i=1}^q
\frac{1}{(2i-1)^2(2i+3)^2}
-\sum_{i=0}^q\frac{(i+2)^2+(i+1)^2}{(2i+1)(2i+5)(2i+3)^2}
\right).
$$

\vspace{8mm}

For the case $i_1=i_2$,
using Theorem 17, we have

\vspace{3mm}
$$
{\sf M}\biggl\{\left(I_{(10)T,t}^{(i_1 i_1)}-
I_{(10)T,t}^{(i_1 i_1)q}
\right)^2\biggr\}=
{\sf M}\biggl\{\left(I_{(01)T,t}^{(i_1 i_1)}-
I_{(01)T,t}^{(i_1 i_1)q}\right)^2\biggr\}=
$$

\vspace{3mm}
\begin{equation}
\label{2007ura1}
=\frac{(T-t)^4}{16}\left(\frac{1}{9}-
\sum_{i=0}^{q}
\frac{1}{(2i+1)(2i+5)(2i+3)^2}
-2\sum_{i=1}^{q}
\frac{1}{(2i-1)^2(2i+3)^2}\right).
\end{equation}

\vspace{6mm}

On the basis of 
the presented 
expansions of 
iterated stochastic integrals we 
can see that increasing of multiplicities of these integrals 
or degree indexes of their weight functions 
leads
to noticeable complication of formulas 
for the mentioned expansions. 

However, increasing of the mentioned parameters leads to increasing 
of orders of smallness with respect to $T-t$ in the mean-square sense 
for iterated stochastic integrals. This leads to a sharp decrease  
of member 
quantities
in expansions of iterated stochastic 
integrals, which are required for achieving the acceptable accuracy
of approximation. In the context of it let us consider the approach 
to approximation of iterated stochastic integrals, which 
provides a possibility to obtain the mean-square approximations of 
the required accuracy without using the 
complex expansions.

Let us analyze the following approximation of triple stochastic integral 
using (\ref{zzz1})

\vspace{3mm}
$$
I_{(000)T,t}^{(i_1i_2i_3)q_1}
=\sum_{j_1,j_2,j_3=0}^{q_1}
C_{j_3j_2j_1}\Biggl(
\zeta_{j_1}^{(i_1)}\zeta_{j_2}^{(i_2)}\zeta_{j_3}^{(i_3)}
-{\bf 1}_{\{i_1=i_2\}}
{\bf 1}_{\{j_1=j_2\}}
\zeta_{j_3}^{(i_3)}-
\Biggr.
$$

\vspace{1mm}
\begin{equation}
\label{sad001}
\Biggl.
-{\bf 1}_{\{i_2=i_3\}}
{\bf 1}_{\{j_2=j_3\}}
\zeta_{j_1}^{(i_1)}-
{\bf 1}_{\{i_1=i_3\}}
{\bf 1}_{\{j_1=j_3\}}
\zeta_{j_2}^{(i_2)}\Biggr),
\end{equation}

\vspace{5mm}
\noindent
where $C_{j_3j_2j_1}$ is defined by (\ref{zzz2}), (\ref{zzz3}).

In particular, from
(\ref{sad001}) when $i_1\ne i_2$, 
$i_2\ne i_3$, $i_1\ne i_3$
we obtain

\vspace{1mm}
\begin{equation}
\label{38}
I_{(000)T,t}^{(i_1i_2i_3)q_1}=
\sum_{j_1,j_2,j_3=0}^{q_1}
C_{j_3j_2j_1}
\zeta_{j_1}^{(i_1)}\zeta_{j_2}^{(i_2)}\zeta_{j_3}^{(i_3)}.
\end{equation}

\vspace{5mm}

Using (\ref{zzz4}), (\ref{zzz5})--(\ref{zzz0}), we get

\vspace{2mm}
$$
{\sf M}\left\{\left(
I_{(000)T,t}^{(i_1i_2 i_3)}-
I_{(000)T,t}^{(i_1i_2 i_3)q_1}\right)^2\right\}=
$$

\vspace{2mm}
\begin{equation}
\label{39}
=
\frac{(T-t)^{3}}{6}
-\sum_{j_1,j_2,j_3=0}^{q_1}
C_{j_3j_2j_1}^2\ \ \ (i_1\ne i_2, i_1\ne i_3, i_2\ne i_3),
\end{equation}

\vspace{8mm}

$$
{\sf M}\left\{\left(
I_{(000)T,t}^{(i_1i_2 i_3)}-
I_{(000)T,t}^{(i_1i_2 i_3)q_1}\right)^2\right\}=
$$

\vspace{2mm}
\begin{equation}
\label{39a}
=
\frac{(T-t)^{3}}{6}-\sum_{j_1,j_2,j_3=0}^{q_1}
C_{j_3j_2j_1}^2
-\sum_{j_1,j_2,j_3=0}^{q_1}
C_{j_2j_3j_1}C_{j_3j_2j_1}\ \ \ (i_1\ne i_2=i_3),
\end{equation}

\vspace{8mm}

$$
{\sf M}\left\{\left(
I_{(000)T,t}^{(i_1i_2 i_3)}-
I_{(000)T,t}^{(i_1i_2 i_3)q_1}\right)^2\right\}=
$$

\vspace{2mm}
\begin{equation}
\label{39b}
=
\frac{(T-t)^{3}}{6}-\sum_{j_1,j_2,j_3=0}^{q_1}
C_{j_3j_2j_1}^2
-\sum_{j_1,j_2,j_3=0}^{q_1}
C_{j_3j_2j_1}C_{j_1j_2j_3}\ \ \ (i_1=i_3\ne i_2),
\end{equation}

\vspace{8mm}

$$
{\sf M}\left\{\left(
I_{(000)T,t}^{(i_1i_2 i_3)}-
I_{(000)T,t}^{(i_1i_2 i_3)q_1}\right)^2\right\}=
$$

\vspace{2mm}
\begin{equation}
\label{39c}
=
\frac{(T-t)^{3}}{6}-\sum_{j_1,j_2,j_3=0}^{q_1}
C_{j_3j_2j_1}^2
-\sum_{j_1,j_2,j_3=0}^{q_1}
C_{j_3j_1j_2}C_{j_3j_2j_1}\ \ \ (i_1=i_2\ne i_3),
\end{equation}

\vspace{8mm}

\begin{equation}
\label{leto1041}
{\sf M}\left\{\left(
I_{(000)T,t}^{(i_1i_2 i_3)}-
I_{(000)T,t}^{(i_1i_2 i_3)q_1}\right)^2\right\}\le
6\left(\frac{(T-t)^{3}}{6}-\sum_{j_1,j_2,j_3=0}^{q_1}
C_{j_3j_2j_1}^2\right)\ \ \ (i_1, i_2, i_3=1,\ldots,m).
\end{equation}

\vspace{7mm}

We may act similarly with more complicated 
iterated stochastic integrals. For example, for the 
approximation of stochastic integral
$I_{(0000)T,t}^{(i_1 i_2 i_3 i_4)}$ 
we may write (see (\ref{zzz10}))

\vspace{3mm}
$$
I_{(0000)T,t}^{(i_1 i_2 i_3 i_4)q_2}=
\sum_{j_1,j_2,j_3,j_4=0}^{q_2}
C_{j_4 j_3 j_2 j_1}\Biggl(
\zeta_{j_1}^{(i_1)}\zeta_{j_2}^{(i_2)}\zeta_{j_3}^{(i_3)}\zeta_{j_4}^{(i_4)}
-\Biggr.
$$
$$
-
{\bf 1}_{\{i_1=i_2\}}
{\bf 1}_{\{j_1=j_2\}}
\zeta_{j_3}^{(i_3)}
\zeta_{j_4}^{(i_4)}
-
{\bf 1}_{\{i_1=i_3\}}
{\bf 1}_{\{j_1=j_3\}}
\zeta_{j_2}^{(i_2)}
\zeta_{j_4}^{(i_4)}-
$$
$$
-
{\bf 1}_{\{i_1=i_4\}}
{\bf 1}_{\{j_1=j_4\}}
\zeta_{j_2}^{(i_2)}
\zeta_{j_3}^{(i_3)}
-
{\bf 1}_{\{i_2=i_3\}}
{\bf 1}_{\{j_2=j_3\}}
\zeta_{j_1}^{(i_1)}
\zeta_{j_4}^{(i_4)}-
$$
$$
-
{\bf 1}_{\{i_2=i_4\}}
{\bf 1}_{\{j_2=j_4\}}
\zeta_{j_1}^{(i_1)}
\zeta_{j_3}^{(i_3)}
-
{\bf 1}_{\{i_3=i_4\}}
{\bf 1}_{\{j_3=j_4\}}
\zeta_{j_1}^{(i_1)}
\zeta_{j_2}^{(i_2)}+
$$

\vspace{-3mm}
$$
+
{\bf 1}_{\{i_1=i_2\}}
{\bf 1}_{\{j_1=j_2\}}
{\bf 1}_{\{i_3=i_4\}}
{\bf 1}_{\{j_3=j_4\}}
+
$$

\vspace{-3mm}
$$
+
{\bf 1}_{\{i_1=i_3\}}
{\bf 1}_{\{j_1=j_3\}}
{\bf 1}_{\{i_2=i_4\}}
{\bf 1}_{\{j_2=j_4\}}+
$$
$$
+\Biggl.
{\bf 1}_{\{i_1=i_4\}}
{\bf 1}_{\{j_1=j_4\}}
{\bf 1}_{\{i_2=i_3\}}
{\bf 1}_{\{j_2=j_3\}}\Biggr),
$$

\vspace{7mm}
\noindent
where $C_{j_4 j_3 j_2 j_1}$ is defined by (\ref{zzz11}), (\ref{zzz12}).
Moreover, according to (\ref{zzz0})

\vspace{2mm}
$$
{\sf M}\left\{\left(
I_{(0000)T,t}^{(i_1i_2 i_3 i_4)}-
I_{(0000)T,t}^{(i_1i_2 i_3 i_4)q_2}\right)^2\right\}\le
24\left(\frac{(T-t)^{4}}{24}-\sum_{j_1,j_2,j_3,j_4=0}^{q_2}
C_{j_4j_3j_2j_1}^2\right)\ \ \ (i_1, i_2, i_3, i_4=1,\ldots,m).
$$

\vspace{5mm}

For pairwise different $i_1, i_2, i_3, i_4=1,\ldots,m$ from 
(\ref{zzz4}) we obtain

\vspace{2mm}
\begin{equation}
\label{r7}
{\sf M}\left\{\left(
I_{(0000)T,t}^{(i_1i_2 i_3i_4)}-
I_{(0000)T,t}^{(i_1i_2 i_3i_4)q_2}\right)^2\right\}=
\frac{(T-t)^{4}}{24}
-\sum_{j_1,j_2,j_3,j_4=0}^{q_2}
C_{j_4j_3j_2j_1}^2.
\end{equation}

\vspace{5mm}

\noindent
\begin{figure}
\begin{center}
\centerline{Table 3.\ Coefficients $\bar C_{3j_2j_1}.$}
\vspace{4mm}
\begin{tabular}{|c|c|c|c|c|c|c|c|c|}
\hline
${}_{j_2} {}^{j_1}$&0&1&2&3&4&5&6\\
\hline
0&$0$&$\frac{2}{105}$&$0$&$-\frac{4}{315}$&$0$&$\frac{2}{693}$&0\\
\hline
1&$\frac{4}{105}$&0&$-\frac{2}{315}$&0&$-\frac{8}{3465}$&0&$\frac{10}{9009}$\\
\hline
2&$\frac{2}{35}$&$-\frac{2}{105}$&$0$&$\frac{4}{3465}$&
$0$&$-\frac{74}{45045}$&0\\
\hline
3&$\frac{2}{315}$&$0$&$-\frac{2}{3465}$&0&
$\frac{16}{45045}$&0&$-\frac{10}{9009}$\\
\hline
4&$-\frac{2}{63}$&$\frac{46}{3465}$&0&$-\frac{32}{45045}$&
0&$\frac{2}{9009}$&0\\
\hline
5&$-\frac{10}{693}$&0&$\frac{38}{9009}$&0&
$-\frac{4}{9009}$&0&$\frac{122}{765765}$\\
\hline
6&$0$&$-\frac{10}{3003}$&$0$&$\frac{20}{9009}$&$0$&$-\frac{226}{765765}$&$0$\\
\hline
\end{tabular}
\end{center}
\vspace{7mm}
\begin{center}
\centerline{Table 4.\ Coefficients $\bar C_{21j_2j_1}.$}
\vspace{4mm}
\begin{tabular}{|c|c|c|c|}
\hline
${}_{j_2} {}^{j_1}$&0&1&2\\
\hline
0&$\frac{2}{21}$&$-\frac{2}{45}$&$\frac{2}{315}$\\
\hline
1&$\frac{2}{315}$&$\frac{2}{315}$&$-\frac{2}{225}$\\
\hline
2&$-\frac{2}{105}$&$\frac{2}{225}$&$\frac{2}{1155}$\\
\hline
\end{tabular}
\end{center}
\vspace{7mm}
\begin{center}
\centerline{Table 5.\ Coefficients $\bar C_{101j_2j_1}.$}
\vspace{4mm}
\begin{tabular}{|c|c|c|}
\hline
${}_{j_2} {}^{j_1}$&0&1\\
\hline
0&$\frac{4}{315}$&$0$\\
\hline
1&$\frac{4}{315}$&$-\frac{8}{945}$\\
\hline
\end{tabular}
\end{center}
\vspace{5mm}
\end{figure}

Using Theorem 17, we can calculate exactly the left-hand
side of (\ref{r7})
for any possible combinations
of $i_1, i_2, i_3, i_4$. These relations were obtained in 
\cite{arxiv-3}, \cite{20}-\cite{12aa-afterxxx}.

In Tables 3--5, we have some examples of exact 
values of the Fourier--Legendre coefficients 
(here and further the Fourier--Legendre coefficients 
have been calculated exactly using DERIVE (computer
algebra system)). Note that in \cite{Kuz-Kuz}, \cite{Mikh-1}
the database with 270,000 exactly
calculated Fourier--Legendre coefficients was described.
This database was used in the software package,
which is written in the Python programming language
for the implementation of the numerical schemes (\ref{al1})-(\ref{al5}),
(\ref{al1x})-(\ref{al5x}).

Assume that $q_1=6$.
Calculating the value of expression 
(\ref{39}) for $q_1=6,$ 
$i_1\ne i_2,$ $i_1\ne i_3,$ $i_3\ne i_2,$
we obtain

$$
{\sf M}\left\{\left(
I_{(000)T,t}^{(i_1i_2 i_3)}-
I_{(000)T,t}^{(i_1i_2 i_3)q_1}\right)^2\right\}\approx
0.01956(T-t)^3.
$$

\vspace{5mm}

Let us choose, for example, $q_2=2.$ 
In the case of pairwise different
$i_1, i_2, i_3, i_4$ we have from (\ref{r7}) the following  
approximate equality

\vspace{2mm}
\begin{equation}
\label{46000}
{\sf M}\left\{\left(
I_{(0000)T,t}^{(i_1i_2i_3 i_4)}-
I_{(0000)T,t}^{(i_1i_2i_3 i_4)q_2}\right)^2\right\}\approx
0.0236084(T-t)^4.
\end{equation}

\vspace{5mm}

Let us analyze the approximations 

\vspace{2mm}
$$
I_{(001)T,t}^{(i_1i_2i_3)q_3},\ \ \
I_{(010)T,t}^{(i_1i_2i_3)q_3},\ \ \
I_{(100)T,t}^{(i_1i_2i_3)q_3},\ \ \
I_{(00000)T,t}^{(i_1i_2i_3i_4 i_5)q_4}
$$

\vspace{5mm}
\noindent
based on the expansions (\ref{sss1})--(\ref{sss4}).

Assume that 
$q_3=2,$ $q_4=1.$  
In the case of pairwise different 
$i_1, \ldots, i_5$ we obtain

\vspace{2mm}
$$
{\sf M}\left\{\left(
I_{(100)T,t}^{(i_1i_2 i_3)}-
I_{(100)T,t}^{(i_1i_2 i_3)q_3}\right)^2\right\}=
\frac{(T-t)^{5}}{60}-\sum_{j_1,j_2,j_3=0}^{2}
\left(C_{j_3j_2j_1}^{100}\right)^2\approx 0.00815429(T-t)^5,
$$

\vspace{4mm}
$$
{\sf M}\left\{\left(
I_{(010)T,t}^{(i_1i_2 i_3)}-
I_{(010)T,t}^{(i_1i_2 i_3)q_3}\right)^2\right\}=
\frac{(T-t)^{5}}{20}-\sum_{j_1,j_2,j_3=0}^{2}
\left(C_{j_3j_2j_1}^{010}\right)^2\approx 0.0173903(T-t)^5,
$$

\vspace{4mm}
$$
{\sf M}\left\{\left(
I_{(001)T,t}^{(i_1i_2 i_3)}-
I_{(001)T,t}^{(i_1i_2 i_3)q_3}\right)^2\right\}=
\frac{(T-t)^5}{10}-\sum_{j_1,j_2,j_3=0}^{2}
\left(C_{j_3j_2j_1}^{001}\right)^2
\approx 0.0252801(T-t)^5,
$$

\vspace{6mm}
$$
{\sf M}\left\{\left(
I_{(00000)T,t}^{(i_1i_2i_3i_4 i_5)}-
I_{(00000)T,t}^{(i_1i_2i_3i_4 i_5)q_4}\right)^2\right\}=
$$

\vspace{2mm}
$$
=
\frac{(T-t)^5}{120}-\sum_{j_1,j_2,j_3,j_4,j_5=0}^{1}
C_{j_5j_4j_3j_2j_1}^2\approx 0.00759105(T-t)^5.
$$

\vspace{7mm}

Note that using (\ref{zzz0}), we can write 

\vspace{2mm}
$$
{\sf M}\left\{\left(
I_{(00000)T,t}^{(i_1i_2 i_3 i_4 i_5)}-
I_{(00000)T,t}^{(i_1i_2 i_3 i_4 i_5)q_4}\right)^2\right\}\le
120\left(\frac{(T-t)^{5}}{120}-\sum_{j_1,j_2,j_3,j_4,j_5=0}^{q_4}
C_{j_5j_4j_3j_2j_1}^2\right),
$$

\vspace{5mm}
\noindent
where $i_1, \ldots, i_5=1,\ldots,m$.

Moreover, from (\ref{zzz0}) we obtain the following useful estimates

\vspace{3mm}

$$
{\sf M}\left\{\left(
I_{(01)T,t}^{(i_1i_2)}-
I_{(01)T,t}^{(i_1i_2)q}\right)^2\right\}\le
2\Biggl(\frac{(T-t)^{4}}{4}-\sum_{j_1,j_2=0}^{q}
\left(C_{j_2j_1}^{01}\right)^2\Biggr),
$$

\vspace{3mm}
$$
{\sf M}\left\{\left(
I_{(10)T,t}^{(i_1i_2)}-
I_{(10)T,t}^{(i_1i_2)q}\right)^2\right\}\le
2\Biggl(\frac{(T-t)^{4}}{12}-\sum_{j_1,j_2=0}^{q}
\left(C_{j_2j_1}^{10}\right)^2\Biggr),
$$

\vspace{3mm}
$$
{\sf M}\left\{\left(
I_{(100)T,t}^{(i_1i_2 i_3)}-
I_{(100)T,t}^{(i_1i_2 i_3)q}\right)^2\right\}\le
6\Biggl(\frac{(T-t)^{5}}{60}-\sum_{j_1,j_2,j_3=0}^{q}
\left(C_{j_3j_2j_1}^{100}\right)^2\Biggr),
$$

\vspace{4mm}

$$
{\sf M}\left\{\left(
I_{(010)T,t}^{(i_1i_2 i_3)}-
I_{(010)T,t}^{(i_1i_2 i_3)q}\right)^2\right\}\le
6\Biggl(\frac{(T-t)^{5}}{20}-\sum_{j_1,j_2,j_3=0}^{q}
\left(C_{j_3j_2j_1}^{010}\right)^2\Biggr),
$$
                          
\vspace{4mm}

$$
{\sf M}\left\{\left(
I_{(001)T,t}^{(i_1i_2 i_3)}-
I_{(001)T,t}^{(i_1i_2 i_3)q}\right)^2\right\}\le
6\Biggl(\frac{(T-t)^5}{10}-\sum_{j_1,j_2,j_3=0}^{q}
\left(C_{j_3j_2j_1}^{001}\right)^2\Biggr),
$$

\vspace{4mm}

$$
{\sf M}\left\{\left(
I_{(20))T,t}^{(i_1i_2)}-
I_{(20)T,t}^{(i_1i_2)q}\right)^2\right\}\le
2\Biggl(\frac{(T-t)^6}{30}-\sum_{j_2,j_1=0}^{q}
\left(C_{j_2j_1}^{20}\right)^2\Biggr),
$$

\vspace{4mm}

$$
{\sf M}\left\{\left(
I_{(11)T,t}^{(i_1i_2)}-
I_{(11)T,t}^{(i_1i_2)q}\right)^2\right\}\le
2\Biggl(\frac{(T-t)^6}{18}-\sum_{j_2,j_1=0}^{q}
\left(C_{j_2j_1}^{11}\right)^2\Biggr),
$$

\vspace{4mm}

$$
{\sf M}\left\{\left(
I_{(02)T,t}^{(i_1i_2)}-
I_{(02)T,t}^{(i_1i_2)q}\right)^2\right\}\le
2\Biggl(\frac{(T-t)^6}{6}-\sum_{j_2,j_1=0}^{q}
\left(C_{j_2j_1}^{02}\right)^2\Biggr),
$$

\vspace{4mm}

$$
{\sf M}\left\{\left(
I_{(1000)T,t}^{(i_1i_2 i_3i_4)}-
I_{(1000)T,t}^{(i_1i_2 i_3i_4)q}\right)^2\right\}\le
24\Biggl(\frac{(T-t)^{6}}{360}-\sum_{j_1,j_2,j_3, j_4=0}^{q}
\left(C_{j_4j_3j_2j_1}^{1000}\right)^2\Biggr),
$$

\vspace{4mm}

$$
{\sf M}\left\{\left(
I_{(0100)T,t}^{(i_1i_2 i_3i_4)}-
I_{(0100)T,t}^{(i_1i_2 i_3i_4)q}\right)^2\right\}\le
24\Biggl(\frac{(T-t)^{6}}{120}-\sum_{j_1,j_2,j_3, j_4=0}^{q}
\left(C_{j_4j_3j_2j_1}^{0100}\right)^2\Biggr),
$$

\vspace{4mm}

$$
{\sf M}\left\{\left(
I_{(0010)T,t}^{(i_1i_2 i_3i_4)}-
I_{(0010)T,t}^{(i_1i_2 i_3 i_4)q}\right)^2\right\}\le
24\Biggl(\frac{(T-t)^6}{60}-\sum_{j_1,j_2,j_3, j_4=0}^{q}
\left(C_{j_4j_3j_2j_1}^{0010}\right)^2\Biggr),
$$

\vspace{4mm}

$$
{\sf M}\left\{\left(
I_{(0001)T,t}^{(i_1i_2 i_3 i_4)}-
I_{(0001)T,t}^{(i_1i_2 i_3 i_4)q}\right)^2\right\}\le
24\Biggl(\frac{(T-t)^6}{36}-\sum_{j_1,j_2,j_3, j_4=0}^{q}
\left(C_{j_4j_3j_2j_1}^{0001}\right)^2\Biggr),
$$

\vspace{2mm}

$$
{\sf M}\left\{\left(
I_{(000000)T,t}^{(i_1 i_2 i_3 i_4 i_5 i_6)}-
I_{(000000)T,t}^{(i_1 i_2 i_3 i_4 i_5 i_6)q}\right)^2\right\}\le
720\left(\frac{(T-t)^{6}}{720}-\sum_{j_1,j_2,j_3,j_4,j_5,j_6=0}^{q}
C_{j_6 j_5 j_4 j_3 j_2 j_1}^2\right).
$$

\vspace{7mm}

In addition, from (\ref{zzz4x}) we get

\vspace{3mm}
$$
{\sf M}\biggl\{\left(I_{(10)T,t}^{(i_1 i_2)}-
I_{(10)T,t}^{(i_1 i_2)q}
\right)^2\biggr\}
=\frac{(T-t)^4}{12}
-\sum_{j_1,j_2=0}^q
\left(C_{j_2j_1}^{10}\right)^2-
\sum_{j_1,j_2=0}^q C_{j_2j_1}^{10}C_{j_1j_2}^{10}\ \ \ (i_1=i_2),
$$

\vspace{5mm}

$$
{\sf M}\biggl\{\left(I_{(10)T,t}^{(i_1 i_2)}-
I_{(10)T,t}^{(i_1 i_2)q}
\right)^2\biggr\}=\frac{(T-t)^4}{12}
-\sum_{j_1,j_2=0}^q
\left(C_{j_2j_1}^{10}\right)^2\ \ \ (i_1\ne i_2),
$$

\vspace{5mm}

$$
{\sf M}\biggl\{\left(I_{(01)T,t}^{(i_1 i_2)}-
I_{(01)T,t}^{(i_1 i_2)q}
\right)^2\biggr\}
=\frac{(T-t)^4}{4}
-\sum_{j_1,j_2=0}^q
\left(C_{j_2j_1}^{01}\right)^2-
\sum_{j_1,j_2=0}^q C_{j_2j_1}^{01}C_{j_1j_2}^{01}\ \ \ (i_1=i_2),
$$

\vspace{5mm}

$$
{\sf M}\biggl\{\left(I_{(01)T,t}^{(i_1 i_2)}-
I_{(01)T,t}^{(i_1 i_2)q}
\right)^2\biggr\}=\frac{(T-t)^4}{4}
-\sum_{j_1,j_2=0}^q
\left(C_{j_2j_1}^{01}\right)^2\ \ \ (i_1\ne i_2),
$$

\vspace{5mm}

$$
{\sf M}\biggl\{\left(I_{(20)T,t}^{(i_1 i_2)}-
I_{(20)T,t}^{(i_1 i_2)q}
\right)^2\biggr\}=\frac{(T-t)^6}{30}
-\sum_{j_1,j_2=0}^q
\left(C_{j_2j_1}^{20}\right)^2-
\sum_{j_1,j_2=0}^q C_{j_2j_1}^{20}C_{j_1j_2}^{20}\ \ \ (i_1=i_2),
$$

\vspace{5mm}

$$
{\sf M}\biggl\{\left(I_{(20)T,t}^{(i_1 i_2)}-
I_{(20)T,t}^{(i_1 i_2)q}
\right)^2\biggr\}=\frac{(T-t)^6}{30}
-\sum_{j_1,j_2=0}^q
\left(C_{j_2j_1}^{20}\right)^2\ \ \ (i_1\ne i_2),
$$

\vspace{5mm}

$$
{\sf M}\biggl\{\left(I_{(11)T,t}^{(i_1 i_2)}-
I_{(11)T,t}^{(i_1 i_2)q}
\right)^2\biggr\}=\frac{(T-t)^6}{18}
-\sum_{j_1,j_2=0}^q
\left(C_{j_2j_1}^{11}\right)^2-
\sum_{j_1,j_2=0}^q C_{j_2j_1}^{11}C_{j_1j_2}^{11}\ \ \ (i_1=i_2),
$$

\vspace{5mm}

$$
{\sf M}\biggl\{\left(I_{(11)T,t}^{(i_1 i_2)}-
I_{(11)T,t}^{(i_1 i_2)q}
\right)^2\biggr\}=\frac{(T-t)^6}{18}
-\sum_{j_1,j_2=0}^q
\left(C_{j_2j_1}^{11}\right)^2\ \ \ (i_1\ne i_2),
$$

\vspace{5mm}

$$
{\sf M}\biggl\{\left(I_{(02)}^{(i_1 i_2)}-
I_{(02)T,t}^{(i_1 i_2)q}
\right)^2\biggr\}=\frac{(T-t)^6}{6}
-\sum_{j_1,j_2=0}^q
\left(C_{j_2j_1}^{02}\right)^2-
\sum_{j_1,j_2=0}^q C_{j_2j_1}^{02}C_{j_1j_2}^{02}\ \ \ (i_1=i_2),
$$

\vspace{5mm}

$$
{\sf M}\biggl\{\left(I_{(02)T,t}^{(i_1 i_2)}-
I_{(02)T,t}^{(i_1 i_2)q}
\right)^2\biggr\}=\frac{(T-t)^6}{6}
-\sum_{j_1,j_2=0}^q
\left(C_{j_2j_1}^{02}\right)^2\ \ \ (i_1\ne i_2),
$$

\vspace{6mm}

Clearly, expansions for iterated Stratonovich stochastic integrals
(see Theorems 9--16) are simpler than expansions for
iterated Ito stochastic integrals (see Theorems 7, 8, and
(\ref{za1})--(\ref{a6})). However, the calculation of the mean-square
approximation error for iterated Stratonovich
stochastic integrals turns out to be much more difficult than for 
iterated Ito stochastic integrals (see Theorem 17 and (\ref{zzz0})).
Below we consider how we can estimate or calculate exactly
(for some particular cases)
the mean-square
approximation error for iterated Stratonovich
stochastic integrals (the development of these results is contained in Chapter 5 of \cite{20aa}
(also see Theorems~10--12 from this paper).

As we mentioned above, on the basis of 
the presented 
approximations of 
iterated Stratonovich stochastic integrals we 
can see that increasing of multiplicities of these integrals 
leads to increasing 
of orders of smallness with respect to $T-t$
in the mean-square sense 
for iterated Stratonovich stochastic integrals
($T-t\ll 1$ 
because the length of integration interval $[t, T]$ 
of the iterated Stratonovich 
stochastic integrals 
plays the role of integration step for the numerical 
methods for Ito SDEs, so $T-t$ is already fairly small).
This leads to a sharp decrease  
of member 
quantities
in the appro\-xi\-ma\-ti\-ons of iterated Stratonovich stochastic 
integrals,
which are required for achieving the acceptable accuracy
of approximation.

From (\ref{fff09}) $(i_1\ne i_2)$ we obtain

$$
{\sf M}\left\{\left(I_{(00)T,t}^{*(i_1 i_2)}-
I_{(00)T,t}^{*(i_1 i_2)q}
\right)^2\right\}=\frac{(T-t)^2}{2}
\sum\limits_{i=q+1}^{\infty}\frac{1}{4i^2-1}\le 
$$

\begin{equation}
\label{teac}
\le \frac{(T-t)^2}{2}\int\limits_{q}^{\infty}
\frac{1}{4x^2-1}dx
=-\frac{(T-t)^2}{8}{\rm ln}\left|
1-\frac{2}{2q+1}\right|\le C_1\frac{(T-t)^2}{q},
\end{equation}

\vspace{3mm}
\noindent
where $C_1$ is a constant.

Since $T-t\ll 1,$ then it is easy to notice that there 
exists such a constant $C_2$ that

\begin{equation}
\label{teac3}
{\sf M}\left\{\left(
I_{(l_1\ldots l_k)T,t}^{*(i_1\ldots i_k)}-
I_{(l_1\ldots l_k)T,t}^{*(i_1\ldots i_k)q}\right)^2\right\}
\le C_2 {\sf M}\left\{\left(I_{(00)T,t}^{*(i_1 i_2)}-
I_{(00)T,t}^{*(i_1 i_2)q}\right)^2\right\},
\end{equation}

\vspace{3mm}
\noindent
where 
$I_{(l_1\ldots l_k)T,t}^{*(i_1\ldots i_k)q}$
is an approximation of the iterated Stratonovich stochastic integral 
$I_{(l_1\ldots l_k)T,t}^{*(i_1\ldots i_k)}.$

From (\ref{teac}) and (\ref{teac3}) we finally obtain

\begin{equation}
\label{teac4}
{\sf M}\left\{\left(
I_{(l_1\ldots l_k)T,t}^{*(i_1\ldots i_k)}-
I_{(l_1\ldots l_k)T,t}^{*(i_1\ldots i_k)q}\right)^2\right\}
\le C \frac{(T-t)^2}{q},
\end{equation}

\vspace{3mm}
\noindent
where constant $C$ does not depends on $T-t$.
Note that, in contrast to the estimate (\ref{teac4}), 
the constant $C$ in Theorems 10--12 does not depend on $q.$

The same idea can be found in \cite{KlPl2} in the framework of 
the method based
on the trigonometric expansion of the
Brownian bridge process.

We can get more information about the numbers $q$ (these
numbers are different for different iterated Stratonovich
stochastic integrals)
using the another approach.
Since for pairwise different $i_1,\ldots,i_k=1,\ldots,m$

\vspace{-1mm}
$$
J^{*}[\psi^{(k)}]_{T,t}=J[\psi^{(k)}]_{T,t}\ \ \ \hbox{w.\ p.\ 1,}
$$

\vspace{4mm}
\noindent
where $J[\psi^{(k)}]_{T,t},$ $J^{*}[\psi^{(k)}]_{T,t}$
are defined by (\ref{ito}) and (\ref{str}) correspondingly,
then 
for pairwise different 
$i_1,\ldots,i_6=1,\ldots,m$ we can write (see (\ref{zzz4}))

\vspace{2mm}

$$
{\sf M}\left\{\left(
I_{(01)T,t}^{*(i_1i_2)}-
I_{(01)T,t}^{*(i_1i_2)q}\right)^2\right\}=
\frac{(T-t)^{4}}{4}-\sum_{j_1,j_2=0}^{q}
\left(C_{j_2j_1}^{01}\right)^2,
$$

\vspace{3mm}
$$
{\sf M}\left\{\left(
I_{(10)T,t}^{*(i_1i_2)}-
I_{(10)T,t}^{*(i_1i_2)q}\right)^2\right\}=
\frac{(T-t)^{4}}{12}-\sum_{j_1,j_2=0}^{q}
\left(C_{j_2j_1}^{10}\right)^2,
$$

\vspace{3mm}

$$
{\sf M}\left\{\left(
I_{(000)T,t}^{*(i_1i_2 i_3)}-
I_{(000)T,t}^{*(i_1i_2 i_3)q}\right)^2\right\}=
\frac{(T-t)^{3}}{6}-\sum_{j_3,j_2,j_1=0}^{q}
C_{j_3j_2j_1}^2,
$$

\vspace{3mm}

$$
{\sf M}\left\{\left(
I_{(0000)T,t}^{*(i_1i_2 i_3 i_4)}-
I_{(0000)T,t}^{*(i_1i_2 i_3 i_4)q}\right)^2\right\}=
\frac{(T-t)^{4}}{24}-\sum_{j_1,j_2,j_3,j_4=0}^{q}
C_{j_4j_3j_2j_1}^2,
$$

\vspace{3mm}

$$
{\sf M}\left\{\left(
I_{(100)T,t}^{*(i_1i_2 i_3)}-
I_{(100)T,t}^{*(i_1i_2 i_3)q}\right)^2\right\}=
\frac{(T-t)^{5}}{60}-\sum_{j_1,j_2,j_3=0}^{q}
\left(C_{j_3j_2j_1}^{100}\right)^2,
$$

\vspace{3mm}

$$
{\sf M}\left\{\left(
I_{(010))T,t}^{*(i_1i_2 i_3)}-
I_{(010)T,t}^{*(i_1i_2 i_3)q}\right)^2\right\}=
\frac{(T-t)^{5}}{20}-\sum_{j_1,j_2,j_3=0}^{q}
\left(C_{j_3j_2j_1}^{010}\right)^2,
$$

\vspace{3mm}

$$
{\sf M}\left\{\left(
I_{(001)T,t}^{*(i_1i_2 i_3)}-
I_{(001)T,t}^{*(i_1i_2 i_3)q}\right)^2\right\}=
\frac{(T-t)^5}{10}-\sum_{j_1,j_2,j_3=0}^{q}
\left(C_{j_3j_2j_1}^{001}\right)^2,
$$

\vspace{3mm}

$$
{\sf M}\left\{\left(
I_{(00000)T,t}^{*(i_1 i_2 i_3 i_4 i_5)}-
I_{(00000)T,t}^{*(i_1 i_2 i_3 i_4 i_5)q}\right)^2\right\}=
\frac{(T-t)^{5}}{120}-\sum_{j_1,j_2,j_3,j_4,j_5=0}^{q}
C_{j_5 i_4 i_3 i_2 j_1}^2,
$$

\vspace{3mm}

$$
{\sf M}\left\{\left(
I_{(20)T,t}^{*(i_1i_2)}-
I_{(20)T,t}^{*(i_1i_2)q}\right)^2\right\}=
\frac{(T-t)^6}{30}-\sum_{j_2,j_1=0}^{q}
\left(C_{j_2j_1}^{20}\right)^2,
$$

\vspace{3mm}

$$
{\sf M}\left\{\left(
I_{(11)T,t}^{*(i_1i_2)}-
I_{(11)T,t}^{*(i_1i_2)q}\right)^2\right\}=
\frac{(T-t)^6}{18}-\sum_{j_2,j_1=0}^{q}
\left(C_{j_2j_1}^{11}\right)^2,
$$

\vspace{3mm}

$$
{\sf M}\left\{\left(
I_{(02)T,t}^{*(i_1i_2)}-
I_{(02)T,t}^{*(i_1i_2)q}\right)^2\right\}=
\frac{(T-t)^6}{6}-\sum_{j_2,j_1=0}^{q}
\left(C_{j_2j_1}^{02}\right)^2,
$$

\vspace{3mm}

$$
{\sf M}\left\{\left(
I_{(1000)T,t}^{*(i_1i_2 i_3i_4)}-
I_{(1000)T,t}^{*(i_1i_2 i_3i_4)q}\right)^2\right\}=
\frac{(T-t)^{6}}{360}-\sum_{j_1,j_2,j_3, j_4=0}^{q}
\left(C_{j_4j_3j_2j_1}^{1000}\right)^2,
$$

\vspace{3mm}

$$
{\sf M}\left\{\left(
I_{(0100)T,t}^{*(i_1i_2 i_3i_4)}-
I_{(0100)T,t}^{*(i_1i_2 i_3i_4)q}\right)^2\right\}=
\frac{(T-t)^{6}}{120}-\sum_{j_1,j_2,j_3, j_4=0}^{q}
\left(C_{j_4j_3j_2j_1}^{0100}\right)^2,
$$

\vspace{3mm}

$$
{\sf M}\left\{\left(
I_{(0010)T,t}^{*(i_1i_2 i_3i_4)}-
I_{(0010)T,t}^{*(i_1i_2 i_3 i_4)q}\right)^2\right\}=
\frac{(T-t)^6}{60}-\sum_{j_1,j_2,j_3, j_4=0}^{q}
\left(C_{j_4j_3j_2j_1}^{0010}\right)^2,
$$

\vspace{3mm}

$$
{\sf M}\left\{\left(
I_{(0001)T,t}^{*(i_1i_2 i_3 i_4)}-
I_{(0001)T,t}^{*(i_1i_2 i_3 i_4)q}\right)^2\right\}=
\frac{(T-t)^6}{36}-\sum_{j_1,j_2,j_3, j_4=0}^{q}
\left(C_{j_4j_3j_2j_1}^{0001}\right)^2,
$$

\vspace{3mm}

$$
{\sf M}\left\{\left(
I_{(000000)T,t}^{*(i_1 i_2 i_3 i_4 i_5 i_6)}-
I_{(000000)T,t}^{*(i_1 i_2 i_3 i_4 i_5 i_6)q}\right)^2\right\}=
\frac{(T-t)^{6}}{720}-\sum_{j_1,j_2,j_3,j_4,j_5,j_6=0}^{q}
C_{j_6 j_5 j_4 j_3 j_2 j_1}^2.
$$

\vspace{9mm}

For example \cite{7} (also see \cite{arxiv-15}-\cite{arxiv-18},
\cite{5-003}, \cite{8}-\cite{12aa-afterxxx}),

\vspace{3mm}

$$
{\sf M}\left\{\left(
I_{(000)T,t}^{*(i_1i_2 i_3)}-
I_{(000)T,t}^{*(i_1i_2 i_3)6}\right)^2\right\}=
\frac{(T-t)^{3}}{6}-\sum_{j_3,j_2,j_1=0}^{6}
C_{j_3j_2j_1}^2
\approx
0.01956000(T-t)^3,
$$

\vspace{3mm}

$$
{\sf M}\left\{\left(
I_{(0000)T,t}^{*(i_1i_2i_3 i_4)}-
I_{(0000)T,t}^{*(i_1i_2i_3 i_4)2}\right)^2\right\}=
\frac{(T-t)^{4}}{24}-\sum_{j_1,j_2,j_3,j_4=0}^{2}
C_{j_4 j_3 j_2 j_1}^2\approx
0.02360840(T-t)^4,
$$

\vspace{3mm}

$$
{\sf M}\left\{\left(
I_{(100)T,t}^{*(i_1i_2 i_3)}-
I_{(100)T,t}^{*(i_1i_2 i_3)2}\right)^2\right\}=
\frac{(T-t)^{5}}{60}-\sum_{j_1,j_2,j_3=0}^{2}
\left(C_{j_3j_2j_1}^{100}\right)^2\approx 0.00815429(T-t)^5,
$$

\vspace{3mm}

$$
{\sf M}\left\{\left(
I_{(010)T,t}^{*(i_1i_2 i_3)}-
I_{(010)T,t}^{*(i_1i_2 i_3)2}\right)^2\right\}=
\frac{(T-t)^{5}}{20}-\sum_{j_1,j_2,j_3=0}^{2}
\left(C_{j_3j_2j_1}^{010}\right)^2\approx 0.0173903(T-t)^5,
$$

\vspace{3mm}

$$
{\sf M}\left\{\left(
I_{(001)T,t}^{*(i_1i_2 i_3)}-
I_{(001)T,t}^{*(i_1i_2 i_3)2}\right)^2\right\}=
\frac{(T-t)^5}{10}-\sum_{j_1,j_2,j_3=0}^{2}
\left(C_{j_3j_2j_1}^{001}\right)^2
\approx 0.0252801(T-t)^5,
$$

\vspace{6mm}
$$
{\sf M}\left\{\left(
I_{(00000)T,t}^{*(i_1i_2i_3i_4 i_5)}-
I_{(00000)T,t}^{*(i_1i_2i_3i_4 i_5)1}\right)^2\right\}=
$$

\vspace{2mm}
$$
=
\frac{(T-t)^5}{120}-\sum_{j_1,j_2,j_3,j_4,j_5=0}^{1}
C_{j_5j_4j_3j_2j_1}^2\approx 0.00759105(T-t)^5.
$$

\vspace{5mm}

Let us consider expansions of the
Ito stochastic integrals $I_{(1)T,t}^{(i_1)},$
$I_{(2)T,t}^{(i_1)}$ based on the approach from
\cite{Mi2} (also see \cite{KlPl2})

\vspace{1mm}
\begin{equation}
\label{ww1}
I_{(1)T,t}^{(i_1)q}=-\frac{{(T-t)}^{3/2}}{2}
\left(\zeta_0^{(i_1)}-\frac{\sqrt{2}}{\pi}\left(\sum_{r=1}^{q}
\frac{1}{r}
\zeta_{2r-1}^{(i_1)}+\sqrt{\alpha_q}\xi_q^{(i_1)}\right)
\right),
\end{equation}

\vspace{4mm}
\begin{equation}
\label{ww2}
I_{(2)T,t}^{(i_1)q}=
(T-t)^{5/2}\left(
\frac{1}{3}\zeta_0^{(i_1)}+\frac{1}{\sqrt{2}\pi^2}
\left(\sum_{r=1}^{q}\frac{1}{r^2}\zeta_{2r}^{(i_1)}+
\sqrt{\beta_q}
\mu_q^{(i_1)}\right)-\right.
$$

\vspace{2mm}
$$
\left.-\frac{1}{\sqrt{2}\pi}\left(\sum_{r=1}^q
\frac{1}{r}\zeta_{2r-1}^{(i_1)}+
\sqrt{\alpha_q}\xi_q^{(i_1)}\right)\right),
\end{equation}

\vspace{6mm}
\noindent
where $\zeta_j^{(i)}$ is defined by the formula (\ref{rr23}),
$\phi_j(\tau)$ is a complete orthonormal system of trigonometric
functions in the space
$L_2([t, T]),$ and
$\zeta_0^{(i)},$ $\zeta_{2r}^{(i)},$
$\zeta_{2r-1}^{(i)},$ $\xi_q^{(i)},$ $\mu_q^{(i)}$ $(r=1,\ldots,q,$\ \
$i=1,\ldots,m)$ are independent 
standard Gaussian random variables, $i_1=1,\ldots,m,$ 

\vspace{2mm}
$$
\xi_q^{(i)}=\frac{1}{\sqrt{\alpha_q}}\sum_{r=q+1}^{\infty}
\frac{1}{r}\zeta_{2r-1}^{(i)},\ \ \ \
\alpha_q=\frac{\pi^2}{6}-\sum_{r=1}^q\frac{1}{r^2},
$$

\vspace{2mm}
$$
\mu_q^{(i)}=\frac{1}{\sqrt{\beta_q}}\sum_{r=q+1}^{\infty}
\frac{1}{r^2}~\zeta_{2r}^{(i)},\ \ \ \
\beta_q=\frac{\pi^4}{90}-\sum_{r=1}^q\frac{1}{r^4}.
$$

\vspace{5mm}

\noindent
It is obvious that (\ref{ww1}), (\ref{ww2})
significantly more complicated compared to
(\ref{4002}), (\ref{4003}).

Another example of obvious advantage of the Legendre polynomials 
over the trigonometric functions (in the framework of the considered 
problem) is the truncated expansion 
of the iterated Stratonovich stochastic integral
$I_{(10)T, t}^{*(i_1 i_2)}$ obtained by Theorem 9, in which
instead of the double Fourier--Legendre series (see (\ref{4004}),
(\ref{4006})) is taken 
the double trigonometric Fourier series

\vspace{1mm}
$$
I_{(10)T,t}^{*(i_1 i_2)q}=-(T-t)^{2}\Biggl(\frac{1}{6}
\zeta_{0}^{(i_1)}\zeta_{0}^{(i_2)}-\frac{1}{2\sqrt{2}\pi}
\sqrt{\alpha_q}\xi_q^{(i_2)}\zeta_0^{(i_1)}+\Biggr.
$$

\vspace{1mm}
$$
+\frac{1}{2\sqrt{2}\pi^2}\sqrt{\beta_q}\Biggl(
\mu_q^{(i_2)}\zeta_0^{(i_1)}-2\mu_q^{(i_1)}\zeta_0^{(i_2)}\Biggr)+
$$

\vspace{1mm}
$$
+\frac{1}{2\sqrt{2}}\sum_{r=1}^{q}
\Biggl(-\frac{1}{\pi r}
\zeta_{2r-1}^{(i_2)}
\zeta_{0}^{(i_1)}+
\frac{1}{\pi^2 r^2}\left(
\zeta_{2r}^{(i_2)}
\zeta_{0}^{(i_1)}-
2\zeta_{2r}^{(i_1)}
\zeta_{0}^{(i_2)}\right)\Biggr)-
$$

\vspace{1mm}
$$
-
\frac{1}{2\pi^2}\sum_{\stackrel{r,l=1}{{}_{r\ne l}}}^{q}
\frac{1}{r^2-l^2}\Biggl(
\zeta_{2r}^{(i_1)}
\zeta_{2l}^{(i_2)}+
\frac{l}{r}
\zeta_{2r-1}^{(i_1)}
\zeta_{2l-1}^{(i_2)}
\Biggr)+
$$

$$
+
\sum_{r=1}^{q}
\Biggl(\frac{1}{4\pi r}\left(
\zeta_{2r}^{(i_1)}
\zeta_{2r-1}^{(i_2)}-
\zeta_{2r-1}^{(i_1)}
\zeta_{2r}^{(i_2)}\right)+
$$

\begin{equation}
\label{944}
+
\Biggl.\Biggl.
\frac{1}{8\pi^2 r^2}\left(
3\zeta_{2r-1}^{(i_1)}
\zeta_{2r-1}^{(i_2)}+
\zeta_{2r}^{(i_2)}
\zeta_{2r}^{(i_1)}\right)\Biggr)\Biggr),
\end{equation}

\vspace{5mm}
\noindent
where the meaning of the notations included 
in (\ref{ww1}), (\ref{ww2}) is preserved.

A deep comparative analysis of the efficiency of application 
of Legendre polynomials and tri\-go\-no\-met\-ric functions to the numerical 
integration of Ito SDEs is given in \cite{arxiv-12}, \cite{5-004}.

\vspace{5mm}

\section{Theorems 7--16 from Point
of View of the Wong--Zakai Approximation}

\vspace{5mm}

The iterated Ito stochastic integrals and solutions
of Ito SDEs are complex and important func\-ti\-o\-nals
from the independent components ${\bf f}_{s}^{(i)},$
$i=1,\ldots,m$ of the multidimensional
Wiener process ${\bf f}_{s},$ $s\in[0, T].$
Let ${\bf f}_{s}^{(i)p},$ $p\in\mathbb{N}$ 
be some approximation of
${\bf f}_{s}^{(i)},$
$i=1,\ldots,m$.
Suppose that 
${\bf f}_{s}^{(i)p}$
converges to
${\bf f}_{s}^{(i)},$
$i=1,\ldots,m$ if $p\to\infty$ in some sense and has
differentiable sample trajectories.

A natural question arises: if we replace 
${\bf f}_{s}^{(i)}$
by ${\bf f}_{s}^{(i)p},$
$i=1,\ldots,m$ in the functionals
mentioned above, will the resulting
functionals converge to the original
functionals from the components 
${\bf f}_{s}^{(i)},$
$i=1,\ldots,m$ of the multidimentional
Wiener process ${\bf f}_{s}$?
The answere to this question is negative 
in the general case. However, 
in the pioneering works of Wong E. and Zakai M. \cite{W-Z-1},
\cite{W-Z-2},
it was shown that under the special conditions and 
for some types of approximations 
of the Wiener process the answere is affirmative
with one peculiarity: the convergence takes place 
to the iterated Stratonovich stochastic integrals
and solutions of Stratonovich SDEs and not to iterated 
Ito stochastic integrals and solutions
of Ito SDEs.
The piecewise 
linear approximation 
as well as the regularization by convolution 
\cite{W-Z-1}-\cite{Watanabe} relate the 
mentioned types of approximations
of the Wiener process. The above approximation 
of stochastic integrals and solutions of SDEs 
is often called the Wong--Zakai approximation.

Let ${\bf w}_{\tau},$ $\tau\in[0, T]$ is a random vector with 
an $m+1$ components: ${\bf w}_{\tau}^{(i)}={\bf f}_{\tau}^{(i)}$ 
for $i=1,\ldots,m$ and 
${\bf w}_{\tau}^{(0)}=\tau,$\ 
${\bf f}_{\tau}^{(i)}$ $(i=1,\ldots,m)$
are independent standard Wiener processes.

It is well known that the following representation 
takes place \cite{Lipt}, \cite{7e}

\vspace{-1mm}
\begin{equation}
\label{um1x}
{\bf w}_{\tau}^{(i)}-{\bf w}_{t}^{(i)}=
\sum_{j=0}^{\infty}\int\limits_t^{\tau}
\phi_j(s)ds\ \zeta_j^{(i)},\ \ \ \zeta_j^{(i)}=
\int\limits_t^T \phi_j(s)d{\bf w}_s^{(i)},
\end{equation}

\vspace{3mm}
\noindent
where $\tau\in[t, T],$ $t\ge 0,$
$\{\phi_j(x)\}_{j=0}^{\infty}$ is an arbitrary complete 
orthonormal system of functions in the space $L_2([t, T]),$ and
$\zeta_j^{(i)}$ are independent standard Gaussian 
random variables for various $i$ or $j.$
Moreover, the series (\ref{um1x}) converges for any $\tau\in [t, T]$
in the mean-square sense.

Let ${\bf w}_{\tau}^{(i)p}-{\bf w}_{t}^{(i)p}$ be 
the mean-square approximation of the process
${\bf w}_{\tau}^{(i)}-{\bf w}_{t}^{(i)},$
which has the following form

\vspace{-5mm}
\begin{equation}
\label{um1xx}
{\bf w}_{\tau}^{(i)p}-{\bf w}_{t}^{(i)p}=
\sum_{j=0}^{p}\int\limits_t^{\tau}
\phi_j(s)ds\ \zeta_j^{(i)}.
\end{equation}

\vspace{3mm}

From (\ref{um1xx}) we obtain

\vspace{-3mm}
\begin{equation}
\label{um1xxx}
d{\bf w}_{\tau}^{(i)p}=
\sum_{j=0}^{p}
\phi_j(\tau)\zeta_j^{(i)} d\tau.
\end{equation}

\vspace{4mm}

Consider the following iterated Riemann--Stieltjes
integral

\vspace{1mm}
\begin{equation}
\label{um1xxxx}
\int\limits_t^T
\psi_k(t_k)\ldots \int\limits_t^{t_2}\psi_1(t_1)
d{\bf w}_{t_1}^{(i_1)p_1}\ldots d{\bf w}_{t_k}^{(i_k)p_k},
\end{equation}

\vspace{4mm}
\noindent
where $p_1,\ldots,p_k\in\mathbb{N},$\ \ $i_1,\ldots,i_k=0,1,\ldots,m,$ 

\begin{equation}
\label{um1xxx1}
d{\bf w}_{\tau}^{(i)p}=
\left\{\begin{matrix}
d{\bf f}_{\tau}^{(i)p}\ &\hbox{\rm for}\ \ \ i=1,\ldots,m\cr\cr\cr
d\tau^p\ &\hbox{\rm for}\ \ \ i=0
\end{matrix}
,\right.
\end{equation}

\vspace{4mm}
\noindent
and $d{\bf f}_{\tau}^{(i)p},$ $d\tau^p$ are defined by the relation (\ref{um1xxx}).

Let us substitute (\ref{um1xxx}) into (\ref{um1xxxx})

\begin{equation}
\label{um1xxxx1}
\int\limits_t^T
\psi_k(t_k)\ldots \int\limits_t^{t_2}\psi_1(t_1)
d{\bf w}_{t_1}^{(i_1)p_1}\ldots d{\bf w}_{t_k}^{(i_k)p_k}=
\sum\limits_{j_1=0}^{p_1}\ldots \sum\limits_{j_k=0}^{p_k}
C_{j_k \ldots j_1}\prod\limits_{l=1}^k \zeta_{j_l}^{(i_l)},
\end{equation}

\vspace{4mm}
\noindent
where 
$$
\zeta_j^{(i)}=\int\limits_t^T \phi_j(s)d{\bf w}_s^{(i)}
$$ 

\vspace{2mm}
\noindent
are independent standard Gaussian random variables for various 
$i$ or $j$ (in the case when $i\ne 0$),
${\bf w}_{s}^{(i)}={\bf f}_{s}^{(i)}$ for
$i=1,\ldots,m$ and 
${\bf w}_{s}^{(0)}=s,$

$$
C_{j_k \ldots j_1}=\int\limits_t^T\psi_k(t_k)\phi_{j_k}(t_k)\ldots
\int\limits_t^{t_2}
\psi_1(t_1)\phi_{j_1}(t_1)
dt_1\ldots dt_k
$$

\vspace{4mm}
\noindent
is the Fourier coefficient.

To best of our knowledge \cite{W-Z-1}-\cite{Watanabe}
the approximations of the Wiener process
in the Wong--Zakai approximation must satisfy fairly strong
restrictions
\cite{Watanabe}
(see Definition 7.1, pp.~480--481).
Moreover, approximations of the Wiener process that are
similar to (\ref{um1xx})
were not considered in \cite{W-Z-1}, \cite{W-Z-2}
(also see \cite{Watanabe}, Theorems 7.1, 7.2).
Therefore, the proof of analogs of Theorems 7.1 and 7.2 \cite{Watanabe}
for approximations of the Wiener 
process based on its series expansion (\ref{um1x})
should be carried out separately.
Thus, the mean-square convergence of the right-hand side
of (\ref{um1xxxx1}) to the iterated Stratonovich stochastic integral 
(\ref{str})
does not follow from the results of the papers
\cite{W-Z-1}, \cite{W-Z-2} (also see \cite{Watanabe},
Theorems 7.1, 7.2).
                                           
From the other hand, Theorems 7--16 from this 
paper can be considered as the proof of the
Wong--Zakai approximation for the iterated 
Stratonovich stochastic integrals (\ref{str}) of multiplicities 1 to 6
based on the approximation (\ref{um1xx}) of the Wiener process.
At that, the Riemann--Stieltjes integrals (\ref{um1xxxx}) converge
(according to Theorems 7--16)
to the appropriate Stratonovich 
stochastic integrals (\ref{str}). Recall that
$\{\phi_j(x)\}_{j=0}^{\infty}$ (see (\ref{um1x}), (\ref{um1xx}))
is a complete 
orthonormal system of Legendre polynomials or 
trigonometric functions 
in the space $L_2([t, T])$ (Theorems 9--13).
The system $\{\phi_j(x)\}_{j=0}^{\infty}$ 
can be arbitrary in Theorems~8, 14--16.

To illustrate the above reasoning, 
consider two examples for the case $k=2,$
$\psi_1(s),$ $\psi_2(s)\equiv 1;$ $i_1, i_2=1,\ldots,m.$

The first example relates to the piecewise linear approximation
of the multidimensional Wiener process (these approximations 
were considered in \cite{W-Z-1}-\cite{Watanabe}).

Let ${\bf b}_{\Delta}^{(i)}(t),$ $t\in[0, T]$ be the piecewise
linear approximation of the $i$th component ${\bf f}_t^{(i)}$
of the multidimensional standard Wiener process ${\bf f}_t,$
$t\in [0, T]$ with independent components
${\bf f}_t^{(i)},$ $i=1,\ldots,m,$ i.e.

\vspace{-2mm}
$$
{\bf b}_{\Delta}^{(i)}(t)={\bf f}_{k\Delta}^{(i)}+
\frac{t-k\Delta}{\Delta}\Delta{\bf f}_{k\Delta}^{(i)},
$$

\vspace{3mm}
\noindent
where 

\vspace{-1mm}
$$
\Delta{\bf f}_{k\Delta}^{(i)}={\bf f}_{(k+1)\Delta}^{(i)}-
{\bf f}_{k\Delta}^{(i)},\ \ \
t\in[k\Delta, (k+1)\Delta),\ \ \ k=0, 1,\ldots, N-1.
$$

\vspace{5mm}

Note that w.~p.~1

\vspace{-1mm}
\begin{equation}
\label{pridum}
\frac{d{\bf b}_{\Delta}^{(i)}}{dt}(t)=
\frac{\Delta{\bf f}_{k\Delta}^{(i)}}{\Delta},\ \ \
t\in[k\Delta, (k+1)\Delta),\ \ \ k=0, 1,\ldots, N-1.
\end{equation}

\vspace{4mm}

Consider the following iterated Riemann--Stieltjes
integral

\vspace{1mm}
$$
\int\limits_0^T
\int\limits_0^{s}
d{\bf b}_{\Delta}^{(i_1)}(\tau)d{\bf b}_{\Delta}^{(i_2)}(s),\ \ \ 
i_1,i_2=1,\ldots,m.
$$

\vspace{5mm}

Using (\ref{pridum})
and the additive property of Riemann--Stieltjes integrals, 
we can write w.~p.~1

\vspace{2mm}
$$
\int\limits_0^T
\int\limits_0^{s}
d{\bf b}_{\Delta}^{(i_1)}(\tau)d{\bf b}_{\Delta}^{(i_2)}(s)=
\int\limits_0^T
\int\limits_0^{s}
\frac{d{\bf b}_{\Delta}^{(i_1)}}{d\tau}(\tau)d\tau
\frac{d {\bf b}_{\Delta}^{(i_2)}}{d s}(s)
ds =
$$

\vspace{3mm}
$$
=
\sum\limits_{l=0}^{N-1}\int\limits_{l\Delta}^{(l+1)\Delta}
\left(
\sum\limits_{q=0}^{l-1}\int\limits_{q\Delta}^{(q+1)\Delta}
\frac{\Delta{\bf f}_{q\Delta}^{(i_1)}}{\Delta}d\tau+
\int\limits_{l\Delta}^{s}
\frac{\Delta{\bf f}_{l\Delta}^{(i_1)}}{\Delta}d\tau\right)
\frac{\Delta{\bf f}_{l\Delta}^{(i_2)}}{\Delta}ds=
$$

\vspace{3mm}
$$
=\sum\limits_{l=0}^{N-1}\sum\limits_{q=0}^{l-1}
\Delta{\bf f}_{q\Delta}^{(i_1)}
\Delta{\bf f}_{l\Delta}^{(i_2)}+
\frac{1}{\Delta^2}\sum\limits_{l=0}^{N-1}
\Delta{\bf f}_{l\Delta}^{(i_1)}
\Delta{\bf f}_{l\Delta}^{(i_2)}
\int\limits_{l\Delta}^{(l+1)\Delta}
\int\limits_{l\Delta}^{s}d\tau ds=
$$

\vspace{3mm}
\begin{equation}
\label{oh-ty}
=\sum\limits_{l=0}^{N-1}\sum\limits_{q=0}^{l-1}
\Delta{\bf f}_{q\Delta}^{(i_1)}
\Delta{\bf f}_{l\Delta}^{(i_2)}+
\frac{1}{2}\sum\limits_{l=0}^{N-1}
\Delta{\bf f}_{l\Delta}^{(i_1)}
\Delta{\bf f}_{l\Delta}^{(i_2)}.
\end{equation}

\vspace{6mm}

Using (\ref{oh-ty}) and the standard relation between Stratonovich and Ito
stochastic integrals, it 
is not difficult to show 
that

\vspace{1mm}
$$
\hbox{\vtop{\offinterlineskip\halign{
\hfil#\hfil\cr
{\rm l.i.m.}\cr
$\stackrel{}{{}_{N\to \infty}}$\cr
}} }
\int\limits_0^T
\int\limits_0^{s}
d{\bf b}_{\Delta}^{(i_1)}(\tau)d{\bf b}_{\Delta}^{(i_2)}(s)=
\int\limits_0^T
\int\limits_0^{s}
d{\bf f}_{\tau}^{(i_1)}d{\bf f}_{s}^{(i_2)}+
\frac{1}{2}{\bf 1}_{\{i_1=i_2\}}\int\limits_0^T ds=
$$

\vspace{3mm}
\begin{equation}
\label{uh-111}
=
{\int\limits_0^{*}}^T
{\int\limits_0^{*}}^s
d{\bf f}_{\tau}^{(i_1)}d{\bf f}_{s}^{(i_2)},
\end{equation}

\vspace{5mm}
\noindent
where $\Delta\to 0$ if $N\to\infty$ ($N\Delta=T$).
Obviously, (\ref{uh-111}) agrees with Theorem 7.1 (see \cite{Watanabe},
p.~486).

The next example relates to the approximation
of the Wiener process based on its series expansion
(\ref{um1x}) for $t=0$, where
$\{\phi_j(x)\}_{j=0}^{\infty}$ 
is a complete 
orthonormal system of Legendre polynomials or 
trigonometric functions 
in the space $L_2([0, T])$.

Consider the following iterated Riemann--Stieltjes
integral

\vspace{-1mm}
\begin{equation}
\label{abcd1}
\int\limits_0^T
\int\limits_0^{s}
d{\bf f}_{\tau}^{(i_1)p}d{\bf f}_{s}^{(i_2)p},\ \ \ 
i_1,i_2=1,\ldots,m,
\end{equation}

\vspace{3mm}
\noindent
where $d{\bf f}_{\tau}^{(i)p}$ is defined by the
relation
(\ref{um1xxx}).

Let us substitute (\ref{um1xxx}) into (\ref{abcd1}) 

\vspace{-1mm}
\begin{equation}
\label{set18}
\int\limits_0^T
\int\limits_0^{s}
d{\bf f}_{\tau}^{(i_1)p}d{\bf f}_{s}^{(i_2)p}=
\sum\limits_{j_1,j_2=0}^p
C_{j_2 j_1} \zeta_{j_1}^{(i_1)}\zeta_{j_2}^{(i_2)},
\end{equation}

\vspace{3mm}
\noindent
where 
$$
C_{j_2 j_1}=
\int\limits_0^T \phi_{j_2}(s)\int\limits_0^s
\phi_{j_1}(\tau)d\tau ds
$$

\vspace{3mm}
\noindent
is the Fourier coefficient; another notations 
are the same as in (\ref{um1xxxx1}).

As we noted above, approximations of the Wiener process that are
similar to (\ref{um1xx})
were not considered in \cite{W-Z-1}, \cite{W-Z-2}
(also see Theorems 7.1, 7.2 in \cite{Watanabe}).
Furthermore, the extension of the results of Theorems 7.1 and 7.2
\cite{Watanabe} to the case under consideration is
not obvious.

However, the authors of the works
\cite{KlPl2}
(Sect.~5.8, pp.~202--204), 
\cite{KPS} (pp.~82-84), \cite{KPW} (pp.~438-439),
\cite{Zapad-9} (pp.~263-264) use 
the Wong--Zakai approximation 
\cite{W-Z-1}-\cite{Watanabe} (without rigorous proof) within the frames
of the approach
\cite{Mi2} based on the series expansion 
of the Brownian bridge process.

On the other hand, we can apply the theory built in Chapters 1 and 2
of the monographs \cite{20aa}-\cite{12aa-afterxxx}. More precisely, 
using 
Theorems 9, 14 from this paper 
we obtain from (\ref{set18}) the desired result

\vspace{-1mm}
$$
\hbox{\vtop{\offinterlineskip\halign{
\hfil#\hfil\cr
{\rm l.i.m.}\cr
$\stackrel{}{{}_{p\to \infty}}$\cr
}} }
\int\limits_0^T
\int\limits_0^{s}
d{\bf f}_{\tau}^{(i_1)p}d{\bf f}_{s}^{(i_2)p}=
\hbox{\vtop{\offinterlineskip\halign{
\hfil#\hfil\cr
{\rm l.i.m.}\cr
$\stackrel{}{{}_{p\to \infty}}$\cr
}} }
\sum\limits_{j_1,j_2=0}^p
C_{j_2 j_1} \zeta_{j_1}^{(i_1)}\zeta_{j_2}^{(i_2)}=
$$

\vspace{2mm}
\begin{equation}
\label{umen-bl}
=
{\int\limits_0^{*}}^T
{\int\limits_0^{*}}^s
d{\bf f}_{\tau}^{(i_1)}d{\bf f}_{s}^{(i_2)}.
\end{equation}

\vspace{4mm}

From the other hand, by Theorem 8
(see (\ref{za2})) for the case
$k=2$ we obtain from (\ref{set18}) the following relation

\vspace{-2mm}
$$
\hbox{\vtop{\offinterlineskip\halign{
\hfil#\hfil\cr
{\rm l.i.m.}\cr
$\stackrel{}{{}_{p\to \infty}}$\cr
}} }
\int\limits_0^T
\int\limits_0^{s}
d{\bf f}_{\tau}^{(i_1)p}d{\bf f}_{s}^{(i_2)p}=
\hbox{\vtop{\offinterlineskip\halign{
\hfil#\hfil\cr
{\rm l.i.m.}\cr
$\stackrel{}{{}_{p\to \infty}}$\cr
}} }
\sum\limits_{j_1,j_2=0}^p
C_{j_2 j_1} \zeta_{j_1}^{(i_1)}\zeta_{j_2}^{(i_2)}=
$$

\vspace{2mm}
$$
=
\hbox{\vtop{\offinterlineskip\halign{
\hfil#\hfil\cr
{\rm l.i.m.}\cr
$\stackrel{}{{}_{p\to \infty}}$\cr
}} }
\sum\limits_{j_1,j_2=0}^p
C_{j_2 j_1} \biggl(\zeta_{j_1}^{(i_1)}\zeta_{j_2}^{(i_2)}-
{\bf 1}_{\{i_1=i_2\}}{\bf 1}_{\{j_1=j_2\}}\biggr)+
{\bf 1}_{\{i_1=i_2\}}\sum\limits_{j_1=0}^{\infty}
C_{j_1 j_1}=
$$

\vspace{2mm}
\begin{equation}
\label{umen-blx}
=
\int\limits_0^T
\int\limits_0^{s}
d{\bf f}_{\tau}^{(i_1)}d{\bf f}_{s}^{(i_2)}+
{\bf 1}_{\{i_1=i_2\}}\sum\limits_{j_1=0}^{\infty}
C_{j_1 j_1}.
\end{equation}

\vspace{5mm}

Since
$$
\sum\limits_{j_1=0}^{\infty}
C_{j_1 j_1}=\frac{1}{2}\sum\limits_{j_1=0}^{\infty}
\left(\int\limits_0^T \phi_j(\tau)d\tau\right)^2
=\frac{1}{2}
\left(\int\limits_0^T \phi_0(\tau)d\tau\right)^2=\frac{1}{2}
\int\limits_0^T ds,
$$

\vspace{5mm}
\noindent
then from (\ref{umen-blx}) and the standard 
relation between Stratonovich and Ito stochastic integrals
we obtain (\ref{umen-bl}).

\end{document}